# The zero loci of ℤ/2 harmonic spinors in dimensions 2, 3 and 4.

Clifford Henry Taubes[†]

Department of Mathematics
Harvard University
Cambridge MA, 02138  USA

chtaubes@math.harvar.edu

ABSTRACT:  Supposing that X is a Riemannian manifold, a ℤ/2 spinor on X is defined by a data set consisting of a closed set in X to be denoted by Z, a real line bundle over X−Z, and a nowhere zero section on X−Z of the tensor product of the real line bundle and a spinor bundle.  The set Z and the spinor are jointly constrained by the following requirement:  The norm of the spinor must extend across Z as a continuous function vanishing on Z.  In particular, the vanishing locus of the norm of the spinor is the complement of the set where the real line bundle is defined, and hence where the spinor is defined.  The ℤ/2 spinor is said to be harmonic when it obeys a first order Dirac equation on X−Z. This monograph analyzes the structure of the set Z for a ℤ/2 harmonic spinor on a manifold of dimension either two, three or four.

[†]Supported in part by the National Science Foundation.

# 1. Introduction

To set the stage for what is to come, let X denote an oriented, Riemannian manifold of dimension either 2, 3 or 4. Suppose in what follows that $\mathbb{S}^+$ and $\mathbb{S}^-$ are vector bundles over X with fiber metric and metric compatible covariant derivatives. These are denoted by $\nabla$. Assume that $\mathcal{D}$ is a first order, elliptic operator mapping $C^\infty(X;\mathbb{S}^+)$ to $C^\infty(X;\mathbb{S}^-)$ that can be written as

$$\mathcal{D} = \mathfrak{S}(\nabla)$$
(1.1)

with $\mathfrak{S}$ denoting the principal symbol of $\mathcal{D}$. This homomorphism from TX to the bundle $\text{Hom}(\mathbb{S}^+;\mathbb{S}^-)$ is assumed to have properties that are listed below in (1.2). The list uses $\langle,\rangle$ to denote the Riemannian metric on TX.

- $\mathfrak{S}^\dagger(v)\mathfrak{S}(u) + \mathfrak{S}^\dagger(u)\mathfrak{S}(v) = -2\langle v,u\rangle$ *at all points* $p \in X$ *and for all* $u, v \in TX|_p$.
- $[\nabla,\mathfrak{S}] = 0$.

(1.2)

Let U denote for the moment a given open set in X. The terminology in what follows has a real line bundle over U being the associated $\mathbb{R}$ bundle to a principal $\{1,-1\}$ bundle over U. Granted this as the definition, then a real line bundle has a canonical fiber metric; and the space of sections of a real line bundle has canonical metric compatible covariant derivative. In particular, this covariant derivative can be used to view $\mathcal{D}$ over U as a map from the space of sections over U of $\mathbb{S}^+\otimes\mathcal{I}$ to the space of sections over U of $\mathbb{S}^-\otimes\mathcal{I}$. By the same token, the covariant derivative for sections of $\mathbb{S}^+$ canonically defines such a covariant derivative for sections over U of $\mathbb{S}^+\otimes\mathcal{I}$. This is denoted also by $\nabla$.

Of interest in this paper are data sets consisting of the following:

- *A continuous non-negative function f with f > 0 somewhere.*
- *Let Z denote $f^{-1}(0)$. A real line bundle $\mathcal{I} \to X-Z$*
- *A section $\nu$ of $\mathbb{S}^+\otimes\mathcal{I}$ over X–Z that obeys $\mathcal{D}\nu = 0$ and is such that*
    i) $|\nu| = f$.
    ii) *The function $|\nabla\nu|^2$ is integrable on X–Z.*

(1.3)

Standard elliptic regularity theorems say that $\nu$ is smooth on X–Z. The third bullet of (1.3) begs the following question:

$$\textit{What can be said about Z?}$$
(1.4)

Four theorems that follow momentarily give something of an answer to this question. when the dimension of X is 2, 3 and 4. These theorems make the following additional assumption about $f$:

*There exists $\varepsilon > 0$ with the following significance: If $p \in X$ and $f(p) = 0$, and if r is small and positive, then $\int_{dist(p,\cdot)<r} f^2 \leq x_p r^{dim(X)+\varepsilon}$ with $x_p$ being independent of r.*

(1.5)

It is conceivable that (1.5) follows from (1.3) but the author has no proof that this is so except when the dimension of X is equal to 2.

The first theorem talks about the case when the dimension of X is two

**Theorem 1.1**: *Suppose that the dimension of X is two. Then the zero locus of $|v|$ is a locally finite set. Supposing that p is a zero of $|v|$, there exists a positive half-integer, $N_0$, with the following significance: Let $D \subset X$ be a disk centered at p with p being the only zero of $|v|$ in its closure. If $N_0$ is an integer, then the bundle $\mathcal{I}$ on D–p is isomorphic to the product $\mathbb{R}$ bundle; and if $N_0$ is half an odd integer, than $\mathcal{I}$ on D–p is isomorphic to the Möbius bundle. In either case, $|v| = c\, dist(\cdot,p)^{N_0} + \cdots$ near p with c being positive number and with the unwritten terms being higher order in the distance to p.*

This theorem is proved in Section 5a. By way of an example when X has dimension 2, take X to be a Riemann surface with its associated complex structure and take $\mathbb{S}^+$ to be $T^{1,0}X$, the $(1, 0)$ part of the complexified cotangent bundle. Take $\mathbb{S}^-$ to be $T^{1,0}X \otimes T^{0,1}X$ and take $\mathcal{D}$ to be the d-bar operator. Then take $v$ to be the square root of a holomorphic, quadratic differential.

The next three theorems describe Z when the dimension of X is either 3 or 4. The first of these theorems introduces the notion of a *point of discontinuity* for the bundle $\mathcal{I}$. A point $p \in Z$ is a point of discontinuity for $\mathcal{I}$ if $\mathcal{I}$ is not isomorphic to the product $\mathbb{R}$ bundle on the complement of Z in any neighborhood of p. The points of discontinuity for $\mathcal{I}$ are the curious and interesting points in Z because if p is not a point of discontinuity for $\mathcal{I}$, then $v$ can be viewed near p as the restriction of an honest, $\mathbb{R}$ valued harmonic section of $\mathbb{S}^+$ to the complement of its zero locus. This harmonic section is defined and $C^\infty$ on the whole of some neighborhood of p. In particular, there is nothing novel about $v$ near p.

**Theorem 1.2**: *Suppose that the dimension of X is three or four. The set of points of discontinuity for $\mathcal{I}$ is the closure of an open subset in Z with the structure of a codimension 2, differentiable submanifold in X*

This theorem is restated and proved as Proposition 10.5 in Section 10d.

The next theorem gives a picture of Z when it is viewed at $\mathcal{O}(1)$ length scales. This theorem introduces the notions of *geodesic arcs* and *geodesic disks* through a given point. Let d denote the dimension of X. Supposing that p is the point in question, a geodesic arc through p is the image in X via a Gaussian coordinate chart map centered at p of the small |x| part of a line through the origin in $\mathbb{R}^d$; and a geodesic disk through p is the corresponding image of the small |x| part of plane through the origin in $\mathbb{R}^d$. The theorem also introduces the notion of a *union of geodesic half-disks with common boundary arc* through p. The definition of this notion requires first the definition of a half-plane through the origin in $\mathbb{R}^d$: Supposing that $\Pi$ is a plane through the origin and $f$ is a non-trivial linear function on $\Pi$, then the $f \geq 0$ part of $\Pi$ is said to be a half-plane through the origin. The line where $f = 0$ is said to be the boundary line of the half plane. A union of geodesic half-disks with common boundary arc is the image via a Gaussian coordinate chart map centered at p of the small |x| part of the union of two or more half-planes all with common boundary line. A geodesic arc or disk or union of half-disks with common boundary arc in a given open ball is said to be *proper* if the set in question is a relatively closed, path connected subset of the ball.

**Theorem 1.3**: *Assume that* X *has dimension three or four. Use* d *to denote this dimension. The set* Z *has Hausdorff dimension at most* d-2. *Moreover,* Z *looks to a first approximation like a stratified d-2 dimensional set in the following sense: Given numbers* $\theta \in (0,1)$ *and* $z > 1$, *there are finite collections* $\mathfrak{U}$ *and* $\mathfrak{V}$ *of the following sort:*
- *Both* $\mathfrak{U}$ *and* $\mathfrak{V}$ *are composed of balls of radius at most* $z^{-1}$.
- *The balls from* $\mathfrak{U}$ *are pairwise disjoint and* $\sum_{B \in \mathfrak{U}} (\text{radius}(B))^\theta < z^{-1}$.
- *The balls in the set* $\mathfrak{V}$ *cover the part of* Z *in* X$-(\cup_{B \in \mathfrak{U}} B)$.
- *Supposing that* B *is a ball from* $\mathfrak{V}$, *let* p *denote its center point and* r *its radius.*
    a) *If* d = 3, *then the part of* Z *in* B *is contained in the radius* $z^{-1}$r *tubular neighborhood of a proper geodesic arc through* p
    b) *If* d = 4, *then the part of* Z *in* B *is contained in the radius* $z^{-1}$r *tubular neighborhood of either a proper geodesic arc through* p *or a proper geodesic disk through* p *or a proper union of half-disks through* p *with common boundary arc.*

The assertion about the Hausdorff dimension of Z is restated and proved as Proposition 7.1 in Section 7a. The remaining assertions of this theorem are restated in the d = 3 case as Proposition 7.2 and in the d = 4 case as Proposition 7.3. These propositions are proved in Sections 7c and 7d respectively.

A paper by De Lellis and Spadaro [DS] about work by Almgren [Al1] on a closely related problem with similar 'points of discontinuity' came to the author's attention near the end of the writing of this manuscript. Almgren, and De Lellis and Spadaro prove the analog of the Hausdorff dimension assertion of Theorem 1.3 in their related context.

The final theorem in this introduction gives a picture of Z at arbitrarilly small length scales near an open and dense subset of Z. This theorem introduces the notion of a codimension 2 Lipshitz graph in $\mathbb{R}^d$ with d being either 3 or 4. Such a graph is defined by data consisting of a Euclidean splitting of $\mathbb{R}^d$ as $\mathbb{R}^2 \times \mathbb{R}^{d-2}$, an open neighborhood $U \subset \mathbb{R}^{d-2}$ about the origin, and a Lipshitz map $\varphi: U \to \mathbb{R}^2$ mapping the origin to the origin. The associated codimension d - 2 Lipshitz graph is the graph of $\varphi$, this being the set of points that can be written with respect to the splitting of $\mathbb{R}^d$ as $\{(\varphi(t), t): t \in U\}$.

**Theorem 1.4**: *Assume that the dimension of* X *is either three or four. Use* d *to again denote this dimension. The set* Z *has a open, dense subset where each point has a neighorhood that is identified by a coordinate chart map with a subset of a codimension 2 Lipshitz graph in* $\mathbb{R}^d$.

This theorem is proved in Section 10b.

The final theorem in this introduction summarize results from a number of the subsequent lemmas and propositions about the behavior of the norm of ν.

**Theorem 1.5**: *Assume that the dimension of* X *is either three or four. The norm of the harmonic section* ν *is uniformly Hölder continuous on* X. *Moreover, if* p ∈ Z, *then there exists* $N_0 \geq \frac{1}{2}$ *and* c > 0 *such that* |ν| *near* p *obeys* |ν| ≤ c dist(·,p)$^{N_0}$. *This number* $N_0$ *is equal to half a positive integer on a dense, open set whose complement is countable if the dimension of* X *is three and has Hausdorff dimension at most 1 if the dimension of* X *is four*.

The author has no examples where the number $N_0$ in this this theorem is not half an integer. It may be that it is always half an integer. A guide for the proof of Theorem 1.5 is given below in Section 1a.

The author explains in [T1] and [T2] how data sets that are described by (1.3) and (1.5) arise from suitably rescaled limits of certain sequences of PSL(2;$\mathbb{C}$) connections on 2, 3 and 4 dimensional manifolds. Similar data sets arise in [GU], [HW] and [Tan]. In the cases of [T1] and [T2], what is denoted by ν is an $\mathcal{I}$ valued harmonic 1-form on X–Z. Theorems that are slightly stronger than Theorems 1.2 -1.4 are stated in [T1] and in the first arXiv version of [T2] for their versions of Z, $\mathcal{I}$ and ν; but the proofs of the latter theorems have a gap due to a silly but none-the-less painful error. This error was pointed out to the author by Thomas Walpuski. The same error appears in both [T1] and in the original version of [T2], the error being in the proof of a lemma to the effect that Z has a unique tangent cone at each of its points. (This is Lemma 8.8 in [T1] and Lemma 3.8 in the original version of [T2].) As explained in what follows, Theorem 1.3 doesn't need unique tangent cones and there are enough points in Z with unique tangent cones to prove Theorems 1.2 and 1.4.

By way of a parenthetical remark, it is likely that all points in Z do infact have unique tangent cones.

**a) Remarks about the proofs of Theorems 1.1-1.5**

A version of the *frequency* function that was first introduced by Almgren [Al2] plays a central role in the proofs of Theorems 1.1-1.4. This is a function on $(0, \infty)$ that can be assigned to each point in Z: If $p \in Z$, then up a multiplicative factor of little consequence, it is defined by the rule

$$r \to r^{d-2} \left( \int_{B_r} |\nabla v|^2 \right) \left( \int_{\partial B_r} |v|^2 \right)^{-1}$$

(1.6)

where $B_r$ denotes the ball of radius r centered at p and $\partial B_r$ denotes the boundary of the closure of this same ball. This function is denoted by $N_{(p)}$.

The function $N_{(p)}$ turns out to be nearly monotonic increasing; and its monotonicity can be used to study the structure of Z and $v$ near p at ever smaller length scales. The point being that if $N_{(p)}$ is monotonic, then it is necessarily bounded as $r \to 0$; and if $N_{(p)}$ is bounded as $r \to 0$, then it follows from the formula in (1.6) that the r dependent rescaling

$$v \to \left( \frac{1}{r^{d-1}} \int_{\partial B_r} |v|^2 \right)^{-1/2} v \ ,$$

(1.7)

when viewed using an r dependent rescaling of Gaussian coordinates on a neighborhood of p, has weak Sobolev $L^2_1$ limits on compact subsets of $\mathbb{R}^d$ as $r \to 0$. These weak limits are non-trivial and they obey (1.3) and (1.5) on $\mathbb{R}^d$. Moreover, because $N_{(p)}$ is monotonic, the weak limits are necessarily homogeneous with respect to the scaling diffeomorphisms of $\mathbb{R}^4$ that multiply the coordinates by any positive number. The limiting process and the properties of these rescaling limits of $v$ can be used to infer properties of $v$ at $\mathcal{O}(r)$ length scales with r small, but non-zero. This is how versions of the frequency function are typically used [DF], [H], [HHL].

Unfortunately, the sequence of $r \to 0$ rescalings of $v$ may have many $L^2_1$ weak limits when p is a point of discontinuity for $\mathcal{I}$. This can't happen if p is not a point of discontinuity for r because $v$ near p is a smooth, $\mathbb{R}$ valued harmonic section of $\mathbb{S}^+$; and as such, it has a non-trivial Taylor's expansion near p. The erroneous lemmas in [T1] and in the original version of [T2] claimed that there was a unique rescaling limit in all cases.

The lack of unique r → 0 rescaling limits of $v$ has no bearing on the proof of Theorem 1.3. The existence of unique r → 0 rescaling limits would, however, imply a stronger assertion, this being that Z is contained in a countable union of d-2 dimensional rectifiable sets.

The proofs of Theorems 1.2 and 1.4 use the fact that there are unique r → 0 rescaling limits of $v$ at certain distinguished points in Z. To say more about this, let $N_{(p)}(0)$ denote the r → 0 limit of the frequency function at p. (This $N_{(p)}(0)$ is the number $N_0$ in Theorem 1.5.) The r → 0 limit of $N_{(p)}$ always exists because $N_{(p)}$ is nearly monotonic. As explained in the upcoming Sections 8 and 9, there is a unique r → 0 rescaling limit of $v$ at p if the function on Z given by the rule p → $N_{(p)}(0)$ is locally constant near p. This is proved using a Taylor's approximation for $v$ when p is not a point of discontinuity for $\mathcal{I}$. In the case of a point of discontinuity for $\mathcal{I}$, it is proved using a matrix version of the frequency function that is defined as follows: A suitable orthonormal frame for TX near p is chosen so as to define the d directional covariant derivatives of $v$ along the frame vectors on a ball centered at p. Denote this set by $\{\nabla_a v\}_{a \in \{1,...,d\}}$ with d being the dimension of X. The matrix function of r has entries $T_{ab}(r)$ for a, b ∈ {1, …, d} given up to an inconsequential factor by the rule

$$r \to r^{d-2} \left( \int_{B_r} \langle \nabla_a v, \nabla_b v \rangle \right) \left( \int_{\partial B_r} |v|^2 \right)^{-1} ,$$

(1.8)

with $\langle\,,\,\rangle$ now denoting the inner product on $\mathbb{S}^+$. This is a non-negative matrix whose trace is the frequency function in (1.6). It therefore has one eigenvalue that is $\mathcal{O}(1)$. It follows from the equation $\mathcal{D}v = 0$ that it has at least two $\mathcal{O}(1)$ eigenvalues. The span of the eigenvectors with very small eigenvalues turns out to be a good proxy at the scale defined by the given value of r for a tangent cone to Z at p. These proxy tangent cones are used in Section 9 to prove an assertion to the affect that the r → 0 rescalings of $v$ at p have a unique limit when $N_{(\cdot)}(0)$ is constant near p. The existence of unique r → 0 rescalings of $v$ are used in Section 10 to prove Theorems 1.2 and 1.4.

What follows is a guide for where to find the proofs of the various assertions in Theorem 1.5: The Hölder continuity of $|v|$ is proved in Section 3e. Meanwhile, the local bound on the Hölder norm near any given point in Z is restated in Lemma 3.3 and proved in Section 3c. The lower bound on what is denoted here by $N_0$ is restated as Lemma 6.4 and the fact that $N_0$ is equal to half a positive integer on an open, dense set is follows from Lemma 6.6. The bound on the Hausdorff dimension of the set where $N_0$ is not a half a positive integer comes from Proposition 7.1.

Many of the arguments in this paper are close kin to arguments in the published version of [T1] and the first arXiv version of [T2], especially in the first four sections. These are the sections that define Almgren's frequency function in the context at hand and the $r \to 0$ scaling limits of $v$.

**b) Examples**

The following is an example of a data set $(Z, \mathcal{I}, v)$ of that obeys (1.3) and (1.5). The ambient manifold is $\mathbb{R}^4$ which is viewed as $\mathbb{C}^2$ with complex coordinates z and w. Take Z to be the union of the $z = 0$ plane and the $w = 0$ plane. Define $\mathcal{I}$ to be the restriction to $\mathbb{C}^2-Z$ of the non-trivial real line bundle over the complement in the ambient $\mathbb{C}^2$ of the $z = 0$ plane. The real part of the 1-form

$$v = \tfrac{3}{2} z^{1/2} w^2 \, dz + 2w z^{3/2} \, dw$$

(1.9)

is a harmonic, $\mathcal{I}$-valued 1-form on $\mathbb{C}^2-Z$ that vanishes on Z. The points of discontinuity for $\mathcal{I}$ are the points in the $z = 0$ plane and the submanifold from Theorem 1.2 is the complement of the origin in the $z = 0$ plane. The function $N_{(\cdot)}(0)$ in this case has value $\tfrac{1}{2}$ on the complement of the origin in the $z = 0$ plane and value 1 on the complement of the origin in the $w = 0$ plane. The value of $N_{(\cdot)}(0)$ at the origin in $\mathbb{C}^2$ is $\tfrac{5}{2}$.

For a second example, take $X = \mathbb{R}^4$ again, but now viewed as $\mathbb{C} \times \mathbb{R}^2$. Use the complex coordinate z for the $\mathbb{C}$ factor and use the real coordinates (s, t) for the $\mathbb{R}^2$ vactor. Take Z to be the union of the $z = 0$ plane and the half planes where both s and the real part of $z^{3/2}$ are zero. Define $\mathcal{I}$ to be the restriction to $\mathbb{R}^4-Z$ of the non-trivial line bundle over the complement of the $z = 0$ plane. The real part of the 1-form

$$v = \tfrac{3}{2} z^{1/2} s \, dz + z^{3/2} \, ds$$

(1.10)

is a harmonic, $\mathcal{I}$-valued 1-form on $\mathbb{R}^4-Z$ that vanishes on Z. The points of discontinuity for $\mathcal{I}$ are the points in the $z = 0$ plane; and the submanifold from Theorem 1.2 is the complement of the origin in the $z = 0$ plane. The function $N_{(\cdot)}(0)$ in this case has value $\tfrac{1}{2}$ on the complement of the origin in the $z = 0$ plane and value 1 on the complement of the origin in the half planes where both s and the real part of $z^{3/2}$ are zero. The value of $N_{(\cdot)}(0)$ on the line where $z = s = 0$ is $\tfrac{3}{2}$.

### c) Table of contents

The following is a table of contents for this paper.

There is also a list of references at the very end.

### d) Notation

The respective inner products on any tensor bundles, on $\mathbb{S}^+$, $\mathbb{S}^-$, their tensor products with $\mathcal{I}$ and likewise with $\mathcal{I}$ and tensor bundles are all denoted by $\langle \, , \, \rangle$. In the cases when $\mathbb{S}^+$ and $\mathbb{S}^-$ are complex bundles with a Hermitian metric, this inner product is the real part of the Hermitian inner product. In particular, the inner product is always a real number. The associated norm is denoted by $|\cdot|$. Covariant derivatives on all of these bundles are denoted by $\nabla$ with it understood that the covariant derivatives on the tensor bundles are defined by the Riemannian metric's Levi-Civita connection.

What is denoted by $c_0$ here is in all cases a number that is greater than 1; its precise value is assumed to be $\mathcal{O}(1)$ and in any case irrelevant to the matters at hand. If X is not compact then the value of $c_0$ in any given instance can depend on the apriori specification of a subset in X with compact closure in the interior of X. This restriction to such a subset will for the most part be implicit in what follows. In any event, the value $c_0$ can be assumed to increase between its successive appearances.

What is denoted by $r_0$ in what follows is a positive number chosen so that the ball of radius $100r_0$ about any given point in X is well inside a Gaussian coordinate chart centered at that point. If X is not compact, then it will go without saying that the definition of $r_0$ (like that of $c_0$) requires the restriction to a fixed open set in X with compact closure in the interior of X.

With regards to points and balls, let $p \in X$ denote a chosen point. Given a number $r \in (0, r_0]$, what is denoted by $B_r$ is the open ball of radius r centered at p. The boundary of the corresponding closed ball is denoted by $\partial B_r$.

This monograph also use χ to denote a chosen, non-increasing function on $\mathbb{R}$ that is equal to 1 on $(-\infty, \frac{1}{4}]$ and equal to 0 on $[\frac{3}{4}, \infty)$.

**e) Acknowledgement**

The author owes a large debt to Thomas Walpuski for his close reading of [T1] and for his patience with the author's mistakes. A similar debt is owed to Ben Mares.

The National Science Foundation supported a part of the writing of this manuscript via the grant DMS-1401192.

## 2. Integrals of $|\nabla v|^2$

This section states and proves various assertions about integrals of $|\nabla v|^2$. These assertions are for the most part consequences of the fact that the function $|v|^2$ obeys the second order differential equation

$$\tfrac{1}{2} d^\dagger d |v|^2 + |\nabla v|^2 + \langle v, \mathcal{R} v \rangle = 0$$

(2.1)

with $\mathcal{R}$ being an endomorphism of $\mathbb{S}^+$. This equation is obtained by using the conditions in (1.2) to write the equation $\mathcal{D}^\dagger \mathcal{D} v = 0$ as

$$\nabla^\dagger \nabla v + \mathcal{R} v = 0$$

(2.2)

and then taking the inner product of both sides of the latter identity with $v$.

**a) The integral of $|\nabla v|^2$ over X**

The function $|\nabla v|^2$ is assumed to have bounded integral on X–Z. This understood, it proves useful to define this function as an $L^1$ function on the whole of X with its integral over X defined to be equal to its integral on X–Z. The following lemma justifies in part this definition:

**Lemma 2.1**: *There exists $\kappa > 0$ with the following significance: Fix $\rho \in (0, \kappa^{-1}]$. Then the integral of $|\nabla v|^2$ over the set in X where $0 < |v| < \rho$ is no greater than $\kappa \rho^{1/\kappa}$.*

*Proof of Lemma 2.1*: To set some notation, suppose that a positive number $\rho$ has been specified. What is denoted below by $\chi_\rho$ is the function on X given by $\chi(1 - \rho^{-1}|v|)$. This function equals 0 where $|v| < \rho$ and it equals 1 where $|v| > 2\rho$. What is denoted below by $\chi'_\rho$ is the function on X that is obtained by evaluating the derivative of $\chi$ at $1 - \rho^{-1}|v|$.

With $\rho$ given, fix for the moment $\varepsilon > 0$ but with $\varepsilon \ll \rho$. Let $f_{\varepsilon,\rho}$ denote the integral over X of the function $(1 - \chi_{\varepsilon\rho})\chi_\varepsilon |\nabla v|^2$. To see about the size of $f_{\varepsilon,\rho}$, multiply both sides of (2.1) by $(1 - \chi_\rho)\chi_\varepsilon$ and then integrate the result over X. Since $\chi'$ is non-positive, an integration by parts leads to an inequality that reads

$$f_{\varepsilon,\rho} < \frac{c_0}{1+c_0} f_{\varepsilon,2\rho} + c_0 \rho^2.$$

(2.3)

Let n be the greatest integer less than $c_0^{-1}|\ln\rho|$. Iterate (2.3) by successively replacing $\rho$ by first $2\rho$, then $2^2\rho$ and so up to $2^n\rho$ leads to a bound on $f_{\varepsilon,\rho}$ by $c_0 \rho^\delta$ with $\delta > c_0^{-1}$. This bound is independent of $\varepsilon$ and taking $\varepsilon$ to zero leads to the bound in the lemma.

The inequality $|d|v|| \leq |\nabla v|$ holds where $|v| > 0$; so it follows that

$$\int_X |d|v||^2 \leq \int_X |\nabla v|^2$$

(2.4)

because $d|v|$ can be defined to be zero on any measurable set where $|v|$ is identically zero. The function $|v|$ is therefore an $L^2_1$ function on X. In particular, the inequalities in the upcoming (2.5) can be invoked if needs be. To set the stage for these inequalities, suppose that $p \in X$ is a given point and that $r \in (0, r_0]$. Use $B_r$ to denote the open, radius r ball centered at p. What is denoted by $f$ in (2.5) is an $L^2_1$ function on $B_r$.

- $(\int_{B_r} |f|^k)^{2/k} \leq c_0 \int_{B_r} (|df|^2 + r^{-2} |f|^2)$ *if* $d = 2$ *and* $k \in [1, \infty)$.
- $(\int_{B_r} |f|^{2d/(d-2)})^{(d-2)/d} \leq c_0 \int_{B_r} (|df|^2 + r^{-2} |f|^2)$ *if* $d > 2$.
- $\int_{B_r} \frac{1}{\mathrm{dist}(\cdot,p)^2} |f|^2 \leq c_0 \int_{B_r} (|df|^2 + r^{-2} |f|^2)$ *if* $d > 2$.
- *If f has compact support in* $B_r$ *or if the integral of f over* $B_r$ *is zero, then the preceding two inequalities hold without the* $r^{-2}|f|^2$ *term in the integrand.*
- $\int_{B_r} f^2 \leq c_0 (r^2 \int_{B_r} |df|^2 + r \int_{\partial B_r} f^2)$.

(2.5)

The first two of these inequalities are Sobolev inequalities and the third is Hardy's inequality. The fifth inequality is proved using much the same sort of integration by parts that proves Hardy's inequality. Given the first three inequalities, then the assertion in the fourth inequality follows from fifth.

The next lemma says something about the integral of $|\nabla v|^2$ over more of X.

**Lemma 2.2**: *Let* $U \subset X$ *denote an open set and let* $V \subset U$ *denote an open subset with compact closure. There exists* $\kappa_{V,U} > 1$ *such that* $\int_V |\nabla v|^2 \leq \kappa_{V,U} \int_U |v|^2$.

*Proof of Lemma 2.2*: Fix a non-negative function that is equal to 1 on V and equal to 0 on the complement of U. Denote this function by $\sigma$. Having fixed $\rho \in (0, 1]$, multiply both sides of (2.1) by $\sigma^2 \chi_\rho$. Integrate the resulting identity by parts to remove a derivative from the resulting $\sigma^2 \chi_\rho d^\dagger d|v|^2$ at the expense of adding two terms, one with integrand $4\sigma \chi_\rho |v| \langle d\sigma, \nabla |v| \rangle$ and the other with the integrand $-2\chi'_\rho \sigma^2 \rho^{-1} |v| ||\nabla|v|||^2$. Granted this integration by parts, the triangle inequality bounds the integral of $|\nabla v|^2$ on V by the sum of two terms. The first is a multiple of the integral of $|v|^2$ on U with the multiplier depending only on the chosen function $\beta$. The second term is $c_0$ times the integral of $|\nabla v|^2$ on the $|v| \leq 2\rho$ part of U. Viewed as a function of $\rho$, the latter integral has limit zero as $\rho$ limits to zero. This is what is said by Lemma 2.1. The terms that remain as $\rho$ limits to zero are bounded and give the inequality that is asserted by the lemma.

**b) The Green's function and $|\nabla v|^2$**

Fix $r \in (0, r_0]$. Supposing that $q \in B_{r/2}$, let $G_q$ denote the Dirichelet Green's function for the Laplacian $d^\dagger d$ on $B_r$ with pole at q. This is a positive function on $B_r - q$ that can be written as

- $G_q = \frac{1}{\omega_d} \frac{1}{\text{dist}(\cdot, q)^{d-2}} + \mathfrak{e}$   *if* $d = 3$ *or* $4$,
- $G_q = -\frac{1}{2\pi} \ln(\frac{1}{r} \text{dist}(\cdot, q)) + \mathfrak{e}$   *if* $d = 2$,

(2.6)

with $\omega_d$ denoting $4\pi$ in the case $d = 3$ and $2\pi^2$ in the case $d = 4$. What is denoted here by $\mathfrak{e}$ is a function on U−p with norm bounded by $c_0$ if $d = 2$ and by $c_0 r^{-1} \text{dist}(\cdot, q)^{3-d}$ if $d = 3$ or if $d = 4$. In all the cases $d = 2, 3$, and $4$, the derivative of $G_q$ obeys $|\nabla G_q| \leq c_0 \text{dist}(\cdot, q)^{1-d}$.

Suppose first that q is not in Z. Let $\sigma$ now denote a smooth, non-negative function that is equal to 1 on V and zero near the boundary of U. Fix $\rho \in (0, c_0^{-1}|v|(p)]$ for the moment and multiply both sides of (2.1) by $\chi_\rho \sigma G_q$ with $\chi_\rho$ as defined in the proof of Lemma 2.1. Integrate the resulting identity over X and then integrate by parts to obtain an equation that reads

$$\tfrac{1}{2}\sigma(q)|v|^2(q) + \int_{B_r}\sigma G_q\,|\nabla v|^2 = -\int_{B_r}\sigma G_q\langle v,\mathcal{R}v\rangle - \int_{B_r}(\tfrac{1}{2}d^\dagger d\sigma\, G_q - \langle d\sigma, dG_q\rangle)|v|^2 + \mathfrak{r}_{q,\rho},$$

(2.7)

with $\mathfrak{r}_{q,\rho}$ being a term whose norm is bounded by $c_0 \operatorname{dist}(Z,q)^{-d+1}\rho^{1/c_0}$. By way of an explanation, this bound on $|\mathfrak{r}_{q,\rho}|$ follows from Lemma 2.1. Taking $\rho \to 0$ in (2.7) leads to the identity

$$\tfrac{1}{2}\sigma(q)|v|^2(q) + \int_{B_r}\sigma G_q\,|\nabla v|^2 = -\int_{B_r}\sigma G_q\langle v,\mathcal{R}v\rangle - \int_{B_r}(\tfrac{1}{2}d^\dagger d\sigma\, G_q - \langle d\sigma, dG_q\rangle)|v|^2.$$

(2.8)

As explained next, this same identity holds when $q \in Z$. To see why, fix for the moment $\varepsilon > 0$ and let $\Theta_\varepsilon$ denote the characteristic function for the ball with center $q$ and radius $\varepsilon$. Let $G_{q,\varepsilon}$ denote the solution to the equation $d^\dagger d G_{q,\varepsilon} + G_{q,\varepsilon} = d\omega_d^{-1}\varepsilon^{-d}\Theta_\varepsilon$. This is a differentiable function whose $\varepsilon \to 0$ limit converges on compact subsets of $B_r - q$ to $G_q$. The function $G_{q,\varepsilon}$ is smooth on the complement of the boundary of the radius $\varepsilon$ ball centered at $q$ with bounded second derivatives near this boundary. Multiply both sides of (2.1) by $\chi_\rho \sigma G_{q,\varepsilon}$ and then integrate both sides of the result. An integration by parts then leads to an inequality that has the form

$$\tfrac{1}{2}d\omega_d^{-1}\varepsilon^{-d}\int_{B_\varepsilon}\sigma|v|^2 + \int_{B_r}\sigma G_{q,\varepsilon}|\nabla v|^2 = -\int_{B_r}\sigma G_{q,\varepsilon}\langle v,\mathcal{R}v\rangle - \int_{B_r}(\tfrac{1}{2}d^\dagger d\sigma\, G_{q,\varepsilon} - \langle d\sigma, dG_{q,\varepsilon}\rangle)|v|^2 + \mathfrak{r}_{q,(\rho,\varepsilon)}.$$

(2.9)

For fixed $\varepsilon$, the term $\mathfrak{r}_{q,(\rho,\varepsilon)}$ has limit zero as $\rho \to 0$. Take this limit first. Since $|v|$ is continuous, the $\varepsilon \to 0$ limit of the left most integral on the left side of (2.9) is equal to $\tfrac{1}{2}\sigma(q)|v|^2(q)$. The explicit integral on the right hand side of (2.9) has an $\varepsilon \to 0$ limit, this being the integral on the right hand side of (2.8). That this is so can be seen using either the inequality in the second the third bullet of (2.5) when $d > 2$ or some $k > 4$ version of the first bullet of (2.5) when $d = 2$. It can also be proved using the fact that $|v|$ is continuous. Granted the preceding, invoke the dominated convergence theorem to see that the $\varepsilon \to 0$ limit of the integral of $G_{q,\varepsilon}|\nabla v|^2$ exists and that limit is the integral over $X$ of $G_q|\nabla v|^2$.

The identity in (2.8) is used to prove the following lemma.

**Lemma 2.3**: *Let $U \subset X_*$ denote an open set and let $V \subset U$ denote an open subset with compact closure. There exists $\kappa_{V,U} > 1$ such that $\sup_V |v| \leq \kappa_{V,U}\int_U |v|^2$.*

***Proof of Lemma 2.3***: Given the bounds in (2.6), the claim follows from the version of (2.8) with $r = r_0$ and $q = p$ since (2.5) with Lemma 2.2 can be used to bound the integrals on the right hand side of (2.8) by a multiple of $\int_{B_r} |v|^2$.

### c) The derivatives of integrals of $|v|^2$

Supposing that $p \in X$, denote by h the function on $[0, r_0)$ given by the rule

$$r \to h(r) = \int_{\partial B_r} |v|^2 \ .$$

(2.10)

As proved momentarily, h is differentiable where $r > 0$ with derivative given by

$$\tfrac{d}{dr} h = \tfrac{d-1}{r} h + 2\int_{B_r} |\nabla v|^2 + 2\int_{B_r} \langle v, \mathcal{R} v \rangle + \int_{\partial B_r} \mathfrak{k} |v|^2 \ ,$$

(2.11)

with $\mathfrak{k}$ obtained from the trace of the second fundamental form of $\partial B_r$ by writing the latter as $\tfrac{d-1}{r} + \mathfrak{k}$. Of particular note is that $\mathfrak{k}$'s norm is no greater than $c_0 r$.

To prove the formula in (2.11), fix $\rho > 0$ for the moment and multiply both sides of (2.1) by $\chi_\rho$. Integrate the result over $B_r$ and then integrate by parts to see that

$$\tfrac{1}{2} \int_{\partial B_r} \tfrac{\partial}{\partial r}(\chi_\rho |v|^2) = \int_{B_r} |\nabla v|^2 + \int_{B_r} \langle v, \mathcal{R} v \rangle + \mathfrak{r}_\rho \ ,$$

(2.12)

with $\mathfrak{r}_\rho$ being a term whose norm is bounded by $c_0$ times the integral of $(|\nabla v|^2 + |v|^2)$ over the subset of $B_r$ where $\rho \leq |v| \leq 2\rho$. This understood, the bound $|\mathfrak{r}_\rho| \leq c_0 \rho^\delta$ with $\delta > c_0^{-1}$ follows from Lemma 2.1. The left hand side of (2.12) can be written as

$$\tfrac{1}{2} \tfrac{d}{dr} (r^{1-d} \int_{\partial B_r} \chi_\rho |v|^2 ) - r^{1-d} \int_{\partial B_r} \chi_\rho \mathfrak{k} |v|^2 \ .$$

(2.13)

The fact that h is differentiable and the formula in (2.11) for its derivative follow from (2.12) and (2.13) by first using difference quotients to define the derivative of h and then approximating any given difference quotient by very small $\rho$ versions of the integral over $\partial B_r$ of $\chi_\rho |v|^2$ is h. This leads to the desired outcome because the $\rho \to 0$ limit of $\mathfrak{r}_\rho$ is zero and because the $\rho \to 0$ limit of $\int_{\partial B_r} \chi_\rho |v|^2$ is h.

### d) The derivative of integrals of $|\nabla v|^2$

Fix $p \in X$ and let H denote the function on $[0, r_0)$ given by the rule

$$r \to H(r) = \int_{B_r} |\nabla v|^2 \ .$$

(2.14)

This function H is bounded and non-decreasing and a continuous and nondecreasing function is differentiable almost everywhere and its derivative is an $L^1$ function. See for example Theorem 3.23 in Rudin's book [R]. In particular, there exists a set of full measure in $(0, r_0]$ that is characterized as follows: If r is in this set, then H is differentiable at r, the function $|\nabla v|^2$ is integrable on $\partial B_r$ and the corresponding integral is the derivative of H at r. The next formula rewrites this integral:

$$\int_{\partial B_r} |\nabla v|^2 = 2 \int_{\partial B_r} |\nabla_r v|^2 + \tfrac{d-2}{r} \int_{B_r} |\nabla v|^2 - \int_{\partial B_r} \langle v, \mathcal{R}^\perp v \rangle + \tfrac{2}{r} \int_{B_r} \langle v, \mathcal{R} v \rangle$$

(2.15)

with $\nabla_r$ denoting the directional covariant derivative of $v$ along the outward pointing unit tangent vectors to the radial geodesics from the point p. Meanwhile, $\mathcal{R}^\perp$ denotes an endomorphism of $\mathbb{S}^+$ over $\partial B_r$ with norm bounded by $c_0$ and $\mathcal{R}$ denotes the endomorphism that appears in (2.1) and (2.2). The proof of (2.15) has six steps.

<u>Step 1</u>: Use $\nabla^\perp$ to denote $\nabla - dr \nabla_r$, this being the part of the covariant derivative that involves only derivatives tangent to the constant radius spheres centered at p. These definitions are such that if $s$ is a section of $\mathbb{S}^+ \otimes \mathcal{I}$ on the complement of Z in $B_r$–p, then

$$|\nabla s|^2 = |\nabla_r s|^2 + |\nabla^\perp s|^2 \ .$$

(2.16)

Let $\mathfrak{S}_r$ denote the homomorphism of $\mathbb{S}^+$ to $\mathbb{S}^-$ on $B_r$–p that is obtained by pairing the principal symbol of $\mathcal{D}$ with the outward pointing tangent vector to the radial geodesics from p and let $\mathfrak{S}^\perp$ denote $\mathfrak{S} - \mathfrak{S}_r \otimes dr$, this being the part of the principal symbol of $\mathcal{D}$ that annihilates this same vector field. Write $\mathcal{D}$ as

$$\mathcal{D} = \mathfrak{S}_r \nabla_r + \mathfrak{S}^\perp(\nabla^\perp) \ .$$

(2.17)

Suppose that $s$ is a section of $\mathbb{S}^+ \otimes \mathcal{I}$ on $\partial B_r - (\partial B_r \cap Z)$ with compact support. Integration by parts and a Bochner-Weitzenboch formula for $\mathfrak{S}^\perp(\nabla^\perp)$ writes

$$\int_{\partial B_r} |\mathfrak{S}(\nabla^\perp)s|^2 = \int_{\partial B_r} |\nabla^\perp s|^2 + \tfrac{d-2}{r} \int_{\partial B_r} \langle s, \mathfrak{S}_r^\dagger \mathfrak{S}(\nabla^\perp)s\rangle + \int_{\partial B_r} \langle s, \mathcal{R}_0^\perp s\rangle$$

(2.18)

with $\mathcal{R}_0^\perp$ being an endomorphism with norm bounded by $c_0$. By way of a parenthetical remark, this endomorphism $\mathcal{R}_0^\perp$ vanishes when the metric on $B_r$ is flat. A point to note when deriving (2.18) is that the volume form on $B_r$ is not $\nabla^\perp$ covariantly constant. This fact accounts for a term on the right hand side in (2.18) that can be written as

$$\tfrac{d-1}{r} \int_{\partial B_r} \langle s, \mathfrak{S}_r^\dagger \mathfrak{S}(\nabla^\perp)s\rangle.$$

(2.19)

The integral of $\langle s, \mathfrak{S}_r^\dagger \mathfrak{S}(\nabla^\perp)s\rangle$ that appears (2.19) also appears in (2.18), but it appears in (2.18) with the multiplicative factor $\tfrac{d-2}{r}$. This is observedly not the factor that multiplies the integral (2.19). The extra multiplicative factor of $-\tfrac{1}{r}$ in (2.18) accounts for the fact that the operator $\nabla^\perp$ does not commute with $\nabla_r$.

<u>Step 2</u>: Step 3 explains why the identity in (2.18) holds for $v$. Granted that this is so, then the identities in (2.16) and (2.18) can be used to write

$$\int_{\partial B_r} |\nabla v|^2 = \int_{\partial B_r} |\nabla_r v|^2 + \int_{\partial B_r} |\mathcal{D}^\perp v|^2 - \tfrac{d-2}{r} \int_{\partial B_r} \langle v, \mathfrak{S}_r^\dagger \mathcal{D}^\perp v\rangle - \int_{\partial B_r} \langle v, \mathcal{R}^\perp v\rangle.$$

(2.20)

Meanwhile, the equation $\mathcal{D}v = 0$ can be written using (2.16) as $\nabla_r v + \mathfrak{S}_r^\dagger \mathcal{D}^\perp v = 0$ because $\mathfrak{S}_r^\dagger \mathfrak{S}_r = 1$. Use this fact to rewrite (2.19) as

$$\int_{\partial B_r} |\nabla v|^2 = 2\int_{\partial B_r} |\nabla_r v|^2 + \tfrac{d-2}{r} \int_{\partial B_r} \langle v, \nabla_r v\rangle - \int_{\partial B_r} \langle v, \mathcal{R}^\perp v\rangle.$$

(2.21)

The formula in (2.15) follows from (2.21) by rewriting the middle integral on the right hand side of (2.21) using an integration by parts with (2.1).

<u>Step 3</u>: The next lemma is used to justify the assertion that the identity in (2.18) holds with $s$ taken equal to $v$. To set the stage, suppose for the moment that M is a given compact, Riemannian manifold. Let $Z \subset M$ denote a given closed subset and let $\mathbb{V}$ denote a vector bundle over M−Z with fiber metric and metric compatible connection. The notation uses $\nabla$ to denote the corresponding covariant derivative.

**Lemma 2.4**: *Suppose that* M, Z *and* $\mathbb{V}$ *are as just described. Let q denote a smooth section of* $\mathbb{L}$ *whose norm extends over* Z *to define a continuous function on* M *with* Z *contained in its zero locus. Assume in addition that* $|\nabla q|$ *defines an* $L^2$ *function on* M–Z. *Given* $\varepsilon > 0$, *there exists a smooth section of* $\mathbb{V}$ *with compact support in* M–Z, *this denoted by* $q_\varepsilon$, *with the integrals of both* $|q - q_\varepsilon|^2$ *and* $|\nabla(q - q_\varepsilon)|^2$ *being less than* $\varepsilon$.

This lemma is proved momentarily

To prove that (2.18) holds with $s = v$, take M in the lemma to be $\partial B_r$, take $\mathbb{V}$ to be the bundle $\mathbb{S}^+ \otimes \mathcal{I}$ and take $q$ to be $v$. Given $\varepsilon > 0$, let $v_\varepsilon$ denote a version of Lemma 2.4's section $q_\varepsilon$ that is obtained using $M = \partial B_r$ and $\mathbb{V} = \mathbb{S}^+ \otimes \mathcal{I}$ and $q = v$. The formula in (2.18) holds with $s = q_\varepsilon$. This understood, then (2.18) must hold for $s = v$ because the $s = v$ versions of both the right and left hand sides of (2.18) differ from the $s = v_\varepsilon$ versions by at most $c_0$ times the integral over $\partial B_r$ of the sum of $|v - v_\varepsilon|^2$ and $|\nabla(v - v_\varepsilon)|^2$.

Step 4: This step contains the

***Proof of Lemma 2.4***: Given $\rho \in (0, c_0^{-1})$, introduce $\tau_\rho$ to denote the function on M given by $\chi(1 - \ln|q|/\ln(\rho))$. Let $q_\rho = \tau_\rho q$. Then

$$\int_{M-Z} |\nabla q_\rho|^2 = \int_{M-Z} \tau_\rho^2 |\nabla q|^2 + \mathfrak{z}_\rho$$

(2.22)

with $\mathfrak{z}_\rho$ being a term with norm at most $c_0 |\ln \rho|^{-1} \int_M |\nabla q|^2$. It follows as a consequence that

$$\lim_{\rho \to 0} \int_{M-Z} |\nabla q_\rho|^2 \leq \int_{M-Z} |\nabla q|^2 .$$

(2.23)

Meanwhile, $\int_{M-Z} |q - q_\rho|^2 \leq c_0 \rho$ so the sequence $\{q_\rho\}_{\rho > 0}$ converges to $q$ in the $L^2$ topology. It follows from this and (2.23) that the equal sign holds in (2.23). This implies in particular that

$$\lim_{\rho \to 0} \int_{\{p \in M : |q(p)| < \rho\}} |\nabla q|^2 = 0 ,$$

(2.24)

and the latter fact implies in turn that

$$\lim_{\rho \to 0} \int_{M-Z} |\nabla(q - q_\rho)|^2 = 0.$$

(2.25)

This last bound implies that $\{q_\rho\}_{\rho>0}$ converges to $q$ in the $L^2_1$ topology on the space of sections of $\mathbb{V}$ on M−Z. Meanwhile, each $q_\rho$ has compact support where $|q| > 0$ and so each $q_\rho$ can be approximated in the $L^2_1$ topology to any desired accuracy by a smooth function with compact support where $|q| > 0$.

### 3. The frequency function

The primary tool that is used subsequently to study Z is an analog of Almgren's frequence function [Al2] that is defined in Section 3a. Some properties of this function are stated in the upcoming Lemmas 3.2-3.4. The final subsection uses the results of Lemmas 3.2 and 3.3 to prove the assertion in Theorem 1.5 to the effect that $|v|$ is Hölder continuous.

**a) The definition of the frequency function**

The definition of the frequency function requires the choice of $p \in X$. With p chosen, two auxilliary functions on $[0, r_0]$ are needed. The first is the function h that is defined in (2.10). The lemma that follows states two relevant facts about this function.

**Lemma 3.1**: *There exists* $\kappa > 1$ *with the following significance: Fix* $p \in X$ *to define the function* h *on* $[0, r_0]$. *This function* h *is strictly positive on* $(0, r_0]$. *Moreover, if* $s \in (0, r_0]$ *and if* $r \in [s, \kappa^{-1}]$, *then* $h(r) \geq (\frac{r}{s})^{d-1} e^{-\kappa(r^2 - s^2)} h(s)$.

This lemma is proved in Section 3b. Take the lemma on faith until then.
Taken together, the various $p \in X$ versions of Lemma 3.1 imply the following:

*The set Z has empty interior .*

(3.1)

The second required function is denoted by $\mathfrak{d}_\diamond$, it is defined below in (3.2). The function M that appears in (3.2) is the same as its namesake in (3.4); it is defined by writing the trace of the second fundamental form of any given $r \in (0, r_0]$ version of $\partial B_r$ as $\frac{d-1}{r} + \mathfrak{k}$. The rule that defines $\mathfrak{d}_\diamond$ is

$$r \to \mathfrak{d}_\diamond(r) = \int_0^r (\frac{1}{h(s)}(\int_{B_s} \langle v, \mathcal{R}v \rangle + \frac{1}{2}\int_{\partial B_s} \mathfrak{k}|v|^2)) ds .$$

(3.2)

Note that $|\mathfrak{d}_\diamond(r)| \leq c_0 r^2$, this being another corollary of Lemma 3.1.

With $h_\diamond$ and $\mathfrak{d}_\diamond$ in hand, define the function $\kappa_\diamond$ on $[0, c_0^{-1}]$ to be the positive square root of $e^{-2\mathfrak{d}_\diamond} r^{1-d} h_\diamond$. It follows from Lemma 3.1 and what was just said about $\mathfrak{d}_\diamond$ that $\kappa_\diamond$ is continous and bounded by $c_0$. It is also the case that $\kappa(0) = \omega_d^{1/2} |v|(p)$ with $\omega_d$ denoting here the area of the radius 1 sphere in $\mathbb{R}^d$. This formula for $\kappa(0)$ holds because $|v|$ is a continuous function.

The desired analog of Almgren's frequency function is denoted by $N$. It is the function on $(0, c_0^{-1}]$ that is defined by the rule

$$r \to N(r) = \frac{1}{r^{d-2} \kappa(r)^2} \int_{B_r - (B_r \cap Z)} |\nabla v|^2 \; .$$

(3.3)

The next lemma states some salient properties of $N$.

**Lemma 3.2**: *There exists $\kappa > 1$ with the following significance: Fix $p \in X$ so as to define the functions $\kappa$ and $N$. The function $N$ is continuous on $(0, r_0]$ and it extends as a continuous function to $[0, r_0]$. Moreover,*

- $\frac{d}{dr} \kappa = \frac{1}{r} N \kappa$ *on* $(0, r_0]$.
- *The function $N$ is almost everywhere differentiable on $(0, r_0]$ and its derivative obeys*

$$\frac{d}{dr} N \geq \frac{2}{r^{d-2} \kappa(r)^2} \int_{\partial B_r} |\nabla_r v - \frac{N}{r} v|^2 \; - \kappa r \; .$$

- *If $s \in (0, r_0]$ and $r \in [0, s]$, then $N(s) \geq N(r) - \kappa(s^2 - r^2)$.*
- *If $|v|(p) > 0$, then $\lim_{r \to 0} N(r) = 0$.*
- *If $|v|(p) = 0$, then $\lim_{r \to 0} N(r) > \kappa^{-1}$.*

Note that identity in Lemma 3.2's first bullet can be integrated to obtain the following assertion: If $s \in (0, r_0]$ and $r \in (0, s]$, then

$$\kappa(s) = \exp\left( \int_r^s \frac{N(t)}{t} dt \right) \kappa(r) \; .$$

(3.4)

Lemma 3.2 is proved in Section 3b also.

The next lemma states some applications of Lemma 3.2, one being an a priori Hölder bound for $|v|$ along $Z$. Keep in mind that Lemma 3.2 asserts in part that $N$ extends as a continuous function to the closed interval $[0, r_0]$.

**Lemma 3.3**: *There exists $\kappa > 1$ with the following significance: Fix a point $p \in Z$ and fix $s \in (0, r_0)$. Then*

- *There exists $z_s \in (0,\kappa)$ such that $\kappa(s) = z_s s^{N(0)}$. Moreover, if $z_s$ is defined by writing $\kappa(s) = z_s s^{N(0)}$, then $\kappa(r) \leq (1+\kappa) z_s r^{N(0)}$ when $r \in (0,s]$.*
- *If $q \in X$ is such that $\mathrm{dist}(p,q) \leq \kappa^{-1} s$, then $|v|(q) \leq \kappa z_s (\mathrm{dist}(p,q))^{N(0)}$.*

This lemma is proved in Section 3c.

Each point in X has its version of the function N. The next lemma talks about the dependence of N on the chosen point. By way of notation, the upcoming lemma uses $N_{(p)}$ to denote the version of N that is defined by a given point $p \in X$. This notation is used subsequently when N's dependence on the point in X plays a central role.

**Lemma 3.4**: *Fix $p \in X$, and let $\{p_k\}_{k=1,2,\ldots}$ denote a sequence in X that converges to p.*
- $\limsup_{k \to \infty} N_{(p_k)}(0) \leq N_{(p)}(0)$.
- *If $\lim_{k \to \infty} N_{(p_k)}(0) = N_{(p)}(0)$, then the following is true: Given $\varepsilon \in (0,1]$ there exists $m_\varepsilon > 1$ such that if $k > m_\varepsilon$ and $r \in [0, m_\varepsilon^{-1}]$, then $|N_{(p_k)}(r) - N_{(p)}(0)| < \varepsilon$.*

Lemma 3.4 is proved in Section 3d.

**b) Proof of Lemmas 3.1 and 3.2**

The proof of these lemmas are intertwined because the formula in the second bullet of Lemma 3.2 is needed to prove the assertion in Lemma 3.1 that $h > 0$ on $(0, r_0]$. The proof of the two lemmas has four parts.

*Part 1*: Fix $p \in X$ to define h and suppose that $s \in (0, r_0]$ is such that $h(s) > 0$. This part proves the assertion in Lemma 3.1 to the effect that $h(r) \geq (\frac{r}{s})^{d-1} e^{-c_0(r^2 - s^2)} h(s)$ when $r \in [s, r_0]$. To start, use the inequality in the fifth bullet of (2.5) to see that

$$\left| \int_{B_r} \langle v, \mathcal{R}v \rangle \right| \leq c_0 r^2 \int_{B_r} |\nabla v|^2 + c_0 r \int_{\partial B_r} |v|^2.$$

(3.5)

Now use the latter inequality to see why (2.11) implies a differential inequality that reads

$$\frac{d}{dr} h \geq \frac{d-1}{r} h - c_0 r h.$$

(3.6)

Integrate (3.6) to obtain the asserted $(\frac{r}{s})^{d-1} e^{-c_0(r^2 - s^2)} h(s)$ lower bound for $h(r)$.

*Part 2*: Suppose that $r_1 \in [0, r_0]$ is such $h > 0$ on $(r_1, r_0]$. If h is positive on $(r_1, r_0]$, then $\kappa$ is also and so the definition in (3.3) for N makes sense at points in $(r_1, r_0]$. The formula in the first bullet of Lemma 3.2 for the derivative of K on $(r_1, r_0]$ follows directly from (3.2) and (2.11).

The fact that N is continuous and almost everywhere differentiable on $(r_1, r_0]$ follows from the fact that K is differentiable where it is non-zero and the fact that the function H in (2.14) is continuous and almost everywhere differentiable. The inequality in the second bullet of Lemma 3.2 for the derivative of N on $(r_1, r_0]$ from the formula in (2.15) and the formula in Lemma 3.2's first bullet for the derivative of K on $(r_1, r_0]$.

If $r \in (r_1, r_0]$ and if $s \in (r_1, r)$, then the inequality in the third bullet of Lemma 3.2 follows by integrating the inequality in this same lemma's second bullet.

*Part 3*: This part of the proof explains why the function h is positive on $(0, r_0]$. To start the explanation, note that if $h(r) = 0$ for some $r \in (0, r_0]$, then h must vanish on the whole of $[0, r]$ because $h(r) \geq (\frac{r}{s})^3 e^{-c_0(r^2 - s^2)} h(s)$. With the preceding understood, suppose for the sake of argument that there is a point in X whose version of h is zero on an interval in $(0, r_1]$. The three steps that follow find nonsense from this assumption.

Step 1: The function $|\nu|$ is not identically zero by assumption; and this implies in turn that there exists $p \in X$ and numbers $r_\Delta \in (0, r_0)$ and $r_\diamond \in [0, r_\Delta)$ such that p's version of h is zero on $[0, r_\diamond]$ and non-zero on $(r_\diamond, r_\Delta]$. Let $r_\ddagger = \frac{1}{2}(r_\diamond + r_\Delta)$ and introduce by way of notation $\delta$ to denote $K(r_\ddagger)$. Keep in mind that $\delta > 0$.

Step 2: Fix $\varepsilon \in (0, 1]$ for the moment. Since the function K is continuous, there exists $r_\varepsilon \in (r_\diamond, r_\ddagger)$ where $K(r_\varepsilon) = \varepsilon \delta$. The formula in (3.4) implies that

$$\int_{r_\varepsilon}^{r_\ddagger} \frac{N(s)}{s} ds = |\ln \varepsilon| \ .$$

(3.7)

This can happen only if there exists $s \in [r_\varepsilon, r_\ddagger]$ where $N(s) \geq |\ln \varepsilon|(\ln(r_\ddagger/r_\diamond))^{-1}$.

Step 3: It follows from the third bullet of Lemma 3.2 that the function N on $[r_\ddagger, r_\Delta]$ is no smaller than $|\ln \varepsilon|(\ln(r_\ddagger/r_\diamond))^{-1} - c_0 r_\Delta^2$. This is nonsense because $\varepsilon$ can be chosen as small as desired and so N must be infinite on the whole of $[r_\ddagger, r_\Delta]$.

*Part 4*: Since h is positive on $(0, r_0]$, so is the function K and so the assertions of the first three bullets of Lemma 3.2 follow from what is said in Part 2. To see about the continuous extension of N to the closed interval $[0, c_0^{-1}]$, note first that $\limsup_{r \to 0} N(r)$ is bounded because if not, then the third bullet of Lemma 3.2 could be invoked to obtain the nonsensical assertion that N is infinite on the whole of $(0, r_0]$.

The assertion of the lemma's third bullet also implies that $\limsup_{r \to 0} N(r)$ is equal to $\liminf_{r \to 0} N(r)$. To see why this is, suppose that these two limits are different so as to generate nonsense. If these limits differ, then there exist numbers $\varepsilon, n > 0$ and sequences $\{r_k\}_{k=1,2,...}$ and $\{s_k\}_{k=1,2,...}$ with the following properties: Both of these sequence have limit zero, the sequence $\{r_k\}_{k \in 1,2,...}$ obeys $N(r_k) > n + \varepsilon$ for all k and the sequence $\{s_k\}_{k \in 1,2,...}$ obeys $N(s_k) < n - \varepsilon$ for all k. This is nonsense because the third bullet of Lemma 3.2 asserts that $N(r) > n + \varepsilon - c_0 r^2$ for r greater than any given $r_k$; and this implies in particular that $N(s_k) > n$ if k is sufficiently large.

The assertion in the fourth bullet of Lemma 3.2 follows from the formula in Lemma 3.2's first bullet for the derivative of K. The assertion in the fifth bullet follows from (1.5). The reason being that if $\lim_{r \to 0} N$ is less than the number $\varepsilon$ in (1.5), then the bound in (1.5) would not be compatible with the formula in (3.4).

### c) Proof of Lemma 3.3

To prove the assertion made by the top bullet of the lemma, first invoke the third bullets of Lemma 3.2 to see that p's version of the function N obeys $N(s) \geq N(0) - c_0 s^2$ when $s \in [0, r_0)$. Use this bound with version of (3.4) that has s replaced by $r_0$ and r replaced by s to see that

$$K(s) \leq s^{N(0)} c_0 \quad \text{when } s \in [0, r_0),$$

(3.8)

This says that $K(s)$ can be written as $z_s s^{N(0)}$ with $z_s \leq c_0$. The bound $N(\cdot) \geq N(0) - c_0 s^2$ on the interval $[r, s]$ with (3.4) lead directly to the bound $K(r) \leq (1 + c_0) z_s r^{N(0)}$.

To prove the lemma's second bullet, assume that q is in the radius $r_0$ centered at p. Since the second bullet's assertion holds if $|v|(q) = 0$, no generality is lost by assuming henceforth that q is not in the set Z. With D denoting the distance between p and q, let $G_q$ denote the Dirichlet Green's function with pole at q for the operator $d^\dagger d$ on the radius $r = 10D$ ball centered at p. Let $\sigma$ denote the function $\chi(D^{-1} \text{dist}(\cdot, p) - 4)$, this being equal to 1 where the distance to p is less than 4D and equal to 0 where the distance to p is greater than 5D. Use this version of $\sigma$ in (2.18). Use these bounds in (2.6) with the fact that the function $K_{(p)}$ is increasing to bound the right hand side of (3.8) by $c_0 (K_{(p)}(5D))^2$. The top bullet of Lemma 3.3 asserts that this in turn is less than $c_0 z_s^2 D^{2N(0)}$ if $D \leq \frac{1}{5} s$.

### d) Proof of Lemma 3.4

The two bullets of the lemma are proved in Parts 1 and 2 that follow.

*Part 1*: This part proves the first bullet of the lemma. If p is not in Z, then the assertion made by the first bullet follows from the fourth bullet of Lemma 3.2. This is also the case if p is in Z but the sequence $\{p_k\}_{k\in\{1,2,...\}}$ is not in Z. Granted the preceding, suppose that p ∈ Z and q ∈ Z has distance $c_0^{-1}$ from p. Let D denote the distance between p and q. Fix r ∈ [4D, $r_0$) and s ∈ [r, $r_0$). Introduce now $K_{(q)}$ to denote q's version of the function K and $K_{(p)}$ to denote p's version. It follows from (3.4) and the third bullet of q's version of Lemma 3.2 that

$$K_{(q)}(r)/K_{(q)}(s) \leq (1+c_0 r^2)(\tfrac{r}{s})^{1+N_{(q)}(0)}.$$

(3.9)

Meanwhile, the fact that $N_{(p)}$ is continuous on [0, $r_0$] has the following implication: Given ε > 0, there exists $r_\varepsilon$ such that $N_{(p)}(r) \leq N_{(p)}(0) + \varepsilon$ if r ∈ [0, $r_\varepsilon$). It follows as a consequence that if D ≤ $\tfrac{1}{100} r_\varepsilon$ and if both r and s are less than $r_\varepsilon$, then

$$K_{(p)}(r)/K_{(p)}(s) \geq (\tfrac{r}{s})^{1+N_{(p)}(0)+c_0\varepsilon}.$$

(3.10)

For any given fixed r ∈ [0, $r_0$) the function on X given by the rule q → $K_{(q)}(r)$ is continous because |v| is continous. It follows that if ε ∈ (0, 1) and r, s ∈ (0, $r_0$] are given, and if q is sufficiently close to p then

$$K_{(p)}(r)/K_{(p)}(s) \leq (1+\varepsilon) K_{(q)}(r)/K_{(q)}(s)$$

(3.11)

This last bound with the bounds in (3.9) and (3.10) have the following as a consequence: If ε ∈ (0, 1) is given, and if r ∈ (0, $r_\varepsilon$] and s ∈ (r, $r_\varepsilon$] are given, and if q is sufficiently close to p, then

$$(\tfrac{r}{s})^{1+N_{(p)}(0)+c_0\varepsilon} \leq (\tfrac{r}{s})^{1+N_{(q)}(0)}.$$

(3.12)

The s = 2r version of this inequality implies that $N_{(q)}(0) \leq N_{(p)}(0) + c_0 \varepsilon$ if q is sufficiently close to p. Since ε can be any positive number, this last observation proves Lemma 3.4.

*Part 2*: To prove the lemma's second bullet, fix k ∈ {1, 2, …} and then invoke the $p_k$ version of Lemma 3.2's third bullet to see that $N_{(p_k)}(r) \geq N_{(p_k)}(0) - c_0 r^2$. Since $|N_{(p_k)}(0) - N_{(p)}(0)| < \tfrac{1}{100}\varepsilon$ if k is sufficiently large, it therefore suffices find z > 1 so that

$$N_{(p_k)}(r) \leq N_{(p)}(0) + \tfrac{1}{100}\varepsilon \;\; if\, k > z\; and\; r \in (0, z^{-1}].$$

(3.13)

The proof (3.13) uses the function F on [0, $r_0$) that is defined the rule whereby F(r) is equal to the integral of $|v|^2$ on the radius r ball centered at p. Fix m > 1 for the moment

and then fix $\rho \in (0, r_0)$ so that $|N_{(p)}(r) - N_{(p)}(0)| < m^{-1}$ on the radius $10\rho$ ball centered at p. The point p versions of (3.4) and Lemma 3.2 imply that

- $F(\rho) \geq (1 - c_0 m^{-1} - c_0 \rho^2) \frac{1}{d+2N(\rho/m)} \rho^d \, m^{2N(\rho/m)} K_{(p)}^2(\rho/m)$ .
- $F(2\rho) \leq F(\rho) + (1 + c_0 \rho^2) \frac{1}{d+2N(\rho/m)} \rho^d \, m^{2N(2\rho)} K_{(p)}^2(\rho/m) (2^{d+2N_{(p)}(2\rho)} - 1)$ .

(3.14)

These inequalities imply in turn that

$$F(2\rho)/F(\rho) \leq 1 + (1 + c_0 m^{-1} + c_0 \rho^2) m^{4/m} (2^{d+2N_{(p)}(2\rho)} - 1)$$

(3.15)

because the difference between $N_{(p)}(2\rho)$ and $N_{(p)}(\rho/m)$ is at most $2m^{-1}$.

With (3.15) in hand, fix $k \in \{1, 2, \ldots\}$ again and define the function $F_k$ on $[0, r_0)$ by setting its value at r to be the integral of $|v|^2$ on the radius r ball centered at $p_k$. Let $D_k$ denote the distance between $p_k$ and p. It follows from the $p_k$ versions of (3.4) and Lemma 3.2 that

$$F_k(2\rho - D_k)/F_i(\rho + D_k) \geq (1 - c_0 \rho)^2 (2^{d+2N_{(p_k)}(\rho)} - 1) .$$

(3.16)

Suppose that k is such that $D_k < c_0^{-1} \rho$ which will be true if k is sufficiently large. Since the ball of radius $2\rho - D_k$ centered at $p_k$ is contained in the ball of radius $2\rho$ centered at p, and the ball of radius $\rho$ centered at p is contained in the ball of radius $\rho + D_k$ centered at $p_i$, the ratio on the left hand side of (3.15) can not be less than the ratio on the left hand side of (3.16). This understood, the inequality in (3.15) is not compatible with that in (3.16) unless

$$N_{(p_k)}(\rho) \leq N_{(p)}(2\rho) + c_0 (m^{-1} \ln m + \rho^2) .$$

(3.17)

The inequality in (3.13) follows if $m$ is equal to $c_0 \varepsilon^{-2}$ and then $\rho$ is chosen to be less than $c_0^{-1} \varepsilon$ and such that $|N_{(p)}(r) - N_{(p)}(0)| \leq c_0^{-1} \varepsilon^2$ on the radius $10\rho$ ball centered at p.

**e) Hölder continuity of |v|**

The two parts that follow prove the assertion in Theorem 1.5 to the effect that $|v|$ is Hölder continuous.

*Part 1*: Differentiating the equation in (2.2) leads to an equation that has the form

$$\nabla^{\dagger}\nabla(\nabla v) + \mathcal{R}_1 \nabla v + \mathcal{R}_0 v = 0$$

(3.18)

where $\mathcal{R}_1$ and $\mathcal{R}_0$ having bounded norms. Given $\delta \in (0,1]$ and let $\beta_\delta$ denote the function $\chi(2 - \delta^{-1}\mathrm{dist}(\cdot,Z))$. This function is equal to 1 where the distance to Z is greater than $2\delta$ and it is equal to zero where the distance to Z is less than $\delta$. Fix a ball $p \in X$ and use $\gamma$ to denote the function $\chi(2r_0^{-1}\mathrm{dist}(\cdot,p) - 1)$. This function is equal to 1 where the distance to p is less than $\frac{1}{2}r_0$ and equal to 0 where the distance to p is greater than $r_0$. The equation in (3.18) leads to the $L^2$ norm bound

$$\int_X \beta_\delta^2 \gamma |(\nabla^{\dagger}\nabla)\nabla v|^2 \leq c_0.$$

(3.19)

An integration by parts then leads from (3.19) to an $L^2_3$ bound for $v$ of the form

$$\int_X \gamma |\nabla^{\otimes 3}(\beta_\delta v)|^2 \leq c_0 \delta^{-6}$$

(3.20)

This bound leads via a standard Sobolev inequality to a $c_0 \delta^{-3}$ bound for the $L^8$ norm $\beta_\delta |\nabla v|$ in the radius $\frac{1}{2}r_0$ ball centered at p. The latter bound leads in turn a $c_0 \delta^{-3}$ bound for the $L^8$ norm of $\beta_\delta d|v|$. This last $c_0 \delta^{-3}$ bound implies via another Sobolev inequality an exponent $\frac{1}{4}$ Hölder norm bound

$$||v(p)| - |v(q)|| \leq c_0 \delta^{-3} \mathrm{dist}(p,q)^{1/4}$$

(3.21)

if both p and q have distance $2\delta$ or more from Z.

*Part 2*: Meanwhile, it follows from the second bullet of Lemma 3.3 and the fifth bullet of Lemma 3.2 that if $\delta \in (0, \frac{1}{100}r_0)$ and if p has distance less than $\delta$ from Z, then

$$|v(p)| \leq c_0 \delta^{\upsilon}$$

(3.22)

with $\upsilon$ positive and independent of both $\delta$ and p. Of course, the analogous inequality holds for $|v(q)|$ if q has distance less than $\delta$ from Z.

Label p and q so that $\mathrm{dist}(p,Z) \leq \mathrm{dist}(q,Z)$. If $\mathrm{dist}(p,q)^{116} \leq \mathrm{dist}(p,Z)$, then it follows from (3.21), that

$$||v(p)| - |v(q)|| \leq c_0 |\mathrm{dist}(p,q)|^{1/16}.$$

(3.23)

Suppose on the other hand that $\text{dist}(p,q)^{1/16} \geq \text{dist}(p,Z)$. Assume in addition that $\text{dist}(q,Z) > 10\,\text{dist}(p,Z)$. This assumption implies that $\text{dist}(q,Z) \leq \frac{10}{9}\,\text{dist}(p,q)$; and since

$$||v(p)|-|v(q)|| \leq |v(p)|+|v(q)|, \tag{3.24}$$

it follows from (3.22) and its $|v(q)|$ analog that

$$||v(p)|-|v(q)|| \leq c_0\,\text{dist}(p,q)^{\upsilon}. \tag{3.25}$$

The last case to consider has $\text{dist}(p,q)^{1/16} \geq \text{dist}(p,Z)$ and $\text{dist}(q,Z) \leq 10\,\text{dist}(p,Z)$. Granted these inequalities, it follows from (3.22), its $|v|(q)$ analog and from (3.24) that

$$||v(p)|-|v(q)|| \leq c_0\,\text{dist}(p,q)^{\upsilon/16}. \tag{3.26}$$

The claimed Hölder continuity of $|v|$ follows from (3.23), (3.25) and (3.26).

## 4. The rescaling limits of $v$

This section looks at the behavior of $|v|$ very near points in Z. What is said here is summarized in Proposition 4.1. Proposition 4.1 looks a lot like Proposition 8.1 in the published version of [T1], and the subsequent proof of Proposition 4.1 borrows heavily from the Sections 8 and 9 in the published version of [T1].

To set the notation for Proposition 4.1, fix $R \in (0, \infty)$ for the moment and define the $\psi_R \colon \mathbb{R}^d \to \mathbb{R}^d$ by the rule $x \to \psi_R(x) = Rx$. A given closed set in $\mathbb{R}^d$ is deemed to be scale invariant when it is mapped to itself by each $R > 0$ version of $\psi_R$. Such a set is either the origin or the cone on a set in the $|x| = 1$ sphere. With regards to scale invariant sets, keep in mind that a real line bundle that is defined on the complement of a scale invariant set is canonically isomorphic to its pull-back by any $R > 0$ version of $\psi_R$.

Keep in mind that the radius $100 r_0$ ball about any given point in X is well inside the domain of a Gaussian coordinate chart centered on the point in question. Suppose in what follows that $\mathrm{T} = (p, \mathrm{L}, \lambda)$ is a data set consisting of a point $p \in \mathrm{X}$, an oriented orthonormal frame at $\mathrm{TX}|_p$, this being $\mathrm{L}$, and a number $\lambda \in (0, r_0]$. The frame $\mathrm{L}$ determines a Gaussian coordinate chart centered at p whose differential at the origin at $\mathbb{R}^4$ sends a fixed fiducial orthonormal frame to $\mathrm{L}$. With the preceding understood, define $\phi_\mathrm{T}$ to be the map from the radius $\lambda^{-1} r_0$ ball about the origin in $\mathbb{R}^d$ to X that composes first $\psi_\lambda$ and then the Gaussian coordinate chart map that is defined by $\mathrm{L}$.

Let $\mathbb{V}^+$ and $\mathbb{V}^-$ denote the respective fibers over the point p of the bundles $\mathbb{S}^+$ and $\mathbb{S}^-$. Parallel transport along the rays from the origin in $\mathbb{R}^d$ is used to identify $\phi_T^*\mathbb{S}^+$ and $\phi_T^*\mathbb{S}^-$ over the $|x| < \lambda^{-1}r_0$ ball with the respective product $\mathbb{V}^+$ and $\mathbb{V}^-$ bundles over this ball. This is done using the $\phi_T$ pull-back of the connections on $\mathbb{S}^+$ and $\mathbb{S}^-$ that define the covariant derivative $\nabla$. The principal symbol of the operator $\mathcal{D}$ at p defines an elliptic operator from $C^\infty(\mathbb{R}^d;\mathbb{V}^+)$ to $C^\infty(\mathbb{R}^d;\mathbb{V}^-)$, this denoted by $\mathcal{D}_0$. A section of the product $\mathbb{V}^+$ bundle over a given open set in $\mathbb{R}^d$ is said to be $\mathcal{D}_0$-*harmonic* when it is annihilated by $\mathcal{D}_0$. The directional covariant derivative on $C^\infty(\mathbb{R}^d-0;\mathbb{V}^+)$ along the outward pointing, unit tangent vector to the rays in $\mathbb{R}^d$ from the origin is denoted in what follows by $\nabla_r$.

With regards to a subsequent convention, the vector spaces $\mathbb{V}^+$ and $\mathbb{V}^-$ will most often not be distinguished by notation or otherwise from the product vector bundles with these spaces as fiber.

Proposition 4.1 has $Z_T$ denoting the $\phi_T$-inverse image of Z and it has $\mathcal{I}_T$ denoting the pull-back via $\phi_T$ of Proposition 7.1's line bundle $\mathcal{I}$. Note that these are defined on the radius $\lambda^{-1}r_0$ ball centered on the origin in $\mathbb{R}^4$. Use $\kappa_{(p)}$ to denote p's version of the function $\kappa$. Proposition 4.1 has $\nu_T$ denoting $\kappa_{(p)}(\lambda)^{-1}\lambda^{-1}\phi_T^*\nu$. This is an $\mathcal{I}_T$ valued 1-form that on the complement of $Z_T$ in $|x| < \lambda^{-1}r_0$ ball centered on the origin in $\mathbb{R}^4$.

An additional convention is used in the proposition, this being the use of $N_{(p)}$ to denote p's version of the function $N$. The $r \to 0$ limit of this function is denoted by $N_{(p)}(0)$. By way of a reminder, Lemma 3.2 asserts that $N_{(p)}$ extends as a continuous function to $[0,r_0]$.

**Proposition 4.1**: *There exists $\kappa > 1$ with the following significance: Let p denote a given point in Z and suppose that $\{T(i)=(p_i,L_i,\lambda_i)\}_{i=1,2,\ldots}$ is a sequence with each element consisting of a point in X, an oriented orthonormal frame at the point and a number in $(0, r_0]$. Assume that $\lim_{i\to\infty} p_i = p$ and that $\lim_{i\to\infty} N_{(p_i)}(0) = N_{(p)}(0)$ and that $\lim_{i\to\infty} \lambda_i = 0$. There is a subsequence in $\{1, 2, \ldots\}$ to be denoted by $\Lambda$ and a data set $(Z_*,\mathcal{I}_*,\nu_*)$ with the properties listed below.*
- *$Z_*$ is a scale invariant set in $\mathbb{R}^d$.*
- *$\mathcal{I}_*$ is a real line bundle defined on the complement of $Z_*$.*
- *$\nu_*$ is an $\mathcal{I}_*$ valued, $\mathcal{D}_0$-harmonic section of $\mathbb{V}^+$ over $\mathbb{R}^d-Z_*$ whose norm extends over $Z_*$ as an $L^2_{1;loc}$ and exponent $\upsilon = \kappa^{-1}$ Holder continuous function with zero locus $Z_*$.*
- $\nabla_r\nu_* - \frac{N_{(p)}(0)}{|x|}\nu_* = 0.$
- *The sequence of data sets $\{(Z_{T(i)},\mathcal{I}_{T(i)}, \lambda_{T(i)}\}_{i\in\Lambda}$ converges to $(Z_*,\mathcal{I}_*,\nu_*)$ in the following sense: Given $\varepsilon \in (0,1)$, there exists $\iota_\varepsilon \geq 1$ such that if $i \in \Lambda$ and $i > \iota_\varepsilon$, then*
  - a) *The functions $|\nu_{T(i)}|$ and $|\nu_*|$ differ where $|x| < \varepsilon^{-1}$ by a function with exponent*

$\upsilon = \kappa^{-1}$ *Holder norm less than $\varepsilon$*
b) *Let $\mathcal{T}_\varepsilon \subset \mathbb{R}^3$ denote the radius $\varepsilon$ tubular neighborhood of $Z_*$. Each $|x| < \varepsilon^{-1}$ point in $Z_{T(i)}$ lies in $\mathcal{T}_\varepsilon$.*
c) *There is an isomorphism between $\mathcal{I}_*$ and $\mathcal{I}_{T(i)}$ over the $|x| < \varepsilon^{-1}$ part of $\mathbb{R}^d - \mathcal{T}_\varepsilon$ that identifies $v_{T(i)}$ as an $\mathcal{I}_*$ valued section of $\mathbb{V}^+$ that differs from $v_*$ by an $\mathcal{I}_*$ valued section of $\mathbb{V}^+$ with any $k < \varepsilon^{-1}$ version of its $C^k$ norm being less than $\varepsilon^{-1}$.*

An important application of this proposition is obtained by fixing a point $p \in Z$ and an orthonormal frame $L$ for $TX|_p$. A sequence of the form $\{T(i) = (p, L, \lambda_i)\}_{i=1,2,...}$ can then be used as the input for Proposition 4.1 if $\lim_{i \to \infty} \lambda_i = 0$.

The function $|v_*|$ is constructed by Lemma 4.5 in Section 4b and the data $Z_*$, $\mathcal{I}_*$ and $v_*$ are constructed in Section 4c. All assertions of the proposition except Item a) of the fifth bullet are proved as part of Lemmas 4.7 in Section 4c. Item a) of the fifth bullet is proved as part of Lemma 4.6 in Section 4b. Section 4a states and proves various lemmas that are used in Sections 4b and 4c.

Section 4d contains a lemma that gives a sufficient condition for two sequences $\{T(i)\}_{i=1,2,...}$ and $\{T(i)'\}_{i=1,2,...}$ of the sort described by the proposition to yield the same limit data set $(Z_*, \mathcal{I}_*, v_*)$. This lemma in Section 4d corrects the statement of an erroneous lemma in the published version of [T1] and in an old version of [T2].

As noted at the outset, the arguments in the proof of Proposition 4.1 borrow much from what is done in Sections 8 and 9 of the published version of [T1]. In particular, the statements and proofs of various results here are very similar to statements and proofs in these sections of this published version.

## a) Integrals of $|v_T|^2$ and $|\nabla v_T|^2$

The convention unless directed otherwise is to define all norms and the integration measure on $\mathbb{R}^4$ using the Euclidean metric. Even so, other metrics lurks in the background. To define these, let $T = (p, L, \lambda)$ denote a data set of the sort that is described at the start of the section. With $T$ given, define $\mathfrak{m}_T$ to be the product of $\lambda^{-2}$ times the $\phi_T$ pull-back of the metric on $X$. The metric $\mathfrak{m}_T$ differs from the Euclidean metric on the $|x| \leq \lambda^{-1} c_0^{-1}$ ball by a tensor with norm bounded by $c_0 \lambda^2 |x|^2$. Moreover, the norm of the derivative of this tensor is bounded by $c_0 \lambda^2 |x|$; and the norms of its derivatives to any given order $k \in \{2, ...\}$ are bounded by a k-dependent multiple of $\lambda^k$.

Any given $T = (p, L, \lambda)$ version of $|v_T|$ is ostensibly defined only on the $|x| \leq \lambda^{-1} r_0$ ball in $\mathbb{R}^d$; but, it is notationally convenient to extend this function to the whole of $\mathbb{R}^d$ by setting the extension equal to zero where $|x| > \lambda^{-1} r_0$. The function $|\nabla v_T|$ is also extended to the whole of $\mathbb{R}^d$ by defining it to be 0 on $Z_T$ and to be zero where $|x| \geq \lambda^{-1} r_0$.

The lemmas in this section refer to a given sequence $\{T(i) = (p_i, L_i, \lambda_i)\}_{i=1,2,...}$ whose elements are triples consisting of a point in Z, an oriented, orthonormal frame at that point and a number $\lambda \in (0, r_0]$. The sequence $\{T(i)\}_{i=1,2,...}$ is chosen so that $\{p_i\}_{i=1,2,...}$ converges to a point $p \in X$ with the additional constraint that $\lim_{i \to \infty} N_{(p_i)}(0) = N_{(p)}(0)$. The corresponding sequence $\{\lambda_i\}_{i=1,2,...}$ converges to 0. The sequence $\{T(i)\}_{i=1,2,...}$ is then used to construct a sequence of data sets $\{(Z_{T(i)}, \mathcal{I}_{T(i)}, v_{T(i)})\}_{i=1,2,...}$ according to the rules set forth prior to the statement of Proposition 4.1.

**Lemma 4.2**: *The sequence $\{v_{T(i)}\}_{i=1,2,...}$ has the properties that are listed in the bullets that follow. These bullets refers to a chosen number $R \in (0, \infty)$.*

- $\{\int_{|x|=R} |v_{T(i)}|^2\}_{i \in \{1,2,...\}}$ *converges as* $i \to \infty$ *with limit* $R^{d-1+2N_{(p)}(0)}$.
- $\{\int_{|x|\leq R} |v_{T(i)}|^2\}_{i \in \{1,2,...\}}$ *converges as* $i \to \infty$ *with limit* $\frac{1}{d+2N_{(p)}(0)} R^{d+2N_{(p)}(0)}$.
- $\{\int_{|x|\leq R} |\nabla v_{T(i)}|^2\}_{i \in \{1,2,...\}}$ *converges as* $i \to \infty$ *with limit* $N_{(p)}(0) R^{d-2+2N_{(p)}(0)}$.

*Moreover, the assertions in these three bullets also hold when the norms, covariant derivative and volume form on $\mathbb{R}^d$ or volume form on the constant $|x|$ spheres are defined for each $i \in \{1, 2, ...\}$ by the corresponding metric $\mathfrak{m}_{T(i)}$.*

This lemma is proved momentarily.

Let $\nabla_r$ denote the directional covariant derivative in the radial direction. The next lemma concerns the $i \to \infty$ limit of integrals on constant $|x|$ spheres in $\mathbb{R}^d$ of $|\nabla_r v_{T(i)}|^2$. As was the case in Section 2d with the integrals over $\partial B_r$ that appear in (2.15), some care must be taken because any given $T = (p, L, \lambda)$ version of the function $|\nabla_r v_T|^2$ need not be square integrable on any given constant $|x|$ sphere. However, the function

$$R \to \int_{|x|=R} |\nabla_r v_T|^2$$

(4.1)

is an $L^1$ function on $(0, \lambda^{-1} r_0)$. The proof that this is so copies an argument in Section 2d. To start, let $H_r: (0, \lambda^{-1} r_0) \to [0, \infty)$ denote the function given by the rule

$$R \to H_r(R) = \int_{|x|\leq R} |\nabla_r v_T|^2 \ .$$

(4.2)

This function is bounded and nondecreasing. The function $H_r$ is also continuous because $|\nabla_r v_T|$ is an $L^2$ function on the $|x| \leq \lambda^{-1} r_0$ ball. It is almost everywhere differentiable because, as noted subsequent to (2.14), a bounded, continuous and increasing function on an interval is almost everywhere differentiable with the derivative being an $L^1$ function. In the case of $H_r$, the derivative is almost everywhere given by (4.1).

**Lemma 4.3**: *The sequence of $L^1$ functions defined by the rule*

$$\{R \to \int_{|x|=R} |\nabla_r v_{T(i)} - \frac{N_{(p)}(0)}{R} v_{T(i)}|^2 \}_{\lambda \in (0, r_0/R)}$$

*is bounded and it converges to zero as $i \to \infty$ on compact subsets of $[0, \infty)$.*

The proof of Lemma 4.3 is given after that of Lemma 4.2.
    The next lemma is not used in subsequent arguments; it is stated because it is an analog to Lemma 8.4 in the published version of [T1]. By way of background, what is said in Section 2d after (2.14) can be repeated to prove that any $T = (p, L, \lambda)$ version of the function on $(0, \lambda^{-1} r_0)$ that is defined by the rule $R \to \int_{|x|=R} |\nabla v_T|^2$ is also an $L^1$ function.

**Lemma 4.4**: *The sequences of $L^1$ functions that are defined by the rules*

$$\{R \to \int_{|x|=R} |\nabla v_{T(i)}|^2 \}_{i=1,2,...} \quad and \quad \{R \to \int_{|x|=R} |\nabla_r v_{T(i)}|^2 \}_{i=1,2,...}$$

*are bounded and converge on compact subsets of $[0, \infty)$ as $i \to \infty$ to the respective limits*

$$N_{(p)}(0)(d - 2 + 2N_{(p)}(0)) R^{d-3+2N_{(p)}(0)} \quad and \quad N_{(p)}(0)^2 R^{d-3+2N_{(p)}(0)}.$$

This lemma is proved after the proof of Lemma 4.4.

*Proof of Lemma 4.2*: What is said at the outset about the difference between the Euclidean metric and any given $T = (p, L, \lambda)$ version of the metric $m_T$ implies that the assertions of the lemma hold using the Euclidean metric to define the norms and integration if and only if they hold with the norms and integration defined at for each index $i \in \{1, 2, ...\}$ by the corresponding metric $m_{T(i)}$.
    To prove the first bullet, fix $i \in \{1, 2,...\}$ and invoke the definition of $\kappa$ to see that

$$\int_{|x|=1} |v_{T(i)}|^2 = 1 + r_i,$$

(4.3)

with $|r_i| \leq c_0 \lambda_i^2$. To exploit this formula, first write the function $N_{(p_i)}$ on $(0, c_0^{-1} r_0]$ as $N_{(p)}(0) + n_i$. The second bullet of Lemma 3.4 can be invoked to prove the following: Given $\varepsilon > 0$, there exists $m_\varepsilon$ such that if $i > m_\varepsilon$ and $r \in (0, m_\varepsilon^{-1}]$, then $|n_i| \leq \varepsilon$. This fact with (3.4) and (4.5) imply that

$$\int_{|x|=R} |v_{T(i)}|^2 = (1+\mathfrak{e}_i(R))\, R^{d-1+2N_{(p)}(0)}$$

(4.4)

with the various $i \in \{1, 2, \ldots\}$ versions of $\mathfrak{e}_i(\cdot)$ defining a sequence of functions on $(0, \infty)$ that converges uniformly to 0 as $i \to \infty$ on compact subsets of $(0, \infty)$.

The second bullet of the lemma follows by integrating the first bullet and invoking the fact that any given point's version of the function $\kappa$ is non-decreasing. The third bullet follows from the first bullet using what was said in the preceding paragraph about the $i \to \infty$ limits of the sequence of functions $\{\mathfrak{n}_i = N_{(p_i)} - N_{(p)}(0)\}_{i=1,2,\ldots}$.

***Proof of Lemma 4.3***: It follows from the fourth bullet of Lemma 4.2 that the indicated sequence has bounded $L^1$ norm on compact subsets of $(0, \infty)$. To see about convergence, fix for the moment an index $i \in \{1, 2,\ldots\}$ and focus for the moment on $T(i) = (p_i, L_i, \lambda_i)$. Assume that $R \in (0, \lambda_i^{-1} r_0)$ is such that the function $N_{(p_i)}$ is differentiable at $\lambda_i R$. The inequality in second bullet of Lemma 3.2 can then be invoked; and a rescaling of this inequality can be written as

$$(r \tfrac{d}{dr} N_{(p_i)})(\lambda_i R) \geq \Big( \int_{|x|=R} |\nabla_r v_{T(i)} - \tfrac{N_{p_i}(\lambda_i R)}{R} v_{T(i)}|^2 - c_0 \lambda_i^2 R^2 \Big) R^{3-d} \Big( \tfrac{K_{p_i}(\lambda_i)}{K_{p_i}(\lambda_i R)} \Big)^2$$

(4.5)

Fix $\varepsilon \in (0, 1)$. If the index $i$ is sufficiently large, with a bound determined by $\varepsilon$ and $R$, then $N(\lambda_i R)$ will differ from $N_p(0)$ by at most $\varepsilon$. That this is so follows from the second bullet of Lemma 3.4. Fix $T > 1$ such that $R \in (T^{-1}, T)$. If $i$ is sufficiently large, with an upper bound given by $\varepsilon$ and $T$, then

$$\Big| \tfrac{K_{p_i}(\lambda_i)}{K_{p_i}(\lambda_i R)} - R^{-N_p(0)} \Big| < \varepsilon$$

(4.6)

also. Thus, if $i$ is sufficiently large (given $\varepsilon$ and $T$), then (4.5) implies that

$$(r \tfrac{d}{dr} N_{(p_i)})(\lambda_i R) \geq c_T^{-1} \Big( \int_{|x|=R} |\nabla_r v_{T(i)} - \tfrac{N_p(0)}{R} v_{T(i)}|^2 - \varepsilon^2 \Big),$$

(4.7)

with $c_T$ denoting a number that is greater than 1 and depends only on $T$.

To continue, let $I$ denote the subset of $(T^{-1}, T)$ where the derivative of $N_{(p_i)}$ is greater than $\varepsilon$. It follows from the second bullet of Lemma 3.2 by integration that

$$\varepsilon \int_I \tfrac{1}{s} ds < 2(\sup_{s \in I} N_{(p_i)}(\lambda_i s) - \inf_{s \in I} N_{(p_i)}(\lambda_i s)) + \varepsilon^2 .$$

(4.8)

when i is sufficiently large with a lower bound determined by $\varepsilon$ and T. The right hand side of (4.10) will be less than $\varepsilon^2$ when i is sufficiently large because of what is said in the second bullet of Lemma 3.4. This implies that the measure of the set I is less than $T\varepsilon$ when i is sufficiently large. Looking at (4.7), it then follows that the measure of the set in $(T^{-1}, T)$ where

$$\int_{|x|=R} |\nabla_r v_{T(i)} - \tfrac{N_p(0)}{R} v_{T(i)}|^2 > \varepsilon$$

(4.9)

is less than $\varepsilon T$ when i is sufficiently large with an upper bound determined by $\varepsilon$ and T. Since $\varepsilon$ and T can be chosen at will subject to $\varepsilon$ being positive and $T > 1$, this last conclusion implies what is asserted by the lemma.

*Proof of Lemma 4.4*: The fact that the sequences have bounded $L^1$ norm on compact intervals in $(0, \infty)$ follows from the third bullet of Lemma 4.2. The assertion that the right most sequence converges and the assertion about the limit follow from Lemma 4.3. The proof of the assertion about the left most sequence uses the fact, as noted previously, that the version of the function $r \to N(r)$ defined by any given point in X is almost everywhere differentiable and that its derivative is an $L^1$ function. It follows directly from the definition in (3.3) that the derivative of N is equal as an $L^1$ function to

$$\tfrac{d}{dr} N = \tfrac{1}{r^{d-2} K(r)^2} \int_{\partial B_r} |\nabla v|^2 - (d - 2 + 2N(r)) N(r) \tfrac{1}{r} .$$

(4.10)

With the preceding understood, fix $R \in (0, \infty)$ and then index $i \in \{1, 2, \ldots\}$ labeling the data set $T(i) = (p_i, L_i, \lambda_i)$ such that $\lambda_i < r_0/R$. The $r = \lambda_i R$ version of (4.10) says that

$$(r \tfrac{d}{dr} N_{(p_i)})(\lambda_i R) = (R^{-d+3} \int_{|x|=R} |\nabla v_{T(i)}|^2 - (d - 2 + 2 N_{(p_i)}(R\lambda_i)) N_{(p_i)}(R\lambda_i) (\tfrac{K_{p_i}(\lambda_i)}{K_{p_i}(\lambda_i R)})^2$$

(4.11)

Let $T \in (1, \infty)$. Fix $\varepsilon > 0$. If i is sufficiently large given $\varepsilon$ and T, and if $R \in (T^{-1}, T)$, then $|N_{(p_i)}(R\lambda_i) - N_p(0)| < \varepsilon$ and (4.6) holds. Meanwhile, the second bullet of Lemma 3.2 says that the left hand side of (4.11) is no smaller than $-c_0 \lambda_i^2 R^2$ and this is less than $\varepsilon$ if i is sufficiently large given T. It follows as a consequence that

$$\int_{|x|=R} |\nabla v_{T(i)}|^2 \; - \; (d-2+2N_p(0))\, N_p(0)\, R^{d-3+2N(0)}$$

(4.12)

is greater than -ε if i is sufficiently large with a lower bound determined by T and ε.

To see about an upper bound for what is written in (4.12), let I again denote the subset of $(T^{-1}, T)$ where the derivative of $N_{(p_i)}$ is greater than ε. As noted in the proof of Lemma 4.3, this set has measure less than εT when i is sufficiently large with a lower bound determined by ε and T. Looking at (4.11), this implies in turn that the set of values for R in $(T^{-1}, T)$ where the expression in (4.12) is greater than ε has measure at most εT when i is sufficiently large. Since ε and T can be chosen at will subject to ε being positive and T > 1, the preceding conclusion implies what is asserted by Lemma 4.4 about its left most sequence.

**b) The section $|v_*|^{-1} v_* \otimes v_*$**

This subsection constructs what turns out to be the section $|v_*|^{-1} v_* \otimes v_*$ of the bundle $\mathbb{V}^+ \otimes \mathbb{V}^+$ that is defined by Proposition 4.1's section $v_*$ of $\mathbb{V}^+ \otimes \mathcal{I}_*$. To set a convention that is used here and in what follows, let $\mathbb{V}$ denote for the moment a vector bundle over $\mathbb{R}^d$. Let $\Xi \subset \{1, 2, \ldots\}$ denote an infinite subsequence and suppose that $\{(\lambda_i, f_i)\}_{i \in \Xi}$ is a sequence with i'th member being a number $\lambda_i \subset (0, r_0)$ and section $f_i$ of $\mathbb{V}$ over the radius $\lambda_i^{-1} r_0$ ball about the origin in $\mathbb{R}^d$. The convention in what follows is to extend $f_i$ to the whole of $\mathbb{R}^d$ by declaring it to be zero where $|x| \geq \lambda_i^{-1} r_0$.

The two lemmas in this section refer to a chosen sequence $\{T(i) = (p_i, L_i, \lambda_i)\}_{i=1,2,\ldots}$ with each member again being a point in Z, an oriented, orthonormal frame at the point and a number in $(0, r_0]$. This sequence is such that $\{p_i\}_{i=1,2,\ldots}$ converges to a point $p \in Z$, and $\{N_{(p_i)}(0)\}_{i=1,2,\ldots}$ converges to $N_{(p)}(0)$, and $\{\lambda_i\}_{i=1,2,\ldots}$ converges to 0. These lemmas are in the respective Parts 1 and 2 of this subsection.

*Part 1*: The first lemma in this section asserts in part the existence of a subsequence of $\{|v_{T(i)}|^{-1} v_{T(i)} \otimes v_{T(i)}\}_{i=1,2,\ldots}$ that converge strongly in the $L^2$ topology on any ball about the origin in $\mathbb{R}^4$. To set the stage, suppose for the moment that $i \in \{1, 2, \ldots\}$ is such that $p_i$ is in the radius $r_0$ ball centered about the limit point $p \in Z$ of the sequence $\{p_i\}_{i=1,2,\ldots}$. Parallel transport along the short geodesic from $p_i$ to p identifies the frame $L_i$ as an orthonormal frame for $TX|_p$. This new incarnation for $L_i$ is still denoted by $L_i$.

**Lemma 4.5**: *There exists $\kappa > 0$ and a subsequence $\Lambda \subset \{1, 2, \ldots\}$ with three salient properties: The first is that the sequence $\{L_i\}_{i \in \Lambda}$ defines a convergent sequence of orthonormal frames at $TX|_p$. The other two properties are listed below.*

- *The sequence of sections $\{|v_{T(i)}|^{-1} v_{T(i)} \otimes v_{T(i)}\}_{i \in \Lambda}$ converges weakly in the $L^2_1$ topology on every ball about the origin in $\mathbb{R}^d$. The limit is an $L^2_1$ and $L^\infty$ section of $\mathbb{V}^+ \otimes \mathbb{V}^+$ on any given radius $R \in (0, \infty)$ ball about the origin with $L^2$ norm equal to the square root of $\frac{1}{d+2N_{(p)}(0)} R^{d+2N_{(p)}(0)}$ and with pointwise norm bounded by $\kappa R^{N_{(p)}(0)}$.*
- *The sequence $\{|\nabla v_{T(i)}|\}_{T(i) \in \Lambda}$ converges weakly in the $L^2$ topology on each such ball.*

Let $\Lambda$ denote one of the sequences from Lemma 4.6. Denote the corresponding limit section of $\{|v_{T(i)}|^{-1} v_{T(i)} \otimes v_{T(i)}\}_{i \in \Lambda}$ by $|v_*|^{-1} v_* \otimes v_*$. Denote the pointwise norm of this section by $|v_*|$. Although the notation indicates that these limits are independent of $\Lambda$, an assertion to this effect is not made.

*Proof of Lemma 4.5*: Fix an infinite subsequence $\Xi \subset \{1, 2, \ldots\}$ so that if $i \in \Xi$, then the frame $L_i$ can be viewed as a frame for $TX|_p$, and so that the corresponding sequence of frames for $TX|_p$ converges. The existence of a subsequence in $\Xi$ such that the corresponding sequence from $\{|v_{T(i)}|^{-1} v_{T(i)} \otimes v_{T(i)}\}_{i \in \Xi}$ converges weakly in the $L^2_{1;\text{loc}}$ topology follows from the bounds given by the second and third bullets of Lemma 4.2. The assertion about the $L^2$ norm of the limit follows from the second bullet in Lemma 4.2 because weak $L^2_1$ convergence on any given ball implies strong $L^2$ convergence on that ball. To see about asserted pointwise bound for the limit section, fix $R \in (0, \infty)$ for the moment and then fix an index $i \in \{1, 2, \ldots\}$. Lemma 3.3 and (4.6) imply that $|v_{T(i)}|$ obeys

$$|v_{T(i)}| \leq c_0 (1 + \mathfrak{r}_i(R)) R^{N_{(p)}(0)} \quad \text{where } |x| \leq R$$

(4.13)

with $\mathfrak{r}_i$ being small when $i$ is large in the sense that $\lim_{i \to 0} \mathfrak{r}_i = 0$ uniformly on compact subsets in $(0, \infty)$. The asserted $c_0 R^{N_{(p)}(0)}$ bound in Lemma 4.5 on the norm of the limit section of the sequence $\{|v_{T(i)}|^{-1} v_{T(i)} \otimes v_{T(i)}\}_{i \in \Lambda}$ follows from (4.13). Lemma 4.2 also implies that the sequence $\{|\nabla v_{T(i)}|\}_{i \in \Xi}$ converges weakly in the $L^2$ topology on any ball in $\mathbb{R}^d$.

The two bullets of Lemma 4.5 follow directly from what is said in the preceding paragraph using a diagonal subsequence argument.

*Part 2*: The next lemma says more about the pointwise convergence of the sequence $\{|v_{T(i)}|^{-1} v_{T(i)} \otimes v_{T(i)}\}_{i \in \Lambda}$ to the limit section $|v_*|^{-1} v_* \otimes v_*$.

**Lemma 4.6**: *There exists a number $\kappa > 10$ that is independent of the sequence $\{T(i)\}_{i=1,2,\ldots}$ and has the following significance: Let $\Lambda \subset \{1, 2, \ldots\}$ denote a sequence as described by Lemma 4.5 and let $|v_*|^{-1} v_* \otimes v_*$ denote the corresponding limit section from the first bullet of Lemma 4.5. The section $|v_*|^{-1} v_* \otimes v_*$ is Holder continuous with exponent $\upsilon = \kappa^{-1}$ and the* *sequence $\{|v_{T(i)}|^{-1} v_{T(i)} \otimes v_{T(i)}\}_{i \in \Lambda}$ converges to $|v_*|^{-1} v_* \otimes v_*$ in the $\upsilon = \kappa^{-1}$ Holder topology on compact subsets of $\mathbb{R}^d$.*

***Proof of Lemma 4.6***: The proof rewrites the arguments from Part 2 of Section 8c in the published version of [T1] that prove the latter's Lemma 8.6. This five steps that follow correspond to the five steps that prove this lemma in the published version of [T1].

<u>Step 1</u>: Fix $i \in \Lambda$ and define $z_i$ by writing $K_{(p_i)}(\lambda_i)$ as $z_i \lambda_i^{N_{(p_i)}(0)}$. Let $Z_{T(i)}$ denote the zero locus of $|v_{T(i)}|$ in the radius $\lambda_i^{-1} r_0$ ball about the origin in $\mathbb{R}^d$. Let $z \in Z_{T(i)}$ denote a given point. This point corresponds to a point, $p_z$, on the zero locus of $v$ in the radius $r_0$ ball centered at p. Suppose that $R > 16$ has been specified, that $\lambda_i < R^{-1} r_0$ and that the point z is in the radius $\frac{1}{8} R$ ball centered at the origin in $\mathbb{R}^d$. Define $z_{p,i} \in (0, c_0]$ to be the infimum of the numbers $z$ such that $K_{(p_z)}(r) \leq z \, r^{N_{(p_z)}(0)}$ for all $r \in (0, \frac{1}{2} R\lambda]$. As explained directly, there exist a purely R-dependent number $c_R$ which is greater than R and has the following significance: If $\lambda_i \leq c_R^{-1} r_0$ then

$$z_{p,i} z_I^{-1} \leq c_0 (R\lambda_i)^{N_{(p_i)}(0) - N_{(p_z)}(0)} \quad .$$
(4.14)

To prove this, use the $p_z$ version of (3.4) with the third bullet of $p_z$'s version of Lemma 3.2 to see that

$$z_{p,i} \leq c_0 (\tfrac{1}{2} R\lambda_i)^{-N_{(p_z)}(0)} K_{(p_z)}(\tfrac{1}{2} R\lambda_i) \quad .$$
(4.15)

Now let $\mathcal{A}$ denote the spherical annulus with inner radius $\frac{1}{2} R\lambda_i$ and outer radius $\frac{9}{16} R\lambda_i$. The fact that $K_{(p_z)}$ is increasing implies that $(K_{(p_z)}(\frac{1}{2} R\lambda_i))^2 (R\lambda_i)^d$ is less than $c_0$ times the integral of $|v|^2$ over $\mathcal{A}$. Meanwhile, this spherical annulus $\mathcal{A}$ is contained in ball of radius $R\lambda_i$ centered at $p_i$ if the index i is large. It follows from Lemma 3.2 that the integral of $|v|^2$ over the latter ball is no greater than $c_0 z_i^2 (R\lambda_i)^{d + 2 N_{(p_i)}(0)}$. The former integral is less than the latter and so $K_{(p_z)}(\frac{1}{2} R\lambda_i)$ is less than $z_i (R\lambda_i)^{N_{(p_i)}(0)}$.

<u>Step 2</u>: Fix $R > 16$, and then an index $i \in \Lambda$ such that $\lambda_i < c_R^{-1} r_0$ with $c_R$ as described in Step 1. Lemma 3.3 and (4.14) have the following implication: Let z denote a point in $Z_{T(i)}$ that lies in the radius $(16 + c_0)^{-1} R$ ball about the origin; and let $x \in \mathbb{R}^d$ denote a point in the concentric ball of twice this radius. Then

$$|v_{T(i)}|(x) \le c_{\diamond R} \, |z-x|^{N_{(p_z)}(0)} \, .$$

(4.16)

with $c_{\diamond R}$ denoting here and subsequently a number greater than 1 that depends only on R. Its value can be assumed to increase between successive appearances. As Lemma 3.2 finds $c_0 > 1$ such that $N_{(p_z)}(0) > c_0^{-1}$, this last bound implies in particular that the section $|v_{T(i)}|^{-1} v_{T(i)} \otimes v_{T(i)}$ is Holder continuous with exponent $\upsilon > c_0^{-1}$ along its zero locus in any given radius ball about the origin in $\mathbb{R}^d$.

<u>Step 3</u>: Fix $i \in \Lambda$ with $\lambda_i < c_R^{-1} r_0$. Given $\delta \in (0, c_0^{-1}]$, let $U_{i,R,\delta}$ denote the set of points in the $|x| \le R$ ball with distance $\delta$ or more from $Z_{T(i)}$. Supposing that $R \ge 1$, let B denote a ball of radius $\frac{1}{4} \delta$ in $U_{i,R,\delta/2}$. Being closed and coclosed, the 1-form $v$ obeys the equation $\nabla^\dagger \nabla v + \text{Ric}(v) = 0$. This equation with Lemma 4.2 implies that

$$\int_B |\nabla \nabla \nabla v_{T(i)}|^2 \le c_0 \delta^{-6} R^{d-2+2N_{(p)}(0)} \, .$$

(4.17)

The preceding inequality with a standard Sobolev inequality can be used to prove that if x and y are points in $U_{i,R,\delta}$, then

$$|v_{T(i)}(x) - v_{T(i)}(y)| \le c_{\diamond R} \delta^{-3} |x-y|^{1/4}$$

(4.18)

It follows from this inequality that $|(|v_{T(i)}|^{-1} v_{T(i)} \otimes v_{T(i)})|_x - (|v_{T(i)}|^{-1} v_{T(i)} \otimes v_{T(i)})|_y|$ is also bounded by $c_{\diamond R} \delta^{-3} |x-y|^{1/4}$ when x and y are both in $U_{i,R,\delta}$.

<u>Step 4</u>: Suppose that x and y are both in the radius R ball about the origin. If either x or y is in $Z_{T(i)}$, then (4.16) implies that

$$|(|v_{T(i)}|^{-1} v_{T(i)} \otimes v_{T(i)})|_x - (|v_{T(i)}|^{-1} v_{T(i)} \otimes v_{T(i)})|_y| \le c_{\diamond R} \, |x-y|^{N_{(p_z)}(0)}$$

(4.19)

for some point $z \in Z_{T(i)}$ because the distance between x and y is at most that between either and the set $Z_{T(i)}$.

Now assume that neither x nor y is in $Z_{T(i)}$. Set $\delta$ to equal the minimum of $\text{dist}(x, Z_{T(i)})$ and $\text{dist}(y, Z_{T(i)})$. If $|x-y| \le \delta^{16}$, then the conclusions of Step 3 imply that

$$|(|v_{T(i)}|^{-1} v_{T(i)} \otimes v_{T(i)})|_x - (|v_{T(i)}|^{-1} v_{T(i)} \otimes v_{T(i)})|_y| \le c_{\diamond R} |x-y|^{1/16} \, .$$

(4.20)

Now suppose that $|x-y| \ge \delta^{16}$ and suppose in addition that $\delta = \text{dist}(x, Z_{T(i)})$. It follows in this case that there exists $z \in Z_{T(i)}$ with $|x-z| \le 2\delta$. Since $|y-z| \le |x-z| + |x-y|$ and $|x-z| \le 2\delta$ and $\delta \le |x-y|^{1/16}$, it then follows that $|y-z| \le 2|x-y|^{1/16} + |x-y|$. Granted this fact, use the bound $|(|\nu_{T(i)}|^{-1}\nu_{T(i)} \otimes \nu_{T(i)})|_x - (|\nu_{T(i)}|^{-1}\nu_{T(i)} \otimes \nu_{T(i)})|_y| \leq \max(|\nu_{T(i)}|(x), |\nu_{T(i)}|(y))$ with (4.16) to see that

$$|(|\nu_{T(i)}|^{-1}\nu_{T(i)} \otimes \nu_{T(i)})|_x - (|\nu_{T(i)}|^{-1}\nu_{T(i)} \otimes \nu_{T(i)})|_y| \leq c_{\diamond R} |x-y|^{N_{z,0}/100}$$

(4.21)

for some $z \in Z_{T(i)}$.

Step 5: The bounds in Step 4 with the lower bound $N_{(\cdot)} > c_0^{-1}$ from the final bullet of Lemma 3.2 imply that the very large i part of the sequence $\{|\nu_{T(i)}|^{-1}\nu_{T(i)} \otimes \nu_{T(i)}\}_{i \in \Lambda}$ is bounded in any $\upsilon < c_0^{-1}$ Holder topology on the radius R ball. This sequence is therefore 'equi-Holder continuous'. As the sequence converges in the $L^2$ topology on the $|x| \leq R$ ball, the arguments that proves the Arzela-Arcoli theorem proves that $\{|\nu_{T(i)}|^{-1}\nu_{T(i)} \otimes \nu_{T(i)}\}_{i \in \Lambda}$ converges in any given exponent $\upsilon < c_0^{-1}$ Holder topology in the $|x| \leq R$ ball. This implies that the section $|\nu_*|\nu_* \otimes \nu_*$ is Holder continuous for any given exponent $\upsilon < c_0^{-1}$.

### c) The construction of $Z_*$, $\mathcal{I}_*$ and $\nu_*$

Supposing that $\{T(i)\}_{i=1,2,...}$ is as described in the beginning of Sections 4a, let $\Lambda$ denote the subsequence of $\{1, 2, ....\}$ that is supplied by Lemma 4.5. The upcoming Lemma 4.7 takes the sequence $\{T(i)\}_{i \in \Lambda}$ as input and gives back a data set $(Z_*, \mathcal{I}_*, \nu_*)$ with $Z_*$ being a closed set in $\mathbb{R}^d$, with $\mathcal{I}_*$ being a real line bundle defined on $\mathbb{R}^d - Z_*$ and with $\nu_*$ being a harmonic section on $\mathbb{R}^d - Z_*$ of the bundle $\mathbb{S}^+ \otimes \mathcal{I}_*$.

With regards to terminology, a set in $\mathbb{R}^d$ is said to be scale invariant when it is the cone on a set in the $|x| = 1$ sphere. Supposing that $R \in (0, \infty)$, let $\psi_R$ to denote the map from $\mathbb{R}^d$ to itself that is defined by the rule $x \to \psi_R(x) = R x$. A set is scale invariant if and only if it is mapped to itself by each $R \in (0, \infty)$ version of $\psi_R$. The lemma also uses $N_0$ to denote $N_{(p)}(0)$ with p being a given point in X.

**Lemma 4.7**: *There exists a number $\kappa > 1$ that is independent of the sequence $\{T(i)\}_{i=1,2,...}$ and is such that the following is true: Fix a subsequence $\Lambda \subset \{1, 2, ...\}$ as described in Lemma 4.5 and let $|\nu_*|$ denote the limit function of the corresponding sequence $\{|\nu_{T(i)}|\}_{i \in \Lambda}$. Use $Z_* \subset \mathbb{R}^d$ to denote the zero locus of $|\nu_*|$.*

- *The set $Z_*$ is a closed, scale invariant set.*
- *There exists a real line bundle $\mathcal{I}_* \to \mathbb{R}^d - Z_*$ and a harmonic, $\mathcal{I}_*$ valued 1-form on $\mathbb{R}^d - Z_*$ with norm equal to $|\nu_*|$. This $\mathcal{I}_*$ valued 1-form is denoted by $\nu_*$.*

- *Supposing that $R \in (0, \infty)$, there is an isometric isomorphism between $\psi_R^* \mathcal{I}_*$ and $\mathcal{I}_*$; and these various $R \in (0, \infty)$ isomorphisms can be used to write $v_*$ as $u_* |x|^{N_{(p)}(0)}$ with $u_*$ being a section of $\mathbb{V}^+$ over the complement of $Z_*$ in the $|x| = 1$ sphere.*
- *The 1-form $v_*$ is such that for each $R > 0$,*

  i) $\int_{|x|=R} |v_*|^2 = R^{d-1+2N_{(p)}(0)}$.

  ii) $\int_{|x|=R} |\nabla v_*|^2 = N_{(p)}(0)(d-2+2N_{(p)}(0)) R^{d-3+2N_{(p)}(0)}$.

- *Fix $\varepsilon \in (0, 1]$ and let $B_{1/\varepsilon}$ denote the $|x| < \varepsilon^{-1}$ ball and introduce $\mathcal{T}_\varepsilon$ to denote the radius $\varepsilon$ tubular neighborhood $Z_* \cap B_{1/\varepsilon}$. There exists $\iota_\varepsilon > 1$ with the following significance: Suppose that $i \in \Lambda$ obeys $i > \iota_\varepsilon$. Then $Z_{T(i)} \cap B_{1/\varepsilon} \subset \mathcal{T}_\varepsilon$ and there is an isomorphism between $\mathcal{I}_*$ and $\mathcal{I}_{T(i)}$ over $B_{1/\varepsilon} - \mathcal{T}_\varepsilon$ that identifies $v_{T(i)}$ as an $\mathcal{I}_*$ valued 1-form that differs from $v_*$ by an $\mathcal{I}_*$ valued 1-form whose $C^k$ norm on $B_{1/\varepsilon} - \mathcal{T}_\varepsilon$ is less than $\varepsilon$.*

***Proof of Lemma 4.7***: The proof has six parts.

*Part 1*: Lemma 4.3 implies that the sequence $\{\partial_r |v_{T(i)}| - \frac{N_{(p)}(0)}{|x|} |v_{T(i)}|\}_{i \in \Lambda}$ converges to zero on any given ball in $\mathbb{R}^d$ in the $L^2$ topology. It follows as a consequence that $|v_*|$ obeys the equation $\partial_r |v_*| - \frac{N_{(p)}(0)}{|x|} |v_*| = 0$ and this implies in turn that $|v_*|$ can be written as $|x|^{N_{(p)}(0)} \hat{u}$ with $\hat{u}$ being a Holder continuous function on the $|x| = 1$ sphere. It follows as a consequence that the zero set of $|v_*|$ is a scale invariant set, it being the union of the origin with the cone on the zero locus of the function $\hat{u}$. The set $Z_*$ is defined to be the zero locus of $|v_*|$.

*Part 2*: The complement in the $|x| = 1$ sphere of its points from $Z_*$ can be covered by a locally finite, countable collection of balls. Moreover, there exists such a cover with the following property: If $B$ is a ball from the cover, then the section $|v_*|^{-1}(v_* \otimes v_*)$ of $\mathbb{V}^+ \otimes \mathbb{V}^+$ on $B$ is nowhere zero. The section $|v_*|^{-1}(v_* \otimes v_*)$ has rank 1 since this is the case for all of the sections in the sequence $\{|v_{T(i)}|^{-1}(v_{T(i)} \otimes v_{T(i)})\}_{i \in \Lambda}$. Since $B$ is contractible, there must exist a Holder continuous section of $\mathbb{V}^+$ on $B$, this denoted by $v_{*B}$, such that $|v_{*B}|^{-1}(v_{*B} \otimes v_{*B})$ and $|v_*|^{-1}(v_* \otimes v_*)$ are equal. Note $v_{*B}$ is unique up to multiplication by $\pm 1$.

Fix a cover with the properties just described and denote it by $\mathfrak{U}$. For each $B \in \mathfrak{U}$, let $C_B \subset \mathbb{R}^d$ denote the cone on $B$. If $n$ is a positive integer and $B \in \mathfrak{U}$, let $C_{B,n}$ denote the $|x| \in (n^{-1}, n)$ part of $C_B$.

*Part 3*: Supposing that $B \subset \mathfrak{U}$ and n is a positive integer, there exists $\iota_{B,n}$ such that if $i \in \Lambda$ and i is greater than $\iota_{B,v}$, then

$$\left| |\nu_{T(i)}|^{-1}(\nu_{T(i)} \otimes \nu_{T(i)}) - |\nu_*|^{-1}(\nu_* \otimes \nu_*) \right| < c_0^{-1} |\nu_*|$$

(4.22)

on the whole of $C_{B,n}$. It follows that the pairing $\langle \nu_{*B}, \nu_{T(i)} \rangle$ gives a nowhere zero section of $\mathcal{I}_{T(i)}$ over B; thus it defines an isomorphism over B between $\mathcal{I}_{T(i)}$ and the product $\mathbb{R}$ bundle.

This isomorphism writes $\nu_{T(i)}$ over $C_{B,n}$ as an honest section of the product $\mathbb{V}^+$ bundle. Denoted this section by $\nu_{T(i)B}$. It follows from Lemma 4.6 that the corresponding sequence of $i \in \Lambda$ and $i > \iota_{B,n}$ versions of $\nu_{T(i)B}$ converge to $\nu_{*B}$ in a Holder topology.

If $i \in \Lambda$ and $i > \iota_{B,n}$, then the equation $\mathcal{D}\nu = 0$ appears as an equation for $\nu_{T(i)B}$ on B that has the form

$$\mathcal{D}_0 \nu_{T(i)B} + \mathfrak{z}_i \nu_{T(i)B} = 0 ,$$

(4.23)

with $\mathfrak{z}_i$ being a first order operator whose coefficients have norms bounded by $c_0 \lambda_i^2 |x|$. Moreover, the norms of the derivatives of these coefficients are bounded by $c_0 \lambda_i^2$ and the norms of the derivatives of any given integer order $k > 1$ are bounded by $c_0 \lambda_i^k$.

Given that the sequence of $i \in \Lambda$ and $i > \iota_{B,n}$ versions of $\nu_{T(i)B}$ converge in the Holder topology, it then follows from (4.23) using standard elliptic estimates that this sequence converges in the $C^\infty$ topology. It follows as a consequence that $\nu_{*B}$ is smooth and that it obeys the equation $\mathcal{D}_0 \nu_{*B} = 0$. Lemma 4.4 implies that $\nabla_r \nu_{*B} = \frac{N_{(p)}(0)}{|x|} \nu_{*B}$.

*Part 4*: Suppose now that B and B´ are in $\mathfrak{U}$ and that they intersect. If the index i is from $\Lambda$ and if it is greater than $\iota_{B,1}$, then there are two identifications on $B \cap B´$ between $\mathcal{I}_{T(i)}$ and the product $\mathbb{R}$ bundle, the one coming from Part 2's isomorphism over B and the other coming from Part 3's isomorphism over B´. These identifications either agree or not, and if not, then one is -1 times the other. Let $I_{i,BB´} = 1$ if they agree and let $I_{i,BB´} = -1$ if they disagree. As explained directly, if i and k are sufficiently large, then $I_{i,BB´} = I_{k,BB´}$. To prove this claim, use the respective i and k identifications between the product $\mathbb{R}$ bundle and $\mathcal{I}_{T(i)}$ and $\mathcal{I}_{T(k)}$ on B to view $\nu_{T(i)}$ and $\nu_{T(k)}$ on $B \cap B´$ as sections of $\mathbb{V}^+$. This understood, then $|\nu_{T(i)B} - \nu_{T(k)B}|$ is very small on $B \cap B´$ because both $\nu_{T(i)B}$ and $\nu_{T(k)B}$ are very close to $\nu_{*B}$. Meanwhile, $|I_{i,BB´} \nu_{T(i)B} - I_{k,BB´} \nu_{T(k)B}|$ is also very small because both $I_{i,BB´} \nu_{T(i)B}$ and $I_{k,BB´} \nu_{T(k)B}$ are very close to $\nu_{*B´}$. Since neither $\nu_{*B}$ nor $\nu_{*B´}$ is zero, it is not possible for both $|\nu_{T(i)B} - \nu_{T(k)B}|$ and $|\nu_{T(i)B} + \nu_{T(k)B}|$ to be simultaneously very small.

It follows from what was just said that for each intersecting pair B, B´ from $\mathfrak{U}$, the sequence $\{I_{i,BB'}: i \in \Lambda$ *and* $i > \iota_B + \iota_{B'}\}$ is constant for i sufficiently large. Use $I_{BB'}$ to denote this constant. This number is either +1 or -1. It also follows from what was said in the preceding paragraph that $\nu_{*B} = I_{BB'} \nu_{*B'}$.

*Part 5*: The collection $\{I_{BB'}\}_{B,B' \in \mathfrak{U}}$ necessarily obey the cocycle constraints that define a $\mathbb{Z}/2$ bundle over $\mathbb{R}^d - Z_*$. This is so because each $I_{BB'}$ is equal to all sufficiently large i versions of $I_{i,BB'}$ and the collection $\{I_{i,BB'}: B, B' \in \mathfrak{U}$ *and* $i > \iota_B + \iota_{B'}\}$ obey the cocycle constraints since this collection defines the vector bundle $\mathcal{I}_{T(i)}$. The $\mathbb{Z}/2$ bundle defined by $\{I_{BB'}\}_{B,B' \in \mathfrak{U}}$ is $\mathcal{I}_*$. This definition using cocycles on the $|x| = 1$ sphere guarantees that any $R > 0$ version of $\psi_R^* \mathcal{I}_*$ is canonically isomorphic to $\mathcal{I}_*$.

The fact that $\nu_{*B} = I_{BB'} \nu_{*B'}$ on $B \cap B'$ for any pair of intersecting balls B, B´ $\in \mathfrak{U}$ implies that the collection $\{\nu_{*B}\}_{B \in \mathfrak{U}}$ defines a section of the bundle $\mathbb{V}^+ \otimes \mathcal{I}_*$ over $\mathbb{R}^d - Z_*$. This section is $\nu_*$. It is harmonic because each $B \in \mathfrak{U}$ version of $\nu_{*B}$ is harmonic. It also obeys $\nabla_r \nu_* = \frac{N_{(p)}(0)}{|x|} \nu_*$ because this equation is obeyed by each $B \in \mathfrak{U}$ version of $\nu_{*B}$. The latter equation implies that $\nu_*$ can be written using the isomorphisms between $\mathcal{I}_*$ and the various $R \in (0, \infty)$ versions of $\psi_R^* \mathcal{I}$ as $|x|^{N_{(p)}(0)} u_*$ with $u_*$ being a section of $\mathbb{V}^+$ over the complement of $Z_*$ in the $|x| = 1$ sphere.

*Part 6*: If $\varepsilon \in (0, 1)$ has been specified, then the fact that all large i versions of any given B and B´ version of $I_{i,BB'}$ are equal to $I_{BB'}$ implies that all sufficiently large i version of $\mathcal{I}_{T(i)}$ are isomorphic to $\mathcal{I}_*$ on the $|x| \in (\varepsilon, \varepsilon^{-1})$ part of the complement of $\mathcal{T}_\varepsilon$. It follows from what is said in Part 5 that there exists such an isomorphism that writes $\nu_{T(i)}$ as $\nu_* + w_i$ with $w_i$ being small when i is large in the sense that any given $k < \varepsilon^{-1}$ version $C^k$ norm less than $\varepsilon$ if i is large. This is what is asserted by the fifth bullet of Lemma 4.7. The two items in the fourth bullet of Lemma 4.7 follow from Lemmas 4.2 and 4.3.

**d) Locally unique limits**

Let $\{T(i) = (p_i, L_i, \lambda_i)\}_{i=1,2,\ldots}$ denote a sequence for use in Lemmas 4.2-4.5. Lemma 4.5 in particular takes this sequence as input and gives back a subsequence $\{T_{(i)}\}_{i \in \Lambda}$ with $\Lambda \subset \{1, 2, \ldots\}$ to which Lemma 4.7 then assigns a limit data set $(Z_*, \mathcal{I}_*, \nu_*)$. Keeping in mind that the only constraint on $\{\lambda_i\}_{i=1,2,\ldots}$ is for its limit to be zero, fix a new sequence $\{\lambda_i'\}_{i \in \Lambda} \subset (0, r_0)$ with limit zero and then use the sequence $\{T'(i) = (p_i, L_i, \lambda_i')\}_{i \in \Lambda}$ as input to Lemmas 4.2-4.5 after relabling its elements consecutively from 1. An assertion is made in the context of [T1] (by Lemmas 8.8 and 9.3 in the published version of [T1]) to the effect that Lemma 4.7's limit set for this new input sequence is invariably the same $(Z_*, \mathcal{I}_*, \nu_*)$. But, as noted by Thomas Walpuski, the proof of Lemmas 8.8 and 9.3 in the published version of [T1] has a fundamental error; and so the claim that there is but one limit data set is in doubt. The next lemma makes a weaker statement.

**Lemma 4.8**: *Let $\{T(i) = (p_i, L_i, \lambda_i)\}_{i=1,2,...}$ denote a sequence that obeys the assumptions in Proposition 4.1 and let $\{T(i)\}_{i \in \Lambda}$ denote a subsequence with $\Lambda \subset \{1, 2, ...\}$ that is then supplied by Lemma 4.5. Use $(Z_*, \mathcal{I}_*, \nu_*)$ to denote the resulting data set from Lemma 4.7. Suppose that $\{L_i'\}_{i \in \Lambda}$ is a sequence with its i'th member being an orthonormal frame for $TX$ at $p_i$. Assume that this sequence converges when viewed as a sequence of frames for $TX|_p$, and that the limit is the same as the limit of the sequence $\{L_i\}_{i \in \Lambda}$. Let $\{\lambda_i'\}_{i \in \Lambda}$ denote a sequence in $(0, r_0)$ with limit zero, and suppose that*

$$\liminf_{i \in \Lambda} |\ln(\lambda_i'/\lambda_i)| \, |N_{(p_i)}(\lambda_i) - N_{(p_i)}(\lambda_i')| = 0 \, .$$

*Then the conclusions of Lemmas 4.5-4.7 hold for $\{T(i)' = (p_i, L_i', \lambda_i')\}_{i \in \Lambda}$, and in particular, Lemma 4.7 produces from $\{T(i)'\}_{i \in \Lambda}$ the same output data set $(Z_*, \mathcal{I}_*, \nu_*)$.*

*Proof of Lemma 4.8*: The proof that the same output data set results when the frame sequence $\{L_i\}_{i \in \Lambda}$ is replaced by $\{L_i'\}_{i \in \Lambda}$ is straightforward and so it is left to the reader. This said, the three parts that follow prove that the assertion of the lemma holds in the case when $\{T(i)' = (p_i, L_i, \lambda_i')\}_{i \in \Lambda}$.

*Part 1*: Suppose that $i \in \Lambda$ and $\lambda \in (0, r_0)$. Define $T(i, \lambda)$ to be the data set $(p_i, L_i, \lambda)$ and define $\nu_{T(i,\lambda)}$ accordingly. Suppose to be that $x$ is a point in $\mathbb{R}^d$ with norm 1. If $x$ is not in $Z_{T(i)}$ then this is so for all points of the form $rx$ with $r$ sufficiently close to 1. With this understood, it follows that assignment $\lambda \to T(i, \lambda)$ for $\lambda$ near $\lambda_i$ is differentiable with the derivative given by the formula

$$\tfrac{\partial}{\partial \lambda} \nu_{T(i,\lambda)} = \lambda^{-1} (\nabla_r \nu_{T(i,\lambda)} - \tfrac{N_{(p_i)}(\lambda)}{|x|} \nu_{T(i,\lambda)}) \, .$$

(4.24)

Granted this formula, it then follows from what is said about the function in (4.2) that the function $|\tfrac{\partial}{\partial \lambda} \nu_{T(i,\lambda)}|^2$ on the $|x| = 1$ sphere in $\mathbb{R}^2$ is integrable for all but a measure zero set of $\lambda$ from the interval $(0, r_0)$. If $\lambda$ is not from this set, then

$$\int_{|x|=1} |\tfrac{\partial}{\partial \lambda} \nu_{T(i,\lambda)}|^2 = \tfrac{1}{\lambda^{d-2} K_{(p_i)}^2(\lambda)} \int_{\mathrm{dist}(p_i, \cdot) = \lambda} |\nabla_r \nu - \tfrac{N_{(p_i)}(\lambda)}{\lambda} \nu|^2 \, .$$

(4.25)

The second bullet of Lemma 3.2 can now be used to write

$$\int_{|x|=1} |\tfrac{\partial}{\partial \lambda} v_{T(i,\lambda)}|^2 = (\tfrac{1}{r} \tfrac{d}{dr} N_{(p_i)})|_{r=\lambda} + \mathfrak{e}_i ,$$

(4.26)

with $\mathfrak{e}_i$ being a function on $(0, r_0)$ with norm bounded by $c_0$.

*Part 2*: Suppose that $\lambda_1$ and $\lambda_2$ are both from $(0, r_0)$ with $\lambda_1 < \lambda_2$. Since

$$\int_{|x|=1} ||v_{T(i,\lambda_2)}| - |v_{T(i,\lambda_1)}||^2 \le (\int_{\lambda_1}^{\lambda_2} (\int_{|x|=1} |\tfrac{\partial}{\partial \lambda} v_{T(i,\lambda)}|^2)^{1/2} d\lambda )^2 ,$$

(4.27)

it follows from (4.26) that

$$\int_{|x|=1} ||v_{T(i,\lambda_1)}| - |v_{T(i,\lambda_2)}||^2 \le (\int_{\lambda_1}^{\lambda_2} (\tfrac{1}{r} |\tfrac{d}{dr} N_{(p_i)}|)^{1/2} dr )^2 + c_0 (\lambda_2 - \lambda_1).$$

(4.28)

The integral on the right hand side of (4.27) at most $\ln(\lambda_2/\lambda_1) |N_{(p_i)}(\lambda_2) - N_{(p_i)}(\lambda_1)|$ and thus

$$\int_{|x|=1} ||v_{T(i,\lambda_1)}| - |v_{T(i,\lambda_2)}||^2 \le \ln(\lambda_2/\lambda_1) |N_{(p_i)}(\lambda_2) - N_{(p_i)}(\lambda_1)| + c_0 (\lambda_2 - \lambda_1) .$$

(4.29)

The inequality in (4.28) can be generalized in the following way: Let $\mathfrak{s}$ for the moment denote a smooth section over the $|x| = 1$ sphere of the bundle $\mathbb{V}^+ \otimes \mathbb{V}^+$. Supposing that $\lambda \in (0, r_0)$, then the assignment $x \to \langle \mathfrak{s}, v_{T(i,\lambda)} \otimes v_{T(i,\lambda)} \rangle |v_{T(i,\lambda)}|^{-1}$ defines a continuous function on $|x| = 1$ sphere that is smooth on the complement of $Z_{T(i,\lambda)}$. With the preceding understood, let $\lambda_1$ and $\lambda_2$ be again from $(0, r_0)$ with $\lambda_1 < \lambda_2$. Then the argument that led to (4.29) can be repeated with only cosmetic changes to prove that

$$\int_{|x|=1} |\langle \mathfrak{s}, |v_{T(i,\lambda_2)}|^{-1} (v_{T(i,\lambda_2)} \otimes v_{T(i,\lambda_2)}) \rangle - \langle \mathfrak{s}, |v_{T(i,\lambda_1)}|^{-1} (v_{T(i,\lambda_1)} \otimes v_{T(i,\lambda_1)}) \rangle|^2 \le$$

$$c_{\mathfrak{s}} (\ln(\lambda_2/\lambda_1) |N_{(p_i)}(\lambda_2) - N_{(p_i)}(\lambda_1)| + (\lambda_2 - \lambda_1))$$

(4.30)

with $c_{\mathfrak{s}}$ denoting a positive number that depends only on the chosen section $\mathfrak{s}$.

*Part 3*: Use the sequence $\{T(i)' = (p_i, L_i, \lambda_i')\}_{i \in \Lambda}$ as input for Lemmas 4.2-4.5 after renumbering successively from 1. Let $\Lambda'$ denote the subsequence supplied by Lemma 4.5 and let $(Z_*', \mathcal{I}_*', \nu_*')$ denote the limit data set supplied by Lemma 4.7.

Suppose that $\varepsilon > 0$ and that $\iota_\varepsilon \geq 1$ is such that the $L^2$ norms on the $|x| = 1$ sphere of the function $|v_{T(i)}| - |v_*|$ is less than $\varepsilon$ when $i \in \Lambda'$ is greater than $\iota_\varepsilon$. Suppose that this is also true for $|v_{T(i)'}| - |v_*'|$ when $i' \in \Lambda$ and $i' > \iota_\varepsilon$. It follows from the assumptions of Lemma 4.8 that

$$|\ln(\lambda_i'/\lambda_i)| \, |N_{(p_i)}(\lambda_i') - N_{(p_i)}(\lambda_i)| < \varepsilon^2 .$$

(4.31)

Granted this bound, then (4.29) bounds the $L^2$ norm of $|v_{T(i)}| - |v_{T(i)'}|$ on the $|x| = 1$ sphere by $c_0 \varepsilon$ and this leads in turn to a $c_0 \varepsilon$ bound for the $L^2$ norm of $|v_*| - |v_*'|$ on the $|x| = 1$ sphere. Since $\varepsilon$ can be chosen at will subject to being positive, it follows that $|v_*| = |v_*'|$ and so $Z_* = Z_*'$. Much the same sort of argument using (4.30) proves that $v_* \otimes v_* = v_*' \otimes v_*'$. This last identification implies that $|v_*|^{-1}|v_*'|^{-1} \langle v_*, v_*' \rangle$ defines an isometric isomorphism between $\mathcal{I}_*$ and $\mathcal{I}_*'$ that identifies $v_*$ with $v_*'$.

## 5. The data sets on $\mathbb{R}^d$

This section describes the data sets $(Z_*, \mathcal{I}_*, v_*)$ from Proposition 4.1 in various special but important cases. To set the stage for this, suppose in this section that $N_0$ is a positive number and that $(Z_*, \mathcal{I}_*, v_*)$ is data set that is described in the upcoming (5.1). The notation has $\nabla_r$ denoting the directional derivative on $\mathbb{R}^d - 0$ along the outward pointing unit vector tangent to the rays from the origin.

- $Z_*$ is a scale invariant set in $\mathbb{R}^d$ and $\mathcal{I}_*$ is a real line bundle on $\mathbb{R}^d - Z_*$.
- $v_*$ is a smooth, $\mathcal{D}_0$-harmonic section of $\mathbb{V}^+ \otimes \mathcal{I}_*$ over $\mathbb{R}^d - Z_*$. Moreover
  a) $|v_*|$ extends to the whole of $\mathbb{R}^d$ as an $L^2_{1;loc}$ and Hölder continuous (for some fixed positive exponent) function with its zero locus being $Z_*$.
  b) $|\nabla v_*|$ extends to the whole of $\mathbb{R}^d$ as an $L^2_{loc}$ function.
  c) $\nabla_r v_* - \frac{N_0}{|x|} v_* = 0$.
  d) $\displaystyle\int_{|x|=1} |v_*|^2 = 1$.
  e) $\displaystyle\int_{|x| \leq 1} |\nabla v_*|^2 = N_0$.

(5.1)

A data set $(Z_*, \mathcal{I}_*, v_*)$ from Proposition 4.1 obeys (5.1) for a suitable choice of $N_0$.

Since the set $Z_*$ in (5.1) is scale invariant, there exists a canonical isomorphism between the bundle $\mathcal{I}_*$ and its pull-backs via the various $R \in (0, \infty)$ versions of the scaling diffeomorphism x → Rx. This fact is used implicitly in what follows to identify $\mathcal{I}_*$ with its restriction to the |x| = 1 sphere in $\mathbb{R}^d$. It then follows from Item c) of the third bullet in (5.1) that this idenfication of $\mathcal{I}_*$ writes

$$v_* = u_* |x|^{N_0}$$

(5.2)

with $u_*$ denoting a section of $\mathbb{V}^+ \otimes \mathcal{I}_*$ over the complement in the |x| = 1 sphere of the set of |x| = 1 points in $Z_*$. This set of |x| = 1 points in $Z_*$ is in some places denoted by $\mathcal{Z}_*$. A key fact in what follows is that the equation $\mathcal{D}_0 v_* = 0$ when written in terms of $u_*$ becomes an eigenvalue equation on the |x| = 1 sphere that has the schematic form

$$\mathfrak{S}_r^\dagger \mathcal{D}_0^\perp u_* + N_0 u_* = 0 \ .$$

(5.3)

What is denoted by $\mathcal{D}_0^\perp$ here is the part of $\mathcal{D}_0$ that involves only derivative tangent to the constant radii sphere centered at the orgin. Meanwhile, $\mathfrak{S}_r^\dagger$ is the adjoint of the homomorphism from $\mathbb{V}^+$ to $\mathbb{V}^-$ that is obtained by evaluating the symbol of $\mathcal{D}_0$ on the unit length vector field on $\mathbb{R}^d$–0 that is tangent to the rays from the origin.

a) **The 2-dimensional case**

The proof of Theorem 1.1 is given at the end of this subsection. The intervening three lemmas and some of what is said in their proofs are used to prove Theorem 1.1. The first of these lemmas says more about the data sets from (5.1) in the case when the dimension is 2.

**Lemma 5.1**: *Let $(Z_*, \mathcal{I}_*, v_*)$ denote a data set from the $\mathbb{R}^2$ version of Lemma 5.1. Then $Z_*$ is a point and $N_0$ is half of a positive integer. The bundle $\mathcal{I}_*$ is isomorphic to the product $\mathbb{R}$ bundle if and only if $N_0$ is a whole integer. Meanwhile, the section $v_*$ can be written as in (5.2) with $|u_*|$ being the constant function $\frac{1}{\sqrt{2\pi}}$.*

This lemma is proved momentarily. The next lemma reasserts a part of Theorem 1.1.

**Lemma 5.2**: *Supposing that X has dimension 2 and that Z, $\mathcal{I}$, and v are as described in (1.3) and (1.5). Then Z is a locally finite set; and if $p \in Z$, then $N_{(p)}(0)$ is half of a positive integer.*

This lemma is also proved momentarily.

The final lemma in this subsection is a uniqueness assertion for the limit data sets that arise in Lemma 4.7 that goes beyond what is asserted by Lemma 4.8. The context here is that X has dimension 2 and that $(Z, \mathcal{I}, \nu)$ is described by (1.3) and (1.5).

**Lemma 5.3**: *Suppose that* $p \in Z$, *that* $L$ *is an oriented, orthonormal frame for* $TX$ *at* $p$ *and that* $\{\lambda_i\}_{i=1,2,\ldots}$ *and* $\{\lambda_i'\}_{i=1,2,\ldots}$ *are two sequences in* $(0, r_0]$ *with limit zero. Use the respective sequences* $\{T(i) = (p, L, \lambda_i)\}_{i=1,2,\ldots}$ *and* $\{T(i)' = (p, L, \lambda_i')\}_{i=1,2,\ldots}$ *as input for Lemma 4.5, and then use corresponding subsequences* $\{T(i) = (p, L, \lambda_i)\}_{i \in \Lambda}$ *and* $\{T(i)' = (p, L, \lambda_i')\}_{i \in \Lambda}$ *as input for Lemma 4.7. Then the resulting limit data sets* $(Z_* = 0, \mathcal{I}_*, \nu_*)$ *and* $(Z_*' = 0, \mathcal{I}_*', \nu_*')$ *have the following property: There is an isometric isomorphism between* $\mathcal{I}_*$ *and* $\mathcal{I}_*'$ *on the complement of the origin identifying* $\nu_*$ *with* $\nu_*'$.

The remainder of this subsection has the proofs of these three lemmas and then the proof of Theorem 1.1

***Proof of Lemma 5.1***: Equation (5.2) writes $\nu_*$ as $u_* |x|^{N_0}$ with $u_*$ a section of $\mathbb{V}^+ \otimes \mathcal{I}_*$ over the $|x| = 1$ circle. The equation $\mathcal{D}_0 \nu_* = 0$ is obeyed if and only if $u_*$ obeys the $d = 2$ version of of (5.3). The $d = 2$ version of the equation in (5.3) is a first ordinary differential equation that can be written schematically as

$$\mathfrak{S}_\perp \nabla^\perp u_* + N_0 u_* = 0$$

(5.4)

with $\nabla^\perp$ being the covariant derivative tangent to the $|x| = 1$ circle and with $\mathfrak{S}_\perp$ being a constant, antisymmetric automorphism of $\mathbb{V}^+$ with square -1. Let $u_+$ and $u_-$ denote the respective orthogonal projections of $u_*$ to the $+i$ and $-i$ eigenspaces of $\mathfrak{S}_\perp$. Equation (5.4) leads to separate equations for $u_+$ and $u_-$. These can be solved directly with the result being as follows: Let $\theta$ denote a Euclidean $\mathbb{R}/2\pi\mathbb{Z}$ parameter for the $|x| = 1$ circle. Then $u_+ = u_{+0} e^{iN_0\theta}$ and $u_- = u_{-0} e^{-iN_0\theta}$ with $u_{+0}$ and $u_{-0}$ being constant.

The writing of $u_*$ as $u_{+0} e^{iN_0\theta} + u_{-0} e^{-iN_0\theta}$ implies that $|u_*|^2 = \frac{1}{2\pi}$, and since this is never zero, it follows that $Z_* = 0$. This depiction of $u_*$ can be used to see that $N_0$ is half an integer by writing $u_* \otimes u_*$ as the sum

$$u_{+0} \otimes u_{+0}\, e^{i2N_0\theta} + u_{-0} \otimes u_{-0}\, e^{-i2N_0\theta} + u_{+0} \otimes u_{-0} + u_{-0} \otimes u_{+0}.$$

(5.5)

and noting that each of the four terms in this sum must be a bonafide section of $\mathbb{V}^+ \otimes \mathbb{V}^+$ on the $|x| = 1$ circle. It also follows from the depiction of $u_*$ as as $u_{+0} e^{iN_0\theta} + u_{-0} e^{-iN_0\theta}$ that $\mathcal{I}_*$ is isomorphic to the product bundle if and only if $N_0$ is an integer.

*Proof of Lemma 5.2*: Suppose that p is a given point in Z. To see that p is isolated, suppose to the contrary that there exists a sequence $\{p_i\}_{i=1,2,\ldots} \subset Z$ with $\lim_{i \to \infty} p_i = p$ and with each point having distance at most $c_0^{-1}$ from p. For each $i \in \{1, 2, \ldots\}$, let $\lambda_i = \text{dist}(p_i, p)$. Supposing that L is a given oriened orthonormal frame for $TX|_p$, then sequence $\{T(i) = (p, L, \lambda_i)\}_{i=1,2,\ldots}$ can be used as input for Proposition 4.1. The fifth bullet of this proposition provides a subsequence $\{T(i) = (p, L, \lambda_i)\}_{i \in \Lambda}$ with $\Lambda$ being a subsequence in $\{1, 2,\ldots\}$ that converges in the manner that is described by the fifth bullet of the proposition. By construction each $i \in \Lambda$ version of $v_{T(i)}$ vanishes at a point on the $|x| = 1$ circle in $\mathbb{R}^2$, this being the point that corresponds via the composition of the rescaling map $x \to \lambda_i x$ and the Gaussian coordinate chart map to the point $p_i$. It therefore follows from the manner of convergence that $v_*$ must vanish at a point on the $|x| = 1$ circle. This conclusion runs afoul of Lemma 5.1.

The assertion in Lemma 5.2 to the effect that each $p \in Z$ version of $N_{(p)}(0)$ is half of a positive integer follows directly from Proposition 4.1 and Lemma 5.1's assertion that the number $N_0$ in any given version of (5.1) is half a positive integer.

*Proof of Lemma 5.3*: The proof has four steps.

 Step 1: Suppose that $p \in Z$, that L is a chosen oriented, orthonormal frame for $TX|_p$ and that $\lambda$ and $\lambda'$ are any two numbers in $(0, r_0)$. Let T denote the data set $(p, L, \lambda)$ and let T' denote $(p, L, \lambda')$. Since p is an isolated point of Z, it follows that if i is any positive integer, then the corresponding versions of $\mathcal{I}_T$ and $\mathcal{I}_{T'}$ must be isomorphic on the complement of the origin in the $|x| \leq 2$ ball. Taking $\lambda = \lambda_i$ and $\lambda' = \lambda_i'$ for any given index $i \in \{1, 2, \ldots\}$ proves that $\mathcal{I}_{T(i)}$ is isomorphic $\mathcal{I}_{T(i)'}$ where $0 < |x| \leq 2$ ball. Meanwhile, taking $T = T(i)$ and $T' = T(i')$ for any pair of indices i and i' proves that the limit bundle $\mathcal{I}_*$ is isomorphic to any large i version of $\mathcal{I}_{T(i)}$ where $0 < |x| \leq 2$. By the same token, the limit bundle $\mathcal{I}_*'$ is isomorphic to any large i version of $\mathcal{I}_{T(i)'}$ where $0 < |x| \leq 2$. It therefore follows that $\mathcal{I}_*$ is isomorphic to $\mathcal{I}_*'$ on $\mathbb{R}^2 - \{0\}$.

 Step 2: Given $r \in (0, r_0)$, set $T(r) = (p, L, \lambda = r)$. Use the isomorphisms from Step 1 to define a function on $(0, r_0)$ to be denoted by Q which is defined by the rule

$$r \to Q(r) = \int_{|x|=1} \langle v_*, v_{T(r)} \rangle \ .$$
(5.6)

Note that $|Q(r)| \le 1 + c_0 r^2$ because $1 - Q(r)$ differs by $c_0 r^2$ from half of the $L^2$ norm on the $|x| = 1$ circle of $|v_* - v_{T(r)}|$. This last observation implies that the isomorphisms between $\mathcal{I}_{T(r)}$ and $\mathcal{I}_*$ can be chosen in a consistant manner so as to have $\lim_{r \to 0} Q(r) = 1$.

Step 3: The function $Q(r)$ is differentiable and it follows from (4.23) and (5.3) that its derivative obeys

$$\tfrac{d}{dr} Q = -\tfrac{1}{r} (N_{(p)}(r) - N_{(p)}(0)) Q + \mathfrak{z} \ ,$$
(5.7)

with the absolute value of $\mathfrak{z}$ being less than $c_0$. To exploit this equation, suppose for the moment that $r \in (0, c_0^{-1}]$ and that $Q(r) > \tfrac{1}{2}$. Fix $r_* \in [0, r)$ such that $Q > \tfrac{1}{4}$ on $(r_*, r]$. What with the third bullet of Lemma 3.2 finding $N_{(p)}(s) - N_{(p)}(0) \ge -c_0 s^2$ for any $s \in (0, c_0^{-1}]$, integration of (5.7) leads to an inequality asserting that

$$Q(s) \ge Q(r) - c_0 r \ .$$
(5.8)

if $s \in (r_*, r)$. It follows from this that if $r < c_0^{-1}$, then $r_* = 0$ and so (5.8) holds on the whole of $(0, r]$.

Step 4: The inequality in (5.8) says that if $r$ is small and $Q(r)$ is very nearly equal to 1, then $Q$ on $(0, r]$ is also nearly equal to 1. In particular, if $\varepsilon \in (0, 1]$ is given and $r$ is a sufficiently small number from $(0, r_0)$, then $Q(r) > 1 - \varepsilon$. Fix such a value for $r$ from the interval $(0, c_0^{-1}\varepsilon]$. Then (5.8) asserts that $Q > 1 - 2\varepsilon$ on $(0, r]$. This implies that the $L^2$ norm of $|v_* - v_{T(s)}|$ on the $|x| \le 1$ sphere is at most $4\varepsilon$ if $s \in (0, r]$. Since $\varepsilon$ can be as small as desired (if positive), this last bound implies what is asserted Lemma 5.3 because it leads to the observation that $\lim_{i \to \infty} |v_* - v_{T(i)'}| = 0$ as can be seen by taking $s$ to be ever larger index $i$ versions of $\lambda_i'$.

*Proof of Theorem 1.1*: Lemma 5.2 reasserts the theorem's claim that the zero locus of $|v|$ is locally finite. Supposing that $p \in X$ is a zero of $|v|$, the theorem also characterizes the behavior of $|v|$ in a small disk centered at p. To see about this characterization of $|v|$, reintroduce the function Q that is defined in (5.6). Step 2 of the proof of Lemma 5.3 explains why Q is greater than $\tfrac{1}{2}$ on an interval of the form $(0, r_*]$ with $r_* \in (0, r_0]$. Since Q is in any event bounded by $1 + c_0 r_*$, and since Lemma 3.2 finds $N_{(p)}(r) > N_{(p)}(0) - c_0 r^2$ it follows from (5.7) that the function defined by the rule $r \to r^{-1}|N_{(p)}(r) - N_{(p)}(0)|$ is integrable on $(0, r_0]$. This last fact and the first bullet of Lemma 3.2 imply in turn that p's version of the function $\kappa$ obeys $\kappa_{(p)}(r) = \kappa_* r^{N_{(p)}(0)} (1 + \mathfrak{e})$ with $\kappa_*$ being a positive constant and with $|\mathfrak{e}|$ having limit zero as $r \to 0$. Theorem 1.1's characterization of $|v|$ follows from this last bound given what Lemma 5.3 says about unique limits, given the manner of convergence dictated by Proposition 4.1 and given what is said in Lemma 5.1 to the effect that $v_*$ can be written as $u_* |x|^{N_{(p)}(0)}$ with $|u_*|$ being constant.

**b) Observations for the d = 3, 4 cases**

The subsequent analysis when $d = 3$ and $d = 4$ make fundamental use of the fact that the data consisting of $Z_*$, $\mathcal{I}_*$ and $v_*$ are described by the $X = \mathbb{R}^d$ version of (1.3) and (1.5). It follows as a consequence that what is said in Sections 2-4 can be invoked to say things about $(Z_*, \mathcal{I}_*, v_*)$.

To set the stage for the applications to come of the $X = \mathbb{R}^d$ versions Sections 2-4, suppose for the moment that $q \in Z_*-0$. follow through on this, suppose for the moment that $q \in Z_*-0$. Let $\kappa_q$ denote the positive function on $(0,\infty)$ whose square is given by the $d = 3$ version of the rule

$$r \to \kappa_q^2(r) = r^{-(d-1)} \int_{|(\cdot)-q|=r} |v_*|^2 \ .$$

(5.9)

The scale invariance implies that $\kappa_{tq}(tr) = \kappa_q(r)$ for any $t \in (0, \infty)$. With $\kappa_q$ in hand, define the function $N_q$ where $\kappa_q > 0$ by the rule:

$$r \to N_q(r) = \frac{2}{r^{d-2}\kappa_q(r)^2} \int_{|(\cdot)-q|\le r} |\nabla v_*|^2 \ .$$

(5.10)

This function obeys $N_{tq}(tr) = N_q(r)$ for any $t \in (0, \infty)$.

Let $r_q$ denote for the function on $\mathbb{R}^d$ giving the distance to the point q. The arguments used for the first and second bullets of Lemma 3.2 can be repeated verbatim to derive the formulas that follow for the derivatives of $\kappa_q$ and $N_q$ where $\kappa_q$ is positive.

- $\frac{d}{dr} \kappa_q = \frac{1}{r} N_q \kappa_q \ .$
- *The function $N_q$ is continuous and almost everywhere differentiable with derivative*

$$\frac{d}{dr} N_q = \frac{2}{r^{d-2}\kappa_q(r)^2} \int_{|(\cdot)-q|=r} |\nabla_{r_q} v_* - \frac{N_q}{r_q} v_*|^2 \ .$$

(5.11)

Granted these formula, then the other arguments for Lemma 3.2 can be repeated to prove that $K_q$ is strictly positive on $(0, \infty)$ and that $N_q$ extends to $[0, \infty)$ as a continuous function. The notation in what follows uses $\tilde{N}_q$ to denote the $r \to 0$ limit of $N_q(r)$. Note that this function on $Z_*$ given by the rule function $q \to \tilde{N}_q$ on $Z_*$ is scale invariant in the sense that $\tilde{N}_{tq} = \tilde{N}_q$ if $q \in Z_*$ and $t \in (0, \infty)$.

There are corresponding versions of Proposition 4.1 for points in $Z_*-0$. To say more about these, suppose that $q$ is such a point. By virtue of the scale invariance, nothing is lost by taking $q$ to have norm 1. Let $L$ denote the oriented orthonormal frame at $q$ that is given by the coordinate basis vectors for $\mathbb{R}^d$. Supposing that $\{\lambda_i\}_{i \in \{1,2...\}}$ is a sequence of positive number with limit zero, the $X = \mathbb{R}^d$ version of Proposition 4.1 can be invoked using as input the sequence $\{T(i) = (q, L, \lambda = \lambda_i)\}_{i \in \{1,2,...\}}$. Let $(Z_{*q}, \mathcal{I}_{*q}, v_{*q})$ denote a limit data set that is supplied by Proposition 4.1. The section $v_{*q}$ of $\mathbb{V}^+ \otimes \mathcal{I}_{*q}$ obeys the three equations

$$\mathcal{D}_0 v_{*q} = 0 \quad and \quad \nabla_r v_{*q} - \frac{\tilde{N}_q}{|x|} v_* = 0 \quad and \quad \nabla_q v_{*q} = 0 .$$

(5.12)

By way of an explanation, the left most equation is the assertion from the second bullet of the $(Z_{*q}, \mathcal{I}_{*q}, v_{*q})$ version of (5.1) and the middle equation is the assertion from Item c) of the second bullet of this same version of (5.1). The right most equation follows from the manner of convergence dictated by Proposition 4.1 given that $v_*$ obeys Item c) of the $(Z_*, \mathcal{I}_*, v_*)$ version of (5.1). Note in this regard that if $\lambda \in (0, \infty)$ and $T = (q, L, \lambda)$, then the map $\phi_T$ that is described at the beginning of Section 4 is the map $x \to \lambda x + q$.

A $(Z_*, \mathcal{I}_*, v_*)$ version of the identity in (4.26) is needed to prove the upcoming Lemma 5.7. To state this version, fix $q \in \mathbb{R}^d - 0$ and $\lambda \in (0, \infty)$; and then define $T(\lambda)$ to be the data set $(q, L, \lambda)$. Let $v_{*T(\lambda)}$ denote $K_q(\lambda)^{-1} \phi_{T(\lambda)}^* v_*$. This is the $v_*$ and $T = T(\lambda)$ version of what is denoted by $v_T$ in Section 4. The corresponding version of (4.26) asserts that

$$\int_{|x|=1} |\tfrac{\partial}{\partial \lambda} v_{*T(\lambda)}|^2 = (\tfrac{1}{r} \tfrac{d}{dr} N_q)|_{r=\lambda} .$$

(5.13)

This equation follows from the $p_i = q$ and $L_i = L$ version of (4.25) and the identity in the second bullet in (5.11)

### c) The case d = 3

Suppose in what follows that $(Z_*, \mathcal{I}_*, v_*)$ are described by some version of (5.1) in the case when the dimension is 3. The next lemma says more about this sort of data set.

**Lemma 5.4**: *Let $(Z_*, \mathcal{I}_*, \nu_*)$ denote a data set as described by the $d = 3$ version of (5.1). The set $Z_*$ is the union of the origin with the cone on a finite set of points in the $|x| = 1$ sphere. Let $\mathcal{Z}_*$ denote this set of $|x| = 1$ points in $Z_*$. If $q \in \mathcal{Z}_*$ then the number $\tilde{N}_q$ is half of a positive integer and it has the following significance:*

- *Write $\nu_* = u_* |x|^{N_0}$ as in (5.2) with $u_*$ being a section of $\mathbb{V}^+ \otimes \mathcal{I}_*$ on the $|x| = 1$ sphere. Supposing that $q \in \mathcal{Z}_*$, then $|u_*|$ near $q$ is bounded by $\kappa_q |q - (\cdot)|^{\tilde{N}_q}$ with $\kappa_q$ being a number that depends only on $q$.*
- *Let $q \in \mathcal{Z}_*$ and let $D$ denote a closed disk in the $|x| = 1$ sphere whose boundary is disjoint from $\mathcal{Z}_*$ and whose only $\mathcal{Z}_*$ point is $q$. The restriction of $\mathcal{I}_*$ to $D-q$ is isomorphic to the product $\mathbb{R}$ bundle if and only if $\tilde{N}_q$ is an integer.*
- *Each $q \in \mathcal{Z}_*$ version of $\tilde{N}_q$ obeys $\tilde{N}_q \leq N_0$.*
- *If there exists $q \in \mathcal{Z}_*$ with $\tilde{N}_q = N_0$, then $Z_*$ is the line through $0$ and $q$. The bundle $\mathcal{I}_*$ in this case is canonically isomorphic to its pull-back via any translation parallel to $Z_*$; and this isomorphism identifies $\nu_*$ with its pull-back by the translation.*

*Proof of Lemma 5.4*: Let $q$ denote an $|x| = 1$ point in $Z_*$ and let $L$ denote the oriented orthonormal frame at $q$ that is given by the coordinate basis vectors for $\mathbb{R}^3$. Supposing that $\{\lambda_i\}_{i \in \{1,2,...\}}$ is a sequence of positive number with limit zero, invoke the $X = \mathbb{R}^3$ version of Proposition 4.1 using as input the sequence $\{T(i) = (q, L, \lambda = \lambda_i)\}_{i \in \{1,2,...\}}$. Let $(Z_{*q}, \mathcal{I}_{*q}, \nu_{*q})$ denote a limit data set that is supplied by Proposition 4.1 and let $\Lambda$ denote a subsequence in $\{1, 2, ...\}$ such that $\{T(i)\}_{i \in \Lambda}$ converges in the manner that is described by Proposition 4.1. Let $(x_1, x_2)$ denote Euclidean coordinates for the 2-dimensional plane through the origin $\mathbb{R}^3$ that is orthogonal to $q$. Use $z$ to denote $x_1 + i x_2$. The middle equation in (5.12) and the right most equation in (5.12) imply that $\nu_{*q}$ can be written as $u_{*q} |z|^{\tilde{N}_q}$ with $u_{*q}$ being a section of $\mathbb{V}^+ \otimes \mathcal{I}_*$ on the complement of the $|z| = 1$ points in the $(x_1, x_2)$ plane. Let $\nabla^\perp$ denote the covariant derivative on the unit radius circle about the origin in this same plane. The left most equation when written using $u_*$ is the equation that is depicted in (5.4). This being the case, the proof of Lemma 5.1 can be repeated with only notational changes to see that $u_{*q}$ is nowhere and that $\tilde{N}_q$ is half a positive integer. Given that $u_{*q}$ is nowhere zero, the argument for Lemma 5.2 can be reused with $u_*$ in (5.2) playing the role of $\nu$ to prove that $Z_*$ has but a finite number of $|x| = 1$ points.

The three parts that follow prove the bulleted items of Lemma 5.4

*Part 1*: The assertion made by the first bullet of Lemma 5.4 restates a part of the $(Z = Z_*, \mathcal{I} = \mathcal{I}_*, \nu = \nu_*)$ version of Lemma 3.3. To prove the lemma's second bullet, fix an $|x| = 1$ point $q \in Z_*$ and reintroduce $Z_{*q}$, $\mathcal{I}_{*q}$ and $v_{*q}$. Letting $z = x_1 + ix_2$ denote the complex coordinate in the plane orthogonal to q, write $v_{*q}$ again as $u_{*q} |z|^{\tilde{N}_q}$. The fact that $u_{*q}$ is nowhere zero implies that $Z_{*q}$ is the line in $\mathbb{R}^d$ through the vector q and the origin. Meanwhile, the $u_{*q}$ version of (5.5) indicates that $\mathcal{I}_{*q}$ is isomorphic to the product $\mathbb{R}$ bundle on the complement of this line if and only if $\tilde{N}_q$ is an integer. The assertion of the second bullet of Lemma 5.4 follows from this observation about $\mathcal{I}_{*q}$ and $\tilde{N}_q$ because of the manner of convergence that is described in Proposition 4.1.

*Part 2*: This part proves the third bullet of Lemma 5.4. To start, fix an $|x| = 1$ point q in $Z_*$ and reintroduce the function $N_q$ from (5.10). Keep in mind that the number $\tilde{N}_q$ is the $r \to 0$ limit of $N_q$. As explained directly, $\lim_{r \to \infty} N_q = N_0$. Since the second bullet of (5.11) says that $N_q$ is non-decreasing, the fact that $N_q(0) = \tilde{N}_q$ and that $\lim_{r \to \infty} N_q = N_0$ require that $\tilde{N}_q \leq N_0$.

To prove that $\lim_{r \to \infty} N_q = N_0$, use Items c) and d) of the third bullet in (5.1) to see that the integral of $|v_*|^2$ on a ball of radius R centered at the origin in equal to $\frac{1}{3+N_0} R^{3+2N_0}$. Meanwhile, if $r \in (1, \infty)$, then the ball of radius r centered at q is contained in the ball of radius r+1 centered at the origin and it contains the ball of radius r - 1 centered at the origin. This implies that the integral of $|v_*|^2$ on the radius r ball centered at q is no less than $\frac{1}{3+2N_0} r^{3+2N_0} (1 - c_0 r^{-1})$ and no greater than $\frac{1}{3+2N_0} r^{3+2N_0} (1 + c_0 r^{-1})$. Since the second bullet of (5.11) asserts that $N_q$ is non-decreasing, these bounds on the integral of $|v_*|^2$ are not compatible with those that come by integrating the equation in the first bullet of (5.11) unless $\lim_{r \to \infty} N_q(r) = N_0$.

*Part 3*: This part proves the fourth bullet of Lemma 5.4. To this end, suppose that q is an $|x| = q$ point with $\tilde{N}_q = N_0$. Since $\lim_{r \to 0} N_q = \tilde{N}_q$ which is $N_0$ in this case, and since $\lim_{r \to \infty} N_q = N_0$ in any case, the function $N_q$ must be constant because the the second bullet of (5.11) implies that $N_q$ is non-decreasing. If $N_q$ is constant then the integrand in the second bullet of (5.11) must be zero and so

$$\nabla_{r_q} v_* - \frac{N_0}{r_q} v_* = 0.$$

(5.14)

If (5.10) holds, then $|v_*|$ can be written as $|v_*|(x) = w_* |x - q|^{N_0}$ with $w_*$ being a function on the radius 1 sphere centered at q. It follows as a consequence that $Z_*$ must be the union of q with a cone on some set in this $|x - q| = 1$ sphere. The only set that contains the ray from the origin through q and is *both* a cone on a finite set in the $|x| = 1$ sphere and a cone on a set in the $|x - q| = 1$ sphere is the line through the origin and q.

Multipling both sides of (5.14) by $r_q$ leads to an equation that reads

$$|x|\nabla_r v_* - \nabla_q v_* - N_0 v_* = 0 ,$$

(5.15)

with $\nabla_q$ denoting the directional derivative along the unit vector q. This equation runs afoul of the equation in Item c) of the third bullet in (5.1) unless $\nabla_q v_* = 0$. The fact that $\nabla_q v_* = 0$ also implies the assertions in the fourth bullet of the lemma to the effect that $Z_*$ is the line through the origin and the point q.

c) **The case d = 4**

The upcoming Lemma 5.5 asserts a dimension 4 analog of Lemma 5.4. This lemma introduces the notion of a *half-plane* through the origin in $\mathbb{R}^4$. A half-plane through the origin in $\mathbb{R}^4$ is defined by a pair consisting of a 2-dimensional linear subspace in $\mathbb{R}^4$ and a non-zero, linear functional on this subspace. Letting $V \subset \mathbb{R}^4$ denote the 2-dimensional subspace and $f$ the linear function, the corresponding half-plane is the $f \geq 0$ part of V. The $f = 0$ part of V is said to be the *boundary line* of the half-plane.

**Lemma 5.5**: *Let $(Z_*, \mathcal{I}_*, v_*)$ denote a data set from the dimension 4 version of (5.1). Use $\mathcal{Z}_*$ to denote the set of $|x| = 1$ points in $Z_*$. The function $\tilde{N}_{(\cdot)}$ on $\mathcal{Z}_*$ is no less than $\frac{1}{2}$ and it has the following properties:*

- *Write $v_* = u_* |x|^{N_0}$ as in (5.2) with $u_*$ being a section of $\mathbb{V}^+ \otimes \mathcal{I}_*$ on the $|x| = 1$ sphere. Fix $q \in \mathcal{Z}_*$ and let $\tilde{N}_q$ denote value of $\tilde{N}_{(\cdot)}$ at q. Then $|u_*|$ near the point q is bounded by $\kappa_q |q - (\cdot)|^{\tilde{N}_q}$ with $\kappa_q$ being a number that depends only on q.*
- *Each $q \in \mathcal{Z}_*$ version of $\tilde{N}_q$ obeys $\tilde{N}_q \leq N_0$.*
- *If there exists $q \in \mathcal{Z}_*$ with $\tilde{N}_q = N_0$, then $Z_*$ is the union of the origin and either a line through the origin; or a finite number of half-planes through the origin with common boundary line, this boundary line being the line through the origin and q.*
  a) *If $Z_*$ is a line through the origin, then the function $\tilde{N}_{(\cdot)}$ is equal to $N_0$ on the two $|x| = 1$ points on this line. The bundle $\mathcal{I}_*$ in this case is isomorphic to the product bundle on $\mathbb{R}^4 - Z_*$ and the number $N_0$ is an integer.*
  b) *If $Z_*$ is a non-empty union of half-planes through the origin, then $\tilde{N}_{(\cdot)} = N_0$ on the two $|x| = 1$ points on their common boundary line. The function $\tilde{N}_{(\cdot)}$ is also locally constant on the complement in $\mathcal{Z}_*$ of these two points and its values on the complement of these two points are half-integers.*
  c) *Suppose that $Z_*$ is a non-empty union of half-planes through the origin and that $q' \in \mathcal{Z}_*$ is in the interior of one such half-plane. Let D denote a closed disk in the $|x| = 1$ sphere whose boundary is disjoint from $\mathcal{Z}_*$ and whose interior intersects $\mathcal{Z}_*$* *transversally at the point* q´ *with* q´ *being its only* $\mathcal{Z}_*$ *point. The bundle* $\mathcal{I}_*$ *on* D−q´ *is isomorphic to the product* $\mathbb{R}$ *bundle if and only if* $Ñ_{(q´)}$ *is an integer.*

d) *The bundle* $\mathcal{I}_*$ *is canonically isomorphic to its pull-back via any translation of* $\mathbb{R}^4$ *by a constant multiple of the vector* q. *This isomorphism identifies* $\nu_*$ *with its pull-back by the translation.*

- *If* $Z_*$ *is a non-empty union of half-planes through the origin and if* $Ñ_{(·)} = N_0$ *on a point in* $\mathcal{Z}_*$ *that is not on their common boundary line, then* $Z_*$ *is a plane through the origin, the function* $Ñ_{(·)} = N_0$ *on the whole of* $\mathcal{Z}_*$. *The number* $N_0$ *is thus half an integer. The bundle* $\mathcal{I}_*$ *in this case is isomorphic to the product* $\mathbb{R}$ *bundle on* $\mathbb{R}^4 − Z_*$ *if and only if* $N_0$ *is an integer.*

*Proof of Lemma 5.5*:  The proof has three parts.

*Part 1*: This part of the proof derives the bound $Ñ_q \geq \frac{1}{2}$. To start, let q denote an $|x| = 1$ point in $Z_*$ and let L denote the oriented orthonormal frame at q that is given by the coordinate basis vectors for $\mathbb{R}^3$. Supposing that $\{\lambda_i\}_{i \in \{1,2...\}}$ is a sequence of positive number with limit zero, invoke the $X = \mathbb{R}^3$ version of Proposition 4.1 using as input the sequence $\{T(i) = (q, L, \lambda = \lambda_i)\}_{i \in \{1,2,...\}}$. Let $(Z_{*q}, \mathcal{I}_{*q}, \nu_{*q})$ denote a limit data set that is supplied by Proposition 4.1. It follows from the right most equation in (5.12) that that $Z_{*q}$ is invariant under translations in the direction of the vector q. This implies that $Z_{*q}$ is foliated by half-planes with a common boundary line, this boundary line being the line through q and the origin in $\mathbb{R}^d$.

Let $V_q$ denote the 3-dimensional hyperplane orthogonal to the vector q. The right most equation in (5.12) implies that $\nu_{*q}$ is determined by its restriction $V_q$ and that it obeys the left and middle equations in (5.12) on $V_q$. This being the case, the restriction of the data set $(Z_{*q}, \mathcal{I}_{*q}, \nu_{*q})$ to $V_q$ obeys the $X = \mathbb{R}^3$ version of (5.1). It follows as a consequence that Lemma 5.4 can be invoked to say more about $(Z_{*q}, \mathcal{I}_{*q}, \nu_{*q})$. In particular, this lemma says first that the restriction of $Z_{*q}$ to $V_q$ is the union of the origin and a finite number of rays through the origin. It associates to each such ray a positive half integer which is no greater than $Ñ_q$, this being the $(Z_{*q}, \mathcal{I}_{*q}, \nu_{*q})$ version of $N_0$. The bound $Ñ_q \geq \frac{1}{2}$ follows from the third bullet of the $(Z_{*q}, \mathcal{I}_{*q}, \nu_{*q})$ version of Lemma 5.4 because the smallest positive half integer is $\frac{1}{2}$.

*Part 2*: The first bullet of Lemma 5.5 follows directly from the $(Z_*, \mathcal{I}_*, \nu_*)$ version of Lemma 3.3. The proof of the second bullet of the lemma is the same but for notation as the proof of the third bullet of Lemma 5.4. To prove the lemma's third bullet, assume henceforth that q is a norm 1 point in $Z_*$ with $\tilde{N}_q = N_0$. If this is so, then the integrand on the right hand side of the second bullet in (5.11) must vanish and so the $d = 4$ versions of (5.13) and (5.14) must hold. The latter equation and the equation in Item c) of the third bullet in (5.1) require that $\nabla_q v_* = 0$ also; and if $\nabla_q v_* = 0$, then $Z_*$ is mapped to itself by any translation of $\mathbb{R}^4$ of the form $x \to x + tq$ with t being constant. The vanishing of $\nabla_q v_*$ also implies that $\mathcal{I}_*$ is canonically isomorphic to its pull-back by these sorts of translations and that this isomorphism identifies $v_*$ with its corresponding pull-back. This is the content of Item d) of Lemma 5.5's third bullet. Note that these isomorphisms between $\mathcal{I}_*$ and its various pull-backs are used implicitly in the subsequent steps of the proof.

Let $V_q \subset \mathbb{R}^4$ denote the orthogonal subspace to q, this being a 3-dimensional linear subspace. The fact that $\nabla_q v_* = 0$ implies that the restriction of $v_*$ to $V_q$ obeys the $d = 3$ version of (5.1). This understood, then what is said by Lemma 5.4 applies, and what is said in this lemma leads directly to Items a)-d) of the third bullet of Lemma 5.4

*Part 3*: Suppose that q and q´ are two points in $\mathcal{Z}_*$ with $\tilde{N}_q = \tilde{N}_{q'} = N_0$ and such that q and q´ are not on a line through the origin. The arguments in Part 2 can be invoked using either q or q´ to see that both $\nabla_q v_* = 0$ and $\nabla_{q'} v_* = 0$. This implies that $Z_*$ is mapped to itself by translations of the form $x \to x + tq + sq'$ for any pair $(t, s) \in \mathbb{R}^2$. Given that $Z_*$ is a union of half planes with a common boundary line, this invariance requires that $Z_*$ be the plane spanned by the vectors q and q´. In fact, the fourth bullet of Lemma 5.4 can be invoked for $v_*$ on $V_q$ to see that $v_*$ on $V_q$ is annihilated by the directional covariant derivative in the direction of the projection to $V_q$ of q´. This appeal to Lemma 5.4 says that $V_q \cap Z_*$ is a line through the origin parallel to this projection of q´ and that $\tilde{N}_{(\cdot)}$ is equal to $N_0$ at each point of $V_q \cap Z_*$. These two facts imply what is asserted by the fourth bullet of Lemma 5.5 because $\nabla_q v_* = 0$.

### d) A compactness lemma

The upcoming Lemma 5.6 asserts that any sequence of data sets that obey (5.1) with an upper bound for the various versions of the number $N_0$ has a subsequence that converges in a sense to a data set that is also described by (5.1).

**Lemma 5.6**: *Suppose that $\{(Z_{*i}, \mathcal{I}_{*i}, v_{*i})\}_{i=1,2,...}$ is a sequence of data sets that are described by (5.1) with an upper bound for the various versions of $N_0$. There exists a data set $(Z_*, \mathcal{I}_*, v_*)$ that is described by (5.1) and there is a subsequence $\Lambda \in \{1, 2, ....\}$ such that the following is true: Fix $\varepsilon \in (0, 1]$. If $i \in \Lambda$ is sufficiently large, then there is an isometric isomorphism between $\mathcal{I}_{*i}$ and $\mathcal{I}_*$ on the complement in the $|x| < \varepsilon^{-1}$ ball of the radius $\varepsilon$ tubular neighborhood of $Z_*$ that allows $v_{*i}$ to be viewed as a section of $\mathbb{V}^+ \otimes \mathcal{I}_*$ on* *this same set. Viewed in this way, $v_{*i}$ differs from $v_*$ on the complement in the $|x| < \varepsilon^{-1}$ ball of the radius $\varepsilon$ tubular neighborhood of $Z_*$ by a section of $\mathbb{V}^+ \otimes \mathcal{I}_*$ with $C^k$ norm less than $\varepsilon$ for all integers $k < \varepsilon^{-1}$.*

*Proof of Lemma 5.6*: The argument is a somewhat simpler version of that used for Proposition 4.1 because one can write each $i \in \{1, 2, \ldots\}$ version of $v_{*i}$ as $u_{*i} |x|^{N_{0,i}}$ with $u_{*i}$ being a section of $\mathbb{V}^+ \otimes \mathcal{I}_{*i}$ on the $|x| = 1$ sphere and with $N_{0,i}$ being the relevant version of the number $N_0$ that appears in (5.1). Items d) and e) in (5.1) lead to a priori $L^2$ bounds on the $|x| = 1$ sphere for $\{u_{*i}\}_{i=1,2,\ldots}$ and for $\{\nabla u_i\}_{i=1,2,\ldots}$. Fix a subsequence $\Xi \subset \{1, 2, \ldots\}$ such that $|u_{*i}|$ and $|\nabla u_{*i}|$ converge weakly in $L^2$ topology. Meanwhile, each $i \in \{1, 2, \ldots\}$ version of $u_{*i}$ obeys the $N_0 = N_{0,i}$ version of (5.3). The latter equation has a corresponding Bochner-Weitzenboch formula that can be used to obtain an a priori $L^\infty$ bound for the sequence $\{u_{*i} \otimes u_{*i}\}_{i=1,2,\ldots}$. The derivation of this bound uses much the same arguments that are used in Section 2b. Given the lower bound of $\frac{1}{2}$ for $N_0$ from Lemmas 5.1, 5.4 and 5.5, the arguments in Section 3e can be used to obtain an a priori bound on some fixed exponent Hölder norm for the sequence $\{u_{*i} \otimes u_{*i}\}_{i=1,2,\ldots}$. These bounds can then be used to find a subsequence $\Lambda \subset \Xi$ such that the sequence $\{u_{*i} \otimes u_{*i}\}_{i \in \Lambda}$ converges in an exponent $\upsilon > 0$ Hölder space of sections of $\mathbb{V}^+ \otimes \mathbb{V}^+$ on the $|x| = 1$ sphere. This convergence can then be used as in Section 4c to define $(Z_*, \mathcal{I}_*, v_*)$.

### e) An isolation lemma

Data sets from (5.1) where $\tilde{N}_{(\cdot)}$ is constant on $Z_*$ play a central role in subsequent sections. The lemma that follows states a sense in which these data sets are isolated.

**Lemma 5.7**: *Given a positive, half-integer $N_\diamond$, there exists $\kappa > 4$ with the following significance: Let $(Z_*, \mathcal{I}_*, v_*)$ denote a data set from (5.1) obeying the conditions below.*
- *The number $N_0$ from (5.1) is less than $N_\diamond + \kappa^{-1}$.*
- *Let q denote an $|x| = 1$ point from $Z_*$. If $\mathcal{I}_*$ is not isomorphic to the product $\mathbb{R}$ bundle on the complement of $Z_*$ in some neighborhood of q, then $\tilde{N}_q \geq N_\diamond$.*

*The number $N_0$ is equal to $N_\diamond$ and the data set $(Z_*, \mathcal{I}_*, v_*)$ is described by either the fourth bullet of Lemma 5.4 or the third and fourth bullets of Lemma 5.5 as the case may be.*

The lemma is proved momentarily.

To see an example of a data set $(Z_*, \mathcal{I}_*, v_*)$ obeying either the fourth bullet of Lemma 5.4 or the third and fourth bullets of Lemma 5.5, let $Z_*$ denote a codimension 2 linear subspace in $\mathbb{R}^d$ and let $\mathcal{I}_*$ denote the product $\mathbb{R}$ bundle on the compliment of this subspace in $\mathbb{R}^d$ if $N_◊$ is an integer or the Möbius bundle line bundle on the complement of this subspace if $N_◊$ is half an odd integer. To exhibit $v_*$, let $(x_1, x_2)$ denote Euclidean coordinates for the orthogonal 2-dimensional vector subspace to $Z_*$. Set $z = x_1 + ix_2$. Let $\mathfrak{S}_1$ and $\mathfrak{S}_2$ denote the components of the principal symbol of $\mathcal{D}_0$ along the respective $x_1$ and $x_2$ directions. The section $v_*$ of $\mathbb{V}^+ \otimes \mathcal{I}_*$ can be written as $v_* = \alpha_+ z^{N_◊} + \alpha_- \bar{z}^{N_◊}$ with $\alpha_+$ and $\alpha_-$ being constant (complex) elements from the respective $+i$ and $-i$ eigenspaces in $\mathbb{V}^+$ for the involution $\mathfrak{S}_1^\dagger \mathfrak{S}_2$.

***Proof of Lemma 5.7***: The proof of the lemma has eleven parts. The implicit assumption has $N_0$ being less than $N_◊ + \frac{1}{2}$ so that $N_◊$ is the greatest half integer that is less than $N_0$. This paper's standard convention with regards to the use of $c_0$ is changed in the the proof so as to allow $c_0$ to depend implicitly on a given upper bound for $N_0$.

*Part 1*: Suppose that $\mathcal{I}_*$ is isomorphic to the product $\mathbb{R}$ bundle on the complement of $Z_*$ in some neighborhood of each $|x| = 1$ point $Z_*$. By way of an example, it follows from Lemma 5.4 that this is the case if $d = 3$ and $N_◊ \in \mathbb{Z}$. Let q denote an $|x| = 1$ point in $Z_*$ and let $B \subset \mathbb{R}^d$ denote a ball centered at q such that $\mathcal{I}_*|_B$ is isomorphic to the product $\mathbb{R}$ bundle on $B - (B \cap Z_*)$. Such an isomorphism writes $v_*$ on $B - (B \cap Z_*)$ as an $\mathbb{R}$ valued harmonic section of the product $B \times \mathbb{V}^+$ whose norm extends over B as a continuous function with zero locus being $B \cap Z_*$. Standard elliptic theory can be invoked to see that $v_*$ extends over $Z_*$ as a harmonic section of the product $\mathbb{V}^+$ bundle. As such, it is real analytic on B and has a Taylor's expansion at each point in $B \cap Z_*$. This Taylor's expansion can be used with the fact that $Z_*$ obeys a first order elliptic equation to see that $Z_* \cap B$ is not locally separating in B.

With the preceding understood, use the fact that the $|x| = 1$ sphere is compact to fix a finite set of balls in the $|x| = 1$ sphere to be denoted by $\mathfrak{U}$ such that if $B \subset \mathfrak{U}$, then $\mathcal{I}_*$ is isomorphic to the product $\mathbb{R}$ bundle on $B - (B \cap Z_*)$. Fix such an isomorphism and denote it by $\iota_B$. If B and B´ are any two balls from $\mathfrak{U}$ that intersect, set $\iota_{BB'}$ to denote the composition $\iota_B(\iota_{B'}^{-1})$. This is a locally constant function on the complement of $Z_*$ in $B \cap B'$ with values in $\{\pm 1\}$. More to the point, it is everywhere either +1 or -1 because $Z_*$ is not locally separating. Granted this, then the collection $\{\iota_{BB'}\}_{B,B' \in \mathfrak{U}}$ unambiguously defines the cocycle data for a real line bundle on the $|x| = 1$ sphere that is isomorphic to $\mathcal{I}_*$ on the complement of $Z_*$. Since the $|x| = 1$ sphere is simply connected, this bundle is the product bundle. This implies that $\mathcal{I}_*$ is isomorphic to the product bundle on $\mathbb{R}^d - Z_*$.

The preceding isomorphism write $v_*$ as an $\mathbb{R}$ valued harmonic section of the product $\mathbb{V}^+$ bundle $\mathbb{R}^d - Z_*$ whose norm extends as zero across $Z_*$. As noted previously, standard regularity theorems can be invoked to see that $v_*$ extends over the whole of $\mathbb{R}^d$ as a harmonic section of the product $\mathbb{V}^+$ bundle. In particular, this extended $\nu_*$ is real analytic and it follows from (5.2) that its coefficients are homogeneous polynomials of degree $N_0$ in the coordinate functions. This implies in particular that $N_0$ is an integer and thus it is equal to $N_\diamond$. Given that $N_0$ is an integer, the remaining assertion of Lemma 5.7 follows from Lemma 5.4 or Lemma 5.5 as the case may be.

*Part 2*: The strategy used to prove Lemma 5.7 when $\mathcal{I}_*$ is not locally isomorphic to the product bundle is modeled on the proof that two eigenvectors of a self-adjoint operator are orthogonal if their eigenvalues are different. In the context of the upcoming proof of Lemma 5.7, the relevant operator is the first order operator $\mathfrak{S}_r^\dagger \mathcal{D}_0^\perp$ that appears in (5.3). One of the relevant eigenvectors in the upcoming proof is the section $u_*$ that is depicted in (5.2) and (5.3); the corresponding eigenvalue is observedly $-N_0$. The proof constructs a second eigenvector of $\mathfrak{S}_r^\dagger \mathcal{D}_0^\perp$ with eigenvalue $-N_\diamond$; and this $-N_\diamond$ eigenvector is very nearly identical to $u_*$ when $N_0 - N_\diamond$ is small. There is a devil here none-the-less because the analog of $Z_*$ for the $-N_\diamond$ eigenvector need not be the same as $Z_*$.

*Part 3*: Assume in this and the remaining parts of the proof that $\mathcal{I}_*$ is not locally isomorphic to the product $\mathbb{R}$ bundle on $\mathbb{R}^d - Z_*$. Define a symmetric, non-negative definite $d \times d$ matrix to be denoted by $T$ with entries $\{T_{ab}\}_{a,b \in \{1,2,\ldots d\}}$ that are given by the rule

$$T_{ab} = \int_{|x|=1} \langle \nabla_a \nu_*, \nabla_b \nu_* \rangle \,.$$

(5.16)

The trace of this matrix is $N_0$. Of interest is its smallest eigenvalue. The next paragraph proves that this smallest eigenvalue is no greater than $c_0 (N_0 - N_\diamond)$.

To prove the claim, let $\hat{u}$ denote an $|x| = 1$ point in $Z_*$ with $\tilde{N}_q \geq N_0$. Such a point exists because of the second bullet's assumption in Lemma 5.7 because $\mathcal{I}_*$ is not locally isomorphic to the product $\mathbb{R}$ bundle on $\mathbb{R}^d - Z_*$. Let $q = 2\hat{u}$. Fix $r \in (0, 1]$ and write the directional covariant derivative of $\nu_*$ on $|x| = r$ sphere in $\mathbb{R}^d$

$$2 \nabla_{\hat{u}} \nu_* = r \nabla_r \nu_* - N_0 \nu_* - r_q \nabla_{r_q} \nu_* + N_q(r_q) \nu_* + (N_0 - N_q(r_q)) \nu_* \,,$$

(5.17)

with $r_q$ being the function defined by the rule $x \to |x - q|$ and with $N_q(\cdot)$ being the function that is defined in (5.10). Let $A$ denote the spherical annulus where $|x - q| \in [1, 3]$. Keeping in mind that $r \nabla_r \nu_* = N_0 \nu_*$, and keeping in mind that $N_q(\cdot) \in [N_\diamond, N_0]$, the identity in (5.17) leads to the bound

$$\langle \hat{u}, T\hat{u} \rangle \leq c_0 \int_A |\nabla_{r_q} v - \tfrac{1}{r_q} N_{(q)} v|^2 + c_0 \sup_{r \in (0,\infty)} (N_0 - N_q(r)).$$

(5.18)

The integral in (5.18) is no greater than $c_0(N_0 - N_\diamond)$. That this is so can be seen by integrating the identity in (5.11). The left most term in (5.18) is also bounded by this same $c_0(N_0 - N_\diamond)$ because the function $N_q(\cdot)$ is non-decreasing and its $r \to 0$ limit is by assumption no smaller than $N_\diamond$.

*Part 4*: Fix $c \in (c_0, c_0^2)$ so that $\langle q, Tq \rangle < c(N_0 - N_\diamond)$ if q is an $|x| = 1$ point in $Z_*$ with $\tilde{N}_q \geq N_\diamond$. Introduce $\Pi$ to denote the span of the eigenvalues of T with eigenvalue less than $c(N_0 - N_\diamond)$. The following dimension bounds are obeyed when $(N_0 - N_\diamond) \leq c_0^{-1}$.

- *The codimension of $\Pi$ is at most* d-2.
- *There are at least two linearly independent eigenvectors of* T *with eigenvalue greater than* $c_0^{-1} N_0$.

(5.19)

The first bullet follows from the second. The second bullet follows from two facts, the first being that the trace of T is equal to $N_0$, and the second being that $\mathcal{D}_0 v_* = 0$. The latter equation bounds the norm of the directional covariant derivative of $v_*$ along any given unit vector by $c_0$ times the norm of the covariant derivative of $v_*$ along the orthogonal d-1 dimensional subspace. It follows from this and the trace(T) = $N_0$ identity that T has at least two linearly independent eigenvectors with eigenvalue greater than $c_0^{-1} N_0$.

The first bullet in (5.19) says that $\Pi$ has codimension 2 if d = 3 and that $\Pi$ can have either codimension 2 or codimension 3 if d = 4. The subsequent parts of the proof assume that $\Pi$ has codimension 2. The case when $\Pi$ has codimension 3 when d = 4 has a different, but not dissimilar argument and will not be given. In any event, it is a consequence of Lemma 5.6 that $\Pi$ can not have codimension 3 in d = 4 if $N_\diamond$ is half an odd integer and $(N_0 - N_\diamond) \leq c_0^{-1}$; and the latter is the relevant case with regards to the forthcoming application of Lemma 5.7 in the proof of Lemma 10.3.

Let q denote an $|x| = 1$ point in $Z_*$ where $\tilde{N}_q \geq N_\diamond$. Write q as $\mathfrak{p}(q) + \Delta \hat{\imath}$ with $\mathfrak{p}$ denoting here and in what follows the orthogonal projection map from $\mathbb{R}^d$ to $\Pi$, with $\Delta$ being non-negative and with $\hat{\imath}$ being a unit length vector that is orthogonal to $\Pi$. Since $\langle q, Tq \rangle < c(N_0 - N_\diamond)$, it follows that $\Delta \leq c_0 c^{1/2}(N_0 - N_\diamond)^{1/2}$. This is to say that all $|x| = 1$ points in $Z_*$ with $\tilde{N}_{(\cdot)} \geq N_\diamond$ have distance at most $c_0 c^{1/2}(N_0 - N_\diamond)^{1/2}$ from $\Pi$. Since the data $(Z_*, \mathcal{I}_*, v_*)$ is scale invariant, it must also be the case that all $|x| \leq 32$ points in $Z_{*=}$ with $\tilde{N}_{(\cdot)} \geq N_\diamond$ have distance at most $32 c_0 c^{1/2}(N_0 - N_\diamond)^{1/2}$ from $\Pi$.

Fix $m \in (1, 32 c_0 c^{1/2})$ so that all $|x| \leq 32$ points in $Z_*$ where $\tilde{N}_{(\cdot)} \geq N_0$ are in the radius $m(N_0 - N_\diamond)^{1/2}$ tubular neighborhood of $\Pi$. It follows from the second bullet's assumption in Lemma 5.7 that the line bundle $\mathcal{I}_*$ is defined on the complement of the radius $m(N_0 - N_\diamond)^{1/2}$ tubular neighborhood of $\Pi$ in the $|x| \leq 32$ sphere in $\mathbb{R}^d$. If $N_\diamond$ is an integer, let $\mathcal{I}_\diamond$ denote the product $\mathbb{R}$ bundle on $\mathbb{R}^d - \Pi$; and if $N_\diamond$ is half an odd integer, let $\mathcal{I}_\diamond$ denote the Möbius line bundle on $\mathbb{R}^d - \Pi$. Supposing that $(N_0 - N_\diamond)^{1/2} \leq c_0^{-1} m^{-1}$, it follows from Lemma 5.6 that $\mathcal{I}_*$ on this part of $\mathbb{R}^d$ is isomorphic to $\mathcal{I}_\diamond$. This isomorphism is used implicitly in what follows.

*Part 5*: Use the function $\chi$ to construct a smooth, non-negative function of the distance from the origin in $\Pi$ that is equal to 1 where the distance is less than 2 and equal to zero where the distance is greater than $\sqrt{5}$. Denote this function by b and let $\langle b \rangle$ denote its integral on $\Pi$. Use $\beta$ to denote the function $\langle b \rangle^{-1} b$. The function $\beta$ has integral 1 on $\Pi$. Let $D \subset \Pi$ denote the disk of radius $\sqrt{5}$, thus containing the support of $\beta$.

Fix Euclidean coordinates $(x_1, x_2)$ for the subspace orthogonal to $\Pi$. Use z to denote the complex coordinate $z = x_1 + i x_2$ and let $\mathcal{A}$ denote the annulus in the $(x_1, x_2)$ plane where $m(N_0 - N_\diamond)^{1/2} \leq |z| \leq 2$. Define a section of $\mathbb{V}^+ \otimes \mathcal{I}_*$ on $\mathcal{A}$ to be denoted by $\underline{v}_*$ by the rule whereby

$$\underline{v}_*(x_1, x_2) = \int_D v_*(x_1, x_2, \cdot) \beta(\cdot) \ .$$

(5.20)

This is an average of $v_*$ taken along the parallels to subspace $\Pi$. Letting u denote a given vector in $\Pi$ with norm less than 2, an appeal to the fundamental theorem of calculus leads to the bound

$$\int_\mathcal{A} |v_*(x_1, x_2, u) - \underline{v}_*(x_1, x_2)|^2 \, dx_1 dx_2 \leq c_0 \int_{|x|<3} |\nabla^\Pi v_*|^2 \ ,$$

(5.21)

with $\nabla^\Pi$ denoting the covariant derivative in directions tangent to $\Pi$. The inequality in (5.21) leads in turn to the observation that

$$\int_\mathcal{A} |v_*(x_1, x_2, u) - \underline{v}_*(x_1, x_2)|^2 \, dx_1 dx_2 \leq c_0 (N_0 - N_\diamond) \ .$$

(5.22)

This bound follows because $\Pi$ is the span of the eigenvectors of the matrix T in (5.16) with small eigenvalue. Note in this regard that the $L^2$ norm of any derivative of $v_*$ over the $|x| \leq 32$ ball in $\mathbb{R}^d$ is bounded by $c_0$ times the corresponding integral over the $|x| \leq 1$ ball because $v_*$ can be written as in (5.2).

*Part 6*: Write the operator $\mathcal{D}_0$ as

$$\mathcal{D}_0 = \mathfrak{S}_1 \nabla_1 + \mathfrak{S}_2 \nabla_2 + \mathfrak{S}_\Pi(\nabla^\Pi)$$

(5.23)

with $\{\mathfrak{S}_k\}_{k=1,2}$ being homomorphisms from $\mathbb{V}^+$ to $\mathbb{V}^-$ and with $\mathfrak{S}_\Pi$ being a homomorphism from $\Pi$ to $\text{End}(\mathbb{V}^+, \mathbb{V}^-)$. Let $\tau$ denote $\mathfrak{S}_1^\dagger \mathfrak{S}_2$; this is a homorphism of $\mathbb{V}^+$ with square $-1$. Let $\underline{v}_+$ and $\underline{v}_-$ denote the respective projections of $\underline{v}_*$ to the $+i$ and $-i$ eigenspaces of $\tau$ on the complexification of $\mathbb{V}^+ \otimes \mathcal{I}_*$. By virtue of the fact that $\mathcal{D}_0 v_* = 0$, these obey a pair of equations on the annulus $\mathcal{A}$ that have the schematic form

$$(\nabla_1 + i\nabla_2)\underline{v}_+ + \mathfrak{r}_+ = 0 \quad \text{and} \quad (\nabla_1 - i\nabla_2)\underline{v}_- + \mathfrak{r}_- = 0 ,$$

(5.24)

with $\mathfrak{r}_+$ and $\mathfrak{r}_-$ being respective complex sections of $\mathbb{V}^+ \otimes \mathcal{I}_*$ on $\mathcal{A}$ with $L^2$ norms bounded by $c_0(N_0 - N_\diamond)^{1/2}$ and $L^\infty$ norms bounded by $c_0$. These bounds are consequences of the fact that $\mathfrak{r}_+$ and $\mathfrak{r}_-$ are averages along $\Pi$ of derivatives of $v_*$ in directions tangent to $\Pi$. Their $L^2$ norms are therefore bounded by $c_0$ times the $L^2$ norm of $|\nabla^\Pi v_*|$ on the $|x| \le 3$ ball in $\mathbb{R}^d$. As noted previously, this in turn is bounded by $c_0(N_0 - N_\diamond)^{1/2}$ because of the definition of $\Pi$ as the span of the small eigenvalues of the matrix T. The $C^0$ norms come about by integrating by parts to replace derivatives on $v_*$ along $\Pi$ with derivatives of the function $\beta$ that appears in (5.20).

Given what was just said about the size of $\mathfrak{r}_+$ and $\mathfrak{r}_-$, standard elliptic theory can be used to find respective complex sections $e_+$ and $e_-$ of $\mathbb{V}^+ \otimes \mathcal{I}_*$ over $\mathcal{A}$ that obey

$$(\nabla_1 + i\nabla_2)e_+ + \mathfrak{r}_+ = 0 \quad \text{and} \quad (\nabla_1 - i\nabla_2)e_- + \mathfrak{r}_- = 0 ,$$

(5.25)

have $L^2_1$ norm bounded by $c_0(N_0 - N_\diamond)^{1/2}$ and $C^0$ norm bounded by $c_0(N_0 - N_\diamond)^{1/2} |\ln((N_0 - N_\diamond)|^{1/2}$.

It follows from (5.24) and (5.25) that $\underline{v}_+ - e_+$ is a holomorphic section of $\mathbb{V}^+ \otimes \mathcal{I}_*$ on the annulus $\mathcal{A}$. As such, it can be written as

$$\underline{v}_+ - e_+ = \alpha_+ z^{N_\diamond} + h_+(z) z^{N_\diamond}$$

(5.26)

with $\alpha_+$ being a constant section of $\mathbb{V}^+ \otimes \mathcal{I}_\diamond$ and with $h_+$ being a convergent Laurent series on $\mathcal{A}$ with no constant term.

To say more about $h_+$, fix a number $\delta$ less than $\frac{1}{1000}$ but greater than $4m(N_0 - N_\diamond)^{1/2}$. Let $\mathcal{A}_\delta \subset \mathcal{A}$ denote the annulus where $\delta \le |z| \le \frac{7}{2}$. Given that $v_*$ can be written as in (5.2) and given what was just said about the $L^2_1$ norm of $e_+$, it follows from the $u = 0$ version of (5.22) that the $L^2$ norm of the function $h_+$ on the annulus $\mathcal{A}_\delta$ is bounded by $c_\delta(N_0 - N_\diamond)^{1/2}$ with $c_\delta$ denoting here and in what follows a number that depends only on $\delta$. Since the $h_+$ is also holomorphic, these $L^2$ bounds lead to the $C^0$ bound $|h_+| \leq c_\delta(N_0 - N_\diamond)^{1/2}$ on the annulus where $2\delta \leq |z| \leq 3$ with this incarnation of $c_\delta$ denoting a second number that depends only on $\delta$.

By way of a convention, $c_\delta$ in what follows will always denote a number of that is greater than 1 and depends only on $\delta$. It's value can be assumed to increase between successive appearances.

By the same token, $\underline{v}_-$ can be written as $\underline{v}_- = \alpha_- \bar{z}^{N_\diamond} + h_-(\bar{z})\bar{z}^{N_\diamond} + e_-$ with $\alpha_-$ being a constant section of $\mathbb{V}^- \otimes \mathcal{I}_\diamond$ and with $h_-$ being a convergent Laurent series with norm bounded by $|h_-| \leq c_\delta(N_0 - N_\diamond)^{1/2}$ on the annulus where $2\delta \leq |z| \leq 3$.

*Part 7*: Define the section $v_\diamond$ of $\mathbb{V}^+ \otimes \mathcal{I}_\diamond$ on the complement of $\Pi$ by the formula

$$v_\diamond = \alpha_+ z^{N_\diamond} + \alpha_- \bar{z}^{N_\diamond}.$$

(5.27)

Let $\chi_\diamond$ denote the function of the coordinate z given by $\chi(4 - m^{-1}(N_0 - N_\diamond)^{-1/2}|z|)$. This function is equal to 1 where the distance to $\Pi$ is greater than $4m(N_0 - N_\diamond)^{1/2}$ and it is equal to zero where the distance to $\Pi$ is less than $3m(N_0 - N_\diamond)^{1/2}$. Of interest here is the integral

$$Q = \int_{|x|=1} \chi_\diamond \langle v_*, v_\diamond \rangle.$$

(5.28)

The three steps that follow prove that $Q > \tfrac{1}{2}$ when $(N_0 - N_\diamond) \leq c_0^{-1}$.

<u>Step 1</u>: Fix for the moment $\delta \in (0, \tfrac{1}{4}]$. The complement of the $\delta < |z| < 1 - \delta$ part of $|x| = 1$ sphere has d-1 dimensional volume less than $c_0 \delta^{1/2}$. Keeping this in mind, the fact that the integral of $|v_*|^2$ over the $|x| = 1$ sphere is equal to 1 implies that

$$\int_{|x|=1 \text{ and } \delta<|z|<(1-\delta)} |v_*|^2 \geq 1 - c_0 \delta^{1/2}$$

(5.29)

because the $v_*$ version of Lemma 2.3 with U being the $|x| \leq 3$ ball in $\mathbb{R}^d$ gives a $c_0$ bound for the pointwise norm of $|v_*|$ on the $|x| \leq \tfrac{8}{3}$ ball.

<u>Step 2</u>: It follows from (5.22) and what is said about $e_+, e_-$ and $h_+$ and $h_-$ that

$$\int_{|x|=1 \text{ and } \delta<|z|<(1-\delta)} |v_* - v_\diamond|^2 \leq c_\delta (N_0 - N_\diamond) .$$

(5.30)

This last bound with (5.29) leads to a lower bound for the $\delta < |z| < 1-\delta$ part of the integral that defines Q:

$$\int_{|x|=1 \text{ and } \delta<|z|<(1-\delta)} \chi_\diamond \langle v_*, v_\diamond \rangle \geq 1 - c_0 \delta^{1/2} - c_\delta (N_0 - N_\diamond)^{1/2} .$$

(5.31)

Note in this regard that $\chi_\diamond = 1$ on the $\delta < |z| < 1-\delta$ part of the $|x| = 1$ sphere.

Step 3: As noted already, there is a $c_0$ bound for $|v_*|$ on the $|x| \leq 2$ ball. Meanwhile, (5.22) or (5.30) leads to a corresponding $c_0$ bound for $|v_\diamond|$ on this ball. These pointwise bounds imply in turn that

$$\int_{|x|=1 \text{ and either } |z|<\delta \text{ or } |z|>(1-\delta)} \chi_\diamond |\langle v_*, v_\diamond \rangle| \leq c_0 \delta^{1/2}$$

(5.32)

because of the $c_0 \delta^{1/2}$ bound on the d-1 dimensional volume of part of the $|x| = 1$ sphere where $|z| < \delta$ or where $|z| > 1-\delta$. The $\delta = c_0^{-1}$ version of (5.31) and (5.32) gives the asserted lower bound for Q when $N_0 - N_\diamond \leq c_0^{-1}$.

*Part 8*: Write $v_*$ as $u_* |x|^{N_0}$ with $u_*$ being viewed here as a section of $\mathbb{V}^+ \otimes \mathcal{I}_\diamond$ on the $|z| > m(N_0 - N_\diamond)$ part of $|x| = 1$ sphere. Likewise write $v_\diamond$ as $u_\diamond |x|^{N_\diamond}$ with $u_\diamond$ being a section of $\mathbb{V}^+ \otimes \mathcal{I}_\diamond$ on the $|z| > 0$ part of the $|x| = 1$ sphere. The section $u_*$ obeys the equation that is depicted in (5.3) and the section $u_\diamond$ obeys the analogous equation with $N_\diamond$ playing the role of $N_0$. These equations with the identities in (1.2) and an integration by parts lead to the writing of $(N_0 - N_\diamond)Q$ as

$$(N_0 - N_\diamond)Q = \int_{|x|=1} \langle v_*, \wp(\nabla^\perp \chi_\diamond) v_\diamond \rangle$$

(5.33)

with $\wp$ denoting here the principal symbol of the operator $\mathfrak{S}_r^\dagger \mathcal{D}_0^\perp$.

The formula in (5.33) leads directly to Lemma 5.7's claim when $N_\diamond \geq 1$. To explain why this is, let $\Delta$ denote for the moment $m(N_0 - N_\diamond)^{1/2}$. The norm of $v_\diamond$ on the support of $d\chi_m$ is no greater than $c_0 \Delta^{N_\diamond}$. Let z denote for the moment an upper bound for the norm of $|v_*|$ on the support of $d\chi_m$ in the $|x| = 1$ sphere. It follows from the fact that $|d^\perp \chi_m| < c_0 \Delta^{-1}$ and the fact that the d-1 dimensional volume of the support of $d\chi_m$ is bounded by $c_0 \Delta^2$ that the absolute value of the integral on the right hand side of (5.33) is bounded by $c_0 z \Delta^{1+N_\Diamond}$. Keeping in mind that $m \leq c_0$ so that $\Delta \leq c_0 (N_0 - N_\Diamond)^{1/2}$ this bound for the right hand side of (5.33) is less than $c_0 z (N_0 - N_\Diamond)^{(1+N_\Diamond)/2}$. Use this bound in (5.33) with the lower bound $Q \geq \frac{1}{2}$ from Part 5 to obtain an inequality that reads

$$(N_0 - N_\Diamond) < c_0 z (N_0 - N_\Diamond)^{(1+N_\Diamond)/2}.$$
(5.34)

Since $z \leq c_0$, this last bound can hold when $N_\Diamond > 1$ and $N_0 - N_\Diamond < c_0^{-1}$ only if $N_0 = N_\Diamond$. As explained in the next paragraph, $z < \frac{1}{1000} c_0^{-1}$ when $N_0 - N_\Diamond < c_0^{-1}$. Granted this bound, then (5.34) also requires that $N_0 = N_\Diamond$ when $N_\Diamond = 1$ and $N_0 - N_\Diamond < c_0^{-1}$. In any event, if $N_0 = N_\Diamond$, then the remaining assertions of Lemma 5.7 follow from Lemma 5.4 or 5.5.

To see about the size of $z$, suppose for the moment that q is a point on the $|x| = 1$ sphere where $d\chi_m \neq 0$. Let $z_q$ denote the norm of $|v_*|$ at q. Fix $\delta \in (m(N_0 - N_\Diamond)^{1/2}, 1]$. Since the function $\kappa_q$ is increasing, it follows that the integral of $|v_*|^2$ on the ball of radius $\delta^{1/2}$ centered at q is greater than $d^{-1} z_q^2 \delta^{d/2}$. On the otherhand, the integral of $|v_\Diamond|^2$ on this same sphere is bounded by $c_0 \delta^{d/2+N_\Diamond}$. What with the bounds on $h_+$, $h_-$ and $e_+$ and $e_-$ where $|z| > \delta$, these bounds for the respective integrals of $|v_*|^2$ and $|v_\Diamond|^2$ will run afoul of the bound in (5.22) unless $z_q^2 \leq c_0(\delta^{N_\Diamond} + \delta^{d/2}) + c_\delta(N_0 - N_\Diamond)$. Supposing that $\varepsilon > 0$, it follows from the latter bound that $z_q$ is less than $\varepsilon$ if $(N_0 - N_\Diamond)$ is less than a number greater than 1 that depends only on $\varepsilon$.

Now suppose that $N_\Diamond$ is half an odd integer, with $N_\Diamond = \frac{1}{2}$ being of particular interest. If $N_\Diamond$ is half an odd integer, then $\mathcal{I}_\Diamond$ is not isomorphic to the product $\mathbb{R}$ bundle on $\mathbb{R}^d - \Pi$, and since $\mathcal{I}_* = \mathcal{I}_\Diamond$ where $|z| \geq m(N_0 - N_\Diamond)^{1/2}$ and $|x| \leq 32$ it follows that $\mathcal{I}_*$ is not isomorphic to the product $\mathbb{R}$ bundle on the $|z| \geq m(N_0 - N_\Diamond)^{1/2}$ part of the $|x| \leq 32$ ball. This in turn has the following implication: Suppose that q a point from $\Pi$ with $|q| \leq 16$. Then there is an $\tilde{N}_{(\cdot)} \geq N_\Diamond$ point from $Z_*$ in the normal plane to $\Pi$ through the point q with distance at most $m(N_0 - N_\Diamond)^{1/2}$ from q. It follows as a consequence that the number z that appears in (5.34) is no larger than $c_0 m^{N_\Diamond} (N_0 - N_\Diamond)^{N_\Diamond/2}$. The right hand side of (5.34) is thus no greater than $c_0 (N_0 - N_\Diamond)^{(1+2N_\Diamond)/2}$. Unfortunately, this bound is not sufficient for the purposes at hand when $N_\Diamond = \frac{1}{2}$ because the exponent is equal to 1, not greater than 1.

As proved momentarily, the absolute value of the right hand side of (5.33) in the case when $N_\Diamond$ is half an odd integer, and in particular, when $N_\Diamond = \frac{1}{2}$, is no greater than

$$c_0 (N_0 - N_\Diamond)^{(1+N_\Diamond)} |\ln(N_0 - N_\Diamond)|^{1/2}.$$
(5.35)

Substituting this expression for the right hand side of (5.34) leads to the desired conclusion that $N_0 = N_\diamond$ if $(N_0 - N_\diamond) \leq c_0^{-1}$.

*Part 9*: This part of and the remaining parts of the proof derive the bound in (5.35) for the right hand side of (5.33). By way of notation, $\Delta$ is used here and in the subsequent parts of the proof to denote the number $m(N_0 - N_\diamond)^{1/2}$.

The first step to deriving (5.35) uses the scaling behavior of $v_*$ and $v_\diamond$ to write the integral on the right hand side of (5.33) as a volume integral. Let $\chi_{\diamond\diamond}$ denote the scaling invariant function $x \to \chi_{\diamond\diamond} = \chi(4 - \Delta^{-1}|x|^{-1}|z|)$. This is the same as $\chi_\diamond$ on the $|x| = 1$ sphere but it has the advantage of being fixed by any $\lambda \in (0, \infty)$ version of the rescaling map $x \to \lambda x$. The homogeneous scaling of $v_*$ and $v_\diamond$ under this same scaling map can be invoked to write the integral on the right hand side of (5.33) as

$$w_\Delta^{-1} \int_{1-\Delta \leq |x| \leq 1} \langle v_*, \wp(\nabla^\perp \chi_{\diamond\diamond}) v_\diamond \rangle ,$$

(5.36)

with $w_\Delta$ denoting $(d + N_0 + N_\diamond)^{-1}(1 - (1-\Delta)^{d + N_0 + N_\diamond})$. Of particular note for what follows is that $w_\Delta \geq c_0^{-1} \Delta$ when $\Delta \leq \frac{1}{10}$. The second step to deriving (5.35) is to use the fact that the symbol $\wp$ is antisymmetric to write (5.36) as

$$w_\Delta^{-1} \int_{1-\Delta \leq |x| \leq 1} \langle v_* - v_\diamond, \wp(\nabla^\perp \chi_{\diamond\diamond}) v_\diamond \rangle .$$

(5.37)

Note with regards to (5.37) that the bundles $\mathcal{I}_*$ and $\mathcal{I}_\diamond$ are canonically isomorphic on the support of the integrand in (5.37) because the support of $d\chi_{\diamond\diamond}$ has distance greater than $3(1 - c_0\Delta)\Delta$ from $\Pi$ on the spherical annulus where $|x|$ is between $1 - \Delta$ and $\Delta$.

Now suppose that $q_\diamond$ is a point in $\Pi$ on the $|x| = 1 - 2\Delta$ sphere and that $q \in Z_*$ is a point in the normal plane to $\Pi$ through $q_\diamond$ with distance less than $\Delta$ from the point $q_\diamond$ and with $\tilde{N}_q \geq N_\diamond$. Let $A_q$ denote the spherical annulus centered at q with inner radius $\Delta$ and outer radius $8\Delta$. The subsequent parts of the proof derive the bound

$$\left| \int_{1-\Delta \leq |x| \leq 1 \text{ and } x \in A_q} \langle v_* - v_\diamond, \wp(\nabla^\perp \chi_{\diamond\diamond}) v_\diamond \rangle \right| \leq c_0 |\ln\Delta|^{1/2} \Delta^{(d+2N_\diamond)} .$$

(5.38)

This bound leads to the bound in (5.35) because the support of the integrand in (5.37) can be covered by a set of $c_0$ annuli of the form $A_q$ if $d = 3$ and a set $c_0 \Delta^{-1}$ such annuli if $d = 4$.

*Part 10*: The bound in (5.36) is obtained in series of steps. The first step is to introduce the pull-back of $v_\diamond$ by the translation $x \to x + q$. This pull-back is denoted by $v_{\Diamond q}$. Its zero-locus is the translation of $\Pi$ so its zero locus passes through the point q. The integral in (5.36) is written as a sum of two terms:

$$\int_{1-\Delta \leq |x| \leq 1 \text{ and } x \in A_q} \langle v_* - v_{\Diamond q}, \wp(\nabla^\perp \chi_{\Diamond\Diamond}) v_\Diamond \rangle + \int_{1-\Delta \leq |x| \leq 1 \text{ and } x \in A_q} \langle v_{\Diamond q} - v_\Diamond, \wp(\nabla^\perp \chi_{\Diamond\Diamond}) v_\Diamond \rangle .$$

(5.39)

The rest of this part of the proof derives a suitable bound for the right most integral in (5.39). To this end, let $q_z$ denote the z coordinate of the point q. This obeys $|q_z| \leq m\Delta$. Let $z'$ denote the complex function $z' = z - q_z$. What is denoted by $v_{\Diamond q}$ can be written as

$$v_{\Diamond q} = \alpha_+ z'^{N_\Diamond} + \alpha_- \bar{z}'^{N_\Diamond} .$$

(5.40)

The integrand of the right most integral in (5.39) is the real part of the sum of two terms, the first being

$$\langle \bar{\alpha}_-, \wp(\nabla^\perp \chi_{\Diamond\Diamond}) \alpha_+ \rangle z'^{N_\Diamond} z^{N_\Diamond} + \langle \bar{\alpha}_+, \wp(\nabla^\perp \chi_{\Diamond\Diamond}) \alpha_- \rangle \bar{z}'^{N_\Diamond} \bar{z}^{N_\Diamond} ,$$

(5.41)

and the second being

$$\langle \bar{\alpha}_+, \wp(\nabla^\perp \chi_{\Diamond\Diamond}) \alpha_+ \rangle \bar{z}'^{N_\Diamond} z^{N_\Diamond} + \langle \bar{\alpha}_-, \wp(\nabla^\perp \chi_{\Diamond\Diamond}) \alpha_- \rangle z'^{N_\Diamond} \bar{z}^{N_\Diamond} .$$

(5.42)

Being that $\wp(\nabla^\perp \chi_{\Diamond\Diamond})$ is anti-Hermitian, the expression in (5.41) is purely imaginary. This being so, it makes no contribution to (5.39). As explained in the next paragraph, the absolute value of each term in (5.42) is bounded by $c_0 \Delta^{2N_\Diamond}$. This leads directly to a $c_0 \Delta^{d+2N_\Diamond}$ bound for the absolute value of the right most integral in (5.39).

The asserted $c_0 \Delta^{2N_\Diamond}$ bound for the absolute value of each term in (5.42) comes about as follows: The factor $\Delta^{2N_\Diamond}$ comes from the factors of z and z' that appear in (5.42). Althought the supremum of $|\nabla^\perp \chi_{\Diamond\Diamond}|$ is greater than $c_0^{-1} \Delta^{-1}$, the contribution of $\wp(\nabla^\perp \chi_{\Diamond\Diamond})$ to the terms in (5.42) is bounded by $c_0$. To prove that this is so, let $\tau$ denote the square -1 automorphism $\mathfrak{S}_1^\dagger \mathfrak{S}_2$ of $\mathbb{V}^+$. The automorphism $\wp(\nabla^\perp \chi_{\Diamond\Diamond})$ can be written as a sum $\wp_+ + \wp_-$ with the commutator of $\tau$ and $\wp_+$ being zero, and with the anticommutator of $\tau$ and $\wp_-$ being zero. The automorphism $\wp_+$ has norm bounded by $c_0$ because the coordinate vectors $\frac{\partial}{\partial x_1}$ and $\frac{\partial}{\partial x_2}$ differ on the support of $\nabla \chi_{\Diamond\Diamond}$ by at most $c_0 \Delta$ from their orthogonal projections to the tangent space of the $|x| = 1$ sphere. The automorphism $\wp_-$ has norm as large as $c_0^{-1} \Delta^{-1}$ but it makes no contribution to the terms in (5.42) because $\alpha_+$ and $\alpha_-$ are complex eigenvectors for the automorphism $\tau$.

*Part 11*:  This part of the proof derives a bound for the left hand integral in (5.39). To start the derivation, let $U_q$ denote the support of the integrand for the left most integral in (5.39). Since the volume to $U_q$ is less than $c_0 \Delta^d$ and the norm of $d\chi_{\diamond\diamond}$ is less than $c_0 \Delta^{-1}$ and $|v_\diamond| \leq c_0 \Delta^{N_\diamond}$, this left most integral can not be larger than

$$c_0 \Delta^{d/2-1+N_\diamond} \left( \int_{U_q} |v_* - v_{\diamond q}|^2 \right)^{1/2} .$$

(5.43)

To see about the integral in (5.43), fix $\lambda \in (0, \infty)$ for the moment and let $\phi_{T(\lambda)}$ denote the combination of scaling and translation given by the rule $x \to \phi_{T(\lambda)}(x) = \lambda x + q$. This map sends the origin to q; and if $r \in (0, \infty)$, and it sends the $|x| = r$ sphere to the radius $\lambda r$ sphere centered at the point q. Of particular interest is the version with $\lambda = 3\Delta$. This version sends the spherical annulus in $\mathbb{R}^d$ where $\frac{1}{3} \leq |x| \leq \frac{8}{3}$ to the set $A_q$. Moreover, the $\phi_{T(3\Delta)}$-inverse image of $U_q$ is contained in the part of the $\frac{1}{3} \leq |x| \leq \frac{8}{3}$ annulus where the coordinate z obeys $1 - c_0\Delta \leq |z| \leq \frac{5}{3}$ and where the distance the $(x_1, x_2)$ plane is greater than $\frac{1}{3}$. Use U to denote the $\phi_{T(3\Delta)}$-inverse image of $U_q$.

The integral in (5.43) is equal to $(3\Delta)^d$ times the integral over U of the pull-back by the map $\phi_{T(3\Delta)}$ of the function $|v_* - v_{\diamond q}|^2$. The $\phi_{T(3\Delta)}$ pull-back of $v_{\diamond q}$ is $(3\Delta)^{N_\diamond} v_{*q}$. Meanwhile, the $\phi_{T(3\Delta)}$ pull-back of $v_*$ can be written in the notation of Section 4 as $K_q(3\Delta) v_{T(3\Delta)}$ with $T(3\Delta)$ being the $\lambda = 3\Delta$ version of the data set $T(\lambda) = (q, L, \lambda)$.

Let $\gamma = K_q(3\Delta)/K_q(1)$. It follows from (5.11) that

$$(3\Delta)^{N_0} \leq \gamma \leq (3\Delta)^{N_\diamond} ,$$

(5.44)

and with this understood, writing $\phi_{T(3\Delta)}^* v_*$ as $K_q(3\Delta) v_{T(3\Delta)}$ leads to the bound

$$\phi_{T(3\Delta)}^*(|v_* - v_{\diamond q}|^2) \leq c_0 \Delta^{N_\diamond} |K_q(1) v_{T(3\Delta)} - v_{\diamond q}|^2 + c_0 \Delta^{N_\diamond} |\ln(\Delta)| .$$

(5.45)

The bound in (5.45) leads in turn to a bound on (5.43) by the expression

$$c_0 \Delta^{d-1+2N_\diamond} \left( \int_U |K_q(1) v_{T(3\Delta)} - v_{\diamond q}|^2 \right)^{1/2} .$$

(5.46)

To say something about (5.46), write $v_{T(1)}$ as $K_q(1)(\phi_{T(1)}^* v_*)$ and so bound (5.46) by

$$c_0 \Delta^{d-1+2N_\diamond} K_q(1) \left( \int_U |v_{T(3\Delta)} - v_{T(1)}|^2 \right)^{1/2} + c_0 \Delta^{d-1+2N_\diamond} \left( \int_U |(\phi_{T(1)}^* v_*) - v_{\diamond q}|^2 \right)^{1/2} .$$

(5.47)

The right most integral in (5.47) is no greater than $c_0 \Delta^2$. This bound follows from (5.30) because $\phi_{T(1)}$ is the rigid translation $x \to x+q$ and $\nu_{\diamond q}$ is the pull-back of $\nu_\diamond$ by this translation. The left most integral in (5.47) is no greater than $c_0 |\ln\Delta| \Delta^2$ as can be seen by integrating the point q's version of the formula in (5.13). These bounds lead to a $c_0 \Delta^{d+2N_\diamond} |\ln\Delta|^{1/2}$ bound for the integral in (5.43) and thus to a $c_0 \Delta^{d+2N_\diamond} |\ln\Delta|^{1/2}$ bound for the left most integral in (5.39).

### 6. Notions of continuity

Assume in what follows that the dimension of X is either 3 or 4. Given $p \in Z$, use $N_{(p)}(0)$ to denote the value at 0 of the function N from p's version of Lemma 3.2. A point $p \in Z$ is said to be *weakly continuous* if there is a sequence $\{p_k\}_{k=1,2,\ldots} \subset Z$ that converges to p and is such that $\lim_{k\to\infty} N_{(p_k)}(0) = N_{(p)}(0)$. By way of a reminder, Lemma 3.3 makes an assertion to the effect that if $\{p_k\}_{k \in \{1,2,\ldots\}} \subset Z$ converges to p, then the lim-sup of $\{N_{(p_k)}(0)\}_{k=1,2,\ldots}$ is no greater than $N_{(p)}(0)$. A point $p \in Z$ is said to be *strongly continuous* when $N_{(\cdot)}(0)$ is constant on a neighborhood of p. A strongly continuous point in Z is a priori a point of continuity for $N_{(\cdot)}(0)$ and a point of continuity for $N_{(\cdot)}(0)$ that is not an isolated point in Z is a priori a weakly continuous point in Z. Neighborhoods in Z of these various sorts of points are described in the subsequent subsections and in Sections 7-10.

The lemma that follows directly makes a formal assertion to the effect that the set of points in Z that are *not* weakly continuous is either finite or countably infinite.

**Lemma 6.1**: *The complement in Z of its set of weakly continuous points is countable.*

*Proof of Lemma 6.1*: Supposing that m denotes a positive integer, let $\mathcal{X}_m \subset Z$ denote the set of points in Z with the property that $N_{(\cdot)}(0) < N_{(p)}(0) - \frac{1}{m}$ on the radius $\frac{1}{m}$ ball centered at p. This definition is such that $\mathcal{X}_m \subset \mathcal{X}_{m+1}$. It is also such that union $\cup_{m=1,2,\ldots} \mathcal{X}_m$ is the complement in Z of the set of weakly continuous points. To see why this is, let p denote for the moment a point in Z that is not weakly continuous. Given the definition of weak continuity and given that p is not weakly continuous, there must exists $\delta > 0$ and $r > 0$ such that $N_{(\cdot)}(0) < N_{(p)}(0) - \delta$ on the radius r ball centered at p. This implies that p is contained in $\mathcal{X}_m$ if m is larger than both $\delta^{-1}$ and $r^{-1}$.

The complement in Z of the set of weakly continuous points is countable if each $m \in \{1, 2, \ldots\}$ version of $\mathcal{X}_m$ is finite; and any given $m \in \{1, 2, \ldots\}$ version of $\mathcal{X}_m$ is finite because the distance between distinct points in $\mathcal{X}_m$ is at least $\frac{1}{m}$. Indeed, were p and p´ in $\mathcal{X}_m$ with $0 < \text{dist}(p,p´) < \frac{1}{m}$, then nonsense would ensue because $N_{(p)}(0)$ would be less than $N_{(p´)}(0) - \frac{1}{m}$ and likewise $N_{(p´)}(0)$ would be less than $N_{(p)}(0) - \frac{1}{m}$.

a) **Uniform convergence near points in Z**

The upcoming Lemma 6.2 plays a central role in the subsequent analysis. The lemma introduces the following terminology: Fix $p \in Z$. Given $r \in (0, r_0)$ and a number $\delta$ which is positive, the lemma has $Z_{r,\delta}(p)$ denoting the $N_{(\cdot)}(0) \geq N_{(p)}(0) - \delta$ part of $B_r \cap Z$. The lemma also uses $M_{(p)}$ denote the largest half integer from the interval $[0, N_{(p)}(0))$.

Lemma 6.2 introduces a second novelty, this being a measure of the colinearity of three points in a given radius $r_0$ ball in X. The definition supposes that $p \in X$ is given and that $p_1$ and $p_2$ are points from the radius $r_0$ ball centered at $p$ with neither being $p$. Define $\Delta_p(p', p'')$ to be the absolute value of the sine of the angle between the respective tangent vectors at $p$ to the short geodesics from $p$ to $p_1$ and from $p$ to $p_2$

The statement of Lemma 6.2 refers to data sets $(Z_*, \mathcal{I}_*, \nu_*)$ from (5.1), and it refers to the function $\tilde{N}_{(\cdot)}$ on the set of $|x| = 1$ points in $Z_*$ that is described in either Lemma 5.4 or 5.5 as the case may be.

**Lemma 6.2**: *Fix $p \in Z$ and then $\varepsilon \in (0, \frac{1}{100})$. There exists $\delta \in (0, \varepsilon)$ and $r \in (0, \varepsilon r_0)$ depending only on $p$ and $\varepsilon$ with the following significance: Suppose that $q \in Z_{r,\delta}(p)$. Let L denote a chosen oriented, orthonormal frame for $TX|_q$ and let $\lambda \in (0, r]$. There is a data set $(Z_*, \mathcal{I}_*, \nu_*)$ from the $N_0 = N_{(p)}(0)$ version of (5.1) with the properties listed in the four bullets given below. The bullets use $\mathcal{T}_\varepsilon$ to denote the $|x| < \varepsilon^{-1}$ part of the radius $\varepsilon$ tubular neighborhood of $Z_*$. They also refer to the data set $(Z_T, \mathcal{I}_T, \nu_T)$ that is constructed using $T = (q, L, \lambda)$ by the rules set forth at the beginning of Section 4.*
- *The $|x| < \varepsilon^{-1}$ part of $Z_T$ is in $\mathcal{T}_\varepsilon$.*
- *There is an isometric isomorphism between $\mathcal{I}_T$ and $\mathcal{I}_*$ on the complement of $\mathcal{T}_\varepsilon$ in the $|x| < \varepsilon^{-1}$ ball that allows $\nu_T$ to be viewed as a section of $\mathbb{V}^+ \otimes \mathcal{I}_*$ on this same complement of $\mathcal{T}_\varepsilon$. Viewed in this way, $\nu_T$ differs from $\nu_*$ on the complement of $\mathcal{T}_\varepsilon$ in the $|x| < \varepsilon^{-1}$ ball by a section of $\mathbb{V}^+ \otimes \mathcal{I}_*$ with $C^k$ norm less than $\varepsilon$ for all integers $k < \varepsilon^{-1}$.*
- *Suppose that $\dim(X) = 3$ and that there is a point $q'$ in Z with $N_{(q')}(0) > M_p + \varepsilon$ with $\text{dist}(q, q') \in [\varepsilon\lambda, \varepsilon^{-1}\lambda]$. Then there is an $|x| = 1$ point in $Z_*$ where $\tilde{N}_{(\cdot)} = N_{(p)}(0)$; thus $Z_*$ is described by the fourth bullet of Lemma 5.4,*
- *Suppose that $\dim(X) = 4$ and there is a point $q' \in Z_{r,\delta}$ with $\text{dist}(q, q') \in [\varepsilon\lambda, \varepsilon^{-1}\lambda]$. Then there is an $|x| = 1$ point in $Z_*$ where $\tilde{N}_{(\cdot)} = N_{(p)}(0)$ and $Z_*$ is described by third bullet of Lemma 5.5. In this case the number of half-planes in $Z_*$ has an upper bound that depends only on $N_{(p)}(0)$.*
- *Suppose again that $\dim(X) = 4$. The set $Z_*$ is described by the fourth bullet of Lemma 5.5 if there are points $q'$ and $q''$ with distance from $q$ between $\varepsilon\lambda$ and $\varepsilon^{-1}\lambda$ from $p$, with $q'$ in $Z_{r,\delta}(p)$, with $q''$ in the $N_{(\cdot)}(0) > M_{(p)} + \varepsilon$ part of Z, and with $\Delta_q(q', q'') > \varepsilon^2$.*

This lemma is proved momentarily. What follows is a useful corollary to Lemma 6.2.

**Lemma 6.3**: *Suppose that $p \in Z$ is weakly continuous point and let $\{p_i\}_{i \in \{1,2,...\}}$ denote a sequence from $Z \cap B_{r_0}$ with limit $p$ and with $\lim_{i \to \infty} N_{(p_i)}(0) = N_{(p)}$. Let $(Z_*, \mathcal{I}_*, \nu_*)$ denote a limit data set from Proposition 4.1 using the input sequence $\{T(i) = (p, L, \lambda_i = \text{dist}(p, p_i))\}$.*

- *If $\dim(X) = 3$, then $(Z_*, \mathcal{I}_*, \nu_*)$ is described by the fourth bullet of Lemma 5.4.*
- *If $\dim(X) = 4$, then $(Z_*, \mathcal{I}_*, \nu_*)$ is described by the third bullet of Lemma 5.5.*

*Proof of Lemma 6.2*: The first two bullets of the lemma are proved in Part 1 of what follows. The third and fourth bullets are proved in Part 2 and the fifth is proved in Part 3.

    *Part 1*: The proof of the first two bullets starts by assuming that they are not true and it then generates nonsense from this assumption. Suppose that there exists $p \in Z$, a number $\varepsilon \in (0, \frac{1}{100})$ and a sequence $\{(p_i, L_i, \lambda_i)\}_{i=1,2,...}$ with each member of this sequence consisting of a point in $Z$, an oriented, orthonormal frame for $TX$ at the point and a number $\lambda_i \in (0, \frac{1}{100} r_0]$. This sequence of data sets has four salient properties, the first being that $\lim_{i \to \infty} p_i = p$; the second being that $\lim_{i \to \infty} \lambda_i = 0$; the third being that $\lim_{i \to \infty} N_{(p_i)}(0) = N_{(p)}(0)$; and the fourth being the following: If $i \in \{1, 2, ...\}$, then there is no $(Z_*, \mathcal{I}_*, \nu_*)$ data set from (5.1) that is described by the first two bullets Lemma 6.2 with $(Z_T, \mathcal{I}_T, \nu_T)$ defined using $T = (p_i, L_i, \lambda_i)$. With regards to the limit of the sequence $\{N_{(p_i)}(0)\}_{i \in \{1,2,...\}}$, keep in mind that if $\{p_i\}_{i \in \{1,2,...\}}$ is any sequence in $Z$ with limit $p$, then no limit of the sequence $\{N_{(p_i)}(0)\}_{i \in \{1,2,...\}}$ can be greater than $N_{(p)}(0)$. This is a consequence of Lemma 3.4.

    Use the sequence $\{T(i) = (p_i, L_i, \lambda_i)\}_{i \in \{1,2,...\}}$ as input for Proposition 4.1. Let $(Z_*, \mathcal{I}_*, \nu_*)$ denote the resulting limit data set and let $\{T(i)\}_{i \in \Lambda}$ denote the corresponding subsequence that is described by the fifth bullet of Proposition 4.1. The first four bullets of Proposition 4.1 say that this data set is described by (5.1) with $N_0$ being $N_{(p)}(0)$. It follows from the manner of convergence dictated by the fifth bullet of Proposition 4.1 that the first bullet and second bullets of Lemma 6.2 are obeyed when $T$ is any very large $i \in \Lambda$ version of $T(i)$. These conclusions constitute the desired nonsense that proves the first two bullets of Lemma 6.2 because the assumption at the outset was made to the effect that the first two bullets of Lemma 6.2 can not be satisfied by any $(Z_*, \mathcal{I}_*, \nu_*)$ data set from the $N_0 = N_{(p)}(0)$ version of (5.1) with $T$ being any $i \in \{1, 2, ...\}$ version of $T(i)$.

    *Part 2*: Suppose for the sake of argument that the third bullet or fourth bullet of Lemma 6.2 is false. If so, there is a number $\varepsilon \in (0, \frac{1}{100})$ and a sequence $\{(p_i, L_i, \lambda_i)\}_{i=1,2,...}$ with each member of this sequence again consisting of a point $p \in Z$, an oriented, orthonormal frame for $TX$ at $p$ and a number $\lambda_i \in (0, \frac{1}{100} r_0]$. There is also a second

sequence of points $\{p_i´\}_{i=1,2,...} \subset Z$; and these sequences have the following properties: The sequence $\{\lambda_i\}_{i=1,2,...}$ has limit zero; both of the sequences $\{p_i\}_{i \in \{1,2,...\}}$ and $\{p_i´\}_{i \in \{1,2,...\}}$ have limit p; and for each $i \in \{1, 2, ..., \}$ the distance between $p_i$ and $p_i´$ is in the interval $[\varepsilon \lambda_i, \varepsilon^{-1} \lambda_i]$. In addition, the sequence $\{ N_{(p_i)}(0) \}_{i \in \{1,2...\}}$ converges to $N_{(p)}(0)$ as does the sequence $\{ N_{(p_i´)}(0) \}_{i \in \{1,2,...\}}$ if the dimension of X is four. If the dimension of X is three, then the latter sequence converges to a number that is greater than $M_p$. And, last but crucial: For each $i \in \{1, 2, ...\}$, there is no $(Z_*, \mathcal{I}_*, \nu_*)$ data set from that $N_0 = N_{(p)}(0)$ version of (5.1) that is described by the first two bullets of Lemma 6.2 and either the third bullet or the fourth bullet (as the case may be) with $(Z_T, \mathcal{I}_T, \nu_T)$ defined using $T = (p_i, L_i, \lambda_i)$. The four steps that follow generate the requisite nonsense from such a sequence.

Step 1: Use the sequence $\{T(i) = (p_i, L_i, \lambda_i)\}_{i \in \{1,2,...\}}$ as input for Proposition 4.1. Let $(Z_*, \mathcal{I}_*, \nu_*)$ denote a limit data set that is supplied by Proposition 4.1 from this input sequence and let $\{T(i)\}_{i \in \Lambda}$ denote the subsequence of the input sequence that is described by the fifth bullet of Proposition 4.1. Nonsense ensues if the set $Z_*$ has an $|x| = 1$ point where the function $\tilde{N}_{(\cdot)}$ is equal to $N_{(p)}(0)$. This is because all three bullets of Lemma 6.2 will be then be satisfied with T being any sufficiently large $i \in \Lambda$ version of $T(i)$; and it was assumed at the outset that these three bullets are not satisfied by any data set from the $N_0 = N_{(p)}(0)$ version of (5.1) with T being any $i \in \Lambda$ version of $T(i)$. Steps 2–4 of the proof exhibit an $|x| = 1$ point in $Z_*$ where $\tilde{N}_{(\cdot)} = N_{(p)}(0)$.

Step 2: Supposing that $i \in \Lambda$, there is a point $x \in \mathbb{R}^d$ with $|x| \in [\varepsilon, \varepsilon^{-1}]$ that is mapped to $p_i´$ by the $T = T(i)$ version of the map $\phi_T$ that is described at the beginning of Section 4. (The map is the composition of first the scaling map $x \to \lambda_i x$ and then the Gaussian coordinate chart map that is defined by the frame L.) Use $x_i$ to denote this point. Note in particular that this is a point where $\nu_{T(i)} = 0$. Let $\Theta \subset \Lambda$ denote a subsequence with the property that $\{x_i\}_{i \in \Theta}$ converges and use $q_\diamond$ to denote the limit point. It follows from Item a) of the fifth bullet of Proposition 4.1 that $q_\diamond$ is a $\nu_* = 0$ point in the spherical annulus where $\varepsilon \leq |x| \leq \varepsilon^{-1}$. Let $q_* = q_\diamond/|q_\diamond|$. This is an $|x| = 1$ point where $\nu_* = 0$.

Step 3: The data set $(Z_*, \mathcal{I}_*, \nu_*)$ is described by (5.1) and so Lemma 5.4 or Lemma 5.5 as the case may be assigns a positive number to each $|x| = 1$ point in $Z_*$ and thus to the point $q_*$. This is the $q = q_*$ version of what these lemmas denote by $\tilde{N}_q$. To see how small this can be, fix $r \in (0, r_0]$ for the moment. The convergence dictated by Proposition 4.1 implies that

$$K_{q_\diamond}(r) = \lim_{i \in \Lambda} K_{(p_i)}(\lambda_i) / K_{(p_i´)}(\lambda_i r) ,$$

(6.1)

Meanwhile, the $p_i$ versions of (3.4) and the third bullet in Lemma 3.2 imply that

$$K_{(p_i')}(\lambda_i, 2r) \geq (1 - c_0 \lambda_i) \, 2^{N_{(p_i')}(0)} \, K_{(p_i')}(\lambda_i r).$$

(6.2)

Let $M_0$ denote the limit of the sequence $\{N_{(p_i')}(0)\}_{i \in \{1,2,\ldots\}}$, this being $N_0$ if the dimension of $X$ is four, and it is larger than the largest half integer less than $N_0$ if the dimension of $X$ is three. Since $K_{q_\diamond} > 0$ on $(0, \infty)$, the formula in (6.1) and the limit in (6.2) imply that

$$K_{q_\diamond}(2r) \geq 2^{M_0} K_{q_\diamond}(r) \;.$$

(6.3)

It then follows from Item c) of the third bullet of (5.1) that this same inequality holds with $q_*$ used in lieu of $q_\diamond$ to define the function K.

Step 4: The $q_*$ version of (6.3) runs afoul of (5.11) unless $N_{q_*}$ is equal to or greater than $M_0$ at some point in $[r, 2r]$. Keeping in mind that the conclusion of the preceding step hold for each small but positive r, and keeping in mind that $N_{q_*}$ is non-decreasing, it follows that $\lim_{r \to 0} N_{q_*}(r)$ is equal to or greater than $M_0$. This implies that $\tilde{N}_{q_*} \geq M_0$ because the function $N_{q_*}$ is continuous on $[0, \infty)$. In the case when X has dimension 3, this implies that $\tilde{N}_{q_*} \geq N_0$ because $\tilde{N}_{q_*}$ is half an integer. The implication $\tilde{N}_{q_*} \geq N_0$ follows directly when X has dimension 4 because $M_0 = N_0$ in this case. Since the second bullets of Lemmas 5.4 and 5.5 assert that $\tilde{N}_{q_*}$ can not be greater than $N_0$, it must equal $N_0$.

*Part 3*: To prove the fifth bullet, assume it to be false to generate nonsense. If the bullet is false, then there exists $\varepsilon > 0$ and sequences $\{(p_i, L_i, \lambda_i)\}_{i=1,2,\ldots}$ and $\{p_i'\}_{i=1,2,\ldots}$ and $\{p_i''\}_{i=1,2,\ldots}$ with the first two sequences as described at the beginning of Part 2 and with the third being a sequence in Z that converges to p with $\{N_{(p_i'')}(0)\}_{i \in \{1,2,\ldots\}}$ converging to a number, $M_0$ that is greater than $M_p$, and such that the following are true: For each index $i \in \{1, 2, \ldots\}$, the distance from $p_i$ to $p_i''$ is between $\varepsilon \lambda_i$ and $\varepsilon^{-1} \lambda_i$; and $\Delta_{p_i}(p_i', p_i'') \geq \varepsilon^2$. Furthermore, there is no $(Z_*, \mathcal{I}_*, \nu_*)$ data set from the $N_0 = N_{(p)}(0)$ version of (5.1) that is described by all the first three bullets and bullets four and five of Lemma 6.2 with T being any given $i \in \{1, 2, \ldots\}$ version of $(p_i, L_i, \lambda_i)$.

An argument much like that used in Part 2 uses the sequences $\{p_i'\}_{i=1,2,\ldots}$ and $\{p_i''\}_{i=1,2,\ldots}$ with what is said about convergence by Proposition 4.1 to obtain a pair of points in the $|x| = 1$ part of $Z_*$ where $\tilde{N}_{(\cdot)} = N_{(p)}(0)$ that lie on distinct lines through the origin. The argument in this case uses these sequences to obtain a two $|x| = 1$ points in $Z_*$ that are on different lines through the origin, one limit point having $\tilde{N}_{(\cdot)} = N_0$ and the other having Ñ$_{(0)}$ ≥ M$_0$ and thus greater than the largest half integer less than N$_0$. Since one of these points has Ñ$_{(·)}$ = N$_0$, the set Z$_*$ is described by the fourth bullet of Lemma 5.5. The other point must be in the interior of one of the half-planes that comprise Z$_*$. Since the function Ñ$_{(·)}$ takes half integer values on the interior of the half planes, its value at this other point must also be N$_0$. This implies that Z$_*$ is described by the fifth bullet of Lemma 5.5. This is nonsense since it runs afoul of the assumptions in the preceding paragraph.

**b) The values of the function N$_{(·)}$(0)**

The lemmas in this subsection summarize some observations about the values of the function N$_{(·)}$(0). The first lemma below asserts a lower bound.

**Lemma 6.4**: *The function N$_{(·)}$(0) on Z is nowhere less than $\frac{1}{2}$.*

*Proof of Lemma 6.4*: The upcoming Lemma 6.7 makes an assertion to the effect that the set of points in Z where N$_{(·)}$(0) is half an integer is dense in Z. Granted this lemma, the assertion follows from the first bullet of Lemma 3.4.

The next lemma implies that the function N$_{(·)}$(0) takes positive integer values on the set of isolated points in Z.

**Lemma 6.5**: *Fix p ∈ Z and suppose that there is a ball B ⊂ X centered at p with an isomorphism between 1 and the product bundle on B−(B∩Z). Then the following is true:*
- *N$_{(p)}$(0) is a positive integer.*
- *The isomorphism identifies ν with the restriction to B−(B∩Z) of a smooth section of $\mathbb{S}^+$ on B that is annihilated by the operator $\mathcal{D}$.*

*Proof of Lemma 6.5*: The isomorphism identifies ν with an $L^2_1$ section of $\mathbb{S}^+$ over B−(B∩Z) that vanishes on Z. This section extends over B∩Z as a continuous section of $\mathbb{S}^+$ on B by declaring it to be zero on B∩Z. The extended section is annihilated by $\mathcal{D}$ with the latter viewed as an operator mapping $L^2_1(B; \mathbb{S}^+|_B)$ to $L^2(B; \mathbb{S}^-|_B)$; and this being the case, standard elliptic regularity theory says that the extended section is a smooth section of $\mathbb{S}^+$ over B in the kernel of $\mathcal{D}$. Granted that this is so, then Taylor's theorem with remainder can be used to prove the assertion in the lemma's first bullet to the effect that N$_{(p)}$(0) is a positive integer.

By way of a reminder, a function point p ∈ Z is said to be strongly continuous if the function N$_{(·)}$(0) is constant on a neighborhood of p.

**Lemma 6.6**: *The function $N_{(\cdot)}(0)$ has half-integer values on the set of strongly continuous points in Z. If the dimension of X is three, then it also has half-integer values on the weakly continuous points.*

*Proof of Lemma 6.6*: If $p \in Z$ is an isolated point in Z, then the assertion follows from Lemma 6.5 Assuming that p is not isolated, then p is in any event weakly continuous. If the dimension of X is three, then the assertions follow from the first bullet of Lemma 6.3 and the fourth bullet of Lemma 5.4.

Now suppose that the dimension of X is four. Let $\{p_i\}_{i \in \{1,2,...\}}$ denote a sequence in Z–p with limit p. Look for a sequence $\{q_i\}_{i \in \{1,2,...\}} \subset Z-p$ that obeys the following:

- *If $i \in \{1, 2, ...\}$, then $\frac{1}{2} \text{dist}(p_i, p) \leq \text{dist}(q_i, p) \leq 2\,\text{dist}(p_i, p)$.*
- $\lim_{i \to \infty} \Delta_p(p_i, q_i) > 0$.

(6.4)

Either there is a sequence $\{q_i\}_{i \in \{1,2,...\}}$ obeying (6.4) or not.

CASE 1: If such a sequence exists, then it follows from the fifth bullet of Lemma 6.2 and the fourth bullet of Lemma 5.5 that $N_{(p)}(0)$ is a half-integer.

CASE 2: Suppose on the other hand that no such sequence exists. Let L denote a chosen, oriented orthonormal frame for TX at p and let $\{p_i\}_{i \in \{1,2,...\}} \subset Z-p$ denote a sequence that converges to p. Use the sequence $\{T(i) = (p, L, \lambda_i = \text{dist}(p_i, p))\}_{i \in \{1,2,...\}}$ as input for Proposition 4.1. Let $(Z_*, \mathcal{I}_*, \nu_*)$ denote a limit data set supplied by Proposition 4.1 and let $\{T(i)\}_{i \in \Lambda}$ denote the corresponding subsequence whose convergence is described by Proposition 4.1. Lemma 6.3 asserts that $(Z_*, \mathcal{I}_*, \nu_*)$ is described by the third bullet of Lemma 5.5. In particular, this implies that $Z_*$ is either a line through the origin or a union of half-planes with a common boundary line. This being the case, then the assumption about the values of the functions $\{\Delta_p(p_i, \cdot)\}_{i \in \{1,2,...\}}$ has the following consequence: Given $\varepsilon \in (0, 1]$, there exists a positive integer $\iota_\varepsilon$ such that if $i \in \Lambda$ and i is greater than $\iota_\varepsilon$, then the bundle $\mathcal{I}_{T(i)}$ is isomorphic to the product bundle on the complement of the radius $\varepsilon$ tubular neighborhood in the $|x| = 1$ sphere of the half-planes in $Z_*$. Given the convergence that is dictated by Proposition 4.1, this implies that $\mathcal{I}_*$ is isomorphic to the product bundle on the complement of $Z_*$. It follows as a consequence that $\nu_*$ can be viewed as a harmonic section of $\mathbb{V}^+$ on $\mathbb{R}^4 - Z_*$ whose norm extends as zero across $Z_*$. Standard elliptic theory can then be invoked to see that $\nu_*$ extends to the whole of $\mathbb{R}^4$ as a harmonic section of $\mathbb{V}^+$. It is therefore smooth and so it real analytic and has a non-trivial Taylor's expansion at the origin. The latter can be used to see that $N_0$ and hence $N_{(p)}(0)$ is an integer.

### c) The $N_{(\cdot)}(0) \in \frac{1}{2}\mathbb{Z}$ points are dense

The next lemma asserts that $N_{(\cdot)}(0)$ is a half of an integer on a dense set in Z.

**Lemma 6.7**: *Suppose that* $p \in Z$ *and that* $N_{(p)}(0)$ *is not half of a positive integer. Then there exists a sequence* $\{p_i\}_{i=1,2,\ldots} \subset Z$ *that converges to p with* $\lim_{i \to \infty} N_{(p_i)}(0)$ *equal to half of a positive integer.*

*Proof of Lemma 6.7*: Suppose that the assertion of the lemma is false so as to generate nonsense. The three parts that follow do this.

   *Part 1*: This part proves the following assertion:

   *Any counter example to Lemma 6.7 is a limit of weakly continuous points in Z.*
(6.5)

To prove (6.5), let p denote a point in Z that violates the assertion of Lemma 6.7 and let B denote a ball in X of radius less than $r_0$ and center p with no points from Z where $N_{(\cdot)}(0)$ is half an integer. Suppose in addition that there are no weakly continuous points from Z in B. Since Z is closed in X, the set Z∩B is closed in B. According to Lemma 6.5, the set Z∩B has no isolated points. This being the case, Z∩B is a *perfect* set. The set Z∩B is therefore uncountable because perfect set in $\mathbb{R}^n$ are uncountable. But this is nonsense because no point in Z∩B is weakly continuous and Lemma 6.1 asserts that the complement of the set of weakly continuous points is countable.

   *Part 2*: The assertion in (6.5) and Lemma 6.6 give the required nonsense when the dimension of X is three. This understood, assume henceforth that the dimension of X is four and let $p \in Z$ denote a weakly continuous counter example to Lemma 6.7. Fix an oriented, orthonormal frame for TX at p and parallel transport it along the geodesic arcs from p to define an orthonormal frame at the points in the radius $r_0$ ball centered at p. Denote this frame by L. Let $\{p_i\}_{i \in \{1,2,\ldots\}}$ denote a sequence of points in Z−p with distance less than $r_0$ from p, with $\lim_{i \to \infty} p_i = p$ and with $\lim_{i \to \infty} N_{(p_i)}(0) = N_{(p)}(0)$. Use the sequence $\{T(i) = (p, L, \text{dist}(p_i, p))\}_{i \in \{1,2,\ldots\}}$ as input for Proposition 4.1. It follows from Lemma 6.3 that any limit data set $(Z_*, \mathcal{I}_*, \nu_*)$ from Proposition 4.1 is described by the third bullet of Lemma 5.5. This implies in particular that $Z_*$ is a union of half-planes with a common boundary line. Note that there must be half-planes in $Z_*$. If $Z_*$ is a line, then $\mathcal{I}_*$ would be isomorphic to the product line bundle on $\mathbb{R}^d - Z_*$. If this is the case, then $\nu_*$ would extend over $Z_*$ as a harmonic section of $\mathbb{V}^+$. This extension would be real analytic and thus have a Taylor's expansion at the origin. The latter could then be used to see that $N_0$ and thus $N_{(p)}(0)$ is an integer.

*Part 3*: Suppose that $\{q_i\}_{i \in \{1,2,...\}}$ is another sequence in Z–p with limit p. If $\lim_{i \to \infty} N_{(q_i)}(0) = N_{(p)}(0)$, then the conditions in (6.4) can not be met for otherwise what is said in CASE 1 of the proof of Lemma 6.6 can be repeated to see that $N_{(p)}(0)$ is half an integer. On the other hand, there must be a sequence $\{q_i\}_{i \in \{1,2,...\}} \subset$ Z–p with limit p that does obey (6.4), for if not, then what is said in CASE 2 of the proof of Lemma 6.6 can be repeated to see that $N_{(p)}(0)$ is an integer. It follows from Lemma 3.4 that any sequence $\{q_i\}_{i \in \{1,2,...\}} \subset$ Z–p with limit p that obeys (6.4) is such that $\lim_{i \to \infty} N_{(q_i)}(0) < N_{(p)}(0)$.

*Part 4*: Let $N_\diamond$ denote the greatest half-integer less than $N_{(p)}(0)$. As explained directly, there is a sequence $\{q_i\}_{i \in \{1,2,...\}} \subset$ Z–p that obeys (6.4) with $\lim_{i \to \infty} N_{(q_i)}(0) \leq N_\diamond$.

To prove this claim, let $\{T(i) = (p, L, \lambda_i = \text{dist}(p, p_i))\}_{i \in \{1,2,...\}}$ denote the sequence from Part 2, let $(Z_*, \mathcal{I}_*, \nu_*)$ denote a limit data set from Proposition 4.1 for the input sequence $\{T(i)\}_{i \in \{1,2,...\}}$, and let $\{T(i)\}_{i \in \Lambda}$ denote the subsequence whose convergence is described in Proposition 4.1. Keep in mind in what follows that the value of Lemma 5.5's function $\tilde{N}_{(\cdot)}$ on the interior of a half-plane in $Z_*$ is a half-integer and it is less than $N_0$ which is $N_{(p)}(0)$.

Denote by $\{q_i\}_{i \in \Lambda}$ a sequence in Z–p with limit p that obeys (6.4). For each index $i \in \Lambda$, let $x_i \in \mathbb{R}^d$ denote the $\phi_{T(i)}$-inverse image of $q_i$. The sequence $\{x_i\}_{i \in \Lambda}$ has a subsequence that converges to a point in the interior of one of the half planes that comprise $Z_*$. Use $\vartheta$ to denote the set of indices in $\Lambda$ that label this subsequence. It follows from Lemma 3.2 and from the manner of convergence dictated by Proposition 4.1 that the lim sup of the numbers $\{N_{(q_i)}(0)\}_{i \in \vartheta}$ is no greater than the value of Lemma 5.5's function $\tilde{N}_{(\cdot)}$ on the interior of the limit half plane. The detailed argument is much like that in Step 3 of the proof of Lemma 6.2 and is it is left to the reader.

*Part 5*: Let p again denote a weakly continuous counter example to Lemma 6.7. Introduce by way of notion M to denote the number

$$M_p = \lim_{r \to 0} (\inf_{q \in B_r \cap Z} N_{(q)}(0)).$$

(6.6)

It follows from the fourth assertion of Lemma 3.2, and from what is said in Parts 1-4 that the number $M_p$ is a postive half-integer.

With the preceding understood, fix $r \in (0, r_0]$, and then fix $r_\diamond \in (0, r]$ such that if q is in $B_{r_\diamond} \cap Z$, then $M_q \geq M_p$. Then, fix $r_1 \in (0, \frac{1}{4} r_\diamond]$ such that the ball of radius $\frac{1}{4} r_1$ centered at p has a weakly continuous point in Z where $N_{(\cdot)}(0) < M_p + \frac{1}{2}$. Let $q_1$ denote a weakly continuous point from Z–p in the radius $\frac{1}{4} r_1$ ball centered at p with this property. Note that $M_{q_1}$ must be equal to $M_p$ because it can't be less than $M_p$ and it has to be a half-integer that is less than $N_{(q_1)}(0)$.

Because $N_{(q_1)}(0) < M_p + \frac{1}{2}$, it follows from Lemma 3.2 that there exists $r_2 \in (0, \frac{1}{4} r_1)$ such that all points in Z with distance $r_2$ or less from $q_1$ have $N_{(\cdot)}(0) < M_p + \frac{1}{2}$. Keeping in mind that $M_{q_1} = M_p$, choose a weakly continuous point in the radius $\frac{1}{4} r_2$ ball centered on the point $q_1$ where $N_{(\cdot)}(0) < M_p + \frac{1}{4}$. Denote such a point by $q_2$. It must also be the case that $M_{q_2} = M_p$.

Continue in this vein to choose a number $r_3$ and a point $q_3$, then a number $r_4$ and a point $q_4$ and so on so as to construct a sequence of numbers $\{r_i\}_{i \in \{1,2,...\}}$ and a sequence of weakly continuous points $\{q_i\}_{i \in \{1,2,...\}}$ in Z and that obey the following:

- If $k > 1$, then $r_k \in (0, \frac{1}{4} r_{k-1}]$
- If $k > 1$, then $N_{(\cdot)}(0) < M_p + 2^{-k}$ *on the radius $r_k$ ball centered on $q_{k-1}$.*
- If $k > 1$, then $\text{dist}(q_k, q_{k-1}) < \frac{1}{4} r_k$.

(6.7)

The sequence $\{q_k\}_{k \in \{1,2,...\}}$ converges to an $N_{(\cdot)}(0) = M_p$ point in $B_r \cap Z$. This is explained in the next part of the proof. Since r was chosen at will, the existence of such a limit constitutes the nonsense that proves Lemma 6.7.

*Part 6*: It follows from the first and third bullets of (6.7) that the sequence $\{q_i\}_{i \in \{1,2,...\}}$ is a Cauchy sequence. Since Z is closed, this sequence has a limit point. Denote the limit point by q. It follows from the first and third bullets of (6.7) that for each index k, the point q is in the radius $r_k$ ball centered on $q_{k-1}$. This being the case, it follows from the second bullet of (6.7) that $N_{(q)}(0) \leq M_p$. Meanwhile, it follows from the first and third bullets of (6.7) that q is in the radius $r_\diamond$ ball centered at p, and this implies that $M_q \geq M_p$. Since $M_q \leq N_{(q)}(0)$, it must be the case that $N_{(q)}(0) = M_p$.

### 7. The large scale structure of Z

This section talks about the large scale structure of Z. The upcoming Proposition 7.1 in Section 7a says that Z has Hausdorff dimension d - 2, and the two propositions in Section 7b say that Z looks to a first approximation like a stratified space with the top stratum having dimension d - 2 and the other strata having smaller, integer dimensions. Proposition 7.1 and two propositions in Section 7b restate the assertions of Theorem 1.3.

### a) The Hausdorff dimension of Z

The proposition that follows gives an upper bound for Z's Hausdorff dimension.

**Proposition 7.1**: *The set Z has Hausdorff dimension at most 1 when the dimension of X is three; and it has Hausdorff dimension at most 2 when the dimension of X is four. The set of points in Z where $N_{(\cdot)}(0)$ is not half an integer has Hausdorff dimension zero in the case when the dimension of X is three, and it has Hausdorff dimension at most 1 in the case when the dimension of X is four.*

*Proof of Proposition 7.1*: The proof of has six parts. The first two parts prove the assertion for the case when X has dimension three and the remaining parts prove the assertion for the case when X has dimension four.

What follows is a reminder about Hausdorff dimension. Supposing that $\mu \geq 0$, a given set K has Hausdorff dimension at most $\mu$ when the following is true. Fix numbers $\upsilon > \mu$ and $\delta > 0$ and $\rho > 0$. Then K has an open cover, this denoted by $\mathfrak{J}$, whose constituent sets have diameter at most $\rho$ and obey $\sum_{U \in \mathfrak{J}_\rho} (\mathrm{diam}(U))^\upsilon < \delta$

*Part 1*: Suppose in what follows that $\dim(X) = 3$. The set Z in this case can be written as a finite disjoint union $(\cup_{m \in \{1,\ldots\}} \mathfrak{Z}_m) \cup Z^{\mathrm{nwc}}$ with $Z^{\mathrm{nwc}}$ denoting the subset of points that are not weakly continuous and with any given $m \in \{1, \ldots\}$ version of $\mathfrak{Z}_m$ being the set of weakly continuous points in Z where $N_{(\cdot)}(0)$ is equal to $\frac{1}{2} m$. Such a decomposition follows from Lemma 6.6. Since it follows from Lemma 6.1 that $Z^{\mathrm{nwc}}$ has Hausdorff dimension zero, it is enough to prove that any given $m \in \{1, 2, \ldots\}$ version of $\mathfrak{Z}_m$ has Hausdorff dimension at most 1.

To this end, fix $\varepsilon \in (0, \frac{1}{1000})$. If $p \in \mathfrak{Z}_m$, then there exists $r_{p,\varepsilon} \in (0, \varepsilon r_0)$ with the following two properties: The conclusions of Lemma 6.2 hold using $r = r_{p,\varepsilon}$; and the $r = r_{p,\varepsilon}$ and $\varepsilon$ version of $Z_{r,\delta}$ contains a point with distance $r_{p,\varepsilon}$ from p.

Supposing that $L$ is a given oriented, orthonormal frame for TX at p, use $T$ to denote the data set $(p, L, \lambda = r_{p,\varepsilon})$. Lemma 6.2 takes this data set as input and then supplies a data set $(Z_*, \mathcal{I}_*, \nu_*)$ with the properties that are laid out in the first three bullets of Lemma 6.2. Note in particular that $Z_*$ is a line through the origin.

*Part 2*: Given $p \in \mathfrak{Z}_m$, let $B(p)$ denote the radius $r_{p,\varepsilon}$ ball centered at p. The set $\mathfrak{Z}_m$ has a locally finite cover by sets from the collection $\{B(p): p \in \mathfrak{Z}_m\}$. It follows as a consequence that it is enough to prove that each $p \in \mathfrak{Z}_m$ version of $\mathfrak{Z}_m \cap B(p)$ has Hausdorff dimension at most 1. It follows from Lemma 6.2 that $\mathfrak{Z}_m \cap B(p)$ can be covered by $10\varepsilon^{-1}$ balls of radius $\varepsilon r_{p,\varepsilon}$. Let $\mathfrak{U}_{p,1}$ denote such a cover. Let B denote a ball from $\mathfrak{U}_{p,1}$. It follows from Lemma 6.2 that $\mathfrak{Z}_m \cap B$ can be covered by $10\varepsilon^{-1}$ balls of radius $\varepsilon^2 r_{p,\varepsilon}$. By way of an explanation, of B is a ball of radius $\varepsilon r_{p,\varepsilon}$ centered at a point in $\phi_T(Z_*) \cap B(p)$, and if B has more than 2 points from $\mathfrak{Z}_m$, then Lemma 6.2 can be invoked to see that $Z \cap B$, and $\mathfrak{Z}_m \cap B$ in particular, has distance at most $\varepsilon^2 r_\varepsilon$ from a line in B.

Supposing that $B \in \mathfrak{U}_{p,1}$, fix a set of $10\varepsilon^{-1}$ balls of radius $\varepsilon^2 r_{p,\varepsilon}$ that cover $\mathfrak{Z}_m \cap B$. Let $\mathfrak{U}_{p,2}$ denote the collection of balls that comprise all of these covers. There are at most $100\varepsilon^{-2}$ balls in $\mathfrak{U}_{p,2}$. By the same token, if $B \subset \mathfrak{U}_{p,2}$, then Lemma 6.2 can be used again to see that $\mathfrak{Z}_m \cap B$ has a cover by at most $10\varepsilon^{-1}$ balls of radius $\varepsilon^3 r_{p,\varepsilon}$. Fix such a cover for each ball in $\mathfrak{U}_{p,2}$ and let $\mathfrak{U}_{p,3}$ denote the collection of balls that comprise all of these covers. This set has at most $(10\varepsilon^{-1})^3$ balls. Continue sequentially in this vein to define covers $\{\mathfrak{U}_{p,k}\}_{k=1,2,\ldots}$ of $\mathfrak{Z}_m \cap B(p)$ with each $k \in \{1, 2, \ldots\}$ version of $\mathfrak{U}_{p,k}$ being a set of at most $(10\varepsilon^{-1})^k$ balls with each ball having radius $\varepsilon^k r_{p,\varepsilon}$.

Fix $\upsilon > 0$. Then fix $k \in \{1, 2, \ldots\}$. Each ball $B \subset \mathfrak{U}_{p,k}$ has diameter at most $\varepsilon^k r_{p,\varepsilon}$ and there are at most $(10\varepsilon^{-1})^k$ such balls. Therefore,

$$\sum_{B \in \mathfrak{U}_{p,k}} (\operatorname{diam}(B))^\upsilon \leq (10\varepsilon^{-1})^k \varepsilon^{k\upsilon} r_{p,\varepsilon}^\upsilon.$$

(7.1)

This sum is no greater than $(10\varepsilon^{\upsilon-1})^k r_{p,\varepsilon}^\upsilon$. In particular, if $\upsilon > 1$ and if $\varepsilon < 10^{-1/(\upsilon-1)}$, then the $k \to \infty$ limit of the left hand side of (7.1) is zero.

*Part 3*: Assume henceforth that X has dimension 4. Let $Z^{wc}$ denote the set of points in Z that are weakly continuous. The complement in Z of $Z^{wc}$ is countable so it has Hausdorff dimension 0. This being the case, the assertion that Z has Hausdorff dimension at most 2 follows if $Z^{wc}$ has Hausdorff dimension at most 2.

To prove this assertion about $Z^{wc}$, it proves useful to write $Z^{wc}$ as the disjoint union $Z^p \cup Z^{np}$ where membership in $Z^p$ is based on the criteria that follows directly. A point p is in $Z^p$ if there exist $\varepsilon > 0$ and a sequence $\{(p_i, p_i´)\}_{i \in \{1,2,\ldots\}}$ of pairs of points in Z obeying the following:

- *Both sequences $\{p_i\}_{i \in \{1,2,\ldots\}}$ and $\{p_i´\}_{i \in \{1,2,\ldots\}}$ are in the radius $r_0$ ball centered at p and both converge to p.*
- $\varepsilon < \operatorname{dist}(p, p_i)/\operatorname{dist}(p, p_i´) < \varepsilon^{-1}$.
- $\Delta_p(p_i, p_i´) > \varepsilon$.
- $\lim_{i \to \infty} N_{(p_i)}(0) = N_{(p)}(0)$ and $\lim_{i \to \infty} N_{(p_i´)}(0) = N_{(p)}(0)$.

(7.2)

It follows from Lemma 6.2 that $N_{(\cdot)}(0)$ takes only half integer values on $Z^p$. This understood, let m denote a given integer and let $\mathfrak{Z}^p_m$ denote the $N_{(\cdot)}(0) = \frac{1}{2} m$ part of $Z^p$.

*Part 4*: Fix a positive integer, m. This part and Part 5 of the proof explain why $\mathfrak{Z}^p_m$ has Hausdorff dimension at most 2. To start, fix $\varepsilon \in (0, \frac{1}{100})$ and fix $p \in \mathfrak{Z}^p_m$. There exists $r_{p,\varepsilon} \in (0, \varepsilon r_0)$ with the following two properties: The conclusions of Lemma 6.2 hold with $r = r_{p,\varepsilon}$; and the $r = r_{p,\varepsilon}$ and $\varepsilon$ version of $Z_{r,\delta}$ contains a pair of points $(p´, p´´)$ with distance between $\varepsilon r_{p,\varepsilon}$ from p and $\varepsilon^{-1} r_{p,\varepsilon}$ from p such that $\Delta_p(p´, p´´) > \varepsilon$.

Let L denote a given oriented orthonormal frame at $TX|_p$ and let T denote the data set $(p, L, \lambda = r_{p,\varepsilon})$. Lemma 6.2 takes this data set as input and then supplies a data set $(Z_*, \mathcal{I}_*, \nu_*)$ with the properties that are laid out in Lemma 6.2. Note in particular that $Z_*$ is a plane through the origin in $\mathbb{R}^4$ and that the whose of the $|x| \leq 10$ part of $Z_T$ has distance at most $\varepsilon$ from this plane.

Let B(p) denote the radius $r_{p,\varepsilon}$ ball centered at p. The set $\mathfrak{Z}^p_m$ has a locally finite cover by balls from the set $\{B(p): p \in \mathfrak{Z}^p_m\}$ and so it is enough to prove that each $p \in \mathfrak{Z}^p_m$ version of $\mathfrak{Z}^p_m \cap B(p)$ has Hausdorff dimension at most 2.

*Part 5*: Suppose that $m \in \{1, 2, \ldots\}$ and that p is a point in $\mathfrak{Z}^p_m$. It follows from what is said by Lemma 6.2 that each point in $\mathfrak{Z}^p_m \cap B(p)$ has distance at most $\varepsilon$ from the $\phi_T$ image of a plane through the origin in $\mathbb{R}^4$. There exists therefore a number $c > 1$ and a collection $\mathfrak{U}_{p,1}$ of at most $c\varepsilon^{-2}$ balls of radius $\varepsilon r_{p,\varepsilon}$ covering $\mathfrak{Z}^p_m \cap B(p)$.

Let $B \in \mathfrak{U}_{p,1}$. It follows from what is said by Lemma 6.2 that $\mathfrak{Z}^p_m \cap B$ has in turn a cover by at most $c\varepsilon^{-2}$ balls of radius $\varepsilon^2 r_{p,\varepsilon}$ because the points in $Z \cap B$ with $N_{(\cdot)}(0) = \frac{1}{2}m$ have distance at most $\varepsilon^2 r_{p,\varepsilon}$ from the image via a Gaussian coordinate chart map of either a line in $\mathbb{R}^4$ or a plane in $\mathbb{R}^4$. Fix such a cover for each ball in $\mathfrak{U}_{p,1}$ and denote the collection of balls from these covers by $\mathfrak{U}_{p,2}$. This is a set of at most $(c\varepsilon^{-2})^2$ balls. Continue sequentially by mimicking what is in Part 2 to define successive covers $\{\mathfrak{U}_{p,k}\}_{k=1,2,\ldots}$ of $\mathfrak{T}_{2\varepsilon,m} \cap B(p)$ with any given $k \in \{1, 2, \ldots\}$ version of $\mathfrak{U}_{p,k}$ containing at most $(c\varepsilon^{-2})^k$ balls of radius $\varepsilon^k r_{p,\varepsilon}$.

Fix $\upsilon > 0$ and then fix $k \in \{1, 2, \ldots\}$. The analog of (6.8) in this case says

$$\sum_{B \in \mathfrak{U}_{p,k}} (\mathrm{diam}(B))^\upsilon \leq (c\varepsilon^{\upsilon-2})^k r_{p,\varepsilon}^\upsilon.$$

(7.3)

It follows from this inequality that $\lim_{k\to\infty} \sum_{B \in \mathfrak{U}_{p,k}} (\mathrm{diam}(B))^\upsilon = 0$ if $\upsilon > 2$ and if $\varepsilon < c^{1/(2-\upsilon)}$.

*Part 6*: This last part of the proof explains why the set of points in $Z^{np}$ has Hausdorff dimension at most 1. To do this, fix positive integers m and k, and then define $\mathfrak{Z}^{np}_{m,k}$ to be the subset of points in $Z^{np}$ where $N_{(\cdot)}(0) \in (m-1+2^{-k}, m]$. It is sufficient to prove that each such $\mathfrak{Z}^{np}_{m,k}$ has Hausdorff dimension at most 1.

Fix $\varepsilon \in (0, \frac{1}{100} 2^{-k})$ and fix $p \in \mathfrak{Z}^{np}_{m,k}$. There exists $r_{p,\varepsilon} \in (0, \varepsilon r_0)$ such that the conclusions of Lemma 6.2 hold for $r = r_{p,\varepsilon}$; and the $r = r_{p,\varepsilon}$ version of $Z_{r,\delta}$ contains a point with distance $r_{p,\varepsilon}$ from p. Supposing that L is a chosen, oriented orthonormal frame for $TX|_p$, let T again denote the data set $(p, L, r_{p,\varepsilon})$. Lemma 6.2 supplies a version of $(Z_*, \mathcal{I}_*, \nu_*)$ where $Z_*$ in this case is either a line through the origin or a non-empty, finite union of half planes through the origin with a common boundary line. Let $Z_{**}$ denote $Z_*$ if $Z_*$ is a line, or the common boundary line of the half planes that comprise $Z_*$ if $Z_*$ is not a line. The key observation here follows from the fifth bullet of Lemma 6.2, this being that all points the points in the radius $r_{p,\varepsilon}$ ball from $\mathfrak{Z}^{np}{}_{m,k}$ have distance at most $\varepsilon$ from the $\phi_T$ image of the line $Z_{**}$. Given that this is so, the arguments from Part 2 can be repeated with only cosmetic changes to prove that $\mathfrak{Z}^p{}_{m,k}$ has Hausdorff dimension at most 1.

**b) The topology of Z**

The two propositions in this section describe what Z looks like when viewed through glasses or a microscope with some finite magnification or resolving power. To set the notation for the first proposition, suppose that X has dimension 3. Let p denote a point in X. A *geodesic arc* through p is the image via a Gaussian coordinate chart map centered at p of the $|x| \le r_0$ part of a line through the origin in $\mathbb{R}^3$. A geodesic arc in a given open ball is said to be *proper* if it is a relatively closed subset of the ball.

**Proposition 7.2**: *Assume that* X *has dimension 3. Given numbers* $\theta \in (0,1)$ *and* $z > 1$ *there are finite sets* $\mathfrak{U}$ *and* $\mathfrak{V}$ *of the following sort:*
- *Both sets are composed of balls of radius at most* $z^{-1}$.
- *The balls from* $\mathfrak{U}$ *are pairwise disjoint and* $\sum_{B \in \mathfrak{U}} (\mathrm{radius}(B))^\theta < z^{-1}$.
- *The balls in the set* $\mathfrak{V}$ *cover the part of* Z *in* X$-(\cup_{B \in \mathfrak{U}} B)$.
- *Supposing that* B *is a ball from* $\mathfrak{V}$, *let* p *denote its center point and* r *its radius. The part of* Z *in* B *is contained in the radius* $z^{-1}r$ *tubular neighborhood of a proper geodesic arc through* p.

This proposition is proved in Section 7c.

The next proposition in this section is a dimension 4 analog of the previous one. This proposition also uses the notion of a geodesic arc through a given point and also the notions of *geodesic disks* and *union of geodesic half-disks with a common boundary arc* through a given point. Supposing that p is the point in question, a geodesic arc through p is the image in X via a Gaussian coordinate chart map centered at p of the $|x| \le r_0$ part of a line through the origin in $\mathbb{R}^4$; and a geodesic disk through p is the corresponding image of the $|x| \le r_0$ part of a plane through the origin in $\mathbb{R}^4$. A union of geodesic half-disk through p with a common boundary line is the $|x| \le r_0$ part of the union of two or more half-planes through the origin in $\mathbb{R}^4$ with a common boundary line. These sets are said to be *proper* in a given open ball if they are relatively closed sets in the ball.

**Proposition 7.3**: *Assume that* X *has dimension 4. Given numbers* $\theta \in (0,1)$ *and* $z > 1$, *there are finite sets* $\mathfrak{U}$ *and* $\mathfrak{V}$ *of the following sort:*
- *All sets are composed of balls of radius at most* $z^{-1}$.

- *The balls from $\mathfrak{U}$ are pairwise disjoint and $\sum_{B \in \mathfrak{U}} (\text{radius}(B))^\theta < z^{-1}$.*
- *The balls in the set $\mathfrak{V}$ cover the part of Z in $X - (\cup_{B \in \mathfrak{U}} B)$.*
- *Supposing that B is a ball from $\mathfrak{V}$, let p denote its center point and r its radius. The part of Z in B is contained in the radius $z^{-1}r$ tubular neighborhood of either a proper geodesic arc through p, a proper geodesic disk through p or a proper union of geodesic half-disks through p with a common boundary arc.*

This proposition is proved in Section 7d.

### c) Proof of Proposition 7.2

The proof of the proposition has five parts. The first two parts supply some facts that are used in the subsequent parts of the proof.

*Part 1*: This part makes various observations about the set of weakly continuous points in Z. According to Lemma 6.6, the function $N_{(\cdot)}(0)$ has only half integer values on this set. This being the case, and supposing that m is a positive integer, let $\mathfrak{Z}_m$ denote the $N_{(\cdot)}(0) = \frac{1}{2} m$ part of the set of weakly continuous points in Z.

Given $\varepsilon \in (0, \frac{1}{100})$ and $p \in \mathfrak{Z}_m$, define $r_{p,\varepsilon} \in (0, \varepsilon r_0)$ as done in Part 1 of Section 7a. There exists $r_p \in (\frac{1}{2} r_{p,\varepsilon}, r_{p,\varepsilon})$ such that the boundary of the ball of radius $r_p$ centered at p contains only weakly continuous points from Z. That this is so follows from the fact that the complement in Z of the set of weakly continuous points is countable but the set of choices for $r_p$ is not. With $r_p$ chosen, use B(p) to denote the ball of radius $r_p$ centered at p. Keep in mind that $r_p$ and B(p) depend implicitly on $\varepsilon$.

The collection $\{B(p): p \in \mathfrak{Z}_m\}$ defines an open cover of $\mathfrak{Z}_m$. What follows directly describes what turns out to be a useful subcover. Given a positive integer b, introduce by way of notation $\mathfrak{Z}_{m,b}$ to denote the set of points in $\mathfrak{Z}_m$ with $r_p \in [2^{-b}\varepsilon, 2^{-b+1}\varepsilon)$. Choose a maximal subset in $\mathfrak{Z}_{m,b}$ to be denoted by $\mathfrak{X}_{m,b}$, with the following property: The distance between any two distinct points in $\mathfrak{X}_{m,b}$ is no smaller than $2^{-b-4}\varepsilon$. The set $\mathfrak{X}_{m,b}$ is by necessity a finite set; and it follows from the definitions that the collection of balls $\mathfrak{V}_{m,b} = \{B(p): p \in \mathfrak{X}_{m,b}\}$ is an open cover of $\mathfrak{Z}_{m,b}$.

Since $\mathfrak{Z}_m = \cup_{b \in \{1,2,\ldots\}} \mathfrak{Z}_{m,b}$, the collection $\cup_{b \in \{1,2,\ldots\}} \mathfrak{V}_{m,b}$ is an open cover of $\mathfrak{Z}_m$. Winnow this set by using the following rule: Supposing that $b \in \{2, 3, \ldots\}$, remove all balls from $\mathfrak{V}_{m,b}$ that lie in the union of the balls from the collection $\cup_{a \in \{1,2,\ldots,b-1\}} \mathfrak{V}_{m,a}$. Use $\mathfrak{V}_{m,b*}$ to denote this possibly smaller set of balls and let $\mathfrak{V}_{m*}$ denote $\cup_{b \in \{1,2,\ldots\}} \mathfrak{V}_{m,b*}$ with it understood that $\mathfrak{V}_{m,1*} = \mathfrak{V}_{m,1}$.

*Part 2*: This part of the proof makes some observations about the complement in Z of the set of weakly continuous points. The first observation concerns the points where $N_{(\cdot)}(0)$ is half an integer that are not weakly continuous. To state the observation, suppose p ∈ Z, that p is not weakly continuous, and that $N_{(p)}(0) = \frac{1}{2} m$ with m being a positive integer. Granted all of this, then there exists $r_p \in (0, \frac{1}{100} r_0)$ such that the ball of radius $r_p$ centered at p contains no points in Z with $N_{(\cdot)}(0) \geq \frac{1}{2} m$. This follows from Lemma 3.4.

The second obervation concerns the points in Z where $N_{(\cdot)}(0)$ is not half an integer. Lemma 6.6 says that no such point is weakly continuous. The promised observation about this set of points is stated below by Lemma 7.4. To set notation, define $Z_m$ to be the subset of Z where $N_{(\cdot)}(0) \geq \frac{1}{2} m$; and given $\upsilon \in (0, \frac{1}{2})$, define $Z_{m,\upsilon}$ to be the subset of $Z_m$ where $N_{(\cdot)}(0) > \frac{1}{2} m + \upsilon$. It follows from Lemma 3.4 that $Z_m$ is a closed subset of Z and it follows from Lemmas 6.6 and 3.4 that any given $\upsilon \in (0, \frac{1}{2})$ version of $Z_{m,\upsilon}$ is also a closed subset of Z.

**Lemma 7.4**: *Fix* $m \in \{1, 2, ...\}$ *and* $\upsilon \in (0, \frac{1}{2})$.
- *Each limit point of* $Z_{m,\upsilon} - Z_{m+1}$ *is in* $Z_{m+1}$.
- *Supposing that* $\delta > 0$, *there are at most a finite number of points in* $Z_{m,\upsilon} - Z_{m+1}$ *with distance greater* $\delta$ *from* $Z_{m+1}$.

*Proof of Lemma 7.4*: Assume that the first bullet is false to derive nonsense. If the first bullet is false, there exists $p \in Z_{m,\upsilon}$ with $N_{(p)}(0) < m+1$ and a sequence $\{p_i\}_{i \in \{1,2,...\}} \subset Z_{m,\upsilon}$ with limit p. Nothing is lost by assuming that this sequence is in the radius $\frac{1}{1000} r_0$ ball centered at p. Fix an oriented, orthonormal frame for TX at p to be denoted by L; and for each $i \in \{1, 2, ...\}$, let T(i) denote the data set $(p, L, \lambda_i = \text{dist}(p, p_i))$. Use the sequence $\{T(i)\}_{i \in \{1,2,...\}}$ as input for Proposition 4.1 and let $(Z_*, \mathcal{I}_*, \nu_*)$ denote the corresponding output. This data set is described by Lemma 5.4. Keep in mind the following fact from Lemma 5.4: The function $\tilde{N}_{(\cdot)}$ on the set of $|x|=1$ points in $Z_*$ has only half-integer values with no value greater than $N_{(p)}(0)$. Since $N_{(p)}(0)$ is not itself half and integer, it follows that $\tilde{N}_{(\cdot)}$ is strictly less than $N_{(p)}(0)$. The tension between $\tilde{N}_{(\cdot)}$ being half an integer and it being less that $N_{(p)}(0)$ leads momentarily to the desired nonsense.

Let $\Lambda \subset \{1, 2, ...\}$ denote the subsequence that is described by Proposition 4.1. It follows from the manner of convergence dictated by Proposition 4.1 that the corresponding subsequence $\{\phi_{T(i)}^{-1}(q_i)\}_{i \in \Lambda}$ has a subsequence that converges to an $|x|=1$ point in $Z_*$. Let $q_\diamond$ denote a limit point of this sequence. Arguments much like those used in Step 3 of the proof of Lemma 6.2 can be used to prove that the value of $\tilde{N}_{(\cdot)}$ at $q_\diamond$ is no smaller than the lim-inf of the sequence $\{N_{(p_i)}(0)\}_{i \in \Lambda}$; and this in turn is no less than $\frac{1}{2} m + \upsilon$. Since $\tilde{N}_{(\cdot)} < N_{(p)}(0)$ in any event, this means that $\tilde{N}_{(\cdot)}$ at $q_\diamond$ is not half an integer. This last conclusion is the desired nonsense.

*Part 3*:  Let N be the positive integer such that $Z_N \neq \emptyset$ but $Z_{N+1} = \emptyset$.  Fix for the moment $\varepsilon_N \in (0, \frac{1}{100} z^{-1})$.  If $\mathfrak{Z}_N \neq \emptyset$, define the collection of balls $\mathfrak{V}_{N*}$ using $\varepsilon = \varepsilon_N$ according to the rules in Part 1.  Let $\mathcal{V}_N$ denote the union of the balls in $\mathfrak{V}_{N*}$.

It is a consequence of Lemmas 7.4 and 3.4 that there are at most a finite set of points in any $\upsilon \in (0, \frac{1}{2})$ version of $Z_{N,\upsilon}$.  This implies in turn that the $N_{(\cdot)}(0) > \frac{1}{2} N$ points in $Z_N - (Z_N \cap \mathcal{V}_N)$ can be labeled as $\{p_1, p_2, \ldots\}$ with the labeling such that $N_{(\cdot)}(0)$ does not increase as the label increases.  Fix a ball of radius at most $\frac{1}{100} 2^{-1} \varepsilon_N^{\;2}$ around $p_1$ so that each point from Z on its boundary is weakly continuous points from Z.  Denote this ball by $B(p_1)$.  Let $n(2)$ denote the smallest positive integer such that $p_{n(2)}$ is not in $B(p_1)$.  Fix a ball of radius no greater than $\frac{1}{100} 2^{-2} \varepsilon_N^{\;2}$ centered on $p_{n(2)}$ so that its closure is disjoint from that of $B(p_1)$ and so that each point from Z on its boundary is weakly continuous.  Denote this ball by $B(p_{n(2)})$.  Let $n(3)$ denote the smallest postive integer such that $p_{n(3)}$ is not in $B(p_1)$ nor in $B(p_{n(2)})$.  Fix a ball of radius no greater than $\frac{1}{100} 2^{-3} \varepsilon_N^{\;2}$ around $p_{n(3)}$ so that its closure is disjoint from the closures of $B(p_1)$ and $B(p_{n(2)})$ and so that each point from Z is weakly continuous.  Continuing in this vein sequentially to define points $\{p_{n(1)=1}, p_{n(2)}, \ldots\}$ and the corresponding set of balls $\{B(p_{n(k)})\}_{k \in \{1,2,\ldots\}}$ with pairwise disjoint closures.  Let $\mathfrak{X}_{N+}$ denote the set $\{p_{n(1)}, \ldots\}$.

As explained directly, the set $\mathfrak{X}_{N+}$ is necessarily a finite set.  To see why, suppose that p is a limit point of this set.  The point p can not be an $N_{(\cdot)}(0) = \frac{1}{2} m$ point because of what is said in the first paragraph of Part 2 and because no point in $\mathfrak{X}_{N+}$ is in $\mathcal{V}_N$.  Since $N_{(p)}(0) > \frac{1}{2} m$, it must be that $p = p_i$ for some index i.  But if $p = p_i$, then there exists a positive integer j such that $p \in B(p_{n(j)})$.  Since the points $\{p_{n(j+1)}, \ldots\}$ are not in $B(p_{n(j)})$, the point p can not be a limit point of the sequence.

Let $\mathfrak{X}_{N0}$ denote the set of $N_{(\cdot)}(0) = \frac{1}{2} m$ points in Z that are not weakly continuous and not in any $p \in \mathfrak{X}_{N+}$ version of B(p).  It follows from what is said in the first paragraph of Part 2 that this is also a finite set.  Let $n_{N0}$ denote the number of elements in $\mathfrak{X}_{N0}$.  Supposing that $p \in \mathfrak{X}_{N0}$, reintroduce the number $r_p$ from the first paragraph of Part 2.  Let B(p) denote a ball of radius at most $(100 n_{N0})^{-1} \varepsilon_N^{\;2} r_p$ centered at p with only weakly continuous points from Z on its boundary.  Note that the closures of B(p) and $\mathcal{V}_N$ are disjoint as are the closures of B(p) and each $p' \in \mathfrak{X}_{N+}$ version of B(p').  Introduce by way of notation $\mathfrak{U}_{N\diamond}$ to denote the collection of balls $\{B(p): p \in \mathfrak{X}_{N0} \cup \mathfrak{X}_{N+}\}$.

*Part 4*:  The collection $\mathfrak{U}_{N\diamond}$ is finite, but this is not necessarily so for $\mathfrak{V}_{N*}$.  Winnow $\mathfrak{V}_{N*}$ by deleting any ball in this set that is entirely inside some ball from $\mathfrak{U}_{N\diamond}$.  Denote the result of this perhaps smaller set by $\mathfrak{V}_{N\diamond}$.  The next paragraph explains why the collection $\mathfrak{V}_{N\diamond}$ is a finite.

Suppose for the sake of argument that $\mathfrak{V}_{N\diamond}$ is not finite to derive nonsense.  Fix an infinite sequence of points $\{q_i\}_{i=1,2,\ldots}$ with each member being the center point of a ball from $\mathfrak{V}_{N\diamond}$.  This sequence must converge to a point in $Z_N$.  Thus, it is contained in a ball from either $\mathfrak{U}_{N*}$ or $\mathfrak{V}_{N*}$. Use B to denote such a ball. Since the radii of the balls from the collection $\{B(q_i)\}_{i \in \{1,2,...\}}$ limit to zero, each sufficiently large i version of $B(q_i)$ lie entirely in B. But this is nonsense, because if B is from $\mathfrak{V}_{N*}$ then these large i versions of $B(q_i)$ would not be in $\mathfrak{V}_{N*}$; and if B is from $\mathfrak{U}_{N*}$, then they would not be in $\mathfrak{U}_{N◊}$.

It follows from the construction that the collections $\mathfrak{V}_{N◊} \cup \mathfrak{U}_{N◊}$ defines an open cover of $Z_N$; and it likewise follows that if $\varepsilon_N \leq c_0^{-1} z^{-1}$ and if $Z = Z_N$, then the collections $\mathfrak{U} = \mathfrak{U}_{N◊}$ and $\mathfrak{V} = \mathfrak{V}_{N◊}$ satisfy the requirements of Proposition 7.2.

*Part 5*: Assuming that Z is bigger than $Z_N$, this part of the proof describes an induction argument that adds balls to the sets $\mathfrak{U}_{N◊}$ and $\mathfrak{V}_{N◊}$ so as to obtain collections that satisfy the requirements of Proposition 10.2. To set up the induction argument, suppose that $m \in \{2, 3, ..., N\}$ and that $\mathfrak{U}_m$ and $\mathfrak{V}_m$ are finite collections of balls such that if $\mathfrak{U}$ is set equal to $\mathfrak{U}_m$ and $\mathfrak{V}$ is set equal to $\mathfrak{V}_m$ then the pair $\mathfrak{U}$ and $\mathfrak{V}$ obey the requirements of Proposition 10.2 with the wording changed so that Z is everywhere replaced by $Z_m$. The following condition is also assumed: Every point in Z on the boundary of the closure of any ball from $\mathfrak{U}_m \cup \mathfrak{V}_m$ is weakly continuous. If $m = N$, then $\mathfrak{U}_{m=N}$ can be taken to be $\mathfrak{U}_{N◊}$ and $\mathfrak{V}_{m=N}$ can be taken to be $\mathfrak{V}_{N◊}$. The two steps that follow explain how to construct sets $\mathfrak{U}_{m-1}$ and $\mathfrak{V}_{m-1}$ so that the $\mathfrak{U} = \mathfrak{U}_{m-1}$ and $\mathfrak{V} = \mathfrak{V}_{m-1}$ version of Proposition 7.2 holds if Z is everywhere replaced by $Z_{m-1}$.

<u>Step 1</u>: Fix $\varepsilon_{m-1} \in (0, \frac{1}{100} z^{-1})$ for the moment. Use $\varepsilon = \varepsilon_{m-1}$ to define the collection of sets $\mathfrak{V}_{(m-1)*}$ using the rules set forth in Part 1. Let $\mathcal{V}_{m-1}$ denote the union of the balls from the three collections $\mathfrak{V}_{(m-1)*}, \mathfrak{U}_m$ and $\mathfrak{V}_m$. It follows from Lemmas 3.4 and 7.4 that there are at most a finite number of points in $Z_{m-1} - (Z_{m-1} \cap \mathcal{V}_{m-1})$ with $N_{(\cdot)}(0) > \frac{1}{2}m + \upsilon$ for any given $\upsilon \in (0, \frac{1}{2})$. This fact implies in turn that the set of $N_{(\cdot)}(0) > \frac{1}{2}(m-1)$ points in $Z_{m-1} - (Z_{m-1} \cap \mathcal{V}_{m-1})$ can be labeled as $\{p_1, p_2, ...\}$ with $N_{(\cdot)}(0)$ being a non-increasing function of the labeling index. A construction that differs only cosmetically from the construction in Part 3 creates a finite set $\mathfrak{X}_{(m-1)+} \subset \{p_1, p_2, ...\}$ and an associated set of balls $\{B(p): p \in \mathfrak{X}_{(m-1)+}\}$ with the following properties: The balls in this set have pairwise disjoint closure and their closures are disjoint from the closures of the balls in $\mathfrak{U}_m$ and $\mathfrak{V}_m$. Meanwhile, the union of the balls $\{B(p): p \in \mathfrak{X}_{(m-1)+}\}$ contains the set $\{p_1, p_2, ...\}$.

Let $\mathfrak{X}_{(m-1)0}$ denote the set of $N_{(\cdot)}(0) = \frac{1}{2}(m-1)$ points in the complement of $\mathcal{V}_{m-1}$ and in the complement of each $p \in \mathfrak{X}_{(m-1)+}$ version of $B(p)$. It follows from the first observation of Part 2 and from Lemma 7.4 that this set is finite. Let $n_{(m-1)0}$ denote the number of elements in $\mathfrak{X}_{(m-1)0}$. Fix $p \in \mathfrak{X}_{(m-1)0}$, reintroduce the number $r_p$ from the first paragraph of Part 2, and let B(p) denote a ball of radius at most $(100 n_{(m-1)0})^{-1} \varepsilon_{m-1}^2 r_p$ centered at p with only weakly continuous points from Z on its boundary. The closures of B(p) and $\mathcal{V}_{m-1}$ are disjoint as are the closures of B(p) and each $p' \in \mathfrak{X}_{(m-1)+}$ version of B(p'). Let $\mathfrak{U}_{(m-1)◊}$ to denote the collection of balls $\{B(p): p \in \mathfrak{X}_{(m-1)0} \cup \mathfrak{X}_{(m-1)+}\}$.

Step 2: Winnow the set $\mathfrak{V}_{(m-1)*}$ by throwing away any ball that lies in the union of the balls from $\mathfrak{U}_{(m-1)\diamond}$ and from $\mathfrak{V}_m$ and from $\mathfrak{U}_m$. Denote the resulting collection by $\mathfrak{U}_{(m-1)\diamond}$. The arguments used in Part 4 can be repeated with only cosmetic changes to prove that $\mathfrak{U}_{(m-1)\diamond}$ is finite.

Set $\mathfrak{U}_{m-1} = \mathfrak{U}_{(m-1)\diamond} \cup \mathfrak{U}_m$ and set $\mathfrak{V}_{m-1} = \mathfrak{V}_{(m-1)\diamond} \cup \mathfrak{V}_m$. If $\varepsilon_{m-1}$ is sufficiently small, then the pair $\mathfrak{U}_{m-1}$ and $\mathfrak{V}_{m-1}$ have what is needed to continue the induction.

**d) Proof of Proposition 7.3**

The proof has five parts. The first part of the proof sets some of the notation.

*Part 1*: Write the set of weakly continuous points in Z as the disjoint union $Z^p \cup Z^{np}$ using the definition of $Z^p$ from Part 3 of Section 7a. As noted in this same part of Section 6d, the function $N_{(\cdot)}(0)$ takes only half integer values on $Z^p$. Given a positive integer m, use $\mathfrak{Z}^p_m$ to denote that $N_{(\cdot)}(0) = \frac{1}{2}$ m part of $Z^p$. Having fixed $\varepsilon \in (0, \frac{1}{100})$ and a point $p \in \mathfrak{Z}^p_m$, reintroduce the number $r_{p,\varepsilon} \in (0, \varepsilon r_0)$ from Part 1 of Section 7a. With $r_{p,\varepsilon}$ in hand, fix $r_p \in (\frac{1}{2} r_{p,\varepsilon}, r_{p,\varepsilon})$ such that every points from Z on the boundary of the radius $r_p$ ball centered on p is weakly continuous. Let B(p) denote the radius $r_p$ ball centered at p.

Replace $\mathfrak{Z}_m$ with $\mathfrak{Z}^p_m$ in the final two paragraphs of Part 1 from Section 7c to define a $\mathfrak{Z}^p_m$ version of the collection $\mathfrak{V}_{m*}$. Denote this analog by $\mathfrak{V}^p_{m*}$. Keep in mind that each ball in this collection is some $p \in \mathfrak{Z}^p_m$ version of B(p) and that this collection gives a cover of $\mathfrak{Z}^p_m$.

Let $\mathfrak{Z}^{np}_m$ denote the set of points in $Z^{np}$ where $N_{(\cdot)}(0) = \frac{1}{2}$ m. Meanwhile, the set of points in $Z^{np}$ where $N_{(\cdot)}(0) \in (\frac{1}{2} m, \frac{1}{2}(m+1))$ is written as a disjoint union of sets labeled by the positive integers. These sets are defined as follows: If k is a positive integer, define $\mathfrak{F}^{np}_{m,k}$ to be the set of points in $Z^{np}$ where $N_{(\cdot)}(0)$ no less than $\frac{1}{2}(m+2^{-k})$ but strictly less than $\frac{1}{2}(m+2^{-k+1})$. Given a point p in $\mathfrak{Z}^{np}_m$ and $\varepsilon \in (0, \frac{1}{100})$ or given a point p in $\mathfrak{Z}^{np}_{m,k}$ and $\varepsilon \in (0, \frac{1}{100} 2^{-k})$, reintroduce the number $r_{p,\varepsilon}$ that is defined in Part 6 of Section 6d. Fix $r_p \in (\frac{1}{2} r_{p,\varepsilon}, r_{p,\varepsilon})$ such that every point from Z on the boundary of the radius $r_p$ ball centered at p is weakly continuous. Let B(p) denote this radius $r_p$ ball centered at p.

Replace $\mathfrak{Z}_m$ by $\mathfrak{Z}^{np}_m$, or by some $k \in \{1, 2,...\}$ version of $\mathfrak{Z}^{np}_{m,k}$ in the final two paragraphs of Part 1 from Section 7c to define the $\mathfrak{Z}^{np}_m$ and $\mathfrak{Z}^{np}_{m,k}$ versions of the collection $\mathfrak{V}_{m*}$. Denote these analogs by $\mathfrak{V}^{np}_{m*}$ and $\mathfrak{V}^{np}_{m,k*}$.

*Part 2*: Supposing that m is a positive integer, let $Z_m$ denote the set of points in Z with $N_{(\cdot)}(0) \geq \frac{1}{2}$ m. Use N to denote the integer such that $Z_N \neq \emptyset$ but $Z_{N+1} = \emptyset$. The sets $\mathfrak{V}$ and $\mathfrak{U}$ will be constructed in an inductive manner by constructing sequentially collections of balls $(\mathfrak{U}_N, \mathfrak{V}_N)$, then $(\mathfrak{U}_{N-1}, \mathfrak{V}_{N-1})$, and so on until $(\mathfrak{U}_1, \mathfrak{V}_1)$ which are the desired pair ($\mathfrak{U}, \mathfrak{V}$). With it understood that $\mathfrak{U}_{N+1} = \emptyset$ and $\mathfrak{V}_{N+1} = \emptyset$, suppose that $m \in \{N+1, N, \ldots, 2\}$ and that $\mathfrak{U}_m$ and $\mathfrak{V}_m$ have been constructed so that the conclusions of Proposition 10.3 hold with $\mathfrak{U} = \mathfrak{U}_m$, $\mathfrak{V} = \mathfrak{V}_m$ and with Z replaced by $Z_m$. Make the assumption in what follows that all points from Z on the boundary of the closure of the union of the balls from $\mathfrak{U}_m$ and $\mathfrak{V}_m$ are weakly continuous. The subsequent parts of the proof give the rules for constructing $\mathfrak{U}_{m-1}$ and $\mathfrak{V}_{m-1}$. As might be expected, the rules are much like those used to define the induction in Section 7c that proved Proposition 7.2.

*Part 3*: Fix $\varepsilon_{m-1} \in (0, \frac{1}{100} z^{-1})$ and define the collections $\mathfrak{V}^p_{(m-1)*}$ and $\mathfrak{V}^{np}_{(m-1)*}$ and each positive integer k version of $\mathfrak{V}^{np}_{m,k*}$ using $\varepsilon_{m-1}$ for the chosen value of the parameter $\varepsilon$. Each of these collections is now winnowed to throw out redundant balls. To start this process, remove all balls from $\mathfrak{V}^p_{(m-1)*}$ that are in the union of the balls from the collections $\mathfrak{V}^{np}_{(m-1)*}$, any positive integer k version of $\mathfrak{V}^{np}_{(m-1),k*}$ and $\mathfrak{U}_m$ and $\mathfrak{V}_m$. Introduce by way of notation $\mathfrak{V}^p_{(m-1)\ddagger}$ to denote the resulting collection of balls. Meanwhile, remove all balls from $\mathfrak{V}^{np}_{(m-1)*}$ that are in the union of the balls from $\mathfrak{V}^p_{(m-1)*}$, any positive integer k version of $\mathfrak{V}^{np}_{(m-1),k*}$ and $\mathfrak{U}_m$ and $\mathfrak{V}_m$. Denote the result of doing so by $\mathfrak{V}^{np}_{(m-1)\ddagger}$. Note that the union of the balls from $\mathfrak{V}^{np}_{(m-1)\ddagger}$ and $\mathfrak{V}^p_{(m-1)\ddagger}$ contains all of the weakly continuous, $N_{(\cdot)}(0) = \frac{1}{2} m$ points in Z that are not in a ball from either $\mathfrak{U}_m$ or $\mathfrak{V}_m$.

Supposing that k is a positive integer, remove all balls from $\mathfrak{V}^{np}_{(m-1),k*}$ that are in the union of the balls from the collections $\mathfrak{V}^p_{(m-1)*}$, $\mathfrak{V}^{np}_{(m-1)*}$, all positive $i \neq k$ versions of $\mathfrak{V}^{np}_{m-1,i*}$ and $\mathfrak{U}_m$ and $\mathfrak{V}_m$. Use $\mathfrak{V}^{np}_{(m-1),k\ddagger}$ to denote the result of doing this winnowing.

*Part 4*: There may be $N_{(\cdot)}(0) \in (\frac{1}{2}(m-1), \frac{1}{2}m)$ points in Z that are not weakly continuous. Use $\mathfrak{Q}_{(m-1)+}$ to denote the points from this set that are not in any ball from any of the collections $\mathfrak{V}^p_{(m-1)\ddagger}$, $\mathfrak{V}^{np}_{(m-1)\ddagger}$, each positive integer k version of $\mathfrak{V}^{np}_{m-1,k\ddagger}$, $\mathfrak{U}_m$ or $\mathfrak{V}_m$. Supposing that $\mathfrak{Q}_{(m-1)+} \neq \emptyset$, then there is a point in $\mathfrak{Q}_{(m-1)+}$ where $N_{(\cdot)}(0)$ must take its supremum value on $\mathfrak{Q}_{(m-1)+}$. To prove this, let $N_1$ denote here the supremum of $N_{(\cdot)}(0)$ on $\mathfrak{Q}_{(m-1)+}$ and let $\{p_1, p_2, \ldots\}$ denote a sequence in $\mathfrak{Q}_{(m-1)+}$ with $N_{(\cdot)}(0)$ not decreasing along the sequence and with $N_1$ being the limit of its values on the sequence. Let p denote a limit point of $\{p_i\}_{i \in \{1,2,\ldots\}}$. The function $N_{(p)}(0)$ is greater than $\frac{1}{2}(m-1)$ because $N_1$ is greater than $\frac{1}{2}(m-1)$ and $N_{(p)}(0)$ is no less than $N_1$. Meanwhile, $N_{(p)}(0)$ must be less than $\frac{1}{2}m$ because p would otherwise be in a ball from $\mathfrak{U}_m$ or $\mathfrak{V}_m$ and no balls from either set has points from $\mathfrak{Q}_{(m-1)+}$. For much the same reason, the point p is not in any positive integer k version of $\mathfrak{F}^{np}_{m-1,k}$. It follows as a consequence that p is in $\mathfrak{Q}_{(m-1)+}$. This implies that $N_1 = N_{(p)}(0)$ by virtue of the definition of $N_1$. This also implies that the sequence $\{p_1, p_2, \ldots\}$ contains but a finite number of distinct points.

Let $p_1$ denote a point in $\mathfrak{Q}_{(m-1)+}$ where $N_{(\cdot)}(0) = N_1$. Let $B(p_1)$ denote a ball of radius at most $\frac{1}{100} 2^{-1} \varepsilon_{m-1}^2$ centered at $p_1$ whose closure is disjoint from the closure of the balls from $\mathfrak{U}_m$ and $\mathfrak{V}_m$ and such that each point from Z on its boundary is weakly continuous.

The argument from the preceding paragraph can be repeated with only cosmetic modifications to prove that $\mathfrak{Q}_{(m-1)+}-(\mathfrak{Q}_{(m-1)+} \cap B(p_1))$ has a point where $N_{(\cdot)}(0)$ takes on its maximum. Fix such a point, to be denoted by $p_2$, and let $B(p_2)$ denote a ball of radius at most $\frac{1}{100} 2^{-2} \varepsilon_{m-1}^2$ centered at $p_2$ whose closure is disjoint from the closure of $B(p_1)$ and from the closures of the balls in $\mathfrak{U}_m$ and $\mathfrak{V}_m$. Require in addition that each point from $Z$ on the boundary of $B(p_2)$ be weakly continuous. Continue in this vein to sequentially define a set of points $\{p_1, p_2, \ldots\}$ in $\mathfrak{Q}_{(m-1)+}$ and an associated set of balls $\{B(p_i)\}_{i \in \{1,2\ldots\}}$ with the following properties:

- *If $i \in \{1, 2, \ldots\}$, then $N_{(p_i)}(0)$ is the maximum of $N_{(\cdot)}(0)$ on the complement in $\mathfrak{Q}_{(m-1)+}$ of its points in $\cup_{1 \le j < i} B(p_j)$.*
- *If $i \in \{1, 2, \ldots\}$, then the radius of $B(p_i)$ is no greater than $\frac{1}{100} 2^{-i} \varepsilon_{m-1}^2$.*
- *The closures of the balls in the collection $\{B(p_i)\}_{i \in \{1,2,\ldots\}}$ are pairwise disjoint and they are disjoint from the closures of the balls in $\mathfrak{U}_m$ and $\mathfrak{V}_m$.*
- *Each point in $Z$ on the boundary of any ball from $\{B(p_j)\}_{i \in \{1,2,\ldots\}}$ is weakly continuous.*

(7.3)

Let $\mathfrak{X}_{(m-1)+}$ denote the set $\{p_1, p_2, \ldots\}$. The set $\mathfrak{X}_{(m-1)+}$ is finite. An argument that differs little from that used in Part 3 of Section 7c proves that this is so.

There may be $N_{(\cdot)}(0) = \frac{1}{2}(m-1)$ points in $Z$ that are not weakly continuous. Let $\mathfrak{X}_{m-1,0}$ denote the set of such points that are not in any ball from $\mathfrak{V}^p_{(m-1)\ddagger}$ or $\mathfrak{V}^{np}_{(m-1)\ddagger}$, or from any positive integer $k$ version of $\mathfrak{V}^{np}_{m-1,k\ddagger}$, or from $\{B(p): p \in \mathfrak{X}_{(m-1)+}\}$ or from $\mathfrak{U}_m$ or from $\mathfrak{V}_m$. This set is finite because all limit points of sequences of $N_{(\cdot)}(0) = \frac{1}{2}(m-1)$ points in $Z$ are in one of the aforementioned balls. Let $n_{m-1}$ denote the number of elements in $\mathfrak{X}_{m-1,0}$. Supposing that $p \in \mathfrak{X}_{m-1,0}$, let $B(p)$ denote a ball of radius at most $(100 n_{m-1})^{-2} \varepsilon_{m-1}^2$ centered at $p$ whose closure is disjoint from the closure of the union of the balls from $\mathfrak{V}^p_{(m-1)\ddagger}$ and $\mathfrak{V}^{np}_{(m-1)\ddagger}$ and from each positive integer $k$ version of $\mathfrak{V}^{np}_{m-1,k\ddagger}$ and from $\{B(p): p \in \mathfrak{X}_{(m-1)+}\}$ and from $\mathfrak{U}_m$ and $\mathfrak{V}_m$. Require of $B(p)$ that each point from $Z$ on its boundary be weakly continuous.

*Part 5*: Let $\mathfrak{U}_{(m-1)\Diamond}$ denote the collection $\{B(p): p \in \mathfrak{X}_{(m-1)+} \cup \mathfrak{X}_{(m-1),0}\}$. Throw away all balls in the collections $\mathfrak{V}^p_{(m-1)\ddagger}$, $\mathfrak{V}^{np}_{(m-1)\ddagger}$ and $\{\mathfrak{V}^{np}_{m-1,k\ddagger}\}_{k \in \{1,2,\ldots\}}$ that lie entirely in any ball from $\mathfrak{U}_{(m-1)\Diamond}$. Denote the result of this winnowing procedure by $\mathfrak{V}^p_{(m-1)\Diamond}$, $\mathfrak{V}^{np}_{(m-1)\Diamond}$ and $\{\mathfrak{V}^{np}_{m-1,k\Diamond}\}_{k \in \{1,2,\ldots\}}$. Let $\mathfrak{V}_{(m-1)\Diamond}$ denote the union of the collections $\mathfrak{V}^p_{(m-1)\Diamond}$, $\mathfrak{V}^{np}_{(m-1)\Diamond}$ and $\{\mathfrak{V}^{np}_{m-1,k\Diamond}\}_{k \in \{1,2,\ldots\}}$. As explained directly, $\mathfrak{V}_{(m-1)\Diamond}$ is a finite collection of balls.

If $\mathfrak{V}_{(m-1)\Diamond}$ is not finite, then at least one of $\mathfrak{V}^p_{(m-1)\Diamond}$, $\mathfrak{V}^{np}_{(m-1)\Diamond}$ or $\cup_{k \in \{1,2,\ldots\}} \mathfrak{V}^{np}_{m-1,k\Diamond}$ is not finite. Suppose that this is the case to generate nonsense. Suppose in particular that $\mathfrak{V}^p_{(m-1)\Diamond}$ contains infinitely many balls. This collection is, by construction, countable since it is obtained by throwing away balls from $\cup_{b \in \{1,2,\ldots\}} \mathfrak{V}^p_{(m-1),b*}$ and each $b \in \{1, 2, \ldots\}$ version of $\mathfrak{V}^p_{(m-1),b*}$ is a finite collection of balls. This understood, let $\{p_1, p_2, \ldots\}$ denote a labeling of the centers of the balls from $\mathfrak{V}^p_{(m-1)\diamond}$. This labeling can be chosen so that the radii of the balls in the corresponding sequence $\{B(p_i)\}_{i \in \{1,2,...\}}$ is not increasing as a function of the index i. Such a labeling can be found because the radii of the balls in any $b \in \{1, 2,...\}$ version of $\mathfrak{V}^p_{(m-1),b*}$ are bounded by $2^{-b+1}\varepsilon_{m-1}$. Let p denote a limit point of the sequence $\{p_1, p_2, ...\}$. This point can not be in a ball from $\mathfrak{V}^p_{(m-1)\diamond}$. If it were, then all sufficiently large i versions of $B(p_i)$ would be in this ball since their radii limit to zero as i gets ever larger; and the event that $B(p_i)$ is in a ball from $\mathfrak{V}^p_{(m-1)\diamond}$ is ruled out by the winnowing process that is described in the final paragraph of Part 1 in Section 10a. The point p can not be in any ball from $\mathfrak{V}^{np}_{(m-1)\diamond}$ or in any ball from any $k \in \{1, 2, ...\}$ version of $\mathfrak{V}^{np}_{m-1,k\diamond}$ because if it were, then all sufficiently large i versions of $B(p_i)$ would be in this same ball; and this event is ruled out by the winnowing process in Part 3. The point p can not be in any ball from $\mathfrak{U}_{(m-1)\diamond}$ because if it were, then all sufficiently large i versions of $B(p_i)$ would be in a ball from this collection; and this event is ruled out by the winnowing process in the preceding paragraph. Meanwhile p can not be in a ball from $\mathfrak{U}_m$ or from $\mathfrak{V}_m$ because this would again imply that all sufficiently large i versions of $B(p_i)$ were in such a ball. This exhausts the possible locations of p and so p can't exist. Much the same sort of argument proves that the collections $\mathfrak{V}^{np}_{(m-1)\diamond}$ and $\cup_{k \in \{1,2,...\}} \mathfrak{V}^{np}_{m-1,k\diamond}$ are finite.

With the preceding understood, let $\mathfrak{U}_{(m-1)} = \mathfrak{U}_{(m-1)\diamond} \cup \mathfrak{U}_m$ and let $\mathfrak{V}_{m-1} = \mathfrak{V}_{(m-1)\diamond} \cup \mathfrak{V}_m$. These are finite collections of balls and, supposing that $\varepsilon_{m-1}$ is sufficiently, they have all the properties that are specified in Part 2 to continue the inductive definition of $\mathfrak{U}$ and $\mathfrak{V}$.

**8. Points of continuity where $N_{(·)}(0)$ is an integer**

A point $p \in Z$ is said to be a *point of continuity* for $N_{(·)}(0)$ when the following is true: If $\{p_i\}_{i \in \{1,2,...\}}$ is a sequence in Z that converges to p, then the sequence of numbers $\{N_{(p_i)}(0)\}_{i \in \{1,2...\}}$ converges to $N_{(p)}(0)$. The proposition that follows describes the bundle $\mathcal{I}$ near a point of continuity for $N_{(·)}(0)$ where $N_{(·)}(0)$ is a positive integer. The subsequent proposition states a consequence.

**Proposition 8.1**: *Let $p \in Z$ denote a point of continuity for $N_{(·)}(0)$ where $N_{(·)}(0)$ is a positive integer. The bundle $\mathcal{I}$ is isomorphic to the product $\mathbb{R}$ bundle on the complement of Z in some ball centered at p.*

It follows as a consequence that $v$ can be viewed as an $\mathbb{R}$-valued harmonic section of $\mathbb{S}^+$ near any point of continuity for $N_{(·)}(0)$ where $N_{(·)}(0)$ is an integer. The proposition that follow momentarily are seen to be consequences of this fact.

The upcoming Proposition 8.2 talks about unique tangent cones. To define this term, fix $p \in Z$ for the moment and let L denote any given oriented, orthonormal frame for $TX|_p$. Parallel transport L along the geodesic arcs from p to view L as an orthonormal frame at each point is a small radius bull centered at p. The set Z is said to have a *unique tangent cone* at p when there exists a data set $(\mathcal{I}_\diamond, Z_\diamond, \nu_\diamond)$ with the following universal property: Let $\{T(i) = (p_i, L, \lambda_i)\}_{i=\{1,2,\ldots\}}$ denote any given sequence of data sets with $\{p_i\}_{i \in \{1,2,\ldots\}}$ in this small ball and converging to p, and with $\{\lambda_i\}_{i \in \{1,2,\ldots\}} \subset (0, c_0^{-1}]$ converging to zero. If such a sequence is used in Proposition 4.1, then any limit data set $(Z_*, \mathcal{I}_*, \nu_*)$ from this proposition has $Z_* = Z_\diamond$ and it has $\mathcal{I}_*$ being isomorphic to $\mathcal{I}_\diamond$ by an isomorphism that identifies $\nu_*$ with $\nu_\diamond$.

**Proposition 8.2**: *The set Z has a unique tangent cones at each point of continuity for $N_{(\cdot)}(0)$ where $N_{(\cdot)}(0)$ is a positive integer*.

Section 8a proves Proposition 8.1 and Section 8b proves Proposition 8.2.

**a) Proof of Proposition 8.1**

Suppose that $B \subset X$ is a ball centered at p. The bundle $\mathcal{I}$ is not isomorphic to the product bundle on $B - (B \cap Z)$ if and only if there exists an embedded loop $\gamma \subset B - (B \cap Z)$ such that $\mathcal{I}$'s restriction to $\gamma$ is not isomorphic to the product $\mathbb{R}$ bundle. Suppose for the sake of argument that there exists such a loop in every ball centered at p. The five parts of the proof that follow generate nonsense from this assumption.

By way of notation, L is used to denote a chosen oriented, orthonormal frame for $TX|_p$. Parallel transport along the short geodesic arcs from p is used to view L as an oriented, orthonormal frame at each point in the radius $r_\varepsilon$ ball centered at p.

*Part 1*: This part of the proof supplies certain facts about the bundle $\mathcal{I}$ near the point p. These facts are used in the subsequent parts of the proof.

Given $r \in (0, \frac{1}{100} r_0]$, let $\mathcal{C}_r$ denote set of maps from $S^1$ to $B_r - (B_r \cap Z)$ that pull-back $\mathcal{I}$ to a bundle that is not isomorphic to the product $\mathbb{R}$ bundle on the circle. Given $\gamma \in \mathcal{C}_r$, let $d(\gamma)$ denote the $\inf_{t \in S^1} \text{dist}(\gamma(t), Z)$ and let $D(r)$ denote the supremum of the numbers that comprise the set $\{d(\gamma) : \gamma \in \mathcal{C}_r\}$. Keep in mind that $D(r) \leq r$. The plan for what follows is to use Lemma 6.2 to obtain a picture of Z near points in $Z \cap B_r$ that can resolve structure on length scales whose size is on the order of $D(r)$. Lemma 8.4 in the next part of the proof makes this precise.

*Part 2*: To set the stage for Lemma 8.4, fix $\varepsilon \in (0, \frac{1}{1000}]$ for the moment. The fact that p is a point of continuity for $N_{(\cdot)}(0)$ can be used to define a number $r_\varepsilon \in (0, \frac{1}{1000} r_0)$ so that $|N_{(\cdot)}(0) - N_{(p)}(0)| < \varepsilon$ on the part of Z in the radius $r_\varepsilon$ ball centered at p, and so that the conclusions of Lemma 6.2 can be invoked for points from Z in this ball.

Fix $r \in (0, r_\varepsilon)$ and suppose that $q \in Z \cap B_r$. Lemma 6.2 can be invoked using as input the data set $T = (q, L, \lambda = D(r))$; and let $(Z_*, \mathcal{I}_*, \nu_*)$ denote the resulting data set. Let $\mathbb{B}_{1/\varepsilon}$ denote the $|x| < \varepsilon^{-1}$ ball in $\mathbb{R}^d$. As in Lemma 6.2, let $\mathcal{T}_\varepsilon$ denote part of the radius $\varepsilon$ tubular neighborhood of $Z_*$ in $\mathbb{B}_{1/\varepsilon}$. Let $\phi_T$ denote the map from the $|x| < \lambda^{-1} r_0$ ball in $\mathbb{R}^d$ to X given by the composition of first the rescaling map $x \to \lambda x$ and then the Gaussian coordinate chart map centered at q that is defined by L.

**Lemma 8.4**: *There exists $\kappa > 1$ that can depend on p but not q and has the following significance: Suppose that $\varepsilon \in (0, \kappa^{-1}]$ and that $r \in (0, r_\varepsilon]$. If the dimension of X is three, then $Z_T \cap \mathbb{B}_{1/\varepsilon}$ has distance $\varepsilon$ or less from a line through the origin; and the bundle $\phi_T^* \mathcal{I}$ is isomorphic to the product $\mathbb{R}$ bundle on the complement in $\mathbb{B}_{1/\varepsilon}$ of the radius $\varepsilon$ tubular neighborhood of this line. If the dimension X is four, then $Z_T \cap \mathbb{B}_{1/\varepsilon}$ has distance at most $\varepsilon$ from either a plane through the origin or a line through the origin, and the bundle $\phi_T^* \mathcal{I}$ is isomorphic to the product $\mathbb{R}$ bundle on the complement in $\mathbb{B}_{1/\varepsilon}$ of the radius $\varepsilon$ tubular neighborhood of this plane or line.*

*Proof of Lemma 8.4*: The proof has three steps. The proof reintroduces $\mathcal{I}_T$ to denote the $\phi_T$ pull-back of $\mathcal{I}$.

Step 1: The dimension of X in this step can be either 3 or 4. If there are no points in Z with distance between $\varepsilon \lambda$ and $\varepsilon^{-1} \lambda$ from q, then $Z_T$ does not intersect the $|x| > \varepsilon$ part of $\mathbb{B}_{1/\varepsilon}$. Assuming this to be the case, then $\mathcal{I}_T$ is defined on the whole of $|x| > \varepsilon$ part of $\mathbb{B}_{1/\varepsilon}$. It is therefore isomorphic on the part of $\mathbb{B}_{1/\varepsilon}$ to the product bundle because this part of $\mathbb{B}_{1/\varepsilon}$ is simply connected.

Step 2: Suppose that the dimension of X is three and that there is a point from Z with distance between $\varepsilon \lambda$ and $\varepsilon^{-1} \lambda$ from q. Granted such a point, then the third bullet of Lemma 6.2 can be invoked and then the fourth bullet of Lemma 5.4. Letting $(Z_*, \mathcal{I}_*, \nu_*)$ denote a data set from Lemma 6.2, the fourth bullet of Lemma 5.4 says that $Z_*$ is a line through the origin. It follows from the fourth bullet of Lemma 5.4 that $\mathcal{I}_*$ is isomorphic to the product $\mathbb{R}$ bundle on $\mathbb{R}^3 - Z_*$ because $N_{(p)}(0)$ is an integer. This implies in turn that $\mathcal{I}_T$ is isomorphic to the product bundle on $\mathbb{B}_{1/\varepsilon} - \mathcal{T}_\varepsilon$.

Step 3: Suppose that the dimension of X is four and that there is a point from Z with distance between $\varepsilon \lambda$ and $\varepsilon^{-1} \lambda$ from q. The fourth bullet of Lemma 6.2 can be invoked and then the third and/or fourth bullets of Lemma 5.5. Supposing that the fourth bullet of Lemma 5.5 describes the data set $(Z_*, \mathcal{I}_*, \nu_*)$ from Lemma 6.2, then $Z_*$ is a plane through the origin and $\mathcal{I}_*$ is isomorphic to the product $\mathbb{R}$ bundle on $\mathbb{R}^d$–$Z_*$. This implies in turn that $\mathcal{I}_T$ is isomorphic to the product $\mathbb{R}$ bundle on $\mathbb{B}_{1/\varepsilon}$–$\mathcal{T}_\varepsilon$.

If the fourth bullet of Lemma 5.5 does not describe the data set $(Z_*, \mathcal{I}_*, \nu_*)$ from Lemma 6.2, then this data set is in any event described by the third bullet of Lemma 5.5. If $Z_*$ is a line through the origin, then the line bundle $\mathcal{I}_*$ is isomorphic to the product $\mathbb{R}$ bundle on $\mathbb{R}^4$–$Z_*$ and so $\mathcal{I}_T$ is isomorphic to the product line bundle on $\mathbb{B}_{1/\varepsilon}$–$\mathcal{T}_\varepsilon$. If $Z_*$ is not a line through the origin, then it is a finite union of half planes through the origin. Even so, it follows from the fourth bullet of Lemma 6.2 that all $\mathbb{B}_{1/\varepsilon}$ points in $Z_T$ have distance at most $\varepsilon$ from a single line through the origin in $\mathbb{R}^4$. The bundle $\mathcal{I}_T$ is then necessarily isomorphic to the product $\mathbb{R}$ bundle on the complement of the radius $\varepsilon$ tubular neighborhood of this line.

*Part 3*: This step proves that $\lim_{r \to 0} r^{-1} D(r) = 0$. To start, suppose to the contrary that there exists $\delta \in (0, \frac{1}{2}]$ such that $\limsup_{r \to 0} r^{-1} D(r) > \delta$. The five steps that follow derive nonsense from this assumption.

<u>Step 1</u>: Fix $\varepsilon \in (0, \frac{1}{1000}\delta)$ and then take r to be a positive number less than $r_\varepsilon$ but chosen so that $r^{-1} D(r) \geq \frac{1}{2}\delta$. With r chosen, let $\gamma: S^1 \to B_r$–$(B_r \cap Z)$ denote an embedded from $\mathcal{C}_r$ with $d(\gamma) > (1-\delta) D(r)$. Let T denote the data set $(p, L, \lambda = D(r))$. Let $\gamma_T$ denote the map from $S^1$ to $\mathbb{R}^d$ that is give by composing first $\gamma$ and then $\phi_T^{-1}$. The image of $\gamma_T$ is a curve in the part of $\mathbb{B}_{1/\varepsilon}$ where $1 - \delta < |x| < \delta^{-1}$.

<u>Step 2</u>: Since each point in $\gamma_T$ has distance at least $(1-\delta)$ from $Z_T$, the image of $\gamma_T$ can be isotoped in the compliment of $Z_T$ by pushing each point inward along the rays through the origin in $\mathbb{R}^4$ a distance at most $\frac{1}{100}\delta$ so that each point of the resulting curve has distance at most $(1 - \frac{99}{100}\delta)$ from $Z_T$ and radial coordinate at most $(1-\delta) r/D(r)$. Use $\gamma_T$ now to denote this new map from $S^1$ to $\mathbb{B}_{1/\varepsilon}$–$Z_T$. The composition $\phi_T(\gamma_T)$ is still in $\mathcal{C}_r$ and all of its points are uniformly far from $\partial B_r$.

<u>Step 3</u>: As noted by Lemma 8.4, the part of $Z_T$ in $\mathbb{B}_{1/\varepsilon}$ is contained in the radius $\varepsilon$ tubular neighborhood of either a line through the origin or a plane through the origin. Let $\Pi$ denote this linear subspace of $\mathbb{R}^d$. At the risk of abusing notation, use $\mathcal{T}_\varepsilon$ now to denote the radius $\varepsilon$ tubular neighborhood of $\Pi$. The bundle $\gamma_T^*\mathcal{I}_T$ would be necessarily isomorphic to the product $\mathbb{R}$ bundle on $S^1$ were all of the image of $\gamma_T$ in $\mathbb{B}_{1/\varepsilon}$–$\mathcal{T}_\varepsilon$. This follows from Lemma 8.4. Therefore $\gamma_T$ must intersect $\mathcal{T}_\varepsilon$. If the image of $\gamma_T$ intersects $\Pi$, then an isotopy of the image $\gamma_T$ that moves no point more than a distance $\varepsilon^3$ from its starting position will give a homotopy of $\gamma_T$ in the complement of $Z_T$ that starts at $\gamma_T$ and ends at a map from $S^1$ into the complement of $Z_T$ in $\mathbb{B}_{1/\epsilon}$ that has the following properties: No point in the image lies on $\Pi$, no point in the image has distance less than $(1 - \frac{98}{100}\delta)$ from $Z_T$, and no point in the image has radial coordinate greater than $(1 - \frac{99}{100}\delta)\,r/D(r)$. Use $\gamma_T$ to now to denote this new map from $S^1$. The pull-back of $\mathcal{I}_T$ by this new version of $\gamma_T$ is isomorphic to the pull-back of $\mathcal{I}_T$ by the old version because the isotopy that made the new version from the old version did not move points on $Z_T$.

Step 4: Let $\mathcal{T}_{2\epsilon}$ denote the radius $2\epsilon$ tubular neighborhood of $\Pi$ and let $\gamma_{T*}$ denote a connected component in $\mathcal{T}_{2\epsilon}$ of the image of $\gamma_T$. Fix $z \in \gamma_{T*}$ for the moment. The distance of z to the origin is greater than $\frac{1}{4}$ and since this is much greater than $\epsilon$, there is a unique closest point in $\Pi$ with the same distance from the origin as z. Use $\mathfrak{p}(z)$ to denote this point. This is a point with distance nearly $\epsilon$ from z and in any event $\mathrm{dist}(z,\mathfrak{p}(z)) \leq 2\epsilon$. Pushing z a distance at most $4\epsilon$ along the great circle in the $|x| = |z|$ sphere through z and $\mathfrak{p}(z)$ results in a point on the $|x| = |z|$ sphere that lies in $\mathcal{T}_{2\epsilon} - \mathcal{T}_{\epsilon}$.

This pushing of the point z away from $\mathfrak{p}(z)$ can be done in a smoothly varying fashion as z varies in $\gamma_{T*}$ so as to define an isotopy of $\gamma_{T*}$ that fixes its endpoints and results in a curve that lies entirely in $\mathcal{T}_{2\epsilon} - \mathcal{T}_{\epsilon}$. This isotopy does not intersect $Z_T$ because points in $Z_T$ have distance at least $(1 - \frac{98}{100})\delta$ from points in $\gamma_{T*}$ and the isotopy moves points at most a distance much less than this. By construction, the isotopy does not change the radial distance of any given point from the origin.

Step 5: Use the isotopies in Step 4 to construct yet another version of the map $\gamma_T$, this being a version whose image lies entirely in $\mathbb{B}_{1/\epsilon} - \mathcal{T}_{\epsilon}$. It follows from Lemma 6.2 that the pull-back of $\mathcal{I}_T$ by the new map is isomorphic to the product bundle. This implies that the pull-back of $\mathcal{I}_T$ by the original version of $\gamma_T$ is likewise isomorphic to the product $\mathbb{R}$ bundle because new version of $\gamma_T$ is homotopic to the old in the complement of $Z_T$ in $\mathbb{B}_{1/\epsilon}$. This is the desired nonsense because the original map $\gamma$ was a map from $\mathcal{C}_r$

*Part 4*: Fix $\delta \in (0, \frac{1}{1000}]$ for the moment and then $\epsilon \in (0, \frac{1}{1000}\delta)$. With $\epsilon$ chosen, fix $r \in (0, \frac{1}{100}r_\epsilon]$ such that $D(r) < \epsilon\delta r$. Then fix $\gamma \in \mathcal{C}_r$ such that $d(\gamma) \geq (1-\delta)D(r)$. The map $\gamma$ will be modified momentarily with the result being a new map from $\mathcal{C}_r$ whose image has distance greater than $D(r)$ from Z. The existence of this new map constitutes the nonsense that proves Proposition 7.1 because $D(r)$ is by definition greater than or equal to the minimum distance to Z on the image of every map from $\mathcal{C}_r$.

To define the modification of $\gamma$, introduce by way of notation $\mathcal{N}$ to denote the set of points in X with distance $4D(r)$ or less from $Z \cap B_r$. Let $\gamma_*$ denote a given component of the part of γ's image in the set $\mathcal{N}$. This part of the proof supplies what is needed in Part 5 to deform $\gamma_*$ so as to give a smoothly parametrized path with the properties that are listed below in (8.1). The list uses $\gamma_\diamond$ to denote this new path.

- *The paths $\gamma_\diamond$ and $\gamma_*$ have the same endpoints on $\partial\mathcal{N}$.*
- *The path $\gamma_\diamond$ is in the ball $B_r$.*
- *The distance to Z from any point in $\gamma_\diamond$ is greater than $(1+\delta)D(r)$.*
- *The bundle $\mathcal{I}$ is isomorphic to the product $\mathbb{R}$ bundle on the loop in X that is obtained by concatenating the paths $\gamma_\diamond$ and $\gamma_*$.*

(8.1)

Each component of $\gamma \cap \mathcal{N}$ will have its own version of $\gamma_\diamond$ and so replacing each component of $\gamma \cap \mathcal{N}$ with its corresponding version of $\gamma_\diamond$ defines a closed path in $B_r$ whose points have distance greater than $(1+\delta)D(r)$ from Z. It follows from the fourth bullet of (8.1) that the bundle $\mathcal{I}$ is not isomorphic to the product bundle on this new loop. Thus a smooth parametrization of this new loop gives a map from $S^1$ to X from the set $\mathcal{C}_r$.

The modifications of $\gamma_*$ that lead to $\gamma_\diamond$ are defined by writing $\gamma_*$ as an end-to-end concatenation of arcs with each constituent arc sitting entirely in the radius $16D(r)$ ball centered at some point in Z. Each successive arc is then modified inside this ball by viewing the arc using a version of the map $\phi_r$ as a subset of the radius 16 ball in $\mathbb{R}^d$. The five steps that follow supply a template that is used in Part 5 to modify the constituent arcs in the concatenation. Step 6 summarizes the important points from Steps 1-5 that are used in the subsequent part of the proof.

Step 1: This step and Step 2 describe a template that will be used to make the successive deformations of arcs from the set $\{\gamma_k\}_{k=0,\ldots,n-1}$. To set the stage for this, suppose that $\Pi$ is a linear subspace in $\mathbb{R}^d$ with codimension no less than 2. Let $\Pi^\perp$ denote the orthogonal subspace to $\Pi$ and write any given point $x \in \mathbb{R}^d$ as $x = y + y^\perp$ with y being the orthogonal projection of x to $\Pi$ and with $y^\perp$ being the orthogonal projection of x to $\Pi^\perp$. Suppose that $y^\perp \neq 0$. The *radial deformation* of x normal to $\Pi$ is the path in $\mathbb{R}^d$ parameterized by $[0, \infty)$ that is defined by the rule

$$t \to x(t) = y + (1+t)y^\perp.$$

(8.2)

A given $t \in [0, \infty)$ point on this path is said to be the *distance* t *radial push-off* of x.

Suppose that $\varepsilon \in (0, \frac{1}{1000})$ and that q is a point in $\mathbb{R}^d$ with distance at most than $\varepsilon$ from $\Pi$ and distance greater than $\frac{1}{2}$ from x. If $t \in [0, \infty)$, then the point x(t) that gives the distance t radial push-off of x obeys $|x(t) - q| \geq |x - q| + (t - \varepsilon)$.

Step 2:  Let σ denote a path in the |x| < 12 ball.  If σ intersects Π, then any sufficiently generic isotopy of σ defines a 1-parameter family of paths that starts with σ and ends in a path that is disjoint from Π.  The desired isotopy can assumed to have the following two properties:  No point along the isotopy is moved where the distance to Π is greater than ε; and no point along the isotopy has distance more than $ε^2$ from its starting position in σ.  Let $σ_1$ denote the end member of such an isotopy.

The radial deformations normal to Π that are depicted in (8.2) can be used to isotope $σ_1$ along a 1- parameter family of paths so that each point in the end member path of the isotopy has distance at least $\frac{3}{2}$ from Π, and such that no point is moved by the isotopy where the distance from Π is greater than $\frac{7}{4}$.  Let $σ_2$ denote the end member path of this second isotopy.

Suppose that $\mathcal{Z}$ is a subset of the radius ε tubular neighborhood of Π with the distance from any given point in σ to $\mathcal{Z}$ being greater than $\frac{1}{2}$.  It then follows from what is said at the end of Step 2 that there are no points from $\mathcal{Z}$ along either of the isotopies that are used to go from σ to $σ_2$.

Step 3:  The template described in Step 2 can be used to deal with the requirements in the third and fourth bullets of (8.1), but it cannot deal with the requirment in the second bullet of (8.1).  This step adds to what is done in Step 2 with the result being a template that can be used to deal with the second, third and fourth bullets of (8.1).

To set stage, let S denote a properly embedded codimension 1 hypersurface in the |x| ≤16 ball that is diffeomorphic to a disk.  If $ε < c_0^{-1}$, then S is looks very much like the intersection of the |x| ≤ 16 ball with a codimension 1 plane.  In this event, the complement of S in the |x| < 16 ball has two connected components, each being homeomorphic to an open ball.  Denote one of these components by $\mathcal{B}$.  Assume in what follows that there is a point in $\mathcal{B}$ with distance at most 4 from the origin since this will be a case in the upcoming application.  By way of a look ahead, the subspace $\mathcal{B}$ in the upcoming applications is the part of the |x| < 16 ball in $\mathbb{R}^d$ that is mapped to the interior of the radius r ball centered at the point p by versions of Section 4's map $\phi_T$.

With $\mathcal{B}$ given, assume that the endpoints of Step 2's path σ is in $\mathcal{B}$ and that the isotopy in Step 2 that leads from σ to $σ_1$ keeps these endpoints in $\mathcal{B}$.  This can be assumed because the movement from σ to $σ_1$ is need only make $σ_1$ disjoint from Π.  The path $σ_2$ from Step 2 may or may not lie entirely in $\mathcal{B}$.  If all of $σ_2$ is in $\mathcal{B}$, then no more needs be done except to change the name of $σ_2$ to $σ_3$.

Let $z_-$ denote the starting point to $σ_2$ with respect to a given orientation.  Assume that $z_- \in \mathcal{B}$ because the applications to come will be using versions of $σ_2$ with this property, but assume that there is a point in $σ_2$ that is not in $\mathcal{B}$.  As explained in the next step, if $ε < c_0^{-1}$, then there is a point in $\mathcal{B}$ with the following properties:  The distance between the point and $\Pi$ is no less than $\frac{15}{8}$, and the distance from the point to the ending point of $\sigma$ is no greater than 4. If the ending point of $\sigma_2$ does not meet these requirements, then fix a point in $\mathcal{B}$ that does. Denote this point by $z_+$ in either case. The next step also explains why there is a path in the $|x| < 16$ part of $\mathcal{B}$ from $z_-$ to $z_+$ whose points have distance no less than $\frac{3}{2}$ from $\Pi$. Denote this path by $\sigma_3$.

Step 4: To see about the existence of $z_+$, let $\sigma_+$ denote the ending point of $\sigma$. The worst case scenario has the ending point of $\sigma_2$ being a distance $\frac{7}{4} + \varepsilon^2$ from the ending point of $\sigma$. To explain, the initial isotopy from $\sigma$ to $\sigma_1$ moved $\sigma_+$ to a point $\sigma_{1+}$ with distance at most $\varepsilon^2$ from $\sigma_+$. The isotopy from $\sigma_1$ to $\sigma_2$ moves $\sigma_{1+}$ to a point a distance at most $\frac{7}{4}$ from $\sigma_{1+}$ with the motion being normal to $\Pi$. The latter point is the ending point of $\sigma_2$. Denote the ending point of $\sigma_2$ by $\sigma_{2+}$. Let $\Pi^\perp$ again denote the orthogonal subspace in $\mathbb{R}^d$ to $\Pi$. Keep in mind that $\sigma_{2+} - \sigma_{1+} \in \Pi^\perp$.

Suppose for the moment that $S$ is a flat, codimension 1 plane. Since $\sigma_{1+}$ is in $\mathcal{B}$ and $\sigma_{2+}$ is not in $\mathcal{B}$, the hypersurface $S$ must bisect the line between $\sigma_{1+}$ and $\sigma_{2+}$. To see the implications, write $\sigma_{+2}$ as $y + y_2^\perp$ with $y \in \Pi$ and $y_2^\perp \in \Pi^\perp$. Since $\sigma_{2+} - \sigma_{1+} \in \Pi^\perp$, what was just said about $S$ implies that $S$ must bisect the line segment between $y + y_2^\perp$ and $y$ because $\sigma_{1+}$ lies on this line segment. This implies that $\mathcal{B}$ intersects any radius $r > 0$ circle or sphere (as the case may be) in the subspace $y + \Pi^\perp$ centered at $y$. In particular, there are points in $\mathcal{B}$ on the $r = \frac{15}{8}$ circle or sphere in this subspace. Any such point will have distance at most 4 from $\sigma_+$ if $\varepsilon < c_0^{-1}$. Keeping in mind that the norms of the extrinsic and intrinsic curvatures of $S$ are less than $c_0\varepsilon$, these considerations for the case when $S$ is a flat plane lead to the desired conclusions when $\varepsilon < c_0^{-1}$ because $S$ will be very near a flat plane codimension 1 plane when $\varepsilon < c_0^{-1}$.

The existence of a path in the $|x| < 16$ part of $\mathcal{B}$ connecting $z_-$ and $z_+$ with distance greater than $\frac{3}{2}$ from $\Pi$ follows when $\varepsilon < c_0^{-1}$ because there such paths exist when $S$ is a codimension 1 plane.

Step 5: Let $\sigma$ again be the path in Step 2 and again use $\sigma_+$ to denote its ending point. Denote by $\sigma_{1+}$, $\sigma_{2+}$ and $z_+$ the respective ending points of the paths $\sigma_1$, $\sigma_2$ and $\sigma_3$. The isotopy from $\sigma$ to $\sigma_1$ defines a path from $\sigma_+$ to $\sigma_{1+}$ of length at most $\varepsilon^2$, and the isotopy from $\sigma_1$ to $\sigma_2$ defines a path of the sort that is depicted in (8.2) from $\sigma_{1+}$ to $\sigma_{2+}$. It follows from the construction of $z_+$ that there is a path from $\sigma_{2+}$ to $z_+$ whose points have distance no less than $\frac{3}{2}$ from $\Pi$. Use $\theta$ to denote the end-to-end concatenation of the path from $\sigma_+$ to $\sigma_{1+}$, then the path from $\sigma_{1+}$ to $\sigma_{2+}$, and then the path from $\sigma_{2+}$ to $z_+$. Denote by $\theta_+$ the part of $\theta$ that goes from $\sigma_+$ to $\sigma_{2+}$ and define $\theta_{++}$ to be the part that goes from $\sigma_{2+}$ to $z_+$.

Suppose that $Z$ is a subset on the $|x| < 16$ ball that is contained in the radius $\varepsilon$ tubular neighborhood of $\Pi$ and that it has distance greater than $\frac{1}{2}$ from $\sigma_+$. If $Z$ has these properties, then all points in $\theta$ have distance greater than $\frac{1}{2} - c_0\varepsilon$ from $Z$.

Step 6: This step summarizes some of the salient features of the constructions in Steps 2-5. The construction starts with a path $\sigma$ whose end points are in $\mathcal{B}$ and Step 2 uses $\sigma$ to define a path $\sigma_2$. What follows are some important observations about $\sigma_2$. The final observation refers to a given set $Z$ that lies in the radius $\varepsilon$ tubular neighborhood of $\Pi$ and whose points have distance $\frac{1}{2}$ or more from all points in $\sigma$.

- *The points on $\sigma_2$ have distance no less than $\frac{3}{2}$ from $\Pi$.*
- *The path $\sigma_2$ is isotopic to $\sigma$; and the isotopy does not move points with distance greater than $\frac{7}{4}$ from $\Pi$. In particular,*
    a) *If the starting point of $\sigma$ has distance greater than $\frac{7}{4}$ from $\Pi$, then the starting point of $\sigma_2$ is the same as that of $\sigma$.*
    b) *If the ending point of $\sigma$ has distance greater than $\frac{7}{4}$ from $\Pi$, then the ending point of $\sigma_2$ is the same as that of $\sigma$.*
- *The isotopy that takes $\sigma$ to $\sigma_2$ does not intersect $Z$.*

(8.3)

Suppose in what follows that $\sigma$ and $\sigma_2$ have the same starting point. The path $\theta_+$ that is described in Step 5 is needed to insure that $\sigma$ and $\sigma_2$ have the same ending point. But, not only does the concatenation of $\sigma_2$ with $\theta_+$ have the same starting and ending points as $\sigma$, this concatenation is isotopic rel its endpoints to $\sigma$ in the complement of the set $Z$. What follows is a direct consequence of the latter fact.

*A real line bundle defined on the complement of $Z$ is isomorphic to the product bundle on the loop that is defined by the end-to-end concatenation of the paths $\sigma$, $\sigma_2$ and $\theta_+$.*

(8.4)

The path $\sigma_3$ from Step 3 is needed when the set $\mathcal{B}$ is not the whole of the $|x| < 16$ ball. If it is the whole of this ball, then $\sigma_3$ is $\sigma_2$. In any event, assume in what follows that the starting point of $\sigma$ is in $\mathcal{B}$ and that it has distance greater than $\frac{7}{4}$ from $\Pi$. The path $\sigma_3$ has the following properties:

- *The points in $\sigma_3$ have distance no less than $\frac{3}{2}$ from $\Pi$.*
- *The whole of $\sigma_3$ is in the set $\mathcal{B}$.*
- *The starting point of $\sigma_3$ is the same as that of $\sigma$.*
- *The ending point of $\sigma_3$ has distance $\frac{15}{8}$ from $\Pi$.*

(8.5)

The path θ from Step 5 is needed to connect the ending point of $\sigma_3$ to the ending point of σ. Since θ is in the complement of $Z$ and $\sigma_3$ is also, then so is their concatenation. If the bundle $\mathcal{I}$ is isomorphic to the product $\mathbb{R}$ bundle on the complement of the radius ε tubular neighborhood of Π, then one has the following:

*A real line bundle defined on the complement of $Z$ is isomorphic to the product bundle on the loop that is defined by the end-to-end concatenation of the paths σ, $\sigma_3$ and θ.*

(8.6)

Indeed, this follows from (8.4) and from the fact that the path $\sigma_2$, the path $\theta_{++}$ from Step 5 and the path $\sigma_3$ end-to-end concatenate to define a closed loop in the complement of the radius 1 tubular neighborhood of Π.

*Part 5*: The first step below explains how $\gamma_*$ is written as a suitable end-to-end concatentation of arcs. The remaining steps describe the modifications of these arcs so as to produce a path $\gamma_\diamond$ that obeys the conditions in (8.1).

<u>Step 1</u>: Parametrize $\gamma_*$ via the interval [0, 1]. The endpoints of $\gamma_*$ are therefore γ(0) and γ(1) with it understood that these points will be equal if $\gamma_* = \gamma$. There is, by assumption, a point in $Z$ with distance 4D(r) from γ(0). Choose such a point and denote it by $p_0$. Use $t_1$ to denote the largest parameter in (0,1] such that γ([0, $t_1$)) is in the radius 8D(r) ball centered at $p_0$. Supposing that $t_1 \neq 1$, fix a point in $Z$ with distance at most 4D(r) from γ($t_1$) and denoted it by $p_1$. Let $t_2$ denote the largest parameter in (0, 1] such that y([$t_1$, $t_2$)) is in the radius 8D(r) ball centered at $p_1$. If $t_2 \neq 1$, fix a point in $Z$ to be denoted by $p_2$ with distance at most 4D(r) from γ($t_2$). Continuing in this vein, define a positive integer n, an increasing set $\{t_0, t_1, \cdots, t_n\} \subset [0,1]$ with $t_0 = 0$ and $t_n = 1$, and a set $\{p_0, p_1, \ldots, p_{n-1}\} \subset Z$ such that the following condition is met:

*Supposing that k ∈ {0, 2, ..., n-1}, then the point $p_k$ has distance at most 4D(r) from γ($t_k$); the point $p_k$ has distance at most 8D(r) from γ([$t_k$, $t_{k+1}$)); and if k < n-1, then $p_k$ has distance exactly 8D(r) from γ($t_{k+1}$).*

(8.7)

The notation in what follows uses $\gamma_k$ to denote a given k ∈ {0, 1, ..., n-1} version of the arc γ([$t_k$, $t_{k+1}$]).

<u>Step 2</u>: Keeping in mind that each of the points in the set $\{p_0, p_1, \ldots, p_{n-1}\}$ is in the radius r ball centered at p, fix k ∈ {0, 1, ..., n-1} and invoke Lemma 6.2 with the input data set $\mathcal{T}$ taken to be $\mathcal{T}(k) = (p_k, L, D(r))$. There is either a line or a plane through the origin in $\mathbb{R}^d$ (and the latter only if d = 4) with the property that its radius ε tubular neighborhood contains all $|x| < \varepsilon^{-1}$ points in the $\tau = \tau(k)$ version of $Z_\tau$. This follows from Lemma 8.4 because there is a point from each $k \in \{0, 1, 2\ldots n-1\}$ version of $Z_{\tau(k)}$ in the spherical annulus where $1 - \delta \leq |x| \leq 13$. That this is so follows from the third condition in (8.7) unless the integer n is 1; and if n = 1, this follows from the fact that $d(\gamma)$ is between $(1-\delta)D(r)$ and $D(r)$. Use $\Pi_k$ in what follows to denote the $|x| < 16$ part of this line or plane.

Suppose the moment that $k \in \{0, \ldots, n-2\}$. The composition of first $\phi_{\tau(k)}$ and then the inverse of $\phi_{\tau(k+1)}$ maps $\Pi_k$ to a submanifold in $\mathbb{R}^d$ that lies in the radius $c_0\varepsilon$ tubular neighborhood of $\Pi_{k+1}$. This follows from three facts: The first is that $\Pi_k$ is a linear space, the second is that the composition of the maps in question differs by at most $\varepsilon$ from an isometry in the $C^2$ topology, and the third is that $\phi_{\tau(k)}(Z_{\tau(k)})$ contains the point $p_k$ and at least one point in Z with distance between $(1-\delta)D(r)$ and $13 D(r)$ from $p_k$. This last observation about the appearance of $((\phi_{\tau(k+1)})^{-1} \circ \phi_{\tau(k)})(\Pi_k)$ is used to infer the following when $\varepsilon$ is less than a constant that depends only on p:

*If a point in the $|x| < 16$ ball has distance greater than $\frac{15}{8}$ from $\Pi_k$, then its image via the composition $((\phi_{\tau(k+1)})^{-1} \circ \phi_{\tau(k)})$ has distance greater than $\frac{7}{4}$ from $\Pi_{k+1}$.*

(8.8)

Assume in what follows that $\varepsilon$ is small so as to be able to invoke (8.8).

Step 3: This step explains how the submanifolds $S$ and $\mathcal{B}$ from Step 3 of Part 4 arise. To this end, fix $k \in \{0, 1, \ldots, n-1\}$. If the ball of radius $16 D(r)$ centered at $p_k$ is entirely inside the radius r ball centered at p, then $\mathcal{B}$ is the whole of the $|x| < 16$ ball and $S$ is not present. If the ball of radius $16 D(r)$ centered at $p_k$ is not entirely inside the radius r ball centered at p, then $\mathcal{B}$ is the part of the $|x| < 16$ ball that is mapped $\phi_{\tau(k)}$ into the radius r ball centered at p. Since $\gamma([t_k, t_{k+1}])$ is in this radius r ball centered at p and since the distance of $\gamma(t_k)$ to $p_k$ is at most $4 D(r)$, there is necessarily a point in $\mathcal{B}$ with distance less than 4 from the origin.

Step 4: Let $\gamma_0$ denote the inverse image via the map $\phi_{\tau(0)}$ of $\gamma([0, t_1])$. The starting point of $\gamma_0$ is on the boundary of $\mathcal{N}$. By construction, it has distance 4 from the origin and Euclidean distance greater than or equal to $4 - c_0\varepsilon$ from the set $Z_{\tau(0)}$. This doesn't preclude the possibility that the starting point of $\gamma_0$ is in the radius $\frac{3}{2}$ tubular neighborhood of $\Pi_0$. Supposing that the starting point of $\gamma_0$ is this close to $\Pi_0$, and supposing that $\varepsilon < c_0^{-1}$, then a path of length less than 2 from the starting point of $\gamma_0$ in a direction perpendicular to $\Pi_0$ connects $\gamma_0$ to a point in $\mathcal{B}$ with distance $\frac{15}{8}$ from the subspace $\Pi_0$. Denote this path by $\tau_s$ and denote its ending point by $z_s$. This whole path has distance at most $\frac{3}{2}$ from $Z_{\tau(0)}$ when $\varepsilon$ is small because the starting point of $\gamma_0$ has distance $4 - c_0\varepsilon$ from $Z_{\tau(0)}$.

Step 5: Let $\gamma_{n-1}$ denote the inverse image via the map $\phi_{T(n-1)}$ of $\gamma([t_{n-1},1])$. The ending point $\gamma_{n-1}$ is on the boundary of $\mathcal{N}$ and so it has Euclidean distance greater than or equal to $4-c_0\varepsilon$ from the set $Z_{T(n-1)}$. This doesn't preclude the possibility that the ending point of $\gamma_{n-1}$ is in the radius $\frac{3}{2}$ tubular neighborhood of $\Pi_{n-1}$. A construction that is identical but for notation to the one just used in Step 4 finds a path of length at most 2 from the ending point of $\gamma_{n-1}$ to a point in $\mathcal{B}$ with distance $\frac{15}{8}$ from $\Pi_{n-1}$ whose points have distance at most $\frac{3}{2}$ from $Z_{T(n-1)}$. Denote this path by $\tau_f$ and its ending point by $z_f$.

Step 6: An initial modification of $\gamma_*$ starts at $\gamma(0)$, traverses the path $\phi_{T(0)}(\tau_0)$ to its ending point, then traverses in reverse back to $\gamma(0)$. Having returned to $\gamma(0)$, the path then traverses $\gamma_*$ to $\gamma(1)$, this being the endpoint of $\gamma_*$. Once at $\gamma(1)$, the path traverses $\phi_{T(n-1)}(\tau_f)$ to its ending point, and finally traverses $\phi_{T(n-1)}(\tau_f)$ in reverse back to $\gamma(1)$. Let $\gamma_\ddagger$ denote the path just described. The concatenation of $\gamma_\ddagger$ with $\gamma_*$ defines a loop in X on which $\mathcal{I}$ is isomorphic to the product line bundle.

Let $\gamma_{\ddagger 0}$ denote the part of the path $\gamma_\ddagger$ that starts at the ending point $\phi_{T(0)}(z_s)$ of $\phi_{T(0)}(\tau_s)$, traverses the path $\phi_{T(0)}(\tau_s)$ in reverse back to $\gamma(0)$, then traverses $\gamma_*$ to its ending point at $\gamma(1)$, and finally traverses $\phi_{T(n-1)}(\tau_f)$ to its ending point at $\phi_{T(n-1)}(z_f)$.

Step 7: The part of $\gamma_{\ddagger 0}$ that starts at $\phi_{T(0)}(z_s)$, runs the path $\phi_{T(0)}(\tau_s)$ in reverse and then traverses $\gamma([0,t_1])$ is in the image via of the map $\phi_{T(0)}$ of the $|x|<12$ ball in $\mathbb{R}^d$. Denote this part of $\gamma_{\ddagger 0}$ by $\gamma_{\ddagger 00}$ if $n>1$. If $n=1$, then $\gamma_{\ddagger 0}$ and $\gamma_{\ddagger 00}$ are the same in what follows. Let $\sigma$ denote the inverse image of $\gamma_{\ddagger 0}$ via $\phi_{T(0)}$. Use this version of $\sigma$ to construct the paths $\sigma_1$, $\sigma_2$ and $\sigma_3$ as done in Steps 2 and 3 of Part 4. This version of $\sigma_3$ is denoted later by $\sigma_{30}$. Note that it begins at the point $z_s$ because $z_s$ has distance greater than $\frac{7}{4}$ from $\Pi_0$. Introduce the corresponding version of the path $\theta$ from Step 5 of Part 4 and denoted it by $\theta_0$. As noted in Step 6 of Part 4, the concatentation of first $\sigma_{30}$ and then $\theta_0$ has the same starting and ending points as $\sigma$. It follows that the $\phi_{T(0)}$ image of this concatenation has the same starting and ending points as $\gamma_{\ddagger 00}$.

Let $\gamma_{\ddagger 1}$ denote the path that is obtained by replacing the $\gamma_{\ddagger 00}$ part of $\gamma_{\ddagger 0}$ by the $\phi_{T(0)}$ image of the concatenation of $\sigma_{30}$ and then $\theta_0$. The path $\gamma_{\ddagger 1}$ has the following properties:

- *The path $\gamma_{\ddagger 1}$ is $\phi_{T(0)}(\sigma_{30})$ if $n=1$, and if $n>1$, it is the concatenation of $\phi_{T(0)}(\sigma_{30})$ and $\phi_{T(0)}(\theta_0)$ and $\gamma([t_1,1])$ and $\phi_{T(n-1)}(\tau_f)$.*
- *The $\phi_{T(0)}(\sigma_{30})$ part of $\gamma_{\ddagger 1}$ is in the radius r ball centered at p. If $n>1$, this is also the case for the $\gamma([t_1,1])$ and $\phi_{T(n-1)}(\tau_f)$ parts of $\gamma_{\ddagger 1}$.*
- *The $\phi_{T(0)}(\sigma_{30})$ part of $\gamma_{\ddagger 1}$ has distance greater than $\frac{5}{4}$ from Z. If $n>1$, this is also the case for the $\phi_{T(n-1)}(\tau_f)$ part of $\gamma_{\ddagger 1}$.*

- *The restriction of $\mathcal{I}$ to the loop defined by the concatenation of $\gamma_{\ddagger 1}$ and $\gamma_{\ddagger 0}$ is isomorphic to the product $\mathbb{R}$ bundle.*

(8.9)

By way of an explanation, the second bullet follows from the second bullet of (8.5), the third bullet follows from the first bullet of (8.5) and (8.8), and the fourth bullet follows from (8.6).

In the case n = 1, set $\gamma_\diamond$ to be the concatenation of the paths $\phi_{T(0)}(\tau_s)$ and $\gamma_{\ddagger 1}$ and $\phi_{T(0)}(\tau_f)$. It follows from (8.9) that this version of $\gamma_\diamond$ obeys all of the conditions in (8.1). This understood, assume in what follows that n > 1.

Step 8: The part of $\gamma_{\ddagger 1}$ that starts at the ending point of $\phi_{T(0)}(\sigma_{30})$ and traverses first $\phi_{T(0)}(\theta_0)$ in reverse and then traverses $\gamma([t_1,t_2])$ lies entirely in the image via the map $\phi_{T(1)}$ of the $|x| < 12$ ball in $\mathbb{R}^d$. Denote this part of $\gamma_{\ddagger 1}$ by $\gamma_{\ddagger 11}$ if n > 2. If n = 2, define $\gamma_{\ddagger 11}$ to be the concatenation of the paths $\phi_{T(0)}(\theta_0)$, $\gamma([t_1,t_2])$ and $\phi_{T(1)}(\tau_f)$.

The $\phi_{T(1)}$-inverse image of this path $\gamma_{\ddagger 11}$ can be used for the path $\sigma$ in Steps 2-5 of Part 4. Note in this regard that its starting point is in $\mathcal{B}$ and has distance greater than $\frac{7}{4}$ from $\Pi_1$. The first of these facts follows from the $\sigma_{30}$ version of the second bullet of (8.5) and the second follows from the $\sigma_{30}$ version of the fourth bullet of (8.5) and from (8.8). Construct from $\sigma$ the corresponding path $\sigma_3$ and denote it by $\sigma_{31}$. This path start at the starting point of $\sigma$. Introduce the $\sigma_{31}$ version of the path $\theta$ from Step 5 of Part 4 and denote this version by $\theta_1$. The concatenation of first $\sigma_{31}$ and then $\theta_1$ has the same starting and ending points as $\sigma$. It follows that the $\phi_{T(1)}$ image of this concatenation has the same starting and ending points as $\gamma_{\ddagger 11}$.

Let $\gamma_{\ddagger 2}$ denote the path that is obtained by replacing the $\gamma_{\ddagger 11}$ part of $\gamma_{\ddagger 1}$ by the $\phi_{T(1)}$ image of the concatenation of $\sigma_{31}$ and then $\theta_1$ run in reverse. The path $\gamma_{\ddagger 2}$ has the following important features:

- *The path $\gamma_{\ddagger 2}$ can be written as the concatenation of paths $\phi_{T(0)}(\sigma_{30})$ and $\phi_{T(1)}(\sigma_{31})$, and if n > 2, also $\phi_{T(1)}(\theta_1)$ and $\gamma([t_2,1])$ and $\phi_{T(n-1)}(\tau_f)$.*
- *The $\phi_{T(0)}(\sigma_{30})$ and $\phi_{T(1)}(\sigma_{31})$ parts of $\gamma_{\ddagger 2}$ are in the radius r ball centered at p. If n>1, this is also the case for the $\gamma([t_2,1])$ and $\phi_{T(n-1)}(\tau_f)$ parts of $\gamma_{\ddagger 2}$.*
- *The $\phi_{T(0)}(\sigma_{30})$ and $\phi_{T(1)}(\sigma_{31})$ parts of $\gamma_{\ddagger 2}$ have distance greater than $\frac{5}{4}$ from Z. If n>2, this is also the case for the $\phi_{T(n-1)}(\tau_f)$ part of $\gamma_{\ddagger 2}$.*
- *The restriction of $\mathcal{I}$ to the loop defined by the concatenation of $\gamma_{\ddagger 2}$ and $\gamma_{\ddagger 1}$ is isomorphic to the product $\mathbb{R}$ bundle.*

(8.10)

The second and third bullets of (8.10) follow from (8.5) and (8.8) in the case of the third bullet. The fourth bullet follows from (8.6).

In the case n = 2, set $\gamma_\diamond$ to be the concatenation of the paths $\phi_{T(0)}(\tau_s)$ and $\gamma_{\ddagger 1}$ and $\phi_{T(1)}(\tau_f)$. It follows from (8.10) that this path obeys all of the conditions in (8.1). Suppose in what follows that n > 2.

Step 9: The part of $\gamma_{\ddagger 2}$ that starts at the ending point of $\psi_{T(1)}(\sigma_{31})$ and traverses first $\phi_{T(1)}(\theta_1)$ in reverse and then traverses $\gamma([t_2, t_3])$ lies entirely in the image via the map $\phi_{T(2)}$ of the $|x| < 12$ ball in $\mathbb{R}^d$. Denote this part of $\gamma_{\ddagger 2}$ by $\gamma_{\ddagger 22}$ if n > 3; and if n = 3 take $\gamma_{\ddagger 22}$ to be the concatenation of $\phi_{T(1)}(\theta_1)$, $\gamma([t_2, t_3])$ and $\phi_{T(2)}(\tau_f)$. What is done in Step 8 can be repeated using $\gamma_{\ddagger 22}$ to construct $\gamma_{\ddagger 3}$, and then continuing in this vein until ending with the construction of $\gamma_{\ddagger n}$. This path can be written as the concatenation of paths of a set of n paths $\{\phi_{T(k)}(\sigma_{3k})\}_{k=0,\ldots,n-1}$ with each $k \in \{0, 1, \ldots, n-1\}$ version of $\sigma_{3k}$ obeying the bullets in (8.5) for a suitable index k dependent version of $\sigma$ in the second bullet. The intermediate paths $\{\gamma_{\ddagger m}\}_{m=1,2,\ldots n-1}$ that lead to $\gamma_{\ddagger n}$ are such that any given $m \in \{1, 2, \ldots, n-1\}$ version can be written as the concatentation of the paths from the set $\{\phi_{T(k)}(\sigma_{3k})\}_{k=0,1,\ldots,m-1}$, a path that can be written as $\phi_{T(m-1)}(\theta_{m-1})$ with $\theta_{m-1}$ being the $\sigma_3 = \sigma_{3(m-1)}$ version of Step 5's path $\theta$, and then the paths $\gamma([t_m, 1])$ and $\phi_{T(n-1)}(\tau_f)$.

Given what is said by the various $\{\sigma_{3k}\}_{k=0,1,\ldots,n-1}$ versions of (8.5) and the corresponding versions of (8.6), and given the preceding description of $\{\gamma_{\ddagger m}\}_{m=1,2,\ldots,n}$, it follows that the concatenation of first $\phi_{T(0)}(\tau_s)$, then $\gamma_{\ddagger n}$ and finally the path $\phi_{T(n)}(\tau_f)$ is the desired path $\gamma_\diamond$ for use in (8.7).

**b) Proof of Proposition 8.2**

Suppose in what follows that $p \in Z$ is a point of continuity for $N_{(\cdot)}(0)$ where $N_{(\cdot)}(0)$ is an integer. According to Proposition 8.1, the bundle $\mathcal{I}$ is isomorphic to the product bundle on a small ball centered at p. Fix such a ball and denote it as B and choose an isomorphism from $\mathcal{I}$ on $B-(B \cap Z)$ to the product $\mathbb{R}$ bundle. The isomorphism identifies $v$ on $B-(B \cap Z)$ as an $\mathbb{R}$ valued harmonic section of $\mathbb{S}^+$. Since Z is the zero locus of $|v|$, it follows that $v$ extends across $B \cap Z$ as a continuous section of $\mathbb{S}^+$ on B whose zero locus is $B \cap Z$. This extention also denoted by $v$. It is annihilated by $\mathcal{D}$ because $\mathcal{D}$ is a first order operator and because the extended $v$ is continuous. It follows as a consequence that $v$ is smooth on B.

A version of Aronszajn's unique continuation [Az] theorem implies that $v$ has a non trivial Taylor's approximation near p. To say what this means, fix an orthonormal frame L for $TX|_p$ and use this frame to define a Gaussian coordinate chart at p. Then there is a non-trivial, $\mathbb{V}^+$-valued, homogeneous polynomial $\wp$ of degree $N_0$ such that the pull-back of $v$ by this chart has the form

$$\nu = \wp + \cdots \tag{8.11}$$

where the unwritten terms are bounded by a constant multiple of $r^{1+N_0}$.

Let $Z_\diamond$ denote the zero locus of $\wp$, let $\mathcal{I}_\diamond$ denote the product $\mathbb{R}$ bundle on $\mathbb{R}^d - Z_*$ and let $\nu_*$ denote $c\,\wp$ with $c$ being the square root of the inverse of the integral of $|\wp|^2$ on the $|x| = 1$ sphere in $\mathbb{R}^d$. It follows from (8.11) that any limit $(Z_*, \mathcal{I}_*, \nu_*)$ from Lemma 4.7 using as input $\{\mathrm{T}(i) = (p_i, L, \lambda_i)\}_{i \in \{1,2,\ldots\}}$ with $\{p_i\}_{i \in \{1,2,\ldots\}} \subset B$ converging to p and with $\{\lambda_i\}_{i \in \{1,2,\ldots\}}$ converging to 0 must be such that $Z_* = Z_\diamond$ and such that $\mathcal{I}_*$ is isomorphic to the product bundle via an isomorphism that identifies $\nu_*$ with $\wp$.

## 9. Strongly continuous points where $N_{(\cdot)}(0) = \frac{1}{2} \mod(\mathbb{Z})$

The proposition that follows makes the formal assertion to the effect that Z has a unique tangent cone at each strongly continuous points in Z where $N_{(\cdot)}(0)$ is half of a positive odd integer.

**Proposition 9.1**: *Let* $p \in Z$ *denote a strongly continuous point where* $N_{(\cdot)}(0) = \frac{1}{2} \mod(\mathbb{Z})$. *Then Z has a unique tangent cone at* p. *The corresponding data set* $(Z_\diamond, \mathcal{I}_\diamond, \nu_\diamond)$ *is described by either the fourth bullet of Lemma 5.4 or the third and fourth bullets of Lemma 5.5 as the case may be.*

The proof of this proposition is split between Section 9a and Section 9b. Supposing that p is a strongly continuous point in Z with $N_{(\cdot)}(0)$ being half an odd integer, Section 9a explains why any data set $(Z_*, \mathcal{I}_*, \nu_*)$ that arises via p's version of Proposition 4.1 is described by the fourth bullet of Lemma 5.4 when the dimension of X is 3, and by the third and fourth bullets of Lemma 5.5 when the dimension of X is 4. Section 9b proves that Z has a unique tangent cone at p.

The arguments in Section 9b invoke the a priori bound from the next proposition for the rate of convergence of the function $r \to N_{(p)}(r)$ to its $r = 0$ value $N_{(p)}(0)$. To set the stage for this proposition, suppose for the moment that $p \in Z$ is a given point. Define the function $\underline{N}_{(p)}$ on $[0, \frac{1}{1000} r_0]$ by the rule

$$r \to \underline{N}_{(p)}(r) = (\sup\nolimits_{p' \in B_{8r}} \sup\nolimits_{0 \leq s \leq 16r} N_{(p')}(s)) - N_{(p)}(0) \,.$$
(9.1)

This function is non-decreasing and it is never less than $-c_0 r^2$.

The upcoming Proposition 9.2 refers to set a $K \subset Z$ of the following sort: Suppose that $B \subset X$ is a ball that is well inside a Gaussian coordinate chart and is such that each point in $B \cap Z$ is strongly continuous. Assume in addition that $N_{(\cdot)}(0)$ on $B \cap Z$ is constant and that its value is half of an odd integer. The set K is the part of Z in the closed, concentric ball to B with half of B's radius.

**Proposition 9.2**: *There exists $\kappa > 1$, and given a set $K \subset Z$ as just described, there exists $\kappa_K > \kappa$; these numbers $\kappa$ and $\kappa_K$ having the following significance: Fix $p \in K$. Supposing that $r \in (0, \kappa_K^{-1})$, then $\underline{N}_{(p)}(r) \leq \kappa \, e^{-\sqrt{|\ln(r)|}/\kappa}$.*

This proposition is proved in Section 9f. Sections 9c-9e supply background for the proof.

### a) The structure of $Z_*$

A point $p \in Z$ is said to be a point of continuity for $N_{(\cdot)}(0)$ when the following condition is met: If $\{p_i\}_{i \in \{1,2,\ldots\}}$ is a sequence in $Z$ that converges with limit $p$, then $\lim_{i \to \infty} N_{(p_i)}(0) = N_{(p)}(0)$. This subsection first states and subsequently proves a lemma that strengthens the conclusions of Lemma 6.2 in the case when $p \in Z$ is a point of continuity for $N_{(\cdot)}(0)$ where $N_{(p)}(0)$ is half of an odd integer. To set the stage for this first lemma, fix an oriented, orthonormal frame for $TX|_p$ and parallel transport the latter along the geodesic arcs from p to define an orthonormal frame for $TX$ at each point in radius $r_0$ ball centered at p. The upcoming lemma uses L to denote this frame.

**Lemma 9.3**: *Suppose that $p \in Z$ is a point of continuity for $N_{(\cdot)}(0)$ with $N_{(p)}(0) = \frac{1}{2} \mod(\mathbb{Z})$. Fix $\varepsilon \in (0, \frac{1}{1000})$ and there exists $r_\varepsilon \in (0, \varepsilon r_0)$ with the following significance: Let q denote a point in Z from the radius $r_\varepsilon$ ball centered at p and let $\lambda \in (0, r_\varepsilon]$. There is a data set $(Z_*, \mathcal{I}_*, \nu_*)$ from the $N_0 = N_{(p)}(0)$ version of (5.1) with the properties listed in the bullets given below. These bullets also refer to the data set $(Z_T, \mathcal{I}_T, \nu_T)$ that is constructed using $T = (q, L, \lambda)$ using the rules at the start of Section 4.*
- *The data set $(Z_*, \mathcal{I}_*, \nu_*)$ is described by the fourth bullet of Lemma 5.4 or the third and fourth bullets of Lemma 5.5 as the case may be.*
- *Suppose that $q \in Z_*$ and $\varepsilon < |q| < \varepsilon^{-1}$. There is a point from $Z_T$ in the normal plane to $Z_*$ through q with distance at most $\varepsilon$ from q.*
- *Each $|x| < \varepsilon^{-1}$ point in $Z_T$ has distance less than $\varepsilon$ from some point in $Z_*$.*
- *Let $\mathcal{T}_\varepsilon$ denote the $|x| < \varepsilon^{-1}$ part of the radius $\varepsilon$ tubular neighborhood of $Z_*$. There is an isometric isomorphism between $\mathcal{I}_T$ and $\mathcal{I}_*$ on the complement of $\mathcal{T}_\varepsilon$ in the $|x| < \varepsilon^{-1}$ ball that allows $\nu_T$ to be viewed as a section of $\mathbb{V}^+ \otimes \mathcal{I}_*$ on this same complement of $\mathcal{T}_\varepsilon$. Viewed in this way, $\nu_T$ differs from $\nu_*$ on the complement of $\mathcal{T}_\varepsilon$ in the $|x| < \varepsilon^{-1}$ ball by a section of $\mathbb{V}^+ \otimes \mathcal{I}_*$ with $C^k$ norm less than $\varepsilon$ for all integers $k < \varepsilon^{-1}$.*

This lemma has the following immediate corollary:

**Lemma 9.4**: *Suppose that $p \in Z$ is a point of continuity for $N_{(\cdot)}(0)$ with $N_{(p)}(0) = \frac{1}{2} \mod(\mathbb{Z})$. Every limit data set $(Z_*, \mathcal{I}_*, \nu_*)$ from p's version of Proposition 4.1 is described by the fourth bullet of Lemma 5.4 when the dimension of X is three, and by the third and fourth bullets of Lemma 5.5 when the dimension of X is four.*

  This next lemma strengthens Lemma 9.3 in the case when p is a strongly continuous point in Z with $N_{(p)}(0)$ being half an odd integer. By way of a reminder, the point p is said to be strongly continuous when $N_{(\cdot)}(0)$ is constant on a neighborhood in Z of p. This upcoming Lemma 9.5 refers to a set K that is defined as follows: Suppose that B denotes a ball in X such that $N_{(\cdot)}(0)$ is constant and half an odd integer on $B \cap Z$. Use $B_\diamond$ to denote the closed ball with B's center and half B's radius; then let K denote $Z \cap B_\diamond$.

**Lemma 9.5**: *Suppose that K is as just described. Fix $\varepsilon \in (0, \frac{1}{100})$ and there exists a K and $\varepsilon$ dependent $r_{K,\varepsilon} \in (0, \varepsilon r_0)$ with the following significance: Let p denote a point in K, let L denote a chosen oriented, orthonormal frame for $TX|_p$, and let $\lambda \in (0, r_{K,\varepsilon}]$. There is a data set $(Z_*, \mathcal{I}_*, \nu_*)$ from the $N_0 = N_{(p)}(0)$ version of (5.1) with the properties listed in the bullets given below. These bullets also refer to the data set $(Z_T, \mathcal{I}_T, \nu_T)$ that is constructed using $T = (p, L, \lambda)$ using the rules at the start of Section 4.*
- *The data set $(Z_*, \mathcal{I}_*, \nu_*)$ is described by the fourth bullet of Lemma 5.4 or the third and fourth bullets of Lemma 5.5 as the case may be.*
- *Suppose that $q \in Z_*$ and $\varepsilon < |q| < \varepsilon^{-1}$. There is a point from $Z_T$ in the normal plane to $Z_*$ through q with distance at most $\varepsilon$ from q.*
- *Each $|x| < \varepsilon^{-1}$ point in $Z_T$ has distance less than $\varepsilon$ from some point in $Z_*$.*
- *Let $\mathcal{T}_\varepsilon$ denote the $|x| < \varepsilon^{-1}$ part of the radius $\varepsilon$ tubular neighborhood of $Z_*$. There is an isometric isomorphism between $\mathcal{I}_T$ and $\mathcal{I}_*$ on the complement of $\mathcal{T}_\varepsilon$ in the $|x| < \varepsilon^{-1}$ ball that allows $\nu_T$ to be viewed as a section of $\mathbb{V}^+ \otimes \mathcal{I}_*$ on this same complement of $\mathcal{T}_\varepsilon$. Viewed in this way, $\nu_T$ differs from $\nu_*$ on the complement of $\mathcal{T}_\varepsilon$ in the $|x| < \varepsilon^{-1}$ ball by a section of $\mathbb{V}^+ \otimes \mathcal{I}_*$ with $C^k$ norm less than $\varepsilon$ for all integers $k < \varepsilon^{-1}$.*

The rest of this subsection proves Lemmas 9.3 and 9.5.

***Proof of Lemma 9.3***: Fix $\varepsilon \in (0, \frac{1}{1000}]$. Since p is a point of continuity for $N_{(\cdot)}(0)$, Lemma 6.2 supplies $\rho_\varepsilon \in (0, \varepsilon r_0)$ such that $N_{(\cdot)}(0)$ differs by at most $\varepsilon$ from $N_{(p)}(0)$ in the radius $\rho_\varepsilon$ ball centered at p and so that the conclusions of the Lemma 6.2 hold for points from Z in this same ball. Thus, if $q \in Z$ and $dist(q, p) < \rho_\varepsilon$ and if $\lambda \in (0, \rho_\varepsilon]$, then Lemma 6.2 takes as input the data set $T = (q, L, \rho_\varepsilon)$ and supplies as output a data set $(Z_*, \mathcal{I}_*, \nu_*)$ that is described by (5.1) such that the $|x| < \varepsilon^{-1}$ part of $Z_T$ is in $\mathcal{T}_\varepsilon$ and such that the fourth bullet of Lemma 9.3 hold. To prove that $r_\varepsilon \in (0, \rho_\varepsilon)$ can be found with a data set $(Z_*, \mathcal{I}_*, \nu_*)$ that obeys all of the conditions of the lemma when $\lambda \in (0, r_\varepsilon]$, assume to the contrary that no such $r_\varepsilon$ exists so as to derive nonsense. If no such $r_\varepsilon$ exists, then there exists a sequence $\{T(i) = (p_i, L_i, \lambda_i)\}_{i \in \{1,2,...\}}$ with the following properties: The sequence $\{p_i\}_{i \in \{1,2,...\}}$ is in $Z \cap B_{\rho_\varepsilon}$ and it converges to p; and the sequence $\{\lambda_i\}_{i \in \{1,2,...\}}$ is in $(0, \rho_\varepsilon]$ and it converges to zero. Moreover, for each index $i \in \{1, 2, ...\}$, there is no set $(Z_*, \mathcal{I}_*, \nu_*)$ from (5.1) that is described by the first bullet of Lemma 9.3 and obeys the $T = T(i)$ version of the second, third and fourth bullets of Lemma 9.3. The four steps that follow generate nonsense from this sequence $\{T(i)\}_{i \in \{1,2,...\}}$.

Step 1: Keeping in mind that $\{p_i\}_{i \in \{1,2,...\}}$ is a sequence in Z that converges to p and that $\{\lambda_i\}_{i \in \{1,2,...\}}$ is a sequence of positive numbers that converges to zero, invoke Proposition 4.1 using the sequence $\{T(i)\}_{i \in \{1,2,...\}}$ as input and let $(Z_*, \mathcal{I}_*, \nu_*)$ denote a limit data set and let $\{T(i)\}_{i \in \Lambda} \subset \{T(i)\}_{i \in \{1,2,...\}}$ denote the corresponding subsequence given by Proposition 4.1. It follows from the manner of convergence dictated by Proposition 4.1 that the third and fourth bullets of Lemma 9.3 are satisfied with T being any sufficiently large $i \in \Lambda$ version of $T(i)$.

Step 2: This step proves the following assertion:

*If $i \in \Lambda$ is sufficiently large, then there is a point in Z with distance $\lambda_i$ from $p_i$.*
(9.2)

To see that this is so, suppose to the contrary that there is an infinite subsequence $\vartheta \subset \Lambda$ with there being no such point if $i \in \vartheta$. To generate nonsense from this assumption, note first that each $i \in \vartheta$ version of the bundle $\mathcal{I}_{T(i)}$ is defined on the whole of the $|x| = 1$ sphere. It is therefore isomorphic to the product $\mathbb{R}$ bundle on the whole of this sphere. The manner of convergence dictated by Proposition 4.1 then implies that the bundle $\mathcal{I}_*$ must be isomorphic to the product $\mathbb{R}$ bundle on $\mathbb{R}^d - Z_*$. Now, if $\mathcal{I}_*$ is the product $\mathbb{R}$ bundle, then $\nu_*$ can be viewed as a $\mathcal{D}_0$-harmonic section of $\mathbb{V}^+$ on $\mathbb{R}^d - Z_*$ that vanishes on $Z_*$. This being the case, standard elliptic regularity theory says that $\nu_*$ extends over $Z_*$ as a smooth section of $\mathbb{V}^+$. In particular, this extended version of $\nu_*$ is real analytic and its Taylor's series centered at the origin can be used to prove that the number $N_0$ that appears in the $(Z_*, \mathcal{I}_*, \nu_*)$ version of (5.1) is an integer. Since $N_0$ is equal to $N_{(p)}(0)$, this is nonsense.

Step 3: What with (9.2) being true, it follows from Lemmas 6.2 and 5.6 that $(Z_*, \mathcal{I}_*, \nu_*)$ is described by the fourth bullet of Lemma 5.4 or the third bullet of Lemma 5.5 as the case may be. This says that $Z_*$ is described by the first bullet of Lemma 9.3 when the dimension of X is three. Suppose for the moment that the dimension of X is four. The paragraphs that follow explains why $Z_*$ contains two $|x| = 1$ points that are not colinear with the origin and are such that $\tilde{N}_{(\cdot)} = N_{(p)}(0)$. Granted that there exists two such points, then $(Z_*, \mathcal{I}_*, \nu_*)$ is described by both the third and fourth bullets of Lemma 5.5. This is the assertion of the first bullet of Lemma 9.3 when the dimension of X is four.

The argument in Step 2 supplies at least one $|x| = 1$ point in $Z_*$ where $\tilde{N}_{(\cdot)} = N_{(p)}(0)$. Denote this point by q. Either the following statement is true or it is false:

*There is an increasing sequence $\{i_n\}_{n=1,2,\ldots} \subset \Lambda$ such that if $n \in \{1, 2, \ldots\}$, then all $|x| = 1$ points in the $T = T(i_n)$ version of $Z_T$ have distance at most $\frac{1}{n}$ from either q or -q.*

(9.3)

If this statement is false, then it follows from Lemmas 6.2 and 5.6 that there are at least two $|x| = 1$ points in $Z_*$ that are not colinear with 0 and have $\tilde{N}_{(\cdot)} = N_{(p)}(0)$.

To prove that (9.3) is always false, suppose for the sake of argument that it is true so as to generate nonsense. If it is true, then it follows directly that if $n \in \{1, 2, \ldots\}$, then the $T = T(i_n)$ version of the bundle $\mathcal{I}_T$ is defined in the complement of two disjoint radius $\frac{1}{n}$ balls in the $|x| = 1$ sphere. This subspace of the sphere is simply connected, and so it follows that the $T = T(i_n)$ version of $\mathcal{I}_T$ is isomorphic to the product bundle on the complement of these two balls. The manner of convergence dictated in Proposition 4.1 implies that the bundle $\mathcal{I}_*$ is isomorphic to the product bundle on the complement of a line through the origin. Such an isomorphism can be used to view $\nu_*$ as a $\mathcal{D}_0$-harmonic section of $\mathbb{V}^+$ on $\mathbb{R}^d - Z_*$ that vanishes on $Z_*$. As in Step 1, this implies in turn that $N_{(p)}(0)$ is an integer, which is nonsense.

Step 4: The subsequent paragraph in this step proves that the second bullet of Lemma 9.3 holds when T is any sufficiently large i version of T(i). The fact that this is so constitutes the desired nonsense that proves Lemma 9.3 because the assumption at the outset posited the following: If $i \in \Lambda$, then the four bullets of the lemma can not all be true using $T = T(i)$ and using any data set $(Z_*, \mathcal{I}_*, \nu_*)$ as described in (5.1).

Assume that there is a subsequence $\vartheta \subset \Lambda$ such that if $i \in \vartheta$, then the second bullet of Lemma 9.3 is not true for $T = T(i)$. To generate nonsense from this assumption, reword this assumption to say the following: For each $i \in \vartheta$, there is a point $q_i \in Z_*$ with $|q_i|$ between $\varepsilon$ and $\varepsilon^{-1}$ and such that the radius $\varepsilon$ disk centered at $q_i$ in the normal plane to $Z_*$ has no points from $Z_{T(i)}$. Now, fix $i \in \vartheta$ and let $\gamma_i$ denote the loop of radius $\frac{1}{2}\varepsilon$ in the normal plane to the linear space $Z_*$ at $q_i$. The bundle $\mathcal{I}_{T(i)}$ is isomorphic to the product $\mathbb{R}$ bundle on this loop. However, if i is sufficiently large, then $\mathcal{I}_*$ will be isomorphic to $\mathcal{I}_{T(i)}$ on the complement of the radius $\varepsilon^2$ tubular neighborhood of $Z_*$ in the $|x| \leq \varepsilon^{-2}$ ball and thus on $\gamma_i$. This is nonsense because it follows from the fourth bullet of Lemma 5.5 and Lemma 5.1 that $\mathcal{I}_*$ is definitely not isomorphic to the product $\mathbb{R}$ bundle such a circle.

***Proof of Lemma 9.5***: If no such $r_{K,\varepsilon}$ exists, then there exists a sequence of data sets $\{T(i) = (p_i, L_i, \lambda_i)\}_{i \in \{1,2,...\}}$ with the following properties: The sequence $\{p_i\}_{i \in \{1,2,...\}}$ is in K and the sequence $\{\lambda_i\}_{i \in \{1,2,...\}}$ converges to zero. Moreover, for each $i \in \{1, 2, ...\}$, there is no set $(Z_*, \mathcal{I}_*, \nu_*)$ from (5.1) that is described by the first bullet of Lemma 9.5 and also satisfies the $T = T(i)$ version of the lemma's second, third and fourth bullets. The next paragraph explains why the existence of the sequence $\{T(i)\}_{i \in \{1,2,...\}}$ leads nonsense.

Since K is compact, there is a subsequence $\Xi \subset \{1, 2, ...\}$ and a point $p \in K$ such that $\{p_i\}_{i \in \vartheta}$ converges to a point p in K. The point p is strongly continuous so it is a point of continuity for $N_{(\cdot)}(0)$. Since $N_{(p)}(0)$ is half an odd integer, the point p's version of Lemma 9.3 can be used with T being any sufficiently large i element from the sequence $\{T(i)\}_{i \in \Xi}$. Fix such $i \in \Xi$ and let $(Z_*, \mathcal{I}_*, \nu_*)$ denote a data set supplied by Lemma 9.3. This data set obeys all of the conditions in Lemma 9.5; and this is nonsense because it was assumed at the outset that no $(Z_*, \mathcal{I}_*, \nu_*)$ obeys all of Lemma 9.5's conditions.

**b) Proof of Proposition 9.1**

Granted that the first assertion of the proposition is true, then the second assertion restates Lemma 9.4. The proof of the first assertion depends on the two lemmas that follow directly. To set the stage for the first lemma, suppose that B is a ball in X that is well inside a Gaussian coordinate chart. Assume that each point in $B \cap Z$ is strongly continuous and that $N_{(\cdot)}(0)$ is constant on $B \cap Z$ with value equal to half of an odd integer. Define K to be the part of Z in the closed ball that is concentric to B with half B's radius.

**Lemma 9.6**: *Suppose that K is as described above. There exists a number $\kappa_K > 1$ with the following significance: Fix $p \in K$. Suppose that $p' \in K$ has $\mathrm{dist}(p, p') < \kappa_K^{-1}$ and that $\lambda$ and $\lambda'$ are two numbers from $(0, \kappa_K^{-1}]$. Let $r$ denote the larger of the numbers $\lambda$ and $\lambda'$. Fix an orthonormal frame for TX at $p'$ to be denoted by L. Use the data sets $T = (p', L, \lambda)$ and $T' = (p', L, \lambda')$ as instructed at the beginning of Section 4 to define $\nu_T$ and $\nu_{T'}$. Let $\mathfrak{s}$ denote a given smooth section of $\mathbb{V}^+ \otimes \mathbb{V}^+$. Then*

$$\int_{|x|=1} |\langle \mathfrak{s}, |\nu_T|^{-1}(\nu_T \otimes \nu_T)\rangle - \langle \mathfrak{s}, |\nu_{T'}|^{-1}(\nu_{T'} \otimes \nu_{T'})\rangle|^2 \leq c_{\mathfrak{s}} \, e^{-\sqrt{|\ln(r)|}/\kappa}$$

*with $c_{\mathfrak{s}}$ being a constant that depends only on $\mathfrak{s}$.*

This lemma is proved momentarily.

To set the stage for the second lemma, suppose that B and K are as before. Let p again denote a point from K. Suppose that $\{p_i\}_{i \in \{1,2,...\}}$ and $\{p_i'\}_{i \in \{1,2,...\}}$ are sequences in K that converge to p. Fix an orthonormal frame at p to be denoted by L and parallel transport this frame along the short geodesic arcs from p to define corresponding frames for TX at each point in the sequence $\{p_i\}_{i \in \{1,2,...\}}$ and each point from the sequence $\{p_i'\}_{i \in \{1,2,...\}}$. The data set sequence $\{T(i) = (p_i, L, \lambda_i = \text{dist}(p_i, p_i'))\}_{i \in \{1,2,...\}}$ can be used as input to Proposition 4.1 as can the sequence $\{T'(i) = (p_i', L, \lambda_i = \text{dist}(p_i, p_i'))\}_{i \in \{1,2,...\}}$. These sequences differ in that they have different index dependent point sets, these being $\{p_i\}_{i \in \{1,2,...\}}$ and $\{p_i'\}_{i \in \{1,2,...\}}$ The following lemma makes an assertion to the effect that the output from Proposition 4.1 for these two sequences agree.

**Lemma 9.7**: *Let $\{T(i)\}_{i \in \{1,2,...\}}$ and $\{T'(i)\}$ be as described above. Let $\Lambda \subset \{1, 2, ...\}$ and $(Z_*, \mathcal{I}_*, \nu_*)$ denote limit data supplied by Proposition 4.1 using the sequence $\{T(i)\}_{i \in \{1,2,...\}}$ for input. In particular, $\{(Z_{T(i)}, \mathcal{I}_{T(i)}, \nu_{T(i)})\}_{i \in \Lambda}$ converges to $(Z_*, \mathcal{I}_*, \nu_*)$ in the manner that is dictated in the fifth bullet of Proposition 4.1. Then $\{(Z_{T'(i)}, \mathcal{I}_{T'(i)}, \nu_{T'(i)})\}_{i \in \Lambda}$ also converges to $(Z_*, \mathcal{I}_*, \nu_*)$ in the manner that is dictated in the fifth bullet of Proposition 4.1.*

This lemma is also proved momentarily.

*Proof of Proposition 9.1*: What with Lemma 9.4, it is only necessary to prove the assertion about the unique tangent cones. This is done in two steps.

<u>Step 1</u>: Suppose that $p \in Z$ is a strongly continuous point with $N_{(p)}(0)$ being half an odd integer. Let $\{p_i\}_{i \in \{1,2,...\}}$ denote a sequence in Z that converges to p and let $\{\lambda_i\}_{i \in \{1,2,...\}}$ and $\{\lambda_i'\}_{i \in \{1,2,...\}}$ denote two sequences in $(0, r_0]$ that converge to 0. Fix an orthonormal frame L at p and parallel transport this frame along the short geodesic arcs from p to define corresponding frames for TX at each point in the sequence $\{p_i\}_{i \in \{1,2,...\}}$. Use the sequence $\{T(i) = (p_i, L, \lambda_i)\}_{i \in \{1,2,...\}}$ as input for Proposition 4.1 and let $(Z_*, \mathcal{I}_*, \nu_*)$ denote a limit data set provided by this same proposition. Let $(Z_*', \mathcal{I}_*', \nu_*')$ denote a limit data set from Proposition 4.1 with the input sequence being $\{T(i) = (p_i, L, \lambda_i')\}_{i \in \{1,2,...\}}$. It follows from Lemma 9.6 that $\nu_* \otimes \nu_* = \nu_*' \otimes \nu_*'$. This implies that $Z_* = Z_*'$ and that $\langle \nu_*', \nu_* \rangle$ defines an isomorphism between $\mathcal{I}_*$ and $\mathcal{I}_*'$ that identifies $\nu_*$ with $\nu_*'$.

<u>Step 2</u>: Let p be as in Step 1. Suppose that $\{p_i\}_{i \in \{1,2,...\}}$ and $\{p_i'\}_{i \in \{1,2,...\}}$ are two sequences in Z that converge to p and suppose that $\{\lambda_i\}_{i \in \{1,2,...\}}$ and $\{\lambda_i'\}_{i \in \{1,2,...\}}$ are two sequences in $(0, r_0]$ that converge to 0. Fix an orthonormal frame L at p and parallel transport this frame along the short geodesic arcs from p to define corresponding frames for TX at each point in the sequence $\{p_i\}_{i \in \{1,2,...\}}$ and each point from the sequence $\{p_i'\}_{i \in \{1,2,...\}}$. It follows from Step 1 that any two limit data set $(Z_*, \mathcal{I}_*, \nu_*)$ from Proposition 4.1 with input sequence $\{T(i) = (p_i, L, \lambda_i)\}_{i \in \{1,2,...\}}$ are identical, and that this one limit data set is the same as the limit data set that is obtained from Proposition 4.1 using the input sequence $\{(p_i, L, \text{dist}(p_i, p_i'))\}$. By the same token, all limit data sets $(Z_*', \mathcal{I}_*', \nu_*')$ from Proposition 4.1 with input sequence $\{T'(i) = (p_i', L, \lambda_i')\}_{i \in \{1,2,...\}}$ are identical, and that this one limit data set is the same as the limit data set that is obtained from Proposition 4.1 using the input sequence $\{(p_i´, L, \text{dist}(p_i, p_i´)\}_{i \in \{1,2,...\}}$.

Meanwhile, it follows from Lemma 9.7 that the respective limit data sets from Proposition 4.1 using as input the sequence $\{(p_i, L, \text{dist}(p_i, p_i´))\}_{i \in \{1,2,...\}}$ and as input the sequence $\{(p_i´, L, \text{dist}(p_i, p_i´)\}_{i \in \{1,2,...\}}$ are identical.

***Proof of Lemma 9.6***: No generality is lost by assuming in what follows that $\lambda > \lambda´$. Let n be the positive integer such that $\lambda´$ is between $2^{-n}\lambda$ and $2^{-n+1}\lambda$. Invoke the formula in (4.30) successively n-1 times for each choice of $k \in \{1, 2, \ldots n-1\}$ using $\lambda_1 = 2^{-k+1}\lambda$ and $\lambda_2 = 2^{-k}\lambda$; and then invoke it one final time with $\lambda_1 = 2^{-n+1}\lambda$ and $\lambda_2 = \lambda´$ to see that

$$(\int_{|x|=1} |\langle \mathfrak{s}, |v_T|^{-1}(v_T \otimes v_T)\rangle - \langle \mathfrak{s}, |v_{T´}|^{-1}(v_{T´} \otimes v_{T´})\rangle|^2 )^{1/2} \leq c_\mathfrak{s} (\sum_{k=1,\ldots,n} |N_{(p´)}(2^{-k}\lambda) - N_0|^{1/2} + \lambda^{1/2})$$

(9.4)

with $c_\mathfrak{s}$ denoting here and subsequently a number that is greater than 1 which depends only on the section $\mathfrak{s}$. Invoke Proposition 9.2 to bound the right hand side of (9.4) by

$$c_\mathfrak{s} ( e^{-\sqrt{|\ln \lambda|}/c_0} \sum_{k=1,\ldots,n} e^{-\sqrt{k}/c_0} + \lambda^{1/2}) .$$

(9.5)

Since the sum indexed by $k \in \{1, \ldots, n\}$ that appears here is bounded by $c_0$, this inequality implies the bound that is asserted by Lemma 9.6.

***Proof of Lemma 9.7***: For each $i \in \Lambda$, the composition of first $\phi_{T(i)}$ and then the inverse of the map $\phi_{T´(i)}$ identifies $Z_{T(i)}$ with $Z_{T´(i)}$. Denote this map by $\psi_i$. The map $\psi_i$ sends the origin of $\mathbb{R}^d$, this being a point in $Z_{T(i)}$, to an $|x| = 1$ point in $Z_{T´(i)}$. Meanwhile, some $|x| = 1$ point in $Z_{T(i)}$ is mapped by $\psi_i$ to the origin in $\mathbb{R}^d$. Denote the latter point by $-q_i$. The map $\psi_i$ differs by $\mathcal{O}(\text{dist}(p_i, p_i´)^2)$ from the translation $x \to x + q_i$. The map $\psi_i$ also identifies $\mathcal{I}_{T(i)}$ with $\mathcal{I}_{T´(i)}$ so as to identify $\nu_{T(i)}$ with $\nu_{T´(i)}$.

It follows from Lemma 9.4 that the data set $(Z_*, \mathcal{I}_*, \nu_*)$ is described by the third of Lemma 5.4 or the third and fourth bullets of Lemma 5.5 as the case may be. This implies that it is nearly invariant with respect to the translation by $q_i$ when i is large since $-q_i$ is very nearly a point in $Z_*$ when i is large. This fact and the fact that $\psi_i$ identifies $\mathcal{I}_{T(i)}$ with $\mathcal{I}_{T´(i)}$ so as to identify $\nu_{T(i)}$ with $\nu_{T´(i)}$ implies what is asserted by the lemma.

### c) The definition of $(Z_{*r}, \mathcal{I}_{*r}, \nu_{*r})$

By way of a look ahead, the proof of Proposition 9.2 has much in common with the proofs of Lemmas 5.3 and 5.7. In particular, many of the constructions that are used in Section 9f's proof of Proposition 9.2 have direct analogs in the proofs of these lemmas. This subsection and the next two supply some of these analogs.

To set the background for this subsection, suppose that $B \subset X$ is a ball that is well inside a Gaussian coordinate chart and is such that $N_{(\cdot)}(0)$ on $B \cap Z$ is constant with its value being half an odd integer. Use K to denote the part of Z in the closed, concentric ball to B with half of the radius of B. Given some small positive $\varepsilon$, let $r_{K,\varepsilon}$ denote the relevant version of $r_{K,\varepsilon}$ from Lemma 9.5. Suppose that $p \in K$, that $L$ is an oriented, orthonormal frame for $TX|_p$ and that $r \in (0, r_{K,\varepsilon})$. The four parts that follow in this subsection construct a canonical $(Z_*, \mathcal{I}_*, \nu_*)$ for the $T = (p, L, \lambda = r)$ version of Lemma 9.3. This canonical data set is denoted below by $(Z_{*r}, \mathcal{I}_{*r}, \nu_{*r})$. The next subsection compares $(Z_{*r}, \mathcal{I}_{*r}, \nu_{*r})$ with the data set $(Z_T, \mathcal{I}_T, \nu_T)$.

*Part 1*: Fix an orthonormal frame for $TX|_p$ to be denoted by $L$ and then use parallel transport along the geodesic arcs through p to view $L$ as a framing of TX on the ball B. Let $\{\nabla_a\}_{a \in \{1,\ldots,d\}}$ denote the covariant derivatives along the frame vectors. Fix a positive number r less than half the distance from p to the boundary of B. With apologies to (5.16), let T(r) denote the symmetric, non-negative definite $d \times d$ matrix with entries $\{T_{ab}(r)\}_{a,b \in \{1,\ldots,d\}}$ given by

$$T_{ab}(r) = \frac{1}{r^{d-2} K_{(p)}(r)^2} \int_{B_r} \langle \nabla_a \nu, \nabla_b \nu \rangle \,.$$

(9.6)

It follows from the fourth bullet of Lemma 9.5 that if $r \in (0, r_{K,\varepsilon})$, then T(r) has two eigenvalues that are $\mathcal{O}(1)$ and $d-2$ eigenvalues that are bounded by $c_0 \varepsilon$. Note that the trace of T(r) is $N_{(p)}(r)$.

As explained directly, the small eigenvalues of T(r) are less than $c_0(\underline{N}_{(p)}(r) + r)$. To prove this bound on the small eigenvalues, let T denote the data set (p,L,r). If û is a length 1 vector in $\mathbb{R}^d$, write the inner product of û with T(r)û as $\langle \hat{u}, T(r) \hat{u} \rangle$. By definition,

$$\langle \hat{u}, T(r) \hat{u} \rangle = \frac{1}{r^{d-2} K_{(p)}(r)^2} \int_{B_r} |\nabla_{\hat{u}} \nu|^2$$

(9.7)

with $\nabla_{\hat{u}}$ denoting the directional covariant derivative in $B_r$ along the unit vector that û defines using the frame $L$. Now suppose that û is a length 1 vector in $\mathbb{R}^d$ and $t \in [2, 4]$ are such that $t\hat{u} \in Z_T$. Keep in mind in this regard what is said by the second and third bullets of Lemma 9.5. Let q denote the point $\phi_r(t\hat{u})$ which is a point in Z. Introduce by way of notation $r_q$ to denote the function on X that gives the distance to q, and let $\nabla_{r_q} \nu$ denote the directional covariant derivative of $\nu$ along the outward pointing, unit tangent vector to the geodesic arcs from q. Supposing that $s \in (0, r]$, a formula for $\nabla_{\hat{u}} v$ on $\partial B_s$ that is analogous to (5.17) writes $\nabla_{\hat{u}} v$ as

$$\operatorname{tr} \nabla_{\hat{u}} v = s \nabla_r v - N_{(p)}(s) v - r_q \nabla_{r_q} v + N_{(q)}(r_q) v + \nabla_{\eta} v + (N_{(p)}(s) - N_{(q)}(r_q)) v \tag{9.8}$$

with $\eta$ being a vector field defined on the complement of q with norm bounded by $c_0 r^2$. This vector field accounts for the fact that the metric may not be flat near p. The formula in (9.8) leads directly to the inequality

$$|\nabla_{\hat{u}} v|^2 \leq c_0 |\nabla_r v - \tfrac{1}{s} N_{(p)}(s) v|^2 + c_0 (\tfrac{r_q}{r})^2 |\nabla_{r_q} v - \tfrac{1}{r_q} N_{(q)}(r_q) v|^2 + c_0 \tfrac{1}{r^2} \underline{N}_{(p)}(\tfrac{1}{3} r)^2 |v|^2 + c_0 r^2 |\nabla v|^2. \tag{9.9}$$

because the definition of $\underline{N}_{(p)}$ in (9.1) implies that $|(N_{(p)}(s) - N_{(q)}(r_q))| \leq 2 \underline{N}_{(p)}(\tfrac{1}{3} r)$. The inequality in (9.9) is now used to bound the right hand side of (9.7) by

$$c_0 \frac{1}{r^{d-2} K_{(p)}(r)^2} \int_{B_r} |\nabla_r v - \tfrac{1}{r} N_{(p)} v|^2 + c_0 \frac{1}{r^{d-2} K_{(p)}(r)^2} \int_{A_r(q)} |\nabla_{r_q} v - \tfrac{1}{r_q} N_{(q)} v|^2 + c_0 \underline{N}_{(p)}(r)^2 + c_0 r, \tag{9.10}$$

where $A_r(q)$ denotes here the spherical annulus centered at the point q with inner radius equal to $(t-1)r$ and outer radius $tr$. Note that the integral of $r^{-2} \underline{N}_{(p)}(r)^2 |v|^2$ from (9.9) is dealt with here by using the fact asserted by Lemma 3.2 to the effect that $K_{(p)}$ is an increasing function. This equation is the analog of (5.18).

Keeping in mind that $K_{(q)}((t-1)r) \leq c_0 K_{(p)}(r)$ because of Lemma 6.2, and keeping in mind that $K_{(p)}$ and $K_{(q)}$ are increasing functions, the version of Lemma 3.2's second bullet for p and then the version with q replacing p can be used to bound (9.10) by

$$c_0 (N_{(p)}(r) - N_0) + c_0 (N_{(q)}(4r) - N_0) + c_0 r. \tag{9.11}$$

Given the definition of $\underline{N}_{(p)}(r)$ in (9.1), this bound for (9.7)'s right hand side implies that

$$\langle \hat{u}, T(r) \hat{u} \rangle \leq c_0 (\underline{N}_{(p)}(\tfrac{1}{3} r) + r). \tag{9.12}$$

Since $T(r)$ is non-negative definite, the inequality in (9.12) implies in turn the claim that the d-2 small eigenvalues of $T(r)$ are bounded by $c_0(\underline{N}_{(p)}(r) + r)$.

Use $\Pi_r$ in what follows to denote the linear subspace in $\mathbb{R}^d$ spanned by the d-2 small eigenvalues of the matrix $T(r)$.

*Part 2*: The constructions in this part of the proof are analogous to those in Parts 5 and 6 of the proof of Lemma 5.7. The constructions in these parts of Lemma 5.7's proof average the section $v_*$ in (5.20). The constructions below do not average $v_T$ but average instead the section $v_T \otimes v_T$. There are some advantages and disadvantages to averaging the latter as opposed to $v_T$, the main advantage being that there is no need to worry about the whereabouts of $Z_T$ when averaging $v_T \otimes v_T$. A construction that averages $v_T$ will work at the expense of some modifications in subsequent parts of the proof.

To start the averaging construction, rotate the coordinates on $\mathbb{R}^d$ so that the resulting Euclidean coordinates have $\Pi_{3r}$ being the subspace where the first two coordinates, $(x_1, x_2)$, are zero. Use the function $\chi$ to construct a smooth, non-negative function of the distance from the origin in $\Pi_{3r}$ that is equal to 1 where the distance to the origin is less than 2 and equal to zero where the distance is greater than $\sqrt{5}$. Denote this function by b and let $\langle b \rangle$ denote its integral on $\Pi_{3r}$. Use $\beta$ to denote the function $\langle b \rangle^{-1} b$. The function $\beta$ has integral 1 on $\Pi_{3r}$. Let $D \subset \Pi_{3r}$ denote the disk of radius $\sqrt{5}$ thus containing the support of $\beta$.

Use $t$ to denote $v_T \otimes v_T$, this being a section of $\mathbb{V}^+ \otimes \mathbb{V}^+$ on the radius 3 ball in $\mathbb{R}^d$. Define a section of $\mathbb{V}^+ \otimes \mathbb{V}^+$ on the radius 2 disk in the $(x_1, x_2)$ plane by averaging $t$ in directions tangent to $\Pi_{3r}$. This section is denoted by $\underline{t}$ and it is defined by the rule

$$\underline{t}(x_1, x_2) = \int_D t(x_1, x_2, \cdot) \beta(\cdot)$$

(9.13)

Of interest in what follows is the norm of the difference between $\underline{t}$ and $t$ on certain d-1 dimensional submanifolds in the $|x| < 3$ part of $\mathbb{R}^d$. To say more, introduce by way of notation $\mathfrak{p}$ to denote the orthogonal projection on $\mathbb{R}^d$ to the $(x_1, x_2)$ plane. Suppose that $S$ is a d-1 dimensional submanifold in the $|x| \leq 2$ part of $\mathbb{R}^d$. Assume that $\mathfrak{p}$ restricts to $S$ as a submersion. Let $\upsilon$ denote an upper bound along $S$ for the norm of the inverse of $\mathfrak{p}$'s differential on T$S$ and let A to denote an upper bound for the d-3 dimensional measures of $\mathfrak{p}$'s fibers in $S$. Let $\nabla^\Pi$ denote the covariant derivative in directions tangent to $\Pi_{3r}$. An appeal to the fundamental theorem of calculus and to (9.12) leads to the bound

$$\int_S |t(\cdot) - \underline{t}(\cdot)|^2 \leq c_0 \upsilon^2 A \int_{|x|<3} |\nabla^\Pi v_T|^2 \leq \chi_0 \upsilon^2 A (\underline{N}_{(\mathfrak{p})}(r) + r).$$

(9.14)

A case that is used in what follows has $S$ being the portion of the $|x| = 1$ sphere in $\mathbb{R}^d$ where the coordinates $(x_1, x_2)$ obey $(x_1^2 + x_2^2)^{1/2} \in [\frac{1}{4}, \frac{3}{4}]$.

*Part 3*: Write the operator $\mathcal{D}_0$ as

$$\mathcal{D}_0 = \mathfrak{S}_k \nabla_k \tag{9.15}$$

with $\{\mathfrak{S}_k\}_{k=1,\ldots,d}$ being homomorphisms from $\mathbb{V}^+$ to $\mathbb{V}^-$. Use $\tau$ to denote $\mathfrak{S}_1^\dagger \mathfrak{S}_2$. As noted in Part 6 of the proof of Lemma 5.7, this is an endomorphism of $\mathbb{V}^+$ with square -1. Let $v_+$ and $v_-$ denote the respective projections of $v_\text{T}$ into the +i and -i eigenspaces of $\tau$ on either $\mathbb{V}^+ \otimes \mathcal{I}_\text{T}$ if $\mathbb{V}^+$ is complex or the complexification of $\mathbb{V}^+ \otimes \mathcal{I}_\text{T}$. The equation $\mathcal{D}v = 0$ when written in terms of $v_+$ and $v_-$ on the $|x| \leq 10$ ball in $\mathbb{R}^d$ leads to a pair of equations that have the form

$$\tau^+_{\Pi}(\nabla^\Pi)v_+ + (\nabla_1 - i\nabla_2)v_- + \mathfrak{r}_- = 0 \quad \text{and} \quad \tau^-_{\Pi}(\nabla^\Pi)v_- + (\nabla_1 + i\nabla_2)v_+ + \mathfrak{r}_+ = 0 \tag{9.16}$$

with $\mathfrak{r}_+$ and $\mathfrak{r}_-$ have norms that obey $|\mathfrak{r}_+| + |\mathfrak{r}_-| \leq c_0 r^2(|\nabla v_\text{T}| + |v_\text{T}|)$ and with $\tau^\pm_\Pi$ denoting certain constant coefficient endomorphisms of either $\mathbb{V}^+$ or its complexification as the case may be. The $\mathfrak{r}_+$ and $\mathfrak{r}_-$ terms account for the fact that the metric $\mathfrak{m}_\text{T}$ need not be the Euclidean metric. Make note that the covariant derivatives in (9.16) are the Euclidean covariant derivatives.

Let $t_+$ denote $v_+ \otimes v_+$ and let $t_-$ denote $v_- \otimes v_-$. The equations in (9.16) lead to equations for $t_+$ and $t_-$ that have the form

$$(\nabla_1 + i\nabla_2)t_+ = \mathfrak{a}_+ \quad \text{and} \quad (\nabla_1 - i\nabla_2)t_- = \mathfrak{a}_-, \tag{9.17}$$

where $\mathfrak{a}_+$ and $\mathfrak{a}_-$ have norms that obey $|\mathfrak{a}_+| + |\mathfrak{a}_-| \leq c_0(|\nabla^\Pi v_\text{T}| + r^2(|\nabla v_\text{T}| + |v_\text{T}|))$.

Now let $\underline{t}_+$ and $\underline{t}_-$ denote the averages of $t_+$ and $t_-$ as defined by (9.13) and let $\underline{\mathfrak{a}}_+$ and $\underline{\mathfrak{a}}_-$ denote the respective averages of $\mathfrak{a}_+$ and $\mathfrak{a}_-$. The equations in (9.17) imply that the averaged versions of $t_+, t_1, \mathfrak{a}_+$ and $\mathfrak{a}_-$ obey the analogous equations:

$$(\nabla_1 + i\nabla_2)\underline{t}_+ = \underline{\mathfrak{a}}_+ \quad \text{and} \quad (\nabla_1 - i\nabla_2)\underline{t}_- = \underline{\mathfrak{a}}_-. \tag{9.18}$$

Keep in mind that these equation have no dependence on the coordinates along $\Pi_{3r}$ and that they are equations on the radius 2 disk about the origin in the $(x_1, x_2)$ plane. It follows from the pointwise bounds on $|\mathfrak{a}_+|$ and $|\mathfrak{a}_-|$ that the $L^2$ norms of $|\underline{\mathfrak{a}}_+|$ and $|\underline{\mathfrak{a}}_-|$ are bounded by $c_0 (\underline{N}_{(p)}(r) + r)^{1/2}$.

Let $\mathcal{V}$ denote $\mathbb{V}^+ \otimes \mathbb{V}^+$ if $\mathbb{V}^+$ is defined over $\mathbb{C}$ or its complexification if not. There are $\mathcal{V}$ valued functions on the radius 2 disk about the origin in the $(x_1, x_2)$ plane, these denoted by $o_+$ and $o_-$, obeying

$$(\nabla_1 + i\nabla_2)o_+ = \underline{\mathfrak{a}}_+ \quad \text{and} \quad (\nabla_1 - i\nabla_2)o_- = \underline{\mathfrak{a}}_-$$

(9.19)

and with $L^2_1$ norms on the $(x_1^2 + x_2^2)^{1/2} \leq \frac{7}{4}$ disk that are bounded by $c_0 \, (\underline{N}_{(p)}(r) + r)^{1/2}$. The fact that $o_+$ and $o_-$ can be found with these $L^2_1$ norm bounds follows from the observation at the end of the preceding paragraph about the $c_0 \, (\underline{N}_{(p)}(r) + r)^{1/2}$ bound on the $L^2$ norms of $\underline{a}_+$ and $\underline{a}_-$.

With $o_+$ and $o_-$ as just described, it follows that $s_+ = \underline{t}_+ - o_+$ is a holomorphic section of $\mathcal{V}$ and that $s_- = \underline{t}_- - o_-$ is an anti-holomorphic section. Being that $s_+$ is holomorphic, it can be written as a convergent power series in the complex variable $z = x_1 + i x_2$. Write this power series as

$$s_+(z) = \alpha_+ z^{2N_0} + \mathit{b}_0(z) + z^{2N_0+1} \mathit{b}_1(z) ,$$

(9.20)

with $\alpha_+ \in \mathcal{V}$, with $\mathit{b}_0$ being a polynomial of degree at most $2N_0 - 1$ and with $\mathit{b}_1(z)$ being convergent power series. There is a similar formula using $\bar{z}$ for $s_-$.

The same construction can be made using the data set $T' = (p, L, \frac{1}{2} r)$. This leads to a holomorphic section of $\mathcal{V}$ to be denoted by $s_+'$. It follows from (9.18) that the $L^2$ norm of the $\mathcal{V}$ valued function $z \to s_+'(z) - 2^{N_0} s_+(\frac{1}{2} z)$ on the $|z| < \frac{7}{4}$ disk is at most $c_0 \, (\underline{N}_{(p)}(r) + r)^{1/2}$. Meanwhile, it follows from (4.30) that the same bound holds for the $L^2$ norm of the $\mathcal{V}$ valued function $s_+(z) - s_+'(z)$. These two bounds imply that the $L^2$ norms of both $\mathit{b}_0$ and $\mathit{b}_1$ on the $|z| < \frac{7}{4}$ disk are bounded by $c_0 \, (\underline{N}_{(p)}(r) + r)^{1/2}$ also. Since $\mathit{b}_0$ and $\mathit{b}_1$ are holomorphic, these $L^2$ bounds lead in turn to $c_0 \, (\underline{N}_{(p)}(r) + r)^{1/2}$ bounds for the $C^1$ norms of $\mathit{b}_0$ and $\mathit{b}_1$ on the $|z| \leq \frac{3}{2}$ disk.

*Part 4*: This part of the subsection defines the data $(Z_{*r}, \mathcal{I}_{*r}, v_{*r})$. To start, take $Z_{*r}$ to be the linear, codimension 2 subspace $\Pi_{3r}$ and take $\mathcal{I}_{*r}$ to be the non-trivial real line bundle on the complement in $\mathbb{R}^d$ of this subspace. It then follows from (9.14) and from what was just said about the norms of $\mathit{b}_0$ and $\mathit{b}_1$ that the element $\alpha_+$ in $\mathcal{V}$ is nearly decomposable in the sense that there is an element $a_{*+}$ in the $+i$ eigenspace of $\mathfrak{S}^{\dagger}_1 \mathfrak{S}_2$ in $\mathbb{V}^+$ or its complexification such that $|a_{*+} \otimes a_{*+} - \alpha_+| \leq c_0 \, (\underline{N}_{(p)}(r) + r)^{1/2}$. Such an element exists because $s_+$ differs from $v_+ \otimes v_+$ by at most $c_0 (\underline{N}_{(p)}(r) + r)^{1/2}$ on a set in the $|z| \leq \frac{3}{2}$ disk with measure greater than $c_0^{-1}$. A corresponding element $a_{*-}$ from the $-i$ eigenspace of $\mathfrak{S}^{\dagger}_1 \mathfrak{S}_2$ in $\mathbb{V}^+$ or its complexification is defined using $s_-$. If $\mathbb{V}^+$ is a vector space over $\mathbb{R}$, take $a_{*-} = \bar{a}_{*+}$. With these elements in hand, define

$$v_{*r} = (1+c)^{-1}(a_{*+} z^{N_0} + a_{*-} \bar{z}^{N_0}) ,$$

(9.21)

with $c$ being a real number with norm at most $c_0\,(\underline{N}_{(p)}(r)+r)^{1/2}$ chosen so that the integral of $|v_{*r}|^2$ on the $|x|=1$ sphere in $\mathbb{R}^3$ is equal to 1. The fact that $c$ has this upper bound follows from (9.14) and the bounds in Part 3 on $o_+$, $b_0$, $b_1$ and on their $t_-$ analogs. Note in this regard that the integral of $|v_T|^2$ on the $|x|\le 1$ ball is no greater than $(d+2N_{(p)}(0)-c_0 r^2)^{-1}$ and no less than $(d+2N_{(p)}(r)+c_0 r^2)^{-1}$ so it differs from the former by at most $c_0(\underline{N}_{(p)}(r)+r)$.

### d) Comparing $(Z_{*r}, \mathcal{I}_{*r}, \nu_{*r})$ and $(Z_T, \mathcal{I}_T, \nu_T)$

This subsection states and then proves two lemmas that compare the data sets $(Z_{*r}, \mathcal{I}_{*r}, \nu_{*r})$ and $(Z_T, \mathcal{I}_T, \nu_T)$. So as to make no mistake about the assumptions, the set K is defined as follows: There is a ball $B \subset X$ such that $N_{(\cdot)}(0)$ is constant and half an odd integer on $B \cap Z$. The set K is the part of Z in the closed ball that is concentric to B with half the radius of B.

The first lemma compares $Z_{*r}$ with $Z_T$.

**Lemma 9.8**: *With the set K given, there exists $\kappa > 1$ with the following significance: Suppose that $p \in K$ and $r \in (0, \kappa^{-1} r_0]$. Let $L$ denote a chosen, oriented orthonormal frame for $TX|_p$ and let $T = (p,L,r)$. Use this data set to define $Z_T \subset \mathbb{R}^d$ as instructed at the beginning of Section 4 and to define $Z_{*r}$ as instructed in Section 4c.*
- *Each point in the $|x| \in (\tfrac{2}{3}, 4)$ part of $Z_T$ has distance at most $\kappa(\underline{N}_{(p)}(r)+r)^{1/2}$ from $Z_{*r}$.*
- *Suppose that q is a point in the $|x| \in (\tfrac{2}{3}, 4)$ part of $Z_{*r}$.*
  a) *There is a point from $Z_T$ in the normal plane to $Z_{*r}$ through q with distance at most $\kappa(\underline{N}_{(p)}(r)+r)^{1/2}$ from q.*
  b) *The distance between any two points from $Z_T$ in the normal plane to $Z_{*r}$ through q is no greater than $\kappa e^{-1/(\kappa \underline{N}_{(p)}(r))}$.*

The next lemma compares $\nu_T$ and $\nu_{*r}$ on the $|x|=1$ sphere in $\mathbb{R}^d$.

**Lemma 9.9**: *With the set K given, there exists $\kappa > 1$ with the following significance: Suppose that $p \in K$ and $r \in (0, \kappa^{-1} r_0]$. Let $L$ denote a chosen, oriented orthonormal frame for $TX|_p$ and let $T = (p,L,r)$. Use this data set to define $Z_T \subset \mathbb{R}^d$ as instructed at the beginning of Section 4 and to define $Z_{*r}$ as instructed in Section 4c. Let $\mathcal{A}$ denote the part of the $|x|=1$ sphere in $\mathbb{R}^d$ with distance greater than $\kappa\,(\underline{N}_{(p)}(r)+r)^{1/2}$ from $Z_{*r}$. There is an isomorphism between $\mathcal{I}_T$ and $\mathcal{I}_{*r}$ on $\mathcal{A}$ identifying $\nu_T$ as a section of $\mathbb{V}^+ \otimes \mathcal{I}_{*r}$ that obeys*

$$\int_{\mathcal{A}} |\nu_T - \nu_{*r}|^2 \le \kappa(\underline{N}_{(p)}(r)+r) .$$

The rest of this subsection contains the proof of Lemma 9.8. The proof of Lemma 9.9 is in the next subsection

*Proof of Lemma 9.8*: The proof of the lemma has three parts.

  *Part 1*: Suppose that $\Pi$ and $\Pi'$ are given d-2 dimensional subspaces in $\mathbb{R}^d$. The *angular distance* between $\Pi$ and $\Pi'$ is defined as follows: Let $\mathfrak{p}$ and $\mathfrak{p}'$ denote the respective orthogonal projections to $\Pi$ and $\Pi'$. The angular distance between $\Pi$ and $\Pi'$ is hereby defined to be the positive square root of

$$1 - \tfrac{1}{d-2}\operatorname{tr}(\mathfrak{p}\,\mathfrak{p}'),$$
(9.22)

with $\operatorname{tr}(\cdot)$ denoting the trace of the indicated endomorphism of $\mathbb{R}^d$. Up to a normalization, (9.22) depicts the square of the Hilbert-Schmid norm of $\mathfrak{p} - \mathfrak{p}'$.

Given $s \in [\tfrac{1}{3}r, 3r]$, define the d-2 dimensional subspace $\Pi_s$ as instructed in Part 1 of the preceding subsection. Of immediate interest is the angular distance between $\Pi_r$ and $\Pi_s$ for $s \in [\tfrac{1}{3}r, 3r]$. This is denoted by $\eth(r,s)$. As proved directly,

$$\eth(r,s) \le c_0 (\underline{N}_{(p)}(r) + r)^{1/2}$$
(9.23)

with it understood that $s \in [\tfrac{1}{3}r, 3r]$. To prove this when $s > r$, let $\hat{u}$ denote a length 1 eigenvector of $T(s)$. Write $\hat{u} = (1+\theta^2)^{-1/2}(\hat{e} + \theta\hat{\imath})$ with $\hat{e}$ being a length 1 vector in $\Pi_r$, with $\hat{\imath}$ being a length 1 vector orthogonal to $\Pi_r$ and with $\theta \ge 0$. It follows from Lemma 9.3 that $\theta$ is small. To see how small, use the definition of $T(s)$ and $T(r)$ to see that

$$s^{d-2}K_{(p)}(s)^2 \langle \hat{u}, T(s)\hat{u}\rangle \ge (1+\theta^2)^{-1/2} r^{d-2} K_{(p)}(r)^2 (\langle \hat{e}, T(r)\hat{e}\rangle + \theta^2 \langle \hat{\imath}, T(r)\hat{\imath}\rangle).$$
(9.24)

Since $\langle \hat{\imath}, T(r)\hat{\imath}\rangle \ge c_0^{-1}$ since $K_{(p)}(s) \le c_0 K_{(p)}(r)$, this inequality is satisfied only if $\theta$ obeys the bound $\theta^2 \le c_0 \langle \hat{u}, T(s)\hat{u}\rangle$. Granted the latter bound, invoke (9.12) with s replacing r to see that $\theta^2 \le c_0 (\underline{N}_{(p)}(\tfrac{1}{3}s) + r)$. Since $s < 3r$, this bound for $\theta^2$ implies the bound that is asserted by (9.23). Much the same argument with the roles of r and s reversed proves (9.23) in the case when $s < r$.

  *Part 2*: Let $\mathcal{C} \subset \mathbb{R}^d$ denote the spherical annulus where $\tfrac{2}{3} \le |x| \le 4$. Use $\Delta_r$ to denote least upper bound of the distances between the points in $Z_r \cap \mathcal{C}$ and the linear subspace $Z_{*r}$. The next paragraph explains why

$$\Delta_r \le c_0 (\underline{N}_{(p)}(r) + r)^{1/2}.$$
(9.25)

The bound in (9.25) gives the first bullet of Lemma 9.8.

To prove this, fix $s \in [\frac{1}{3}r, 3r]$ for the moment. It is a consequence of (9.23) that the $\frac{1}{2} \le |x| \le 4$ part of $\Pi_s$ is in the radius $c_0 (\underline{N}_{(p)}(r) + r)^{1/2}$ tubular neighborhood of $Z_{*r}$. This being the case, then (9.25) follows from the following assertion: If $s \in [\frac{1}{3}r, r]$, then no point in the $|x| \in [2s, 4s]$ part of $Z_T$ has distance greater than $c_0(\underline{N}_{(p)}(r) + r)^{1/2}$ from $\Pi_s$. To prove this last assertion, suppose that $\hat{u}$ is a length 1 vector in $\mathbb{R}^d$ and that $t \in [2s, 4s]$ are such that $t\hat{u} \in Z_T$. Write $\hat{u}$ as $(1+\Delta^2)^{-1/2}(\hat{e} + \Delta\hat{\imath})$ with $\hat{e}$ being a length 1 vector tangent to the linear space $\Pi_s$, with $\hat{\imath}$ being a length 1 vector orthogonal to $\Pi_s$ and with $\Delta$ being a non-negative number. The bound in the version of (9.12) with s replacing r leads to a bound on $\Delta$ by $c_0(\underline{N}_{(p)}(\frac{1}{3}s) + s)^{1/2}$ and this is less than $c_0(\underline{N}_{(p)}(r) + r)^{1/2}$.

*Part 3*: This part proves Item a) of the second bullet of Lemma 9.8. Supposing that Item a) is false, then what is said in the first bullet of the lemma implies that $\mathcal{I}_T$ is isomorphic to the product $\mathbb{R}$ bundle on the complement in the $|x| = 1$ sphere of the radius $c_0(\underline{N}_{(p)}(r)+r)^{1/2}$ tubular neighborhood of $Z_{*r}$. Because $N_{(p)}(0)$ is half an odd integer, this runs afoul of what is said in Lemma 9.5.

*Part 4*: This part proves Item b) of the second bullet of Lemma 9.8. To this end, let $\hat{u}$ denote a length 1 unit vector in $\mathbb{R}^d$ and t a number from $[\frac{2}{3}, 4]$ such that $t\hat{u}$ is in $Z_T$. Let $\Pi^\perp$ denote the 2-dimensional normal vector space to $\Pi_{3r}$. The plane in $\mathbb{R}^d$ orthogonal to $\Pi_{3r}$ through the point $t\hat{u}$ is $(1+\Delta^2)^{-1}t\hat{e} + \Pi^\perp$. There may be other points from $Z_T$ in this plane. Even so, the points from $Z_T$ in this plane must lie in a disk with radius D obeying

$$D \le e^{-1/(c_0 \underline{N}_{(p)}(r))} .$$
(9.26)

The bound in (9.26) is derived as follows: Let q denote $\phi_T(t\hat{u})$ which is again in Z. Apply Lemma 9.5 using for the input data the triple $T_1 = (q, L, r)$. Any version of $Z_*$ that is supplied by Lemma 9.5 using the input data $(q, L, r)$ must have angular distance less than $c_0 \varepsilon$ from any version of $Z_*$ that Lemma 9.5 assigns to $T = (p, L, r)$. Suppose for the sake of argument that there are two points from $Z_T$ in the plane $(1+\Delta^2)^{-1}t\hat{e} + \Pi^\perp$ with distance less than 1 from $t\hat{u}$. Applying Lemma 9.5 a second time with the input data set $(q, L, Dr)$ begets a version of $Z_*$ that is nearly orthogonal to p's version of $Z_*$. These two versions of $Z_*$ can differ by an $\mathcal{O}(1)$ angle only if D obeys (9.27) with $N_{(q)}(r)$ replacing $\underline{N}_{(p)}$. That this is so follows from (4.30). Meanwhile, $N_{(q)}(r)$ is no smaller than $\underline{N}_{(p)}(r)$, this being a tautological consequence of the definition of $\underline{N}_{(\cdot)}$ in (9.1).

### e) Proof of Lemma 9.9

This is a long proof with eight parts.

*Part 1*: Let W denote the part of $\mathbb{R}^d$ where $|x| \in [\frac{3}{4}, \frac{5}{4}]$ and where $|z| \in [\frac{1}{32}, 1]$. Write $v_T$ again as $v_+ + v_-$ with $v_+$ and $v_-$ being the respective parts of $v$ in the $+i$ and $-i$ eigenspaces of the endomorphism $\tau$ from Part 3 in the previous subsection. Use $v_{*r+}$ to denote the part of $v_{*r}$ in the $+i$ eigenspace of $\tau$. As can be seen from (9.21), $v_{*r+}$ is given by $(1+c)^{-1/2} a_{*+} z^{N_0}$. The construction of $v_{*r+}$ is such that

$$|v_+ \otimes v_+ - v_{*r+} \otimes v_{*r+}| \leq c_0 (c + |\flat_0| + |\flat_1| + |t - \underline{t}| + |o_+|)$$

(9.27)

on W. There is an analogous inequality for $|v_- \otimes v_- - v_{*r-} \otimes v_{*r-}|$. Use (9.14) with the bounds in Part 3 of the preceding subsection on the norms of $o_+$, $c$, $\flat_0$, $\flat_1$, and with the analogous bounds on the norms of their $v_-$ analogs to see that

$$\int_W |v_T \otimes v_T - v_{*r} \otimes v_{*r}|^2 \leq c_0 (\underline{N}_{(p)}(r) + r) .$$

(9.28)

Now suppose that $\varepsilon \in (0, 1]$ is such that r is less than the corresponding $r_{K,\varepsilon}$ from Lemma 9.5. Let $(Z_*, \mathcal{I}_*, v_*)$ be a data set that is supplied by the $T = (p, L, r)$ version of Lemma 9.5. It follows from (9.28) and from what is said in Lemma 9.5 about $v_*$ being close to $v_T$ that

$$\int_W |v_* \otimes v_* - v_{*r} \otimes v_{*r}|^2 \leq c_0 \varepsilon .$$

(9.29)

This last inequality leads directly to a $c_0 \varepsilon$ bound on the angular distance between the linear subspaces $Z_*$ and $Z_{*r}$. It follows as a consequence that if $\varepsilon < c_0^{-1}$, then $|v_*| > 0$ on W and moreover $|\langle v_*, v_{*r} \rangle| > (1 - \frac{1}{10,000})|v_*||v_{*r}|$ on W. This implies in turn that

$$2|\langle v_T, v_{*r} \rangle|^2 > (1 - \frac{1}{1000})(|v_T|^2 + |v_{*r}|^2)$$

(9.30)

on the set W. In particular, $\langle v_T, v_* \rangle$ defines an isomorphism between $\mathcal{I}_T$ and $\mathcal{I}_{*r}$ on W.

Lemma 9.8 says in part that $|v_T| > 0$ on the $|z| > c_0 (\underline{N}_{(p)}(r) + r)^{1/2}$ part of the spherical annulus where $|x| \in (\frac{2}{3}, 4)$. This being the case, it follows that the isomorphism just defined extends in a unique fashion to give an isomorphism between $\mathcal{I}_T$ and $\mathcal{I}_{*r}$ on this same $|z| > c_0(\underline{N}_{(p)}(r) + r)^{1/2}$ part of the $|x| \in (\frac{2}{3}, 4)$ spherical annulus. The latter isomorphism is used implicitly in the subsequent parts of the proof to view $v_T$ on this region as a section of $\mathbb{V}^+ \otimes \mathcal{I}_{*r}$.

*Part 2*: This part proves that the integral of $|v_T - v_{*r}|^2$ on the $|z| > \frac{1}{16}$ part of the $|x| = 1$ sphere is bounded by $c_0(\underline{N}_{(p)}(r) + r)$. To this end, let $v_+$ again denote the part of $v_T$ in the $+i$ eigenspace of $\tau$ and let $v_{*r+}$ denote the corresponding part of $v_{*r}$. Write $v_+$ on W as $v_+ = \sigma v_{*r+} + q_+$ with $\sigma$ being a $\mathbb{C}$-valued function and with $q_+$ being orthogonal to $v_{*r+}$. Then write

$$v_+ \otimes v_+ = \sigma^2 v_{*r+} \otimes v_{*r+} + \sigma v_{*r+} \otimes q_+ + \sigma q_+ \otimes v_{*r+} + q_+ \otimes q_+ . \tag{9.31}$$

The four terms on the right hand side of (9.31) are pairwise orthogonal in $\mathcal{V}$. Thus,

$$|v_+ \otimes v_+ - v_{*r+} \otimes v_{*r+}| \geq c_0^{-1}(|1 - \sigma^2| + |\sigma||q_+| + |q_+|^2) . \tag{9.32}$$

Write $\sigma$ as $1 + \epsilon$. It follows from (9.30) that $|\epsilon| < \frac{1}{100}$. This being so, (9.29) implies that

$$\int_W (|\epsilon|^2 + |q_+|^2) < c_0(\underline{N}_{(p)}(r) + r) . \tag{9.33}$$

The latter bound and its $v_-$ analog imply that

$$\int_W |v_T - v_{*r}|^2 \leq c_0(\underline{N}_{(p)}(r) + r) . \tag{9.34}$$

The bound in (9.34) requires there be $\rho \in [\frac{7}{8}, \frac{9}{8}]$ with the integral of $|v_T - v_{*r}|^2$ on the $|z| > \frac{1}{16}$ part of the $|x| = \rho$ sphere being less than $c_0(\underline{N}_{(p)}(r) + r)$. Granted this, then an appeal to (4.30) proves that the integral of $|v_T - v_{*r}|^2$ on the $|z| > \frac{1}{16}$ part of the $|x| = 1$ sphere is also bounded by $c_0(\underline{N}_{(p)}(r) + r)$. Part 3 derives an analogous bound for the integral of $|v_T - v_{*r}|^2$ over other sorts of d-1 dimensional surfaces.

*Part 3*: Let S denote a d-1 dimensional surface in the $|x| \leq \frac{3}{2}$ part of $\mathbb{R}^d$ with the projection map to the $(x_1, x_2)$ plane being a submersion on S. This projection map is again denoted by $\mathfrak{p}$. Let z again denote $x_1 + ix_2$; and assume that $|z| > \frac{1}{16}$ on S. Use $\upsilon$ and A to denote respective upper bounds for the norm of the inverse of $\mathfrak{p}$'s differential along S and for the d-3 dimensional measures of the fibers of $\mathfrak{p}$ in S. This part of the proof of Lemma 9.9 explains why

$$\int_S |v_T - v_{*r}|^2 \le c_0 \upsilon^2 A(\underline{N}_{(p)}(r)+r) \ .$$

(9.35)

Granted what was said previously about the norms of $c$, $h_0$ and $h_1$, granted what was said about the $L^2_1$ norm of $o_+$, the bound in (9.14) leads to the integral bounds

$$\int_S |v_+ \otimes v_+ - v_{*r+} \otimes v_{*r+}|^2 \le c_0 \upsilon^2 A(\underline{N}_{(p)}(r)+r) \ .$$

(9.36)

There are analogous integral bounds for $|v_- \otimes v_- - v_{*r-} \otimes v_{*r-}|$. Write $v_+$ as $(1+\varepsilon)v_{*r+} + q_+$ as done in Part 2. Keeping in mind that $|\varepsilon| \le \frac{1}{100}$, the bound in (9.36) implies the $L^2$ bounds

$$\int_S (|\varepsilon|^2 + |q_+|^2) < c_0 \upsilon^2 A(\underline{N}_{(p)}(r)+r)$$

(9.37)

This bound and its $v_-$ analog lead directly to the bound in (9.35).

*Part 4*: Let $S_0$ denote the part of the $|x| = 1$ sphere in $\mathbb{R}^d$ where $|z| \in [\frac{1}{4}, \frac{3}{4}]$. Use (9.35) with $S = S_0$ to see that

$$\int_{S_0} |v_T - v_{*r}|^2 \le c_0 (\underline{N}_{(p)}(r)+r) \ .$$

(9.38)

Fix $m > c_0$ with the lower bound chosen so that $|v_T|$ is positive on the $|z| \ge m(\underline{N}_{(p)}(r)+r)^{1/2}$ part of the spherical annulus $|x| \in [\frac{7}{8}, \frac{9}{8}]$. Lemma 9.8 guarantees that this will be the case for any $m > c_0$. Having chosen $m$, introduce M to denote the greatest integer such that $2^{-M-2} \ge m(\underline{N}_{(p)}(r)+r)^{1/2}$. For each $k \in \{1, 2, \ldots, M\}$, let $S_k$ denote the part of the $|x| = 1$ sphere in $\mathbb{R}^d$ where $|z| \in [\frac{1}{2^{k+2}}, \frac{3}{2^{k+2}}]$. The remaining parts of the proof explain why

$$\int_{S_k} |v_T - v_{*r}|^2 \le c_0 k\, 2^{-2(N_0+1)k} (\underline{N}_{(p)}(r)+r)$$

(9.39)

for each $k \in \{1, 2, \ldots M\}$ if $m > c_0$. Lemma 9.9 follows from these bounds, the bound in (9.38) and the bound asserted at the end of Part 2 for the integral of $|v_T - v_{*r}|^2$ on the $|z| \ge \frac{1}{2}$ part of the $|x| = 1$ sphere.

The proof of (9.39) exploits the fact that each point in $S_k$ has distance on the order of $2^{-k}$ from an $|x| = 1$ point in $Z_T$. In the $d = 3$ case, the whole of $S_k$ is contained in two balls of radius $2^{-k}$, each with center at a point in $Z_T$. In the $d = 4$ case, $S_k$ is in the union of $c_0 2^k$ balls of radius $2^{-k}$ centered at points in $Z_T$. In either case, these points in $Z_T$ have distance at most $c_0 (\underline{N}_{(p)}(r) + r)^{1/2}$ from some $|x| = 1$ point in $Z_{*r}$. That this is so follows from Lemma 9.8.

Suppose that $q_*$ is an $|x| = 1$ points from $Z_T$. If the dimension of X is three, let $B_*$ to denote the radius $2^{-k}$ ball centered at $q_*$. If the dimension of X is four, then define $B_*$ as follows: The coordinates on $\mathbb{R}^4$ are $(x_1, x_2, x_3, x_4)$ with $Z_*$ being the $(x_3, x_4)$ plane. Let t denote the angular $\mathbb{R}/2\pi\mathbb{Z}$ coordinate in the $(x_3, x_4)$ plane. Take $B_*$ to be the part of the radius $2^{-k}$ ball centered at $q_*$ where $|t - t(q_*)| \leq 2^{-k-2}$. The plan for what follows is to bound the contribution to left hand side of (9.39) from any given version of $B_*$ by

$$\int_{S_k \cap B_*} |\nu_T - \nu_{*r}|^2 \leq c_0 k \, 2^{-(2N_0 + d - 1)k} (\underline{N}_{(p)}(r) + r).$$

(9.40)

The bound in (9.40) implies the bound in (9.39) because $S_k$ is contained in the union of at most $c_0 2^{(d-3)k}$ versions of $B_*$.

By way of a look ahead, the bound in (9.40) is obtained by mapping $S_k$ to a surface in the $|x| \leq 1$ ball of the sort that appears in (9.35). For example, the translation and rescaling given by the rule $x \to 2^{-k}(x - q_*)$ maps $S_k$ to a surface in the $|x| \leq 1$ ball of the sort that appears in (9.35). This map is very close to the map that is used, the latter being the composition of first $\phi_T$ and then the inverse of the analog of $\phi_T$ that is defined by the data set $(\phi_T(q_*), L, 2^{-k}r)$. This composition and the map defined by the rule $x \to 2^{-k}(x - q_*)$ are identical when the metric is Euclidean near p.

The mapping of $S_k$ rewrites the integrand in (9.40) also. The affect of this rewriting is dealt with using (4.30). The appeal to (4.30) is the source of the factor of k that appears on the right hand side of (9.40). After controlling for the fact that the integrand is also rewritten, the bound in (9.40) is seen as a consequence of (9.35).

Note that the arguments for (9.40) have antecedents in Part 11 of the proof of Lemma 5.7.

*Part 5*: Let q denote $\phi_T(q_*)$ and let $T(k) = (q, L, 2^{-k}r)$. Use $\theta_k$ to denote the composition of first $\phi_{T(k)}$ and then the inverse of the mapping $\phi_T$. With $\theta_k$ understood, let $\mathcal{A}_k$ denote the inverse image via $\theta_k$ of the part of $S_k$ in the radius $2^{-k}$ ball centered at $q_*$. As noted previously, this is a surface that can be used as a version of $S$ in (9.36). Use $\theta_k$ to write the integral on the left hand side of (9.40) as

$$\int_{\mathcal{A}_k} \theta_k^*(|\nu_T - \nu_{*r}|^2) J_k$$

(9.41)

What is denoted by $J_k$ in (9.41) is the Jacobian that arises when comparing the induced measure on $\mathcal{A}_k$ with its pull-back by $\theta_k$. This Jacobian $J_k$ can be written as

$$J_k = (1 + \mathfrak{r}) 2^{-(d-1)k},$$

(9.42)

with $\mathfrak{r}$ obeying $|\mathfrak{r}| \leq c_0 r$. This term $\mathfrak{r}$ is zero when the metric near p is Euclidean.

The factor $\theta_k^*|\nu_T - \nu_{*r}|^2$ in the integrand of (9.41) is bounded by

$$(1 + c_0 r^2) |\theta_k^* \nu_T - \theta_k^* \nu_{*r}|^2 .$$

(9.43)

Keeping in mind that $\nu_{*r}$ is invariant under translations tangent to $Z_{*r}$ and that it is homogenous with respect to positive real number rescalings of the coordinate z, it follows from the fact that $q_*$ has distance at most $c_0(\underline{N}_{(p)}(r)+r)^{1/2}$ from $Z_{*r}$ that

$$|\theta_k^* \nu_{*r} - 2^{-N_0 k} \nu_{*r}| \leq c_0 (\underline{N}_{(p)}(r)+r)^{1/2} 2^{-N_0 k} .$$

(9.44)

The comparison between $\theta_k^* \nu_T$ and $\nu_T$ is not so straightforward. In any event, it follows from the definitions that

$$\theta_k^* \nu_T = \frac{K_{(q)}(2^{-k} r)}{K_{(p)}(r)} \nu_{T(k)} .$$

(9.45)

Let $\mu_k$ denote $K_{(q)}(r)/K_{(p)}(r)$ so as to write the ratio in (9.45) as $\frac{K_{(q)}(2^{-k} r)}{K_{(q)}(r)} \mu_K$. Since the first and third bullets of Lemma 3.2 imply that

$$(1 - c_0 k (\underline{N}_{(p)}(r)+r)) 2^{-N_0 k} \leq \frac{K_{(q)}(2^{-k} r)}{K_{(q)}(r)} \leq (1 + c_0 k r^2) 2^{-N_0 k} ,$$

(9.46)

the bounds in (9.42)-(9.44) lead to the inequality

$$\int_{\mathcal{A}_k} \theta_k^*(|\nu_T - \nu_{*r}|^2) J_k \leq 2^{-(2N_0 + d-1)k} \left( \int_{\mathcal{A}_k} |\mu_k \nu_{T(k)} - \nu_{*r}|^2 + c_0 (\underline{N}_{(p)}(r)+r) \right).$$

(9.47)

The next two parts of the proof analyzes the integral on the right hand side of (9.47).

*Part 6*: A bound for the integral that appears on the right hand side of (9.47) is derived using the sequence of inequalities that are presented in the five steps that follow.

Step 1: The first set of inequalities involves the T(k) version of $v_{*r}$. This is denoted by $v_{*r(k)}$ and it is defined by using T(k) in the previous subsection instead of T. The T(k) version of (9.36) is applied using $S = \mathcal{A}_k$ with the result being the bound

$$\int_{\mathcal{A}_k} |v_{T(k)} - v_{*r(k)}|^2 \leq c_0 (\underline{N}_{(p)}(r) + r).$$

(9.48)

The integral that appears on the right hand side of (9.47) can be bounded using the inequality in (9.48) by

$$\int_{\mathcal{A}_k} |\mu_k v_{T(k)} - v_{*r}|^2 \leq 2 \int_{\mathcal{A}_k} |\mu_k v_{*r(k)} - v_{*r}|^2 + c_0 (\underline{N}_{(p)}(r) + r).$$

(9.49)

A key point to keep in mind is that both $v_{*r(k)}$ and $v_{*r}$ are given explicitly by the respective T(k) and T versions of (9.21). In particular, a bound on $|\mu_k v_{*r(k)} - v_{*r}|$ on any reasonable open set in the $|x| \leq 1$ ball can be used to bound the integral that appears on the right hand side of (9.49).

Step 2: As it turns out, a good open set to use is the part of the $|x| = 1$ sphere where $|z| \in [\frac{1}{4}, \frac{3}{4}]$. This set was denoted previously by $S_0$. If it turns out that

$$\int_{S_0} |\mu_k v_{*r(k)} - v_{*r}|^2 \leq c_0 k (\underline{N}_{(p)}(r) + r),$$

(9.50)

then the T(k) and T version of (9.21) can be used to see that the right hand side of (9.49) is also bounded by $c_0 k (\underline{N}_{(p)}(r) + r)$. This with (9.47) leads directly to the bound in (9.40).

By way of a look ahead, the bound in (9.50) is derived from a sequence of comparisons on $S_0$: The first compares $v_{T(k)}$ and $v_{*r(k)}$, the second compares what is very nearly the $x \to x - q_*$ translate of $v_T$ with $v_{T(k)}$, and the third compares this translate of $v_T$ with $v_{*r}$.

Step 3: A second appeal to the T(k) version of (9.36) is used to start the derivation of (9.50). This new appeal to (9.35) takes $S = S_0$ with $S_0$ and it gives the bound

$$\int_{S_0} |v_{T(k)} - v_{*r(k)}|^2 \leq c_0 (\underline{N}_{(p)}(r) + r).$$

(9.51)

This is the promised comparison between $v_{T(k)}$ and $v_{*r(k)}$ on $S_0$.

<u>Step 4</u>: Let T´ denote the data set (q,L,r). The corresponding $v_{T'}$ is very nearly the $x \to x - q_*$ translate of $v_T$ because it is very nearly the pull-back of $v_T$ by the composition of first the map $\phi_{T'}$ and then the inverse of the map $\phi_T$. More is said about $v_T - v_{T'}$ in the subsequent steps. What follows directly compares $v_{T'}$ with $v_{T(k)}$ on $S_0$.

Being that q is in Z, it has a version of (4.30) which allows a comparison between $v_{T(k)}$ and $v_{T'}$ on the $|x| = 1$ sphere with the result being the bound

$$\int_{|x|=1} |v_{T(k)} - v_{T'}|^2 \leq c_0 k (\underline{N}_{(p)}(r) + r).$$

(9.52)

Since $S_0$ is an open set in the $|x| = 1$ sphere, the right hand side of (9.52) also bounds the integral of $|v_{T(k)} - v_{T'}|^2$ over $S_0$.

<u>Step 5</u>: Let $\theta$ denote the composition of first $\phi_{T'}$ and then the inverse of $\phi_T$. The definitions from the Section 4.1 of $v_T$ and $v_{T'}$ are such that

$$v_{T'} = \frac{K_{(p)}(r)}{K_{(q)}(r)} \theta^* v_T.$$

(9.53)

Since the ratio in (9.53) is $\mu_k^{-1}$, it follows that

$$\int_{S_0} |\mu_k v_{T'} - v_{*r}|^2 \leq c_0 \int_{S_0} |\theta^* v_T - v_{*r}|^2 + c_0 (\underline{N}_{(p)}(r) + r).$$

(9.54)

The integral on the on the left hand side of (9.54) can be written as

$$\int_{\theta(S_0)} |v_T - (\theta^{-1})^* v_{*r}|^2 J_\theta$$

(9.55)

with $J_\theta$ denoting the relevant Jacobian for the integration measure. This differs from 1 by at most $c_0 r$ because $\theta_0$ differs by $\mathcal{O}(r^2)$ from a translation. Keeping in mind that this translation is very nearly a translation by a unit vector along $Z_{*r}$, the formula in (9.21) can be used to bound $|\theta^{-1}v_{*r} - v_{*r}|$ by $c_0(\underline{N}_{(p)}(r)+r)^{1/2}$. This same fact about $\theta$ implies that there is an $S = \theta(S_0)$ version of (9.35) that bounds the integral of $|v_r - v_{*r}|^2$ on $\theta(S_0)$ by $c_0(\underline{N}_{(p)}(r)+r)$.

The preceding bounds lead to a $c_0(\underline{N}_{(p)}(r)+r)$ bound for (9.55) and thus a $c_0(\underline{N}_{(p)}(r)+r)$ for the left hand side of (9.54). The latter bound and those in (9.51) and (9.52) give the bound that is asserted by (9.50).

**f) Proof of Proposition 9.2**

The six parts of this subsection supply a proof of Proposition 9.2. The implicit assumption in what follows is that r is small enough to define the data set $(Z_{*r}, \mathcal{I}_{*r}, v_{*r})$ using the rules set forth in Section 9c. This will be so if $r < r_{*K}$ with $r_{*K}$ being a positive number that depends only on K.

*Part 1*: According to Lemma 9.8, all $|x| = 1$ points in $Z_{T(r)}$ have distance at most $c_0(\underline{N}_{(p)}(r)+r)^{1/2}$ from $Z_{*r}$. Lemma 9.9 describes an identification between $\mathcal{I}_{T(r)}$ and $\mathcal{I}_{*r}$ on this part of the $|x| = 1$ sphere. Fix this identification in what follows.

Suppose that $s \in (0, r]$ and assume for the moment that all $|x| = 1$ points in $Z_{T(s)}$ have distance at most $\frac{1}{1000}$ from $Z_{*r}$. It follows from Lemmas 3.3 and 4.8 that this is so if

$$s \geq r e^{-1/(c_0 \underline{N}_{(p)}(r))}.$$
(9.56)

It follows as a consequence that the identification between $\mathcal{I}_{T(r)}$ and $\mathcal{I}_{*r}$ on the part of the $|x| = 1$ sphere with distance greater than $\frac{1}{1000}$ from $Z_{*r}$ extends to give an identification between $\mathcal{I}_{T(s)}$ and $\mathcal{I}_{*r}$ on this same part of the $|x| = 1$ sphere if s obeys (9.56). Let $S_0$ denote the part of the $|x| = 1$ sphere where the coordinate z obeys $|z| \in [\frac{1}{4}, \frac{3}{4}]$. This part of the sphere has distance greater than $\frac{1}{1000}$ from $Z_{*r}$ so the identification between $\mathcal{I}_{T(s)}$ and $\mathcal{I}_{*r}$ is defined on $S_0$ if $s \in (0, r]$ obeys (9.56).

Supposing that $s \in (0, r]$ and that (9.56) holds, the identification given in the preceding paragraph between the bundles $\mathcal{I}_{T(s)}$ and $\mathcal{I}_{*r}$ extends under certain circumstances to a larger part of the $|x| = 1$ sphere. In particular, if $x > 1000$ is given and if all $|x| = 1$ points in $Z_{T(s)}$ have distance less than $x^{-1}$ from $Z_{*r}$, then the identification of line bundles extends to the part of the $|x| = 1$ sphere where the distance to $Z_{*r}$ is greater than $x^{-1}$.

Given that s obeys (9.56), these line bundle identifications are assumed henceforth without further comment. They are used implicitly to view $v_{T(s)}$ as a section of $\mathbb{V}^+ \otimes \mathcal{I}_{*r}$ over the appropriate part of the $|x| = 1$ sphere.

*Part 2*: Suppose that $s \in (0, r]$ obeys (9.56). As explained directly, the square of the distance from any $|x| = 1$ point in $Z_{T(s)}$ to $Z_{*r}$ is no greater than

$$c_0 \int_{S_0} |v_{T(s)} - v_{*r}|^2 + c_0 (\underline{N}_{(p)}(r) + r).$$

(9.57)

To prove this, invoke the T(s) version of (9.35) to obtain the bound

$$\int_{S_0} |v_{T(s)} - v_{*s}|^2 \leq c_0 (\underline{N}_{(p)}(r) + r),$$

(9.58)

and then use this bound with the triangle inequality to bound the integral of $|v_{*s} - v_{*r}|^2$ over $S_0$ by the expression in (9.57). Then use the depictions of $v_{*s}$ and $v_{*r}$ the respective T(s) and T(r) versions of (9.21) to infer that (9.57) bounds the square of the angular distance between $Z_{*r}$ and $Z_{*s}$. This means that the distance between any given $|x| = 1$ point in $Z_{T(s)}$ and $Z_{T(r)}$ is bounded by the number that is given in (9.57). Meanwhile, the T(s) version of Lemma 9.8 asserts in part that the distance between any $|x| = 1$ point in $Z_{T(s)}$ and $Z_{*s}$ is at most $c_0(\underline{N}_{(p)}(r) + r)^{1/2}$.

*Part 3*: Fix $m > 1$ and use $\chi$ to construct a non-negative function on the $|x| = 1$ sphere that equals one where the distance to $Z_{*r}$ is greater than $2m(\underline{N}_{(p)}(r) + r)^{1/2}$ and equals zero where the distance to $Z_{*r}$ is less than $m(\underline{N}_{(p)}(r) + r)^{1/2}$. This function can and should be constructed so that the norm of its derivative is at most $c_0 m^{-1}(\underline{N}_{(p)}(r) + r)^{-1/2}$. Denote this function by $\sigma$.

If $m > c_0$, then Lemma 9.8 can be invoked to see that all $|x| = 1$ points in $Z_{T(r)}$ have distance at most $m^{1/2} (\underline{N}_{(p)}(r) + r)^{1/2}$ from $Z_{*r}$. Such a choice guarantees that $|v_{T(r)}|$ is positive on the support of $\sigma$. If $s \in (0, r)$ and if $s$ is sufficiently close to $r$, then $|v_{T(s)}|$ will also be positive on the support of $\sigma$. This is because the function on the product of $(0, r]$ and the $|x| = 1$ sphere given at $s \in (0,r]$ by $|v_{T(s)}|$ is continuous.

The fact that this function is continuous has the following related consequence: Supposing that $m > c_0$, then $|v_{T(s)}| > 0$ on the support of $\sigma$ for all points s in an interval of the form $[s_*, r]$ for $s_* < r$ provided that

$$\int_{|x|=1} \sigma |v_{T(s)} - v_{*r}|^2 < m^{1/2} (\underline{N}_{(p)}(r) + r)$$

(9.59)

for all points $s \in [s_*, r]$.

*Part 4*: Supposing that $s \in [s_*, r]$, let $\mathfrak{m}_{T(s)}$ denote as in Section 4 the product of $s^{-2}$ times the $\phi_{T(s)}$ pull-back of the metric from X. This metric differs from the Euclidean metric by a term whose norm is bounded by $c_0 s^2$ and whose first and second derivatives are likewise bounded by $c_0 s^2$. The inner product on $\mathbb{V}^+$ that is defined by $\mathfrak{m}_{T(s)}$ is denoted in what follows by $\langle \cdot, \cdot \rangle_s$, the norm is denoted by $|\cdot|_s$ and the integration measure on the $|x| = 1$ sphere defined by $\mathfrak{m}_{T(s)}$ is denoted by $\omega_s$.

If $m > c_0$, then what was said in the Part 3 about the significance of (9.59) is true whether the norm and integration measure on the $|x| = 1$ sphere are defined by the Euclidean metric on the $|x| \leq 1$ ball or the metric $\mathfrak{m}_{T(s)}$. This is to say that the function $|v_{T(s)}|$ is positive on the support of $\sigma$ if $s \in [s_*, r]$ provided that

$$\int_{|x|=1} \sigma |v_{T(s)} - v_{*r}|_s^2 \omega_s \leq m^{1/2} (\underline{N}_{(p)}(r) + r)$$

(9.60)

for all $s \in [s_*, r]$. To prove that this is so, note first that it gives a $c_0 m^{1/2} (\underline{N}_{(p)}(r) + r)$ upper bound for the integral of $|v_{*s} - v_{*r}|^2$ on the $|z| \geq \frac{1}{100}$ part of the $|x| = 1$ sphere. The latter integral bound leads in turn to a $c_0 m^{1/4} (\underline{N}_{(p)}(r) + r)^{1/2}$ bound on the angular distance between $Z_{*s}$ and $Z_{*r}$. Meanwhile, the version of Lemma 9.8 with s replacing r asserts in part that all $|x| = 1$ points in $Z_{T(s)}$ have distance at most $c_0(\underline{N}_{(p)}(s) + s)$ from $Z_{*s}$.

*Part 5*: Suppose that $s_* < r$ and that (9.60) holds for all $s \in [s_*, r]$ so that $|v_{T(s)}| > 0$ on the support of $\sigma$. Use Q to denote the function on $[s_*, r]$ that is given by the rule

$$s \to Q(s) = \int_{|x|=1} \sigma \langle v_{T(s)}, v_{*r} \rangle_s \omega_s \ .$$

(9.61)

Looking at (9.60) finds that

$$Q(s) = 1 - \tfrac{1}{2} \int_{|x|=1} \sigma |v_{T(s)} - v_{*r}|_s^2 \omega_s - q$$

(9.62)

with q obeying $|q| \leq c_0 m (\underline{N}_{(p)}(r) + r)^{1+N_0}$. This is because the integral of $|v_{T(s)}|_s^2 \omega_s$ on the $|x| = 1$ sphere is 1, as is the integral of $|v_{*r}|^2$ with the standard integration measure is 1. Meanwhile, the integrals of both on the support of $1 - \sigma$ are les than $m (\underline{N}_{(p)}(r) + r)^{1+N_0}$. On the other hand, these same bounds for the integrals of $|v_{T(s)}|_s^2 \omega_s$ and $|v_{*r}|^2$ imply that

$$Q(s) \leq 1 + c_0(\underline{N}_{(p)}(r) + r) \ .$$

(9.63)

Given the bound in (9.63), there is little room for $Q(\cdot)$ to increase as s decreases in $[s_*, r]$.

*Part 6*: The derivative of $Q$ has the schematic form

$$s \frac{d}{ds} Q = \int_{|x|=1} \sigma \langle \tfrac{\partial}{\partial r} v_{T(s)}, v_{*r} \rangle_s \omega_s - N_{(p)}(s) Q + \mathfrak{e}_0 s,$$

(9.64)

with $\mathfrak{e}_0$ being a function of s with norm bounded by $c_0$. This function $\mathfrak{e}_0$ accounts for the fact that the inner product $\langle\,,\,\rangle_s$ and integration measure $\omega_s$ are defined using the metric $\mathfrak{m}_{T(s)}$ and not the Euclidean metric.

The identity in (9.64) is exploited by first using the equation $\mathcal{D} v = 0$ to write $\tfrac{\partial}{\partial r} v_{T(s)}$ as $-\tfrac{1}{s} \mathfrak{S}_r^\dagger \mathcal{D}_s^\perp v_{T(s)}$ with $\mathcal{D}_s^\perp$ involving only derivatives tangent to the $|x|=1$ sphere. The operator $\mathfrak{S}_r^\dagger \mathcal{D}_s^\perp$ is self-adjoint with respect to the $L^2$ inner product on the $|x|=1$ sphere as defined by the metric $\mathfrak{m}_{T(s)}$. (The metric $\mathfrak{m}_{T(s)}$ is used for precisely this reason). Use this self-adjoint property of $\mathcal{D}_s^\perp$ to write the equation in (9.64) as

$$s \frac{d}{ds} Q = - \int_{|x|=1} \sigma \langle v_{T(s)}, \mathfrak{S}_r^\dagger \mathcal{D}_s^\perp v_{*r} \rangle_s \omega_s + \int_{|x|=1} \langle v_{T(s)}, \wp_s(\nabla \sigma) v_{*r} \rangle_s \omega_s - N_{(p)}(s) Q + \mathfrak{e}_0 s,$$

(9.65)

with $\wp_s$ denoting the symbol homomorphism for the operator $\mathfrak{S}_r^\dagger \mathcal{D}_s^\perp$.

An important point to be made now is that the operator $\mathfrak{S}_r^\dagger \mathcal{D}_s^\perp$ can be written as $\mathfrak{S}_r^\dagger \mathcal{D}_0^\perp + \mathfrak{r}_1(\nabla^\perp) + \mathfrak{r}_0$ with $\mathcal{D}_0^\perp$ being the operator that appears in (5.3), with $\mathfrak{r}_1$ being an endomorphism with norm bounded by $c_0 s^2$ and $\mathfrak{r}_0$ being an endomorphism with norm bounded by $c_0 s$. Since $\mathfrak{S}_r^\dagger \mathcal{D}_0^\perp v_{*r} = -N_0 v_{*r}$, the equation in (9.65) can be written as

$$s \frac{d}{ds} Q = -(N_{(p)}(s) - N_0) Q + \mathfrak{w}(s) + \mathfrak{e}_1,$$

(9.66)

with $\mathfrak{e}_1$ being a function with norm bounded by $c_0 s$ and with $\mathfrak{w}$ being short hand for

$$\mathfrak{w}(s) = \int_{|x|=1} \langle v_{T(s)}, \wp_s(\nabla \sigma) v_{*r} \rangle_s \omega_s .$$

(9.67)

As it turns out the norm of $\mathfrak{w}$ obeys

$$|\mathfrak{w}| \leq c_0 |\ln(\underline{N}_{(p)}(r) + r)|\, (m(\underline{N}_{(p)}(r) + r))^{1+N_0} .$$

(9.68)

Parts 9 and 10 explain why this is so.

*Part 7*: Assume in what follows that $m < (\underline{N}_{(p)}(r)+r)^{-1/100}$. Granted this, the bound in (9.68) and (9.60) and (9.62) imply that Q on $[s_*,r]$ obeys the differential inequality

$$s \frac{d}{ds} Q \leq -(N_{(p)}(s) - N_0) + c_0 (\underline{N}_{(p)}(r)+r)^{5/4} + \mathfrak{e}_1 s \,.$$

(9.69)

Suppose that $x$ is less than $(\underline{N}_{(p)}(r)+r)^{-1/8}$. Given this constraint, if $x$ is none-the-less greater than $c_0$ and if $s_\diamond \in [s_*,r]$ is such that

$$N_{(p)}(s) - N_0 > x^{-1}(\underline{N}_{(p)}(r)+r) + xs$$

(9.70)

on $[s_\diamond, r]$, then (9.69) with (9.68) lead to the bound

$$s \frac{d}{ds} Q \leq -\tfrac{1}{2} x^{-1}(\underline{N}_{(p)}(r)+r)$$

(9.71)

on this interval $[s_\diamond, r]$. It follows as a consequence of (9.62) that (9.60) is not near to being violated on $[s, r]$ because Q is increasing as s decreases. This is to say that if $s_*$ is the smallest s such that (9.60) holds on $[s_*, r]$, and if (9.71) holds on $[s_\diamond, r]$, then $s_\diamond > s_*$. Keeping in mind (9.63) and (9.60), if (9.71) holds on $[s_\diamond,r]$, then

$$\ln(\tfrac{r}{s_\diamond}) \leq c_0 m^{1/2} x.$$

(9.72)

This can be said differently as follows: Fix $m \in [c_0, c_0^2]$ so that the all of preceding constructions can be made. The number $x$ in what follows will be taken to large enough so that (9.71) is a consequence of (9.70), but in any event less than $c_0$. By way of notation, let $\lambda_x = e^{c_0 \sqrt{m}\, x}$. Assume that $r < c_0^{-1} r_{*K}$ with $c_0$ chosen so that the preceding constructions can be made with the chosen value of $m$. Then

$$N_{(p)}(r) - N_0 \leq \tfrac{1}{x} \underline{N}_{(p)}(\lambda_x r) + c_0 r$$

(9.73)

Indeed, (9.73) follows from (9.72) by changing the notation by the rule whereby $s_\diamond \to r$ and $r \to \lambda_x r$. The fact that function $\underline{N}_{(p)}(\cdot)$ is increasing is also used.

*Part 8*: With (9.73) in hand, there is one more thing needed for the proof of Proposition 9.2, this being the formula below in (9.74) bounding $\underline{N}_{(p)}(\cdot)$ in terms of $N_{(p)}(\cdot)$. The number R that appears in the formula is assumed to be larger than 100. The assumption on r is that it be from $(0, R^{-1} r_0]$.

$$\underline{N}_{(p)}(r) \le (N_{(p)}(R\,r) - N_0) + c_0 \tfrac{1}{R} + c_0 R^2 r^2$$

(9.74)

Part 11 of the proof derives this formula.

Use (9.74) with r replaced by $\lambda_x r$ and with $R = R_1$ to go from (9.73) to

$$N_{(p)}(r) - N_0 \le \tfrac{1}{x}(N_{(p)}(\lambda_x R_1 r) - N_0) + c_0 \tfrac{1}{xR_1} + c_0(r + \tfrac{1}{x}\lambda_x^2 R_1^2 r^2).$$

(9.75)

This can be iterated now to rewrite $N_{(p)}(\lambda_x R_1 r)$ on the right hand side of (9.75) provided that $\lambda_x R_1 r < r_{*K}$. Using a different version of R, a number $R_2 > 100$, this iteration leads to the bound

$$N_{(p)}(r) - N_0 \le \tfrac{1}{x^2}(N_{(p)}(\lambda_x^2 R_2 R_1 r) - N_0) + c_0(\tfrac{1}{x^2 R_2} + \tfrac{1}{xR_1})$$
$$+ c_0((1 + \tfrac{1}{x}\lambda_x R_1)r + (\tfrac{1}{x}\lambda_x^2 R_1^2 + \tfrac{1}{x^2}\lambda_x^4 R_2^2 R_1^2)r^2).$$

(9.76)

Additional iterations can be made along the following lines: Suppose that L is an integer greater than 1 and that $\{R_1, \ldots, R_L\} \subset (100, \infty)$ is a sequence that obeys

$$\lambda_x^L R_L \cdots R_1 r \le r_{*K}.$$

(9.77)

Then the L iterations of (9.75) can be made resulting in the bound

$$N_{(p)}(r) - N_0 \le \tfrac{1}{x^L}(N_{(p)}(\lambda_x^L R_L \cdots R_1 r) + c_0(\tfrac{1}{x^L R_L} + \cdots + \tfrac{1}{xR_1}) + c_0(w_1 r + w_2 r^2)$$

(9.78)

with $w_1$ and $w_2$ being the respective sums

$$w_1 = 1 + \tfrac{1}{x}\lambda_x R_1 + \cdots + \tfrac{1}{x^L}(\lambda_x^L R_L \cdots R_1) \quad and \quad w_2 = \tfrac{1}{x}\lambda_x^2 R_1^2 + \cdots + \tfrac{1}{x^L}(\lambda_x^{2L} R_L^2 \cdots R_1^2).$$

(9.79)

To go from theory to practice, take $R_k$ for $k \in \{1, 2, \ldots, L\}$ to be $x^{L-k+1}$. The constraint in (9.77) is obeyed if

$$c_0 xL + \tfrac{1}{2}(L^2 + L) \le \tfrac{1}{\ln x} |\ln(r/r_{*K})| - c_0.$$

(9.80)

In particular, this is guaranteed if L is the greatest integer less than

$$10^{-4} (\tfrac{1}{\ln x} |\ln(r/r_{*K})|)^{1/2} - c_0 x^2$$

(9.81)

assuming r is such that (9.81) is greater than 1. With L chosen in this way, the sums $w_1$ and $w_2$ in (9.79) are bounded respectively by

$$c_0 \, e^{L^2 \ln x + c_0 L x} \quad and \quad c_0 \, e^{2L^2 \ln x + c_0 L x} \,, \tag{9.82}$$

which in turn are bounded respectively by $c_0 r^{-1/4}$ and $c_0 r^{-1/2}$ given that $x < c_0$. These bounds can be used in (9.78) to see that $N_{(p)}(r) - N_0 \leq c_0 \, e^{-\sqrt{|\ln r|}/c_0}$. Apply this bound with r replaced by 16r and with p replaced by any given point in the ball $B_{8r}$ to obtain the bound for $\underline{N}_{(p)}$ that is asserted by Proposition 9.2.

*Part 9*: To bound $\mathfrak{w}$ as indicated in (9.68), use the recipe in Section 9c to define $\nu_{*s}$ using $T(s)$ in lieu of $T(r)$. With $\nu_{*s}$ in hand, write $\mathfrak{w}(s)$ as the sum

$$\mathfrak{w}(s) = \int_{|x|=1} \langle (\nu_{T(s)} - \nu_{*s}), \wp_s(\nabla\sigma)\nu_{*r} \rangle_s + \int_{|x|=1} \langle \nu_{*s}, \wp_s(\nabla\sigma)\nu_{*r} \rangle_s \,. \tag{9.83}$$

The absolute value of the left most integral in (9.83) is no greater than

$$c_0 \, (m(\underline{N}_{(p)}(r) + r))^{N_0/2} \, \left( \int_{S_*} |\nu_{T(s)} - \nu_{*s}|^2 \right)^{1/2} \tag{9.84}$$

with $S_*$ denoting the set in the $|x| = 1$ sphere where the distance to $Z_{*r}$ is greater than $m(\underline{N}_{(p)}(r) + r)^{1/2}$ but less than $2m(\underline{N}_{(p)}(r) + r)^{1/2}$. The point being that this set contains the support of $\nabla\sigma$. The $T(s)$ version of (9.39) can be invoked to see that

$$\int_{S_*} |\nu_{T(s)} - \nu_{*s}|^2 \leq c_0 \, |\ln(\underline{N}_{(p)}(r) + r)| \, (m(\underline{N}_{(p)}(r) + r))^{2+N_0} \,. \tag{9.85}$$

Use this bound in (9.84) to bound the left most integral in (9.84) by the expression on the right hand side of (9.68).

*Part 10*: The task for this part of the proof is to supply a useful bound for the right most integral in (9.83). The arguments here are much like those in Part 10 of the proof of Lemma 5.7. The first point to make in this regard is that the endomorphism $\wp_s(\nabla\sigma)$ is antisymmetric with respect to the inner product $\langle \cdot, \cdot \rangle_s$; and as a consequence the right most integral in (9.83) is equal to

$$\int_{|x|=1} \langle v_{*s} - v_{*r}, \wp_s(\nabla\sigma) v_{*r} \rangle_s .$$

(9.86)

The second point to note is that (9.61) and (9.63) lead to the bound

$$\int_{S_0} |v_{*s} - v_{*r}|^2 \leq c_0 m^{1/2} (\underline{N}_{(p)}(r) + r) ;$$

(9.87)

and this in turn has the following implication: Let $v_{*r+}$ and $v_{*s+}$ denote the respective parts of $v_{*r}$ and $v_{*s}$ in the $+i$ eigenspace of the endomorphism $\tau$. Write $v_{*r+}$ as in (9.21) and then $v_{*s+}$ can be written as

$$v_{*s+} = v_{*r} + (1+c)^{-1/2}(a_{*+} z_s^{N_0} + a_{*-} \overline{z}_s^{N_0}) + (1+c)^{-1/2}(\omega_+ a_{*+} z^{N_0} + \omega_+ a_{*-} \overline{z}^{N_0}) + \mathfrak{z} ,$$

(9.88)

with the notation as follows: What are denoted by $\omega_+$ and $\omega_-$ indicate endomorphisms with norm bounded by $c_0 m^{1/2} (\underline{N}_{(p)}(r) + r)^{1/2}$. What is denoted by $z_s$ is a rotated version of the complex function $z$ to account for the fact that $Z_{*s}$ need not equal $Z_{*r}$. In any event, $|z_s - z| \leq c_0 m^{1/2} (\underline{N}_{(p)}(r) + r)^{1/2}$. Meanwhile, $\mathfrak{z}$ is an element in $\mathbb{V}^+ \otimes \mathcal{I}_{*r}$ obeying

$$|\mathfrak{z}| \leq c_0 m (\underline{N}_{(p)}(r) + r) |z|^{N_0 - 1}$$

(9.89)

on the support of $\nabla\sigma$. Since $|z|$ is between $m(\underline{N}_{(p)}(r) + r)^{1/2}$ and $2m(\underline{N}_{(p)}(r) + r)^{1/2}$ on the support of $\nabla\sigma$, this is less than $(m(\underline{N}_{(p)}(r) + r))^{(1+N_0)/2}$.

These bounds for $|\mathfrak{z}|$ and the bounds for $|\omega_+|$ and $|\omega_-|$ can be used to see that the contributions to (9.86) from $\mathfrak{z}$ and $(1+c)^{-1/2}(\omega_+ a_{*+} z^{N_0} + \omega_+ a_{*-} \overline{z}^{N_0})$ are no greater than

$$c_0 (m(\underline{N}_{(p)}(r) + r))^{1+N_0} .$$

(9.90)

Keep in mind in this regard that $|\nabla\sigma| \leq c_0 (m(\underline{N}_{(p)}(r) + r))^{-1/2}$ and that the volume of the support of $\nabla\sigma$ is bounded by $c_0 m (\underline{N}_{(p)}(r) + r)$.

The left most term on the right hand side of (9.88) is $v_{*r}$ and this term makes no contribution to (9.86) because the endomorphism $\wp_s(\nabla\sigma)$ is antisymmetric with respect to the inner product $\langle , \rangle_s$. This antisymmetry also implies that the contribution to (9.86) from $(1+c)^{-1/2}(a_{*+} z'^{N_0} + a_{*-} \overline{z}'^{N_0})$ is given by the integral of

$$(\langle a_{*+}, \wp_s(\nabla\sigma) a_{*+}\rangle_s - \langle a_{*-}, \wp_s(\nabla\sigma) a_{*-}\rangle_s)(\overline{z}'^{N_0} z^{N_0} - \overline{z}^{N_0} z'^{N_0}) .$$

(9.91)

The contribution of this to (9.86) is bounded by (9.90) also. This is so because the endomorphism $\wp_s(\nabla\sigma)$ almost anti-commutes with the endomorphism $\tau$ and $a_{*+}$ and $a_{*-}$ are eigenvectors for $\tau$. In particular, the norm of $\tau\wp_s(\nabla\sigma) + \wp_s(\nabla\sigma)\tau$ is bounded by $c_0$ because the orthogonal projection of $\nabla\sigma$ tangent to $Z_{*r}$ has norm bounded by $c_0$, and because $\wp_s$ differs from its Euclidean metric analog by at most $c_0 r^2$.

*Part 11*: This part of the proof derives the inequality that is depicted in (9.74). Fix $T > 80$. If $q \in B_{8r}$, then the ball of radius $Tr$ centered at $q$ is contained in the ball of radius $(T+8)r$ centered at p. This implies that

$$N_{(q)}(Tr) \leq (1 + \tfrac{8}{T})^2 \left(\frac{K_{(p)}((T+8)r)}{K_{(q)}(Tr)}\right)^2 N_{(p)}((T+8)r).$$

(9.92)

The next task is to compare the two versions of $K_{(\cdot)}$ that appear in (9.92). To do this, use the fact that $K_{(q)}$ is an increasing function to see that

$$\int_{(T-17)r}^{Tr} \tfrac{1}{s} K_{(q)}(s)\, ds \leq |\ln(1 - \tfrac{17}{T})| K_{(q)}(Tr).$$

(9.93)

Let $A_q$ denote the spherical annulus centered at q with inner radius $(T-17)r$ and with outer radius $Tr$. The integral on the left hand side of (9.93) can be written as

$$\int_{A_q} \tfrac{1}{\text{dist}(\cdot,q)^d} |\nu|^2 \ .$$

(9.94)

Let $A_p$ denote the spherical annulus centered at p with inner radius $(T-9)r$ and with outer radius $(T-8)r$. Since the annulus $A_p$ is inside the annulus $A_q$, the integral in (9.99) is no smaller the integral over $A_p$ of the same integrand. Since $\text{dist}(\cdot,q) \leq \text{dist}(\cdot,p) + 8$ and since $T \geq 80$, the integral in (9.94) is no smaller than

$$(1 - \tfrac{10}{T}) \int_{A_q} \tfrac{1}{\text{dist}(\cdot,p)^d} |\nu|^2 \ ,$$

(9.95)

which in turn is no smaller than $(1 - \tfrac{10}{T})|\ln(1 - \tfrac{8}{T})| K_{(p)}((T-9)r)$. These inequalities starting from (9.93) give the lower bound

$$(1 + \tfrac{c_0}{T})^{-1} K_{(p)}((T-9)r) \leq K_{(q)}(Tr) \ .$$

(9.96)

Since $N_{(p)}(\cdot)$ is bounded, appeal to p's version of (3.4) leads from (9.96) to

$$(1+ \tfrac{c_0}{T} )^{-1} K_{(p)}((T+8)r) \le K_{(q)}(Rr)$$

(9.97)

Use this bound in (9.92) to see that

$$N_{(q)}(Tr) \le (1+ \tfrac{c_0}{T} ) N_{(p)}((T+8)r) \ .$$

(9.98)

Now invoke q's version of the third bullet of Lemma 3.4 to deduce that

$$\underline{N}_{(p)}(r) \le (N_{(p)}((T+8)r) - N_0) + \tfrac{c_0}{T} N_0 + c_0 (T+8)^2 r^2 \ .$$

(9.99)

Use R to denote $T+8$ to obtain the bound that is asserted by (9.72).

## 10. The structure of Z at small length scales

Sections 10a and 10b characterize certain neighborhoods in Z of the points in Z where the function $N_{(\cdot)}(0)$ is continuous. These characterizations are summarized in Proposition 10.1 and Lemmas 10.2 and 10.3. Section 10c says somethings about the points in Z where $N_{(\cdot)}(0)$ is not continuous. Section 10d talks about the points of discontinuity for the bundle $\mathcal{I}$. Theorem 1.4 is proved in Section 10b as a consequence of Propositions 10.1 and Lemma 10.3. Meanwhile, Proposition 10.5 in Section 10d restates Theorem 1.2.

### a) Neighborhoods of the strongly continuous points in Z

A codimension 2, Lipshitz graph through the origin in $\mathbb{R}^d$ is determined by an orthogonal, product decomposition of $\mathbb{R}^d$ as $\mathbb{R}^2 \times \mathbb{R}^{d-2}$, an open ball $D \subset \mathbb{R}^{d-2}$ centered at the origin and a Lipshitz map $\varphi \colon D \to \mathbb{R}^2$ with $\varphi(0) = 0$. The map $\varphi$ is Lipshitz if there exists $\kappa > 0$ such that $|\varphi(x) - \varphi(y)| \le \kappa |x-y|$ for all $x, y \in D$. The corresponding Lipshitz graph is the graph of $\varphi$ with respect to the given splitting of $\mathbb{R}^d$ as $\mathbb{R}^2 \times \mathbb{R}^{d-2}$, the graph being the set of pairs $\{(\varphi(z), z)\}_{z \in D}$. Let p denote a point in Z. A neighborhood in Z of p is said to be contained in a *codimension 2 Lipshitz graph* when a Gaussian coordinate chart centered at p identifies the neighborhood as a subset of a codimension 2 Lipshitz graph through the origin in $\mathbb{R}^{d-}$. A neighborhood in Z of p is said to be a *codimension 2 differentiable graph* if it is the image by a Gaussian coordinate chart centered at p of the graph of a $C^1$ map from D to $\mathbb{R}^2$ taking 0 to 0.

**Proposition 10.1**: *Let* $p \in Z$ *denote a strongly continuous point*.
- *A neighborhood of p in Z is contained in codimension 2 Lipshitz graph*.

- *If a neighborhood of p in a codimension 2 Lipshitz graph is contained in Z, then a neighborhood of p in Z is contained in a codimension 2 differentiable graph.*
- *If $N_{(p)}(0)$ is half of an odd integer, then a neighborhood of p in Z is a codimension 2 differentiable graph.*

*Proof of Proposition 10.1*: The assertion of the proposition is automatic in the case when p is isolated, so assume that this is not the case. Fix an oriented, orthonormal frame for $TX|_p$, this denoted by L. With L fixed, either Proposition 9.2 or Proposition 9.1 as the case may be assigns to p a unique data set as described in (5.1), this denoted by $(Z_\diamond, \mathcal{I}_\diamond, \nu_\diamond)$, with the following property: Any limit data $(Z_*, \mathcal{I}_*, \nu_*)$ set from the (p, L) version of Proposition 4.1 has $Z_* = Z_\diamond$ and it has $\mathcal{I}_*$ being isomorphic to $\mathcal{I}_\diamond$ by an isomorphism that identifies $\nu_*$ and $\nu_\diamond$. This set $Z_\diamond$ is always a line through the origin if the X has dimension 3, and it is either a line or or a plane through the origin if the dimension of X is 4. It is always a plane if the dimension of X is 4 and $N_{(p)}(0)$ is half an odd integer.

The following lemma strengthens the assertion of the proposition's first bullet in the case when the X has dimension 4 and $Z_\diamond$ is a line through the origin.

**Lemma 10.2**: *Suppose that $\dim(X) = 4$ and that $p \in Z$ is a strongly continuous point with the corresponding set $Z_\diamond$ being a line through the origin. There is a Gaussian coordinate chart centered at p having the following property: With $\mathbb{R}^4$ given the orthogonal splitting as $\mathbb{R} \times \mathbb{R}^3$, the Gaussian coordinate chart map identifies a neighborhood of Z near p with the graph of a Lipshitz map from an interval containing 0 in $\mathbb{R}$ to $\mathbb{R}^3$ sending 0 to 0.*

Parts 1-3 of what follows prove the first bullet of the proposition with the assumption in the dimension 4 case that $Z_\diamond$ is a plane through the origin. But for notation, these same arguments prove Lemma 10.2 when $Z_\diamond$ is a line through the origin. Part 4 of what follows proves the second bullet of the proposition. The third bullet follows from the second bullet and Lemma 9.3.

*Part 1*: Fix $r_p \in (0, c_0^{-1}]$ so that the ball of radius $r_p$ centered at p is well inside a Gaussian coordinate chart for a neighborhood of p in X. Use the chosen frame L to define a Gaussian coordinate system centered at p and let $\eta: B_{r_p} \to \mathbb{R}^d$ denote the inverse to the Gaussian coordinate chart map. Given $\delta \in (0,1]$ and $r \in (0, \frac{1}{100} r_p]$, define $Z_{\delta,r}$ to be the image via $\eta$ of the part of $Z \cap B_r$ where $N_{(\cdot)}(0) > N_{(p)}(0) - \delta$. Write $\mathbb{R}^d$ as an orthogonal product $\mathbb{R}^2 \times \mathbb{R}^{d-2}$ with $\{0\} \times \mathbb{R}^{d-2}$ being the linear space $Z_\diamond$. Use $\mathfrak{p}^\perp$ to denote the projection map $\mathbb{R}^4 \to \mathbb{R}^{d-2}$ onto the right most factor in the splitting $\mathbb{R}^2 \times \mathbb{R}^{d-2}$. Writing any given $x \in \mathbb{R}^4$ with respect to the splitting $\mathbb{R}^2 \times \mathbb{R}^{d-2}$ as $x = (z,w)$, then w is $\mathfrak{p}^\perp(x)$.

*There exists $\delta \in (0, 1]$ such that the restriction of $\mathfrak{p}^\perp$ to*
*a neighborhood in $Z_{\delta,r}$ of 0 defines a 1-1 map.*

(10.1)

The proof of (10.1) assumes it to be false so as to generate nonsense. If (10.1) is false, then there exists a sequence $\{j_n\}_{n \in \{1,2,...\}}$ whose n'th member for $n \in \{1, 2, ...\}$ is a pair $(x_n, y_n)$ with $x_n$ and $y_n$ being elements in the $\delta = 1/n$ version of $Z_{\delta,r}$ that are mapped to the same point by $\mathfrak{p}$. In addition, both $\{|x_n|\}_{n \in \{1,2,...\}}$ and $\{|y_n|\}_{n \in \{1,2,...\}}$ are decreasing sequences with limit 0 as n gets ever larger. To obtain nonsense from the sequence $\{j_n\}_{n \in \{1,2,...\}}$, invoke Proposition 4.1 using for input the sequences $\{p_n = \eta^{-1}(x_n)\}_{n \in \{1,2,...\}}$ and $\{\lambda_n = |x_n - y_n|\}_{n \in \{1,2,...\}}$. Let $(Z_*, \mathcal{I}_*, \nu_*)$ denote a limit data set that is then supplied by Proposition 4.1 and let $\Lambda$ denote the corresponding subsequence of $\{1, 2, ...\}$. The manner of convergence that is described by Proposition 4.1 implies that the sequence of points $\{q_n = \lambda_n^{-1}(x_n - y_n)\}_{n \in \Lambda}$ contains a subsequence that converges to a point on the $|x| = 1$ sphere in $\mathbb{R}^d$. Denote this limit point by $q_*$. Since $\mathfrak{p}^\perp(x_n) = \mathfrak{p}^\perp(y_n)$ holds for every index $n \in \{1,2,...\}$, the point $q_*$ is not in $\{0\} \times \mathbb{R}^{d-2}$, and since the latter plane is $Z_\diamond$, the point $q_*$ is not in $Z_\diamond$. This last conclusion is nonsense for the following reason: Proposition 4.1 requires that $\nu_*|_{q_*} = 0$ because $\eta^{-1}(y_n)$ is in $Z$ for each $n \in \{1, 2, ...\}$. This means that $q_*$ is in $Z_*$. But as noted previously, $Z_* = Z_\diamond$.

*Part 2*: Given what is said in (10.1), there exists $\delta \in (0, 1]$ and $r \in (0, r_\mathfrak{p})$ such that the restriction of $\mathfrak{p}^\perp$ to $Z_{\delta,r}$ is 1-1. The map $\mathfrak{p}^\perp$ on $Z_{\delta,r}$ has the following additional property:

*Given $\varepsilon \in (0, 1)$, there exists $\delta \in (0,1)$ and a disk $D_\varepsilon \subset \mathbb{R}^{d-2}$ about the origin such that if x and y are points in $(\mathfrak{p}^\perp)^{-1}(D_\varepsilon) \cap Z_{\delta,r}$, then $|x - y| \leq (1+\varepsilon) |\mathfrak{p}^\perp(x) - \mathfrak{p}^\perp(y)|$.*

(10.2)

The proof of this assertion assumes that there exists $\varepsilon \in (0,1)$ and no $\delta$ and disk $D_\varepsilon$ with the desired properties. This assumption is used to derive nonsense. If there is no $\delta$ and associated disk $D_\varepsilon$ as described, then there is a sequence $\{j_n\}_{n \in \{1,2,...\}}$ whose n'th member for $n \in \{1, 2, ...\}$ is a pair $(x_n, y_n)$ with both in the $\delta = 1/n$ version of $Z_{\delta,r}$, with $|\mathfrak{p}^\perp(x_n)|$ and $|\mathfrak{p}^\perp(y_n)|$ being less than $\frac{1}{n}$, and with $|x_n - y_n|$ greater than $(1 + \varepsilon)|\mathfrak{p}^\perp(x_n) - \mathfrak{p}^\perp(y_n)|$. To generate nonsense from such a sequence, invoke Proposition 4.1 with the input sequences being $\{p_n = \eta^{-1}(x_n)\}_{n \in \{1,2,...\}}$ and $\{\lambda_n = |x_n - y_n|\}_{n \in \{1,2,...\}}$. Let $(Z_*, \mathcal{I}_*, \nu_*)$ denote a limit data set that is then supplied by Proposition 4.1 and let $\Lambda$ denote the corresponding subsequence of $\{1, 2, ...\}$. The manner of convergence that is described by Proposition 4.1 again implies that the sequence $\{q_n = \lambda_n^{-1}(x_n - y_n)\}_{n \in \Lambda}$ has a subsequence that converges to a point on the $|x| = 1$ sphere. Let $q_*$ denote the limit point of such a subsequence.

Write this limit point $q_*$ as $(z_*, w_*)$ with respect to the splitting $\mathbb{R}^4 = \mathbb{R}^2 \times \mathbb{R}^{d-2}$. The norm of $z_*$ obeys $|z_*|^2 \geq (2\varepsilon + \varepsilon^2)|w_*|^2$ because it is assumed that $|x_n - y_n| \geq (1+\varepsilon)|\mathfrak{p}^\perp(x_n) - \mathfrak{p}^\perp(y_n)|$ for each index $n \in \{1, 2, \ldots\}$. Since $|z_*|^2 + |w_*|^2 = 1$, the inequality $|z_*|^2 \geq (2\varepsilon + \varepsilon^2)|w_*|^2$ implies that $z_*$ is non-zero. If $z_*$ is non-zero, then $q_*$ is not in $\{0\} \times \mathbb{R}^{d-2}$ and thus not in $Z_\diamond$. But this is nonsense because $\eta^{-1}(y_n)$ is in $Z$ for each $n \in \{1, 2, \ldots\}$; and granted that this is so, then Proposition 4.1 requires $\nu_*|_{q_*} = 0$. This is to say that $q_* \in Z_*$ which is $Z_\diamond$.

*Part 3*: Fix $\delta \in (0, 1)$ and $\varepsilon \in (0, 1]$ so that (10.1) and (10.2) hold. Let $D_\varepsilon \subset \mathbb{R}^{d-2}$ be as described in (5.2). Given $w \in D_\varepsilon \cap \mathfrak{p}^\perp(Z_{\delta,r})$, write the point in $(\mathfrak{p}^\perp)^{-1}(w) \cap Z_{\delta,r}$ with respect to the splitting $\mathbb{R}^4 = \mathbb{R}^2 \times \mathbb{R}^{d-2}$ as $(\varphi_Z(w), w)$. It follows from (10.2) that the map $\varphi_Z(\cdot)$ is such that $|\varphi_Z(w) - \varphi_Z(w')| \leq 2\varepsilon(1+\varepsilon)|w - w'|$ for any two points $w$ and $w'$ in $D_\varepsilon \cap \mathfrak{p}^\perp(Z_{\delta,r})$. The Kirszbraun theorem [K] (see, e.g. [V]) asserts that there is a Lipshitz map from $D_\varepsilon$ to $\mathbb{R}^2$, this denoted by $\varphi$, with the following properties:

- $\varphi = \varphi_Z$ on $\mathfrak{p}^\perp(Z_{\delta,r}) \cap D_\varepsilon$.
- *If $w$ and $w'$ are any two points in $D_\varepsilon$, then $|\varphi(w) - \varphi(w')| \leq 2\varepsilon(1+\varepsilon)|w - w'|$.*

(10.3)

A neighborhood of $p$ in $Z$ is contained in a strong Lipshitz graph because $\varphi$ is a Lipshitz map with a first order zero at the origin and a neighborhood in $Z_{\delta,r}$ of 0 is contained in the graph of $\varphi$.

*Part 4*: This part of the proof verifies that a neighborhood of $p$ in $Z$ is a codimensional 2 differentiable graph if it is a 2-dimensional Lipshitz graph. To see that this is so, let $D \subset \mathbb{R}^{d-k}$ denote a disk about the origin in the domain of the map $\varphi_Z$ from Part 3. It follows from (10.3) that $\varphi_Z$ is differentiable at the origin in $D$ and so the graph of $D$ has a tangent plane at the origin. Since each point in $Z$ near $p$ is strongly continuous, the analysis that leads to (10.3) can be applied with $p$ replaced by any point near $p$ in $Z$ to see that the $\varphi_Z$-graph of $D$ has a tangent plane at each point. Proposition 4.1 and either Proposition 9.2 or 9.1 can be used to prove that these tangent planes vary continuously on the $\varphi_Z$-graph of $D$. It follows from this and (10.3) that $\varphi_Z$ has continuous first derivatives near 0 in $D$ and so $p$ has a neighborhood that is contained in a $C^1$ graph.

**b) Points of continuity for $N_{(\cdot)}(0)$**

By way of a reminder, a point $p \in Z$ is said to be a point of continuity for $N_{(\cdot)}(0)$ when the following condition is met: If $\{p_i\}_{i \in \{1,2,\ldots\}}$ is a sequence in $Z$ that converges with limit $p$, then $\{N_{(p_i)}(0)\}_{i \in \{1,2,\ldots\}}$ converges to $N_{(p)}(0)$. A strongly continuous point is a point of continuity for $N_{(\cdot)}(0)$. The next lemma says that these are the only points of continuity.

**Lemma 10.3**: *The function $N_{(\cdot)}(0)$ is constant on some neighborhood of each of its points of continuity*.

It follows from this lemma that Proposition 10.1 is describing neighborhoods of the points of continuity for the function $N_{(\cdot)}(0)$.

This lemma is proved momentarily. Accept it as true for the moment.

*Proof of Theorem 1.4*: The set of strongly continuous points in Z is an open subset of Z. Supposing that this set is also dense in Z, then the assertion made by Theorem 1.4 follows from the characterization in Proposition 10.1. To see why the set of strongly continuous points in Z is dense in Z, suppose that p is a given point in Z. Given $r \in (0, r_0]$, define $o_{(p)}(r)$ to be the smallest half integer value of $N_{(\cdot)}(0)$ on $Z \cap B_r$. It follows from Lemmas 6.5 and 6.7 that this is a well defined function of r. It is also non-increasing; and it is less than $N_{(p)}(0)$ unless p is strongly continuous. Let $o_{(p)}(0)$ denote the $r \to 0$ limit of the numbers $\{o_{(p)}(r): r \in (0, r_0)\}$. It is a consequence of the definition of $o_{(p)}(0)$ that there exists a sequence $\{p_i\}_{i \in \{1,2,\ldots\}} \subset Z$ with limit p and with $N_{(\cdot)}(0) = o_{(p)}(0)$ on each member. Meanwhile, it follows from Lemma 6.7 that each sufficiently large i version of $p_i$ is a point of continuity for $N_{(\cdot)}(0)$; and Lemma 10.3 says that these points of continuity are strongly continuous.

*Proof of Lemma 10.3*: It follows from Lemma 6.7 that $N_{(\cdot)}(0)$ must be half an integer at each if its points of continuity. Suppose first that $p \in Z$ is a point of continuity for $N_{(\cdot)}(0)$ with $N_{(p)}(0)$ being an integer. It follows from Proposition 8.1 that the bundle $\mathcal{I}$ is isomorphic to the product $\mathbb{R}$ bundle on some ball centered at p. Let B denote such a ball. The identification between $\mathcal{I}$ and the product $\mathbb{R}$ bundle on B−Z writes $\nu$ as an $\mathbb{R}$-valued harmonic section of $\mathbb{S}^+$ on B−Z. Since the norm of $\nu$ is zero on Z, the $\mathbb{R}$-valued incarnation of $\nu$ extends over $B \cap Z$ to define a smooth harmonic section of $\mathbb{S}^+$ on B. This extension of $\nu$ has a non-trivial Taylor's expansion at each point in B of the form that is depicted in (8.11). Such an expansion can be used to see that $N_{(\cdot)}(0)$ has only integer values on $B \cap Z$ and so it must be constant on some neighborhood of any point of continuity for $N_{(\cdot)}(0)$ on $B \cap Z$.

The proof of the assertion in Lemma 10.3 for the case when $N_{(p)}(0)$ is half an odd integer has four parts.

*Part 1*: The proof starts with the following claim: A point p from Z is a point of continuity for $N_{(\cdot)}(0)$ if and only if p has a neighborhood in Z where $N_{(\cdot)}(0) \geq N_{(p)}(0)$. The 'if' implication follows from Lemma 3.3. To prove the 'only if' part of this claim, suppose for the sake of argument that it is false so as to generate nonsense. If it is false, then there is a point of continuity for $N_{(\cdot)}(0)$ in Z, this denoted by p, and a sequence $\{p_i\}_{i=1,2,...} \subset Z$ converging to p with $N_{(\cdot)}(0)$ less than $N_{(p)}(0)$ on each member. Keep in mind that the $\lim_{i\to\infty} N_{(p_i)}(0)$ is none-the-less equal to $N_{(p)}(0)$ because p is a point of continuity for $N_{(\cdot)}(0)$. It is also the case that $N_{(p)}(0)$ is half an integer. It follows from these last two observations that there is no loss when assuming that each $i \in \{1, 2, ...\}$ version of $N_{(p_i)}(0)$ is *not* equal to half of an integer. With this last assumption understood, fix for the moment $i \in \{1, 2, ...\}$. According to Lemma 6.7, there is a sequence $\{p_{i,k}\}_{k=1,2,...} \subset Z$ that converges to $p_i$ on which $N_{(\cdot)}(0)$ converges to half of a positive integer. Since $N_{(p_i)}(0) < N_{(p)}(0)$ it follows from Lemma 3.3 that this half of a positive integer must be less than $N_{(p)}(0)$ and so it is less than $N_{(p)}(0) - \frac{1}{4}$. It follows from this that there exists a 'diagonal' subsequence $\{p'_m\}_{m=1,2,...} \subset \{\{p_{i,k}\}_{k=1,2,...}\}_{i=1,2,...}$ that converges to p and such that $\lim_{m\to\infty} N_{(p'_m)}(0) < N_{(p)}(0) - \frac{1}{4}$. This existence of the latter sequence is nonsense because p is assumed to be a point of continuity for the function $N_{(\cdot)}(0)$.

*Part 2*: Suppose that p is a point of continuity for $N_{(\cdot)}(0)$ in Z with $N_{(p)}(0)$ being half of an odd integer. Assume for the sake of an argument that $N_{(\cdot)}(0)$ is not constant on some neighborhood of p. There is in this case a sequence $\{p_i\}_{i\in\{1,2,...\}} \subset Z$ in the radius $r_0$ ball centered at p with limit p and with $\{N_{(p_i)}(0)\}_{i\in\{1,2,...\}}$ converging to $N_{(p)}(0)$. Because of what is said in Part 1, nothing is lost by assuming that $N_{(p_i)}(0) > N_{(p)}(0)$ for each index i.

This sequence is used in the subsequent parts of the proof to generate nonsense thus proving the lemma.

*Part 3*: Choose an oriented, orthonormal frame for TX at p and parallel transport this frame along the geodesic arcs from p to define oriented, orthonormal frames for TX at the points in the radius $r_0$ ball centered at p. Use L to denote this frame. Fix the index i to be large. Having done so, choose a sequence $\{\lambda_k\}_{k\in\{1,2,...\}}$ in $(0, r_0)$ with limit zero and use the sequence $\{T(k) = (p_i, L, \lambda_k)\}_{k\in\{1,2,...\}}$ as input for Proposition 4.1. Let $(Z_{*i}, \mathcal{I}_{*i}, \nu_{*i})$ denote a limit data set from (5.1) that is supplied by Proposition 4.1 for this input sequence. Let $\{T(k)\}_{k\in\Lambda}$ with $\Lambda \subset \{1, 2, ...\}$ denote the subsequence that converges in the manner described by Proposition 4.1.

The desired nonsense comes from an appeal to the version of Lemma 5.7 with $N_\diamond$ being $N_{(p)}(0)$ and with $(Z_*, \mathcal{I}_*, \nu_*)$ being a sufficiently large i version of $(Z_{*i}, \mathcal{I}_{*i}, \nu_{*i})$. The assertion of this version of Lemma 5.7 to the affect that $N_0 = N_\diamond$ constitutes the desired nonsense because the $(Z_{*i}, \mathcal{I}_{*i}, \nu_{*i})$ version of $N_0$ is $N_{(p_i)}(0)$ and this is, by assumption, greater than $N_{(p)}(0)$.

With regards to an $N_\diamond = N_{(p)}(0)$ and $(Z_*, \mathcal{I}_*, \nu_*) = (Z_{*i}, \mathcal{I}_{*i}, \nu_{*i})$ version of Lemma 5.7, note that the requirements of the lemma's first bullet are met when i is large because $N_0$ is $N_{(p_i)}(0)$ and the $i \to \infty$ limit of sequence $\{N_{(p_i)}(0)\}_{i \in \{1,2,\ldots\}}$ is $N_{(p)}(0)$ which is $N_\diamond$. The next part of the proof verifies that the requirement of the second bullet of Lemma 5.7 is met when the index i is large.

*Part 4*: Suppose that $\varepsilon \in (0, \frac{1}{1000})$ and that $i \in \{1, 2, \ldots\}$ is such that the distance from p to $p_i$ is less than Lemma 9.3's number $r_\varepsilon$. Let q denote an $|x| = 1$ point in $Z_{*i}$ and suppose $\mathcal{I}_*$ is not isomorphic to the product $\mathbb{R}$ bundle on the complement of $Z_{*i}$ in any neighborhood of q. Fix $\rho \in (0, \frac{1}{100}\varepsilon)$ and let $B_\rho$ denote the radius $\rho$ ball centered on q. It follows from the manner of convergence that is dictated by Proposition 4.1 that if $k \in \Lambda$ is large, then there must exist points from $Z_{T(k)}$ in $B_\rho$. Indeed, if $B_\rho$ lacked points from $Z_{T(k)}$, then the bundle $\mathcal{I}_{T(k)}$ would be isomorphic to the product $\mathbb{R}$ bundle on $B_\rho$. If this were the case for a subsequence $\vartheta \subset \Lambda$, then the manner of convergence that is dictated by Proposition 4.1 could be used to construct an isomorphism between $\mathcal{I}_{*i}$ and the product $\mathbb{R}$ bundle on the complement in $B_\rho$ of $B_\rho \cap Z_{*i}$.

It follows from what was said in the preceding paragraph that there exists a sequence $\{q_k\}_{k \in \Lambda}$ that converges to q with the index k member being in $Z_{T(k)}$. Keeping in mind that value of $N_{(\cdot)}(0)$ on $\phi_{T(k)}(q_k)$ is no smaller than $N_{(p)}(0)$, cosmetic modifications of what is said in Step 3 of the proof of Lemma 6.2 can be used to prove that $\tilde{N}_q \geq N_{(p)}(0)$. Since $N_{(p)}(0)$ is $N_\diamond$ the inequality $\tilde{N}_q \geq N_{(p)}(0)$ verifies that the conditions stated by the second bullet of the $N_\diamond = N_{(p)}(0)$ and $(Z_*, \mathcal{I}_*, \nu_*) = (Z_{*i}, \mathcal{I}_{*i}, \nu_{*i})$ version of Lemma 5.7 are met when i is large.

c) **Points in Z where $N_{(\cdot)}(0)$ is not continuous**

Let $Z^{nc}$ to denote the set of points in Z where $N_{(\cdot)}(0)$ is not continuous. The upcoming Proposition 10.4 makes an assertion to the effect that $Z^{nc}$ has no part that is dense in a codimension 2 Lipshitz graph. This implies in particular that any subset of $Z^{nc}$ with finite d-2 dimensional Hausdorff measure is either nonrectifiable or full of d-2 dimensional holes. This is a curious result because of what is said in Lemma 6.8 about Z and thus $Z^{nc}$ having Hausdorff dimension at most d-2.

**Proposition 10.4**: *No part of $Z^{nc}$ is dense in a codimension 2 Lipshitz graph*.

*Proof of Proposition 10.4*: The proof has five steps.

Step 1: It follows from Proposition 7.1 that $Z^{nc}$ can be written as a finite union $(\cup_{m \in \{2,3,\ldots\}} \mathfrak{Z}_m) \cup Z'$ with $Z'$ having Hausdorff dimension no greater than d-3 and with $N_{(\cdot)}(0)$ on any given $m \in \{2, 3, \ldots\}$ version of $\mathfrak{Z}_m$ being equal to $\frac{1}{2}$m. The label m must be greater than 1 because Lemma 6.7 implies that $N_{(\cdot)}(0)$ is continuous near each point where it equals $\frac{1}{2}$.

Fix $p \in \mathfrak{Z}_m$ and let L denote a chosen, oriented orthonormal frame for $TX|_p$. Use $T = (p, L, \frac{1}{1000} r_0)$ to define the $(Z_T, \mathcal{I}_T, \nu_{*T})$ using the rules that are laid out at the start of Section 4. Fix an orthogonal decomposition of $\mathbb{R}^d$ as $\mathbb{R}^2 \times \mathbb{R}^{d-2}$, a ball $D \subset \mathbb{R}^{d-2}$ of radius less than 1 and a Lipshitz map $\varphi: D \to \mathbb{R}^2$ sending 0 to 0. Let $\psi: D \to \mathbb{R}^d$ denote the map that is defined by the rule $z \to (\varphi(z), z)$, this being the tautological parametrization of the graph of $\varphi$. Suppose that $Z$ is a given set in $\mathfrak{Z}_m$ with $p \in Z$ and with $\phi_T^{-1}(Z) \subset \psi(D)$. Assume for the sake of argument that for each $\varepsilon \in (0, 1)$, there is an open set $\mathcal{U}_\varepsilon \subset D$ obeying the following:

- *Each point in $\mathcal{U}_\varepsilon$ has distance less than $\varepsilon$ from 0.*
- *There is a dense set $\mathcal{U}_{\Diamond\varepsilon} \subset \mathcal{U}_\varepsilon$ with $\psi(\mathcal{U}_{\Diamond\varepsilon}) \subset \phi_T^{-1}(Z)$.*

(10.4)

The existence of small $\varepsilon$ versions of $\mathcal{U}_\varepsilon$ are used momentarily to generate nonsense.

<u>Step 2</u>: Since Z is closed and $\mathcal{U}_{\Diamond\varepsilon}$ is dense in $\mathcal{U}_\varepsilon$ and $\psi(\mathcal{U}_{\Diamond\varepsilon})$ is in $Z_T$, the whole of $\psi(\mathcal{U}_\varepsilon)$ is contained in $Z_T$. This understood, let $U_\varepsilon \subset Z_T$ denote $\psi(\mathcal{U}_\varepsilon)$ and let $U_{\Diamond\varepsilon} \subset U_\varepsilon$ denote the subset of $U_\varepsilon$ in $\phi_T^{-1}(Z)$. It follows from Lemma 3.2 that the only points in $\phi_T(U_\varepsilon)$ where $N_{(\cdot)}(0)$ is half an integer are points where $N_{(\cdot)}(0) = \frac{1}{2} m$. At the same time, each point in $\phi_T(U_\varepsilon)$ is the limit of a sequence of points in Z where $N_{(\cdot)}(0)$ is half an integer strictly less than m. Keeping the preceding in mind, let q denote an interior point in $U_\varepsilon$ from $U_{\Diamond\varepsilon}$ and let $\{p_n\}_{n \in \{1,2,\ldots\}}$ denote a sequence in Z that converges to $\phi_T(q)$ with the value of $N_{(\cdot)}(0)$ on each element being half of some integer that is strictly less than m. For each index $n \in \{1, 2, \ldots\}$, let $x_n$ denote the $\phi_T^{-1}(p_n)$.

Let $\mathfrak{p}^\perp: \mathbb{R}^d \to \mathbb{R}^{d-2}$ again denote the orthogonal projection. If the index i is large, then the splitting of $\mathbb{R}^d$ as $\mathbb{R}^2 \times \mathbb{R}^{d-2}$ can be used to write $x_n$ as $(w_n + \varphi(\mathfrak{p}^\perp(x_n)), \mathfrak{p}^\perp(x_n))$. Note in particular that $w_n \neq 0$ because $x_n$ is not in the image of $\psi$.

<u>Step 3</u>: If the index n is large, then there is a sequence $\{q_{n,k}\}_{k \in \{1,2,\ldots\}} \subset U_{\Diamond\varepsilon}$ with $\lim_{k \to \infty} \mathfrak{p}^\perp(q_{n,k}) = \mathfrak{p}^\perp(x_n)$. This is because $\mathfrak{p}^\perp(U_\varepsilon)$ is an open set in $\mathbb{R}^{d-2}$ and $U_{\Diamond\varepsilon}$ is dense in $U_\varepsilon$. Write the vector $q_{n,k} - x_n$ as

$$(w_n + \varphi(\mathfrak{p}^\perp(q_{n,k})) - \varphi(\mathfrak{p}^\perp(x_n)), \mathfrak{p}^\perp(q_{n,k}) - \mathfrak{p}(x_n)) .$$

(10.5)

Since $\varphi$ is continuous, choose for each index n, an index k(n) such that both

$$|\mathfrak{p}^\perp(q_{n,k(n)}) - \mathfrak{p}^\perp(x_n)| < \varepsilon |w_n| \quad \text{and} \quad |\varphi(\mathfrak{p}^\perp(q_{n,k(n)})) - \varphi(\mathfrak{p}^\perp(x_n))| \leq \varepsilon |w_n| \ .$$
(10.6)

Fix a large index n so that $|x_n - q| < \varepsilon$, then use $\lambda_\varepsilon$ to denote $|w_n|$, use $\hat{\mathfrak{p}}_\varepsilon$ to denote the corresponding version of $p_n$, and use $q_\varepsilon$ to denote the corresponding versin of $q_{n,k(n)}$.

Because $\mathfrak{p}^\perp(U_{\diamond\varepsilon})$ is dense in $\mathfrak{p}^\perp(U_\varepsilon)$, a point in $U_{\diamond\varepsilon}$ to be denoted by $q_\varepsilon{'}$ obeys

$$|\mathfrak{p}^\perp(q_\varepsilon{'}) - \mathfrak{p}^\perp(q_\varepsilon)| \in [\tfrac{1}{2}\lambda_\varepsilon, 2\lambda_\varepsilon] \ .$$
(10.7)

In the case when $\dim(X) = 4$, a second point in $U_{\diamond\varepsilon}$ can be found, this one to be denoted by $q_\varepsilon{''}$, that obeys (10.7) and is such that $\langle \mathfrak{p}(q_\varepsilon{'}), \mathfrak{p}(q_\varepsilon{''}) \rangle \leq \varepsilon \lambda_\varepsilon^2$.

Step 4: The sequence $\{\phi_T(q_{\varepsilon=1/i})\}_{i \in \{1,2,\ldots\}}$ is a sequence in $\mathfrak{Z}_m$ that converges to q. Let i denote a given positive integer and let T(i) denote the data set $(\phi_T(q_{1/i}), L, \lambda_{\varepsilon=1/i})$. The sequence indexed by $\{1, 2, \ldots\}$ whose i'th member is the $\phi_{T(i)}$-inverse image of $\phi_T(\hat{\mathfrak{p}}_{\varepsilon=1/i})$ is a bounded sequence in $\mathbb{R}^d$ with norm bounded aways from zero. That this is so follows from (10.5) and (10.6). Meanwhile, the sequence indexed by $\{1, 2, \ldots\}$ with i'th member being the $\phi_{T(i)}$-inverse image of the point $\phi_T(q_{\varepsilon=1/i}{'})$ is also bounded and bounded away from zero; this being a consequence of (10.7). Such is also the case when $d = 4$ for the sequence whose i'th member is the $\phi_{T(i)}$-inverse image of $\phi_T(q_{\varepsilon=1/i}{''})$.

With the preceding understood, choose a subsequence $\Theta \subset \{1, 2, \ldots\}$ so that the first two sequences converge, and when $d = 4$, so that the third sequence converges also. Let $\hat{p}$ denote the limit of the first sequence, let $q'$ denote the limit of the second sequence and when $d = 4$, let $q''$ denote the limit of the third. It follows from (10.6) that $\mathfrak{p}^\perp(\hat{p}) = 0$; and it follows from (10.7) that $\mathfrak{p}^\perp(q') \neq 0$. In the case when $d = 4$, it likewise follows that $\mathfrak{p}^\perp(q'') \neq 0$; and it follows from what is said in the last sentence of Part 3 that $\mathfrak{p}^\perp(q')$ is orthogonal to $\mathfrak{p}^\perp(q'')$.

Step 5: Use sequence $\{T(i)\}_{i \in \Theta}$ as input for Proposition 4.1 after renumbering $\Theta$ consecutively from 1. The corresponding limit data set $(Z_*, \mathfrak{I}_*, \nu_*)$ is the line spanned by $q'$ in the case when the dimension of X is three and it is the plane spanned by the vectors $q'$ and $q''$ in the case when the dimension of X is four. That this is so follows from Lemma 6.2 given the manner of convergence that is dictated by Proposition 4.1 Meanwhile, $\hat{p}$ is a point in $Z_*$ with distance between $\tfrac{1}{2}$ and 2 from the origin. This last conclusion is the desired nonsense that proves the proposition. It is nonsense in the case $d = 3$ because $\mathfrak{p}^\perp(\hat{p}) = 0$ and $\mathfrak{p}^\perp(q') \neq 0$; and it is also nonsense in the case $d = 4$ because $\mathfrak{p}^\perp(q') \neq 0$, $\mathfrak{p}^\perp(q'') \neq 0$ and $\mathfrak{p}^\perp(q')$ is orthogonal to $\mathfrak{p}^\perp(q'')$.

### d) Points of discontinuity for $I$

A point in Z is said to be a *point of discontinuity for $I$* if the line bundle $I$ is not isomorphic to the product $\mathbb{R}$ bundle on the complement of Z in any neighborhood of the point. The following proposition characterizes in part the points of discontinuity for $I$.

**Proposition 10.5**: *The set of points of discontinuity for $I$ is the closure of an open set in Z that has the structure of a codimension 2, differentiable submanifold. This open set is the set of points of continuity of $N_{(\cdot)}(0)$ where its value is half of an odd integer.*

*Proof of Proposition 10.5*: Let $C$ denote for the moment the set of points of continuity for $N_{(\cdot)}(0)$ where its value is half an odd integer. It follows from Lemma 10.4 and Proposition 10.1 that $C$ has the structure of a differentiable, codimension 2 submanifold. It follows meanwhile from Lemma 10.4, Lemma 6.2 and Proposition 9.1 that each point in $C$ is a point of discontinuity for $I$. The proof that $C$ is a dense in the set of point discontinuity of $I$ has four parts.

*Part 1*: Write Z as the nested union $Z_{\leq 1} \subset Z_{\leq 2} \subset \cdots \subset Z_{\leq N} = Z$ where any integer m version of $Z_{\leq m}$ denotes the set of points in Z where $N_{(\cdot)}(0) \leq \frac{1}{2} m$ and where N is the least integer upper bound for the function $2N_{(\cdot)}(0)$ on Z. The proof that $C$ is dense in the set of points of discontinuity for $I$ starts with the observation that $Z_{\leq 1}$ is the set of $N_{(\cdot)}(0) = \frac{1}{2}$ points in Z and so it is part of $C$. In particular, $C \cap Z_{\leq 1}$ is dense in the set of points of discontinuity for $I$ in $Z_{\leq 1}$. Given the preceding, assume that $m \in \{1, 2, \ldots, N-1\}$ and that $C \cap Z_{\leq m}$ is known to be dense in the set of points of discontinuity for $I$ in $Z_{\leq m}$. The subsequent parts of the proof uses this assumption to prove that $C \cap Z_{\leq m+1}$ is dense in the set of points of discontinuity for $I$ in $Z_{\leq m+1}$.

*Part 2*: Suppose that $p \in Z_{\leq m+1}$ is a point of discontinuity for $I$ that is not in $C$. The point p can not be a point of continuity for $N_{(\cdot)}(0)$ so there is a sequence $\{p_i\}_{i \in \{1,2,\ldots\}}$ of points in Z with limit p and such that $\lim_{i \to \infty} N_{(p_i)}(0) \leq N_{(p)}(0)$. It then follows from Lemma 6.7 that there is a sequence of this sort with each $i \in \{1, 2, \ldots\}$ version of $N_{(p_i)}(0)$ being half an integer. Either such a sequence exists with each $i \in \{1, 2, \ldots\}$ version of $N_{(p_i)}(0)$ being half an odd integer or not. Suppose that they are half of odd integers. Proposition 4.1 can be used to see that each element of this sequence is a point of discontinuity for $I$. This understood, it follows from the induction step that there is a sequence of points in $C$ with limit p also.

*Part* 3: The following lemma is needed to deal with the event that p has a neighborhood in Z with no points where $N_{(·)}(0)$ is less than $N_{(p)}(0)$ and half of an odd integer.

**Lemma 10.6**: *Suppose that* U *is an open set in* X *with the following two properties*:
- *The first cohomology with* $\mathbb{Z}/2$ *coefficients is zero.*
- *Given any point* $p \in U$, *there exists a ball centered at* p, *this denoted by* B, *such that* $\mathcal{I}$ *on* $B-(B \cap Z)$ *is isomorphic to the product* $\mathbb{R}$ *bundle.*

*Then* $\mathcal{I}$ *is isomorphic to the product* $\mathbb{R}$ *bundle on* $U-(U \cap Z)$.

*Proof of Lemma 10.6*: If $p \in U \cap Z$, let $B_\diamond(p)$ denote a ball centered at p of radius less than $\frac{1}{100} r_0$ such that $\mathcal{I}$ on $B_\diamond(p)-(B_\diamond(p) \cap Z)$ is isomorphic to the product $\mathbb{R}$ bundle. Let B(p) denote the concentric ball with half the radius of the ball $B_\diamond(p)$. The collection $\{B(p): p \in U\}$ forms an open cover of U. Let $\mathfrak{U}$ denote a locally finite subcover. If $B \in \mathfrak{U}$, let $\iota_B$ denote an isomorphism from $\mathcal{I}$ on $B-(B \cap Z)$ to the product $\mathbb{R}$ bundle on $B-(B \cap Z)$. Given balls B and B´ from $\mathfrak{U}$, the composition $\iota_B \iota_{B'}^{-1}$ defines an automorphism of the product $\mathbb{R}$ bundle on the complement in $B \cap B'$ of the set Z. Denote this by $\varepsilon_{B,B'}$. This is a locally constant function on the complement of Z in $B \cap B'$ with values in $\{\pm 1\}$. The subsequent paragraph explains why $\varepsilon_{B,B'}$ is in fact contant. Assuming that this is so, then the constant +1 or -1 value of $\varepsilon_{B,B'}$ can be assigned to the whole of $B \cap B'$. The collection $\{\varepsilon_{B,B'}: B, B' \in \mathfrak{U}\}$ is then the cocycle data for a real line bundle on U because the cocycle condition $\varepsilon_{B,B'} \varepsilon_{B',B''} \varepsilon_{B'',B} = 1$ is obeyed when B, B´ and B´´ are any three sets in $\mathfrak{U}$ with a point in common. This real line bundle is a priori isomorphic to the product $\mathbb{R}$ bundle because real line bundles on U are classified by the $\mathbb{Z}/2$ first cohomology. It is also canonically isomorphic to $\mathcal{I}$ on $U-(U \cap Z)$. Thus, $\mathcal{I}$ is isomorphic to the product $\mathbb{R}$ bundle on $U-(U \cap Z)$.

This paragraph explains why the $\pm 1$ valued function $\varepsilon_{B,B'}$ is in fact constant. The bundle $\mathcal{I}$ is assumed to be isomorphic to the product $\mathbb{R}$ bundle on $B_\diamond-(B_\diamond \cap Z)$. Fix such an isomorphism and use it to view $\nu$ as an $\mathbb{R}$ valued harmonic 1-form on $B_\diamond$. The zero locus of $\nu$ on $B_\diamond$ is the critical point set of a harmonic function and this understood, Taylor's theorem with remainder can be used to see that Z does not disconnect any convex set in $B_\diamond$. This implies in particular that the complement of Z in $B \cap B'$ is path connected. Since $\varepsilon_{B,B'}$ is locally constant on this complement, it must have the same value everywhere.

*Part* 4: Suppose that p is a point in $Z_{\leq m+1}$ with a neighborhood that has no points where $N_{(·)}(0)$ is less than $N_{(p)}(0)$ and half of an odd integer. Let B denote a ball of radius less than $\frac{1}{100} r_0$ centered at p with this property. The first observation is that $N_{(\cdot)}(0)$ is an integer on all points in $B \cap Z$ where $N_{(\cdot)}(0)$ is less than $\frac{1}{2}(m+1)$. The argument for this has four steps.

  Step 1: If $q \in B \cap Z$ and $N_{(q)}(0) \leq \frac{1}{2} m$, then $\mathcal{I}$ is isomorphic to the product $\mathbb{R}$ bundle on the complement of Z in a neighborhood of q. This understood, it then follows from Lemma 6.6 that $N_{(\cdot)}(0)$ at q is an integer. Suppose next $q \in B \cap Z$ and assume, for the sake of argument, that $N_{(q)}(0) \in (\frac{1}{2} m, \frac{1}{2}(m+1))$. This assumption is used to generate nonsense.

Suppose first that $\dim(X) = 3$. In this case, q is not weakly continuous. Since the set of points in Z that are not weakly continuous is countable, there is a point in $B \cap Z$, this denoted by q´ with $N_{(q')}(0) \in (\frac{1}{2} m, \frac{1}{2}(m+1))$ with a neighborhood where $N_{(\cdot)}(0) \leq \frac{1}{2} m$ on all points from Z–q´. Let B´ denote a ball centered at q´ with this property. It follows from the induction hypothesis and the U = B´–q´ version of Lemma 10.4 that $\mathcal{I}$ is isomorphic to the product $\mathbb{R}$ bundle on the complement of Z in B´. This being the case, Lemma 6.6 says that $N_{(q')}(0)$ is an integer. This is nonsense because it runs afoul of the assumptions.

  Step 2: The argument in Step 1 to rule out points in $B \cap Z$ where $N_{(\cdot)}$ is strictly between $\frac{1}{2} m$ and $\frac{1}{2}(m+1)$ works when X has dimension 4 if there are no weakly continuous points in $B \cap Z$ with $N_{(\cdot)}(0)$ strictly between $\frac{1}{2} m$ and $\frac{1}{2}(m+1)$. To prove that no weakly continuous points have $N_{(\cdot)}(0)$ in this range, suppose to the contrary that there were such a point. Denote it by q. Since q is weakly continuous, there exists a sequence $\{p_i\}_{i \in \{1,2,\ldots\}} \subset B \cap Z$ with limit q and with $\lim_{i \to \infty} N_{(p_i)}(0) = N_{(q)}(0)$. Fix an oriented, orthonormal frame for TX at q to be denoted by L and for each $i \in \{1, 2, \ldots\}$, define the data set T(i) to be $(q, L, \lambda_i = \text{dist}(q, p_i))$. Use q instead of p in Lemma 6.2 to see that a data set $(Z_*, \mathcal{I}_*, \nu_*)$ that is supplied by Proposition 4.1 using q instead of p and the input sequence $\{T(i)\}_{i \in \{1,2,\ldots\}}$ is such that $Z_*$ is either a line through the origin in $\mathbb{R}^4$ or a non-empty, finite union of half planes through the origin in $\mathbb{R}^4$ with a common boundary line.

  Step 3: The case when $Z_*$ is a line can be ruled because Item a) of the third bullet of Lemma 5.5 implies that $N_{(q)}(0)$ is an integer if $Z_*$ is a line. If $Z_*$ is a union of half planes through the origin, then the function $Ň_{(\cdot)}$ on the interior of the set of $|x| = 1$ points in at least one of these half planes must be half of an odd integer. If not, then $\mathcal{I}_*$ would be isomorphic to the product $\mathbb{R}$ bundle on $\mathbb{R}^4 – Z_*$, such an isomorphism would identify $\nu_*$ as an $\mathbb{R}$-valued 1-form, and then a Taylor's expansion of $\nu_*$ based at the origin proves that $N_{(q)}(0)$ is an integer.

Let $\Pi \subset Z_*$ denote a half plane with $\tilde{N}_{(\cdot)}(0)$ being half an odd integer on the set of its $|x| = 1$ points that are not on the common boundary line. Let U denote a very small radius ball centered on an $|x| = 1$ point in the interior of $\Pi$. Choose the radius so that all points from $Z_*$ in the closure of U are in the interior of $\Pi$. Let $\{T(i)\}_{i \in \Lambda}$ denote the subsequence that is supplied by q's version of Proposition 4.1. It follows from the manner of convergence that is dictated by Proposition 4.1 that all sufficiently large $i \in \Lambda$ versions of $\mathcal{I}_{T(i)}$ are not isomorphic to the product $\mathbb{R}$ bundle on $U-(U \cap Z_{T(i)})$. It follows as a consequence of Lemma 10.4 that $\phi_{T(i)}(U \cap Z_{T(i)})$ contains an uncountable set of points of discontinuity for $\mathcal{I}$ if the index i is sufficiently large.

Meanwhile, it follows from the fifth bullet of Lemma 6.2 that the function $N_{(\cdot)}(0)$ on $\phi_{T(i)}(U \cap Z_{T(i)})$ is no greater than $\frac{1}{2}(N_{(q)}(0)+m)$ if i is sufficiently large. Fix an index i with this property, and then choose a weakly continuous point $q' \in \phi_{T(i)}(U \cap Z_{T(i)})$ where $N_{(\cdot)}(0) \leq \frac{1}{2}(N_{(q)}(0)+m)$ which is a point of discontinuity for $\mathcal{I}$. Note that it follows from the induction hypothesis that $N_{(q')} > \frac{1}{2} m$.

Step 4: Let $q_1$ denote $q'$. Repeat the constructions in Steps 2 and 3 using $q_1$ in lieu of q and let $q_2$ denote the new version of what Step 3 denotes by $q'$. Then return to Step 3 using $q_2$ in lieu of q and let $q_3$ denote the new version of $q'$. Continue in this vein to define a sequence $\{q_1, q_2, \ldots\}$ with the following properties:

- *Each point in $\{q_1, q_2, \ldots\}$ is a point of discontinuity for $\mathcal{I}$.*
- *Fix $\varepsilon > 0$, and there exists a ball $B_\varepsilon \subset B$ such that all sufficiently large i versions of $q_i$ are in $B_\varepsilon$ and $N_{(\cdot)}(0) < \frac{1}{2}m+\varepsilon$ on $B_\varepsilon \cap Z$.*

(10.8)

It follows from the second bullet of (10.8) that the sequence $\{q_i\}_{i \in \{1,2,\ldots\}}$ has a limit point in B to be denoted by $q_*$ with $N_{(q_*)}(0) = \frac{1}{2}m$. The induction hypothesis says that $\mathcal{I}_*$ is isomorphic to the product $\mathbb{R}$ bundle on the complement of Z in a neighborhood of $q_*$. This last observation is nonsensical because it runs afoul of the first bullet of (10.8).

*Part 5*: It follows from what is said in Part 4 that $N_{(p)}(0) = \frac{1}{2}(m+1)$. If p is not weakly continuous, then $N_{(\cdot)}(0) \leq \frac{1}{2}m$ on the points from Z in a neighborhood of p. This too follows from what is said in Part 4. Lemma 10.6 in this case implies that $\mathcal{I}$ is isomorphic to the product $\mathbb{R}$ bundle on the complement of Z in some neighborhood of p.

Supposing that p is weakly continuous, there are two cases to consider; they are distinguished by the parity of the integer $m+1$. Assume first that $m+1$ is even. The argument that proves Proposition 8.1 can be repeated with only cosmetic changes to prove that $\mathcal{I}$ is isomorphic to the product $\mathbb{R}$ bundle near Z. The cosmetic changes amount to the following: Define $Z_{(m+1)}$ to be the part of Z where $N_{(\cdot)}(0) \geq \frac{1}{2}(m+1)$ and replace Z in the proof by the set of points in $Z_{m+1}$. The homotopies that are made in the various parts of Section 8a's proof of Proposition 8.1 can ignore the rest of Z because the bundle $\mathcal{I}$ is isomorphic to the product $\mathbb{R}$ bundle on some neighborhood of each point in $Z \cap B_r$ where $N_{(\cdot)}(0) < \frac{1}{2}(m+1)$. No more is said about these changes except with regards to Lemma 8.4 because this lemma plays such a central role in the proof of Proposition 8.1.

The role played by Lemma 8.4 is played here by the upcoming Lemma 10.7. The background for Lemma 10.7 is much like that for Lemma 8.4 as described at the beginning of Part 2 in Section 8a. By way of a summary, $\varepsilon$ is chosen from $(0, \frac{1}{1000})$. With $\varepsilon$ chosen, then it follows from Lemma 3.4 that there exists $r_\varepsilon \in (0, \frac{1}{1000} r_0)$ so that $N_{(\cdot)}(0) \leq \frac{1}{2}(m+1) + \varepsilon$ on the ball of radius $r_\varepsilon$ centered at p, and so that conclusions of Lemma 6.2 can be invoked for points from Z in this ball. If $r \in (0, r_\varepsilon)$, the number $D(r)$ is defined just as in Part 1 of the proof of Proposition 8.1. Fix such r and then $q \in Z_{m+1} \cap B_r$. Lemma 10.7 uses T to denote the data set $(q, L, \lambda = D(r))$ with L being an oriented orthonormal frame for $TX|_q$ that is obtained from a chosen frame at $TX|_p$ by parallel transport along the short geodesic arc between the two points. The lemma uses $Z_{T,m+1}$ to denote the part of $Z_T$ that is mapped by $\phi_T$ to $Z_{m+1}$. The rest of the notation in Lemma 10.7 is the same as that in Lemma 8.4.

**Lemma 10.7**: *There exists $\kappa > 1$ that can depend on p but not q and has the following significance: Suppose that $\varepsilon \in (0, \kappa^{-1}]$ and that $r \in (0, r_\varepsilon]$. If the dimension of X is three, then $Z_{T,m+1} \cap \mathbb{B}_{1/\varepsilon}$ has distance $\varepsilon$ or less from a line through the origin; and the bundle $\phi_T^* \mathcal{I}$ is isomorphic to the product $\mathbb{R}$ bundle on the complement in $\mathbb{B}_{1/\varepsilon}$ of the radius $\varepsilon$ tubular neighborhood of this line. If the dimension X is four, then $Z_{T,m+1} \cap \mathbb{B}_{1/\varepsilon}$ has distance at most $\varepsilon$ from either a plane through the origin or a line through the origin, and the bundle $\phi_T^* \mathcal{I}$ is isomorphic to the product $\mathbb{R}$ bundle on the complement in $\mathbb{B}_{1/\varepsilon}$ of the radius $\varepsilon$ tubular neighborhood of this plane or line.*

The proof of this lemma differs only in notation from the proof of Lemma 8.4.

*Part 6*: The eight steps that follow prove that if $m+1$ is an odd integer, then p is in the $N_{(\cdot)}(0) = \frac{1}{2}(m+1)$ part of C given the following assumptions: The point p is weakly continuous, the value of $N_{(p)}(0)$ is $\frac{1}{2}(m+1)$, and there is no point near p where $N_{(\cdot)}(0)$ is half an odd integer that is strictly less than $m+1$.

Step 1: If the dimension of X is four, then p, being weakly continuous, can be in $Z^p$ or in $Z^{np}$. In the latter case, an almost verbatim repetition of the arguments from Steps 2-4 of Part 4 will prove that $\mathcal{I}$ is isomorphic to the product $\mathbb{R}$ bundle on the complement of Z in some neighborhood of p. This being the case, assume henceforth that p is in $Z^P$ when the dim(X) = 4. Since p is weakly continuous, there exists a sequence $\{p_i\}_{i \in \{1,2,...\}}$ in the radius $r_0$ ball centered at p with limit p and with $\lim_{i \to \infty} N_{(p_i)}(0) = \frac{1}{2}(m+1)$. When the dim(X) = 4, there exists a second sequence with this property, $\{p_i'\}_{i \in \{1,2,...\}}$, and such that the sequence $\{\Delta_{(p)}(p_i, p_i')\}_{i \in \{1,2,...\}}$ is bounded away from zero and the sequence $\{\ln(\text{dist}(p,p_i)/\text{dist}(p,p_i'))\}_{i \in \{1,2,...\}}$ is bounded from above and below. Fix some large index i and set $r = \text{dist}(p, p_i)$. Parallel transport a chosen an oriented, orthonormal frame for $TX|_p$ along the geodesic rays from p to define a frame for TX on the radius $r_0$ ball centered at p. Use this frame to define the matrix T(r) using the formula in (9.6). If the index i is large, and thus r is small, this matrix will have a d-2 linearly independent eigenvectors with very small eigenvalues and then 2 linearly independent eigenvectors with $\mathcal{O}(1)$ eigenvalues. The argument from Part 1 of Section 9c can be used to prove this. Let $Z_{*r} \subset \mathbb{R}^d$ denote d-2 dimensional subspace that is spanned by the eigenvectors of T(r) with very small eigenvalue.

Step 2: Supposing that the index i is large and so r is small, set $T = (p, L, r)$. It follows from Lemma 6.2 that the whole of the $|x| \leq 4$ part of $Z_T$ is very close to $Z_{*r}$. The following analog of Lemma 9.8 says more in this regard. The notation has $Z_{T,m+1}$ denoting the set of points in $Z_T$ that are mapped by $\phi_T$ to the $N_{(\cdot)}(0) \geq \frac{1}{2}(m+1)$ part of Z.

**Lemma 10.8**: *There exists $\kappa > 1$ with the following significance: Suppose that the index $i \in \{1, 2, ...\}$ is large so that $r = \text{dist}(p, p_i)$ is small. Use $T = (p, L, r)$ to define $Z_T \subset \mathbb{R}^d$ as instructed at the beginning of Section 4 and to define $Z_{*r}$ as instructed above*
- *The points in $Z_{T,m+1}$ where $|x| \in (\frac{2}{3}, 4)$ have distance at most $\kappa(\underline{N}_{(p)}(r) + r)^{1/2}$ from $Z_{*r}$.*
- *Suppose that q is a point in the $|x| \in (\frac{2}{3}, 4)$ part of $Z_{*r}$.*
  a) *There is a point from $Z_{T,m+1}$ in the normal plane to $Z_{*r}$ through q with distance at most $\kappa(\underline{N}_{(p)}(r) + r)^{1/2}$ from q.*
  b) *The distance between any two points from $Z_T$ in the normal plane to $Z_{*r}$ through q is no greater than $\kappa e^{-1/(\kappa \underline{N}_{(p)}(r))}$.*

*Proof of Lemma 10.8*: The proof is essentially the same as that of Lemma 9.8 except that Item a) of the second bullet needs the induction hypothesis to conclude that the points of discontinuity for $\mathcal{I}_T$ in the normal plane to $Z_{*r}$ at q must come from $Z_{T,m+1}$.

Step 3: Fix a large index i that can be used in Lemma 10.8. Let q and q´ denote orthogonal points points in the $|x| = \frac{3}{4}$ part of $Z_{*r}$ and let $p_{i1}$ and $p_{i1}'$ denote the $\phi_T$ images of points from $Z_{T,m+1}$ in the respective normal planes to $Z_{*r}$ at q and q´. Repeat step 2 using $p_{i1}$ in lieu of $p_i$ to obtain a version of the two bullets from Lemma 10.8 that has the

data set T replaced by $(p, L, \frac{3}{4} r)$, the subspace $Z_{*r}$ replaced by the span of the eigenvectors of $T(\frac{3}{4} r)$ with very small eigenvalue, and in general, r is replaced everywhere by $\frac{3}{4} r$.

Step 4: Repeat Step 3 using the new version of Lemma 10.8's two bullets to obtain a second version of the two bullets of Lemma 10.8; this new version has T being $(p, L, (\frac{3}{4})^2 r)$ and $Z_{*r}$ replaced by the spanby the span of the eigenvectors of $T((\frac{3}{4})^2 r)$ with very small eigenvalue, in general, r replaced everywhere by $(\frac{3}{4})^2 r$. Continue in this vein iterating Step 3 to see that the assertions of the two bullets in Lemma 10.8 hold with r being any positive number that is less than some large i version of dist$(p, p_i)$, with T being $(p, L, r)$ and with $Z_{*r}$ being the d-2 dimensional span of the eigenvectors of $T(r)$ with small very small eigenvalues.

Step 5: The results in Step 4 imply that there exists $r_* \in (0, \frac{1}{100} r_0)$ such that if r´ is positive and less than $r_*$, then there is a point from $Z_{m+1}$ on the boundary the radius r´ ball centered at p that is a point of discontinuity for $\mathcal{I}$. Fix r´ $\in (0, r_*)$ and let p´ denote such a point. The analysis in Steps 1-4 can be repeated with p´ used in lieu of p to obtain a p´ version of Lemma 10.8. This is stated below as Lemma 10.9. The lemma views L as an orthonormal frame for TX at p´ that is obtained from its namesake at p by parallel transport along the short geodesic segment between p and p´.

**Lemma 10.9**: *There exists* $\kappa > 1$ *and* $r_* \in (0, \frac{1}{100} r_0)$ *with the following significance: Let* p´ *denote a point of discontinuity of* $\mathcal{I}$ *from* $Z_{m+1}$ *in the radius* $r_*$ *ball centered at* p. *Supposing that* $r \in (0, r_*)$, *let* $T_{p'}(r)$ *denote the* p´ *version of the matrix* $T(r)$ *in (9.6)*.
- *The matrix* $T_{p'}(r)$ *has d-2 linearly independent eigenvectors with eigenvalue less than* $\kappa(\underline{N}_{(p')}(r) + r)$. *Let* $Z'_{*,r}$ *denote the span of these eigenvectors*.
- *The points in* $Z_{T',m+1}$ *where* $|x| \in (\frac{2}{3}, 4)$ *have distance at most* $\kappa(\underline{N}_{(p')}(r) + r)^{1/2}$ *from* $Z'_{*,r}$.
- *Suppose that* q *is a point in the* $|x| \in (\frac{2}{3}, 4)$ *part of* $Z'_{*,r}$. *Let* $T' = (p', L, r)$.
  a) *There is a point from* $Z_{T',m+1}$ *in the normal plane to* $Z'_{*,r}$ *through* q *with distance at most* $\kappa(\underline{N}_{(p')}(r) + r)^{1/2}$ *from* q.
  b) *The distance between any two points from* $Z_{T'}$ *in the normal plane to* $Z'_{*,r}$ *through* q *is no greater than* $\kappa e^{-1/(\kappa \underline{N}_{(p')}(r))}$.

Step 6: With Lemma 10.9 in hand, the arguments that prove Proposition 9.2 can be repeated with the only significant modification being the replacement of K in the proof by the part of $Z_{m+1}$ in the radius $\frac{1}{4} r_*$ ball centered at p. These arguments lead to the following lemma.

**Lemma 10.10**: *There exists* $r_* \in (0, \frac{1}{100} r_0)$ *and* $\kappa > 1$ *with the following significance: Let K denote the set of points in* $Z_{m+1}$ *from the radius* $\frac{1}{4} r_*$ *ball centered at* p. *Supposing that* $p' \in K$ *and* $r \in (0, \kappa^{-1} r_*)$, *then* $\underline{N}_{(p')}(r) \leq \kappa\, e^{-\sqrt{|\ln(r)|}/\kappa}$.

By way of a remark about Proposition 9.2: A close reading of the proof of Proposition 9.2 finds that the key inputs are summarized by the three bullets in Lemma 10.9.

Given Lemma 10.10, then the arguments that are used to prove Proposition 9.1 can be repeated with only cosmetic changes to prove the next lemma.

**Lemma 10.11**: *There exists* $r_* \in (0, \frac{1}{100} r_0)$ *with the following significance: The set* Z *has a unique tangent cone at each point of discontinuity for* $\mathfrak{I}$ *from the part of* $Z_{m+1}$ *in the radius* $\frac{1}{4} r_*$ *ball centered at* p.

Step 7: The existence of unique tangent cones on the part of $Z_{m+1}$ that is described by Lemma 10.11 can be used with the same arguments prove Proposition 10.1 to prove the following assertion:

*The set of points of discontinuity for $\mathfrak{I}$ from the part of $Z_{m+1}$ in the radius $\frac{1}{8} r_*$ ball centered at p is a codimension 2, differentiable submanifold of* X.

(10.9)

Note that the Item a) of the second bullet of Lemma 10.9 implies that the set in question is not just part of a differentiable submanifold, it is the whole of the submanifold.

Step 8: Given (10.9) and the unique tangent cone assertion in Lemma 10.11, then the arguments that lead to Proposition 10.4 can be repeated to see that there are no points in the radius $\frac{1}{16} r_*$ ball centered at p where $N_{(\cdot)}(0)$ is less than $m+1$. This implies that p is strongly continuous and thus p is in the set *C*.